\documentclass[12pt]{amsart}
\usepackage{amsfonts,amssymb,amsmath,amsthm,amstext,color}

\usepackage[colorlinks, linkcolor=blue, citecolor=blue, urlcolor=blue,
pagebackref,hypertexnames=false]{hyperref}
\usepackage{enumerate}
\usepackage{array}
\usepackage{geometry}
\usepackage{mathrsfs}
\usepackage{txfonts}
\usepackage{tikz}
\usetikzlibrary{arrows.meta, positioning, decorations.pathreplacing}
\usepackage{url}

\usepackage{comment}
\allowdisplaybreaks
\numberwithin{equation}{section}
\newtheorem{mainthm}{Theorem}

\newtheorem{theorem}{Theorem}[section]
\newtheorem{proposition}[theorem]{Proposition}
\newtheorem{lemma}[theorem]{Lemma}
\newtheorem{corollary}[theorem]{Corollary}
\theoremstyle{definition}

\theoremstyle{remark}
\newtheorem{remark}[theorem]{Remark}

\title[Inequalities for flag square functions]
{Chang--Wilson--Wolff and Hoeffding inequalities for flag square functions}

\author{Yongsheng Han}
\address{Department of Mathematics, Auburn University, AL 36849-5310, USA}
\email{hanyong@auburn.edu}

\author{Ji Li}
\address{School of Mathematical and Physical Sciences, Macquarie University, NSW 2109, Australia}
\email{ji.li@mq.edu.au}

 \author{Brett D. Wick}
\address{Department of Mathematics, Washington University in St.~Louis, St.~Louis, MO 63130-4899, USA}
\email{wick@math.wustl.edu}

\author{Liangchuan Wu}
\address{School of Mathematical Sciences, Anhui University, Hefei 230601, P.R.~China}
\email{wuliangchuan@ahu.edu.cn}


\begin{document}

\begin{abstract}
We prove sharp local and global reverse $L^p$ inequalities for the Poisson area function $S_F^{(k)}$ of a $k$-step flag on $\mathbb R^N$. For a product $R=E_1\times\cdots\times E_k$ of axis-parallel cubes in the coordinate blocks, set $Q_R=\prod_{j=1}^k(\mathrm{Id}-P_j)$, where $P_j$ averages over $E_j$ in the $j$th block. If $f\in L^u(\mathbb R^N)$ for some $1\le u<\infty$, then for $2\le p<\infty$,
$$
\|Q_Rf\|_{L^p(R)}\le Cp^{k/2}\|S_F^{(k)}(f)\|_{L^p(R)},
$$
with all scales included in the area. When $S_F^{(k)}(f)\le1$ on $R$, we obtain a Chang--Wilson--Wolff type estimate $|R|^{-1}\int_R\exp(c|Q_Rf|^{2/k})\,dx\le2$. The factor $p^{k/2}$ and the power $2/k$ are sharp even on the unit cube. Globally, $f\in L^p$ if and only if $S_F^{(k)}(f)\in L^p$, with optimal constants of order $p^{k/2}$ in both directions. For $f-f_R$, the sum of the terms centred in $\ell$ blocks satisfies a local estimate with the areas of the corresponding coordinate averages and the sharp factor $p^{\ell/2}$.

For $k=2$, every $L^2$ function with bounded flag area belongs to the closure of $L^\infty$ in flag $\mathrm{BMO}$. Every $L^2$ function in this closure satisfies $\sup_R|R|^{-1}\int_R\exp(\lambda|Q_Rf|^{1/2})\,dx<\infty$ for every $\lambda>0$. The power $1/2$ is sharp on a fixed flag rectangle.

We split the reconstruction kernels according to their coordinate cancellations. In blocks without cancellation, decay in the scale ratios gives summability over side lengths. In the remaining blocks, we iterate a dyadic martingale estimate using the intermediate square functions of Li, Pipher and Wu. Each such block contributes $\sqrt p$ to the local Hoeffding inequality for finite sums of atoms.
\end{abstract}

\maketitle
\tableofcontents


\section{Introduction}\label{s:intro}


For $f\in L^2(\mathbb R^d)$, let $u(y,t)=e^{-t\sqrt{-\Delta}}f(y)$ be its Poisson extension.
The Lusin area integral with the full gradient is
$$
S_{\mathrm{cl}}^{(d)}(f)(x)^2
=\int_0^\infty\int_{|x-y|<t}
\left |t\nabla_{y,t}u(y,t) \right |^2\,\frac{dy}{t^d}\frac{dt}{t},
$$
where $|x-y|$ is the Euclidean distance and $\nabla_{y,t}=(\nabla_y,\partial_t)$.
Stein asked how much local integrability follows from boundedness of the area integral (see \cite[p.~219]{CWW}).
If $S_{\mathrm{cl}}^{(d)}(f)\in L^\infty$, then $f\in\mathrm{BMO}$, and the John--Nirenberg inequality
yields exponential integrability of $|f-f_Q|$ on every cube $Q$, where $f_Q=|Q|^{-1}\int_Q f$.
Chang, Wilson and Wolff~\cite[Theorem~3.2]{CWW} proved exponential integrability of $|f-f_Q|^2$
under a bound on the Poisson area formed with the spatial gradient.
Since $|\nabla_yu|\leq|\nabla_{y,t}u|$, for $0<\|S_{\mathrm{cl}}^{(d)}(f)\|_\infty<\infty$ their theorem implies
$$
\sup_Q\frac1{|Q|}\int_Q
\exp\bigg(c_d\frac{|f(x)-f_Q|^2}
{\|S_{\mathrm{cl}}^{(d)}(f)\|_\infty^2}\bigg)\,dx
\leq C_d,
$$
where $\|\cdot\|_\infty$ denotes the essential supremum on $\mathbb R^d$ and $c_d,C_d>0$ depend only on $d$.
Hoeffding's inequality~\cite{Hoeffding} establishes a Gaussian tail bound for the deviation of a sum
of independent bounded random variables from its mean. There the sum of the squared lengths
of their ranges takes the place of $\|S_{\mathrm{cl}}^{(d)}(f)\|_\infty^2$.
For dyadic martingale differences, \eqref{e:scalar-square} bounds the $L^p$ norm of their sum by
$C\sqrt p$ times the $L^p$ norm of their pointwise square function, $2\leq p<\infty$.

The exponential-square estimate implies the moment bound
$$
\bigg(\frac1{|Q|}\int_Q|f-f_Q|^p\,dx\bigg)^{1/p}
\leq C\sqrt p\,\|S_{\mathrm{cl}}^{(d)}(f)\|_\infty,\qquad 2\leq p<\infty.
$$
For $f\in L^p(\mathbb R^d)$, the global reverse inequality has the same factor~\cite{CWW,FP,Wilson08},
$$
\|f\|_{L^p(\mathbb R^d)}
\leq C\sqrt p\,\|S_{\mathrm{cl}}^{(d)}(f)\|_{L^p(\mathbb R^d)},
\qquad 2\leq p<\infty.
$$
For the full gradient, the examples in Proposition~\ref{p:flag-reverse-sharp} with $k=1$
show that the factor $\sqrt p$ in both moment estimates and the exponential power $2$ are sharp.
For weighted estimates see~\cite{ChanilloWheeden,Wilson89,Wilson08,CDWY}, and
for semigroups that need not preserve constants see~\cite{WuJGA}.
Related results on spaces of homogeneous type include exponential-square integrability~\cite{WuYan}
and dyadic good-$\lambda$ inequalities~\cite{CLOW}.

For double dyadic martingales on $[0,1]^2$, Pipher~\cite{Pipher86}
proved that the best constant in the reverse $L^p$ square-function
estimate has order $p$. Her proof uses the fact that a martingale
difference $d_I$ takes opposite values on the two halves of a dyadic
interval $I$, so $|d_I(x)|=\|d_I\|_{L^\infty(I)}$ for $x\in I$.
On a higher-dimensional dyadic cube, the values of a martingale
difference on the children sum to zero and can have different
magnitudes. In a product, the child on which the supremum in the
first variable is attained may therefore depend on the second
variable. This is where a direct iteration of the one-parameter
estimate breaks down.

Li, Pipher and Wu~\cite{LiPipherWu} used the pointwise values of
the differences to form their intermediate square functions.
With the first variable fixed, the sum inside the $\ell^2$ norm
is still a martingale sum in the second variable. The one-parameter
estimate can thus be applied in each variable in turn, with a
factor $\sqrt p$ at each step. For $k$ parameters the resulting
factor is $p^{k/2}$. In particular, if a martingale sum $F$ on a
product of dyadic cubes has bounded square function, then
$|F|^{2/k}$ is exponentially integrable on that product.%
\footnote{Pipher's sharpness example is a product of two normalized
Rademacher sums. The same exponents occur in the moment and tail
bounds for $L^2$-normalized polynomials of degree at most $k$ in
independent Rademacher variables~\cite[Theorems~9.21 and~9.23]{ODonnell}.}

In the product setting, each Poisson average acts on a separate coordinate block.
The flag averages studied here act on nested subspaces. Already for a two-step flag,
both scale parameters enter the average in the second coordinate block.

\subsection{The two-step flag and its Poisson area}\label{s:intro-background}

Flag kernels arose in the study of singular integrals on the Heisenberg group~\cite[Section~1]{NRSW1}.
M\"uller, Ricci and Stein~\cite{MRS,MRS2} proved that the joint Marcinkiewicz multipliers of the sub-Laplacian
and the central derivative have flag kernels and are bounded on $L^p$ for $1<p<\infty$.
In the Euclidean model of~\cite{NRS,HLLW}, the two-step flag on $\mathbb R^{n+m}$, with $n,m\geq1$, is
$$
\{0\}\subset\{0\}\times\mathbb R^m\subset\mathbb R^{n+m}.
$$
Let $P_t^{(q)}$ denote the Poisson kernel on $\mathbb R^q$. The average at scale $t$ acts on
$(x,y)\in\mathbb R^n\times\mathbb R^m$, and the one at scale $s$ acts on $y$ alone.
Their composition has kernel
$$
P_{t,s}(x,y)=(P_t^{(n+m)}*_y P_s^{(m)})(x,y)
=\int_{\mathbb R^m}P_t^{(n+m)}(x,y-v)P_s^{(m)}(v)\,dv,
$$
where $*_y$ denotes convolution in $y$. The Poisson extension of $f$ is
$$
U_f(x,y,t,s)=(P_{t,s}*f)(x,y)
=e^{-t\sqrt{-\Delta_{x,y}}}e^{-s\sqrt{-\Delta_y}}f(x,y),
\qquad t,s>0,
$$
as in~\cite[(1.7)--(1.8)]{HLLW}. The function $U_f$ is harmonic in $(x,y,t)$ and in $(y,s)$.
The integral in $y$ is $\int_{\mathbb R^m}P_{t,s}(x,y)\,dy=P_t^{(n)}(x)$.
Integrating in $x$ and using the semigroup identity, we obtain
$$
\int_{\mathbb R^n}P_{t,s}(x,y)\,dx=(P_t^{(m)}*P_s^{(m)})(y)=P_{t+s}^{(m)}(y).
$$
Thus the $x$-scale is $t$ and the $y$-scale is $t+s$. Table~\ref{tab:two-scales} compares these
with the scales for product areas and Zygmund dilations~\cite{Zygmund,HanLiLinTanZygmund}.

\begin{table}[htbp]
\caption{Spatial scales for two independent parameters $t,s>0$.}\label{tab:two-scales}
\centering
\small
\begin{tabular}{@{}l@{\qquad}l@{\qquad}l@{}}
\hline
Setting & Spatial scales & For fixed $t$ \\ 
\hline 
Product & $(t,s)$ & Varying $s$ changes only the second scale.\\[1mm]
Zygmund & $(t,s,ts)$ & The third scale is $t$ times the second.\\[1mm]
Two-step flag & $(t,t+s)$ & For $s\leq t$, both scales are comparable to $t$.\\
              &           & For $s>t$, the second scale is comparable to $s$.\\
\hline
\end{tabular}
\end{table}

The parameters $t,s$ vary independently, but the $x$-scale $t$ is always at most the $y$-scale $t+s$.
A product $I\times J$ of cubes is called a flag rectangle if $\ell(I)\leq\ell(J)$, where $\ell(I)$ denotes the side length of $I$.
Partial convolution also affects cancellation. In the formula that reconstructs $f$ from its Poisson
derivatives, the kernel has integral zero in $y$ for each $x$. Its integral in $x$ can be a nonzero
function of $y$. This distinction determines where we can use a martingale estimate in the local proof.

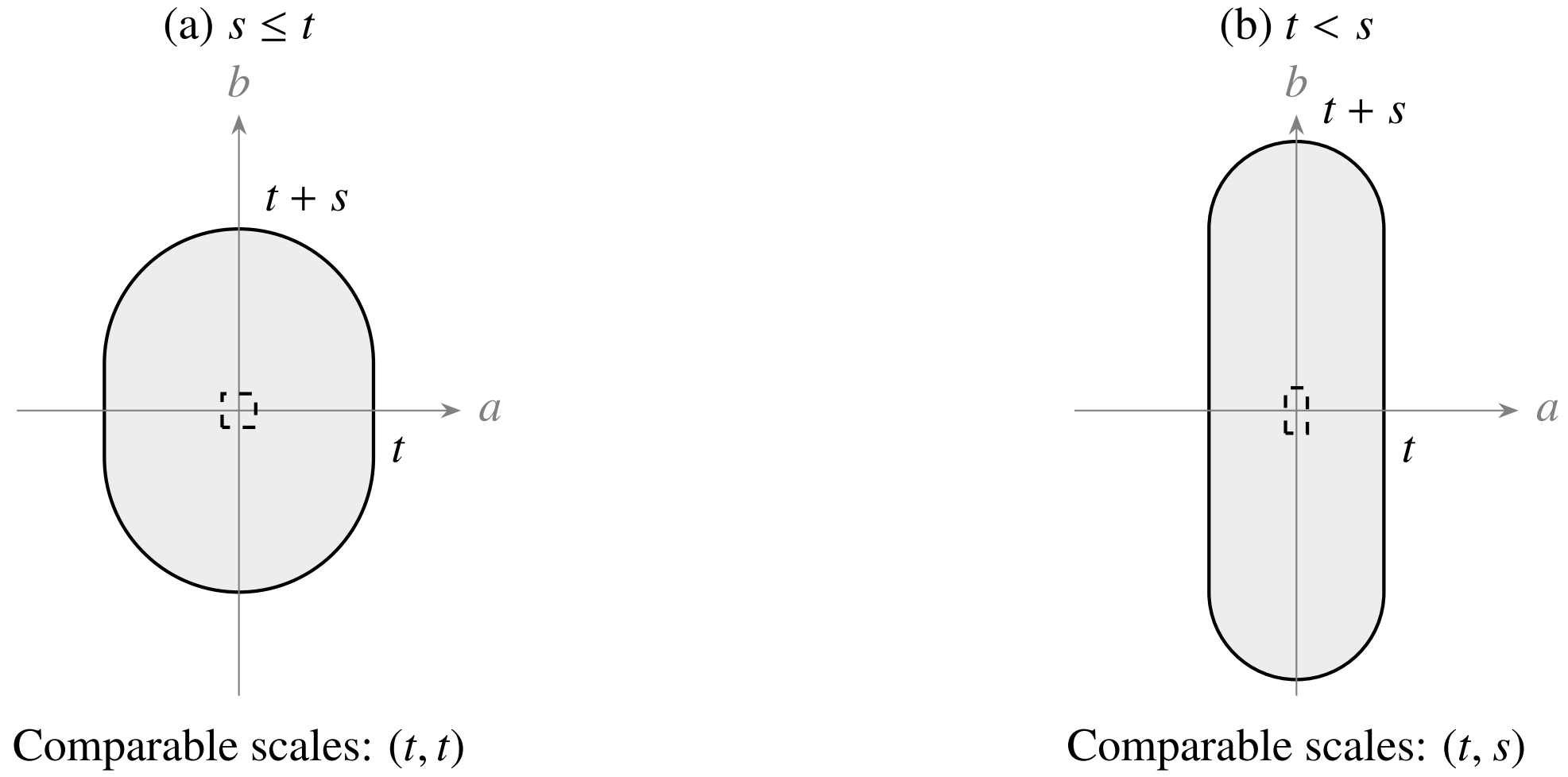
\begin{figure}[bp]
\centering
\begin{minipage}[t]{0.48\textwidth}
\centering
\begin{tikzpicture}[x=1.12cm,y=1.12cm,>=Stealth,
every node/.style={font=\small}]
\path[fill=black!7,draw=black,thick]
(-1,-.35)--(-1,.35) arc[start angle=180,end angle=0,radius=1]
--(1,-.35) arc[start angle=0,end angle=-180,radius=1]--cycle;
\draw[->,gray] (-1.65,0)--(1.65,0) node[right] {$a$};
\draw[->,gray] (0,-2.12)--(0,2.2) node[above] {$b$};
\draw[dashed,thick] (-.125,-.125) rectangle (.125,.125);
\node[anchor=north west] at (1,-.06) {$t$};
\node[anchor=south west] at (.06,1.35) {$t+s$};
\node at (0,2.85) {(a) $s\leq t$};
\node at (0,-2.52) {Comparable scales: $(t,t)$};
\end{tikzpicture}
\end{minipage}\hfill
\begin{minipage}[t]{0.48\textwidth}
\centering
\begin{tikzpicture}[x=1.12cm,y=1.12cm,>=Stealth,
every node/.style={font=\small}]
\path[fill=black!7,draw=black,thick]
(-.65,-1.35)--(-.65,1.35) arc[start angle=180,end angle=0,radius=.65]
--(.65,-1.35) arc[start angle=0,end angle=-180,radius=.65]--cycle;
\draw[->,gray] (-1.65,0)--(1.65,0) node[right] {$a$};
\draw[->,gray] (0,-2.12)--(0,2.2) node[above] {$b$};
\draw[dashed,thick] (-.08125,-.16875) rectangle (.08125,.16875);
\node[anchor=north west] at (.65,-.06) {$t$};
\node[anchor=south west] at (.06,2) {$t+s$};
\node at (0,2.85) {(b) $t<s$};
\node at (0,-2.52) {Comparable scales: $(t,s)$};
\end{tikzpicture}
\end{minipage}
\caption{Horizontal and vertical displacements $(a,b)$ in the averaging region at scales $t,s$ for $n=m=1$.
The shaded region is the union of the disks of radius $t$ centred on the vertical segment from $(0,-s)$ to $(0,s)$.
The dashed rectangle has half-width $t/8$ and half-height $\max\{t,s\}/8$.}\label{fig:flag-cones}
\end{figure}

Let $B_1$ and $B_2$ be the open Euclidean unit balls in $\mathbb R^{n+m}$ and $\{0\}\times\mathbb R^m$, respectively.
For fixed $a_2\in B_2$, the points $(x,y)-ta_1-sa_2$ with $a_1\in B_1$ fill the ball of radius $t$
about $(x,y)-sa_2$. Figure~\ref{fig:flag-cones} shows the union of these balls for $n=m=1$.

Write $\nabla_1=(\partial_t,\nabla_x,\nabla_y)$ and $\nabla_2=(\partial_s,\nabla_y)$.
With the Euclidean norm on the mixed derivatives $t\nabla_1s\nabla_2U_f$, the flag Poisson area
of~\cite[Definition~1.3.1]{HLLW} is
$$
S_F(f)(x,y)^2
=\int_0^\infty\int_0^\infty\int_{B_1\times B_2}
\big|t\nabla_1s\nabla_2U_f\big((x,y)-ta_1-sa_2,t,s\big)\big|^2
\,da_1\,da_2\,\frac{dt}{t}\frac{ds}{s}.
$$
The equivalent formula in the spatial variables $(u,v)$ is \eqref{e:poisson-area-weight}.

Nagel, Ricci and Stein~\cite{NRS} developed the Euclidean theory of flag kernels, and Nagel, Ricci,
Stein and Wainger~\cite{NRSW1,NRSW2} extended it to homogeneous groups.
Partial convolution also occurs for Euclidean flag singular integrals. For multipliers
$m_1(\xi,\eta)$ and $m_2(\eta)$ satisfying the Mikhlin estimates of all orders, the kernel of
$m_1(\xi,\eta)m_2(\eta)$ is $K_1*_yK_2$ in the sense of distributions, where $K_1,K_2$ are
the respective kernels. For example, $R_i^{(1)}R_j^{(2)}$ is the composition of a Riesz transform
on $\mathbb R^{n+m}$ and one acting in $y$ on $\mathbb R^m$ (Subsection~\ref{s:hardy-estimates}).
Han and Lu~\cite{HL} constructed the Euclidean flag Hardy spaces.
Han, Lee, Li and Wick~\cite[Sections~1.2--1.3]{HLLW} characterized the flag Hardy space $H_F^1$ by the flag Poisson area,
by maximal functions and by flag Riesz transforms. For the weak $L\log L$ endpoint see~\cite{HLLWEndpoint},
and for flag Hardy spaces on the Heisenberg group see~\cite{HLS}.

\subsection{Local oscillation and the main questions}\label{s:intro-question}

Let $R=I\times J$ be a product of axis-parallel cubes in $\mathbb R^n$ and $\mathbb R^m$,
and let $f_R=|R|^{-1}\int_R f$. The two coordinate averages are
$$
f_J(x)=\frac1{|J|}\int_J f(x,y)\,dy,\qquad
f_I(y)=\frac1{|I|}\int_I f(x,y)\,dx.
$$
With no restriction on the side lengths of $I$ and $J$, which part of $f-f_R$ can we estimate by the flag area on $R$ itself?
The derivative $s\nabla_2$ annihilates the Poisson extension of every function of $x$ alone, so $S_F(f)$ does not
detect oscillation in $x$ that is independent of $y$.
The flag area controls the part of $f$ whose average vanishes in each variable separately.
On $R$ this doubly centred part is
$$
Q_Rf=f-f_J-f_I+f_R.
$$
For example, if $f(x,y)=g(x)h(y)$, then $Q_Rf=(g-g_I)(h-h_J)$ on $R$, where $g_I$ and $h_J$ are the means on $I$ and $J$.
For $f\in L^u(\mathbb R^{n+m})$ with $1\leq u<\infty$, the two-step estimate is
$$
\|Q_Rf\|_{L^p(R)}\leq Cp\|S_F(f)\|_{L^p(R)},\qquad 2\leq p<\infty.
$$
Both norms are taken on the same $R$, and the area function at each point of $R$ includes all scales.
Thus an $L^\infty$ bound for $S_F(f)$ on $R$ alone suffices for exponential integrability of $|Q_Rf|$ throughout $R$.
To recover the full oscillation on $R$, we also need the oscillations of the two coordinate averages, since
$$
f-f_R=Q_Rf+(f_J-f_R)+(f_I-f_R).
$$

These observations lead to four questions.
\begin{itemize}
\item \textbf{Local and global estimates} (Theorem~\ref{t:concentration}).
For a flag of $k$ steps, how does the best constant in the local reverse estimate grow with $p$?
What are the sharp orders in the two global $L^p$ inequalities between $f$ and its area function?

\smallskip

\item \textbf{The full oscillation} (Theorem~\ref{t:coordinate}).
Which area functions control the remaining terms in $f-f_R$?
How does the power of $p$ depend on the number of coordinate blocks left after averaging?

\smallskip

\item \textbf{Approximation by bounded functions} (Theorem~\ref{t:approximation}).
For a two-step flag, suppose $f\in L^2$ has bounded area on the whole space.
Can one approximate $f$ in flag $\mathrm{BMO}$ by a bounded function, and how large an $L^\infty$ norm suffices for a given error?

\smallskip

\item \textbf{Exponential integrability in the closure} (Theorem~\ref{t:classes}).
For $L^2$ functions in the closure of $L^\infty$ in flag $\mathrm{BMO}$, which powers of $|Q_Rf|$ are exponentially integrable uniformly over all rectangles?
Does a uniform bound for the averages of $\exp(c|Q_Rf|)$, for some $c>0$, characterize this closure?
\end{itemize}

\subsection{Finite flags and sharp concentration}\label{s:intro-flags}

The flags of Nagel, Ricci, Stein and Wainger have arbitrary length~\cite{NRSW1}.
For a flag of $k\geq1$ steps, let $X_j=\mathbb R^{d_j}$, $d_j\geq1$, and define
$$
V_j=X_j\oplus\cdots\oplus X_k,\qquad 1\leq j\leq k.
$$
Thus $D_j=\dim V_j=d_j+\cdots+d_k$ and $N=D_1$. We identify each $V_j$ with the last $k-j+1$ coordinate spaces in $\mathbb R^N$, so that
$V_1\supset\cdots\supset V_k$. The Poisson semigroup at scale $t_j$ acts on $V_j$, with
Laplacian $\Delta_{V_j}=\sum_{i=j}^k\Delta_{X_i}$. For $\mathbf t=(t_1,\ldots,t_k)\in(0,\infty)^k$, let
$$
U_f(x,\mathbf t)=\prod_{j=1}^k e^{-t_j\sqrt{-\Delta_{V_j}}}f(x).
$$
Here $\nabla_j=(\partial_{t_j},\nabla_{V_j})$ denotes the time and spatial gradient on the $j$th subspace.
Each component of $\mathcal G_{\mathbf t}f=(\prod_{j=1}^k t_j\nabla_j)U_f$ uses one entry of each $\nabla_j$.
We take the Euclidean norm of this array.

At scales $\mathbf t$, the spatial displacement from $x$ is a sum $t_1a_1+\cdots+t_ka_k$,
where $a_j$ ranges over the Euclidean unit ball $B_j=B_{V_j}(0,1)$ in $V_j$. In these variables the flag Poisson area
is defined by
$$
S_F^{(k)}(f)(x)^2
=\int_{(0,\infty)^k}\int_{B_1\times\cdots\times B_k}
 \bigg|\bigg(\prod_{j=1}^k t_j\nabla_j\bigg)U_f
 \bigg(x-\sum_{j=1}^kt_ja_j,\mathbf t\bigg)\bigg|^2
 \,da_1\cdots da_k\prod_{j=1}^k\frac{dt_j}{t_j},
$$
where $da_j$ is Lebesgue measure on $V_j$. For $k=1$ this coincides with $S_{\mathrm{cl}}^{(d_1)}$.
For $k=2$ it coincides with the area $S_F$ of Subsection~\ref{s:intro-background}.
Formula \eqref{e:flag-area} writes the integrand for general $k$ as a sum over the components of the mixed gradient.
See also~\cite[Section~8.6]{NRSW1} for convolution square functions along finite flags
and~\cite[pp.~549--551]{HLW} for the Euclidean discrete construction.

Let $R=E_1\times\cdots\times E_k$, where each $E_j$ is an axis-parallel cube in $X_j$. For
$x=(x_1,\ldots,x_k)\in R$, define the average in block $j$ by
$$
P_jf(x)=\frac1{|E_j|}\int_{E_j}f(x_1,\ldots,x_{j-1},u_j,x_{j+1},\ldots,x_k)\,du_j.
$$
Subtracting the mean in each block in turn, we obtain
$$
Q_Rf=\prod_{j=1}^k(I-P_j)f\quad\text{on }R.
$$
Here $I$ is the identity. Subsection~\ref{s:flag-local} extends each centred difference by zero in $x_j$ outside $E_j$.
Only the first $j$ displacements act in $X_j$, so the projection of the averaging region onto $X_j$
is a ball of radius $t_1+\cdots+t_j$, comparable to $\max_{i\leq j}t_i$. These radii increase with $j$.
Accordingly, we call $R$ a flag rectangle if $\ell(E_1)\leq\cdots\leq\ell(E_k)$. The local estimates in
the following theorem hold for every ordering of the side lengths, since the localization argument applies
separately in each coordinate (Subsection~\ref{s:intro-method}).

\begin{mainthm}\label{t:concentration}
Let $f\in L^u(\mathbb R^N)$ for some $1\leq u<\infty$.
\begin{itemize}
\item[(i)] \textbf{Local reverse estimate and Chang--Wilson--Wolff inequality.}
For every product $R=E_1\times\cdots\times E_k$ of axis-parallel cubes
$E_j\subset\mathbb R^{d_j}$,
$$
\|Q_Rf\|_{L^p(R)}
\leq Cp^{k/2}\|S_F^{(k)}(f)\|_{L^p(R)},\qquad 2\leq p<\infty.
$$
If $0<\|S_F^{(k)}(f)\|_{L^\infty(R)}<\infty$, then
$$
\frac1{|R|}\int_R
\exp\bigg[c\, \bigg(\frac{|Q_Rf|}{\|S_F^{(k)}(f)\|_{L^\infty(R)}}\bigg)^{2/k}\bigg]\,dx\leq2.
$$
The order $p^{k/2}$ and the exponential power $2/k$ are sharp even for $R=[-1/2,1/2]^N$.

\smallskip

\item[(ii)] \textbf{Global $L^p$ estimates.}
For every $2\leq p<\infty$, $f\in L^p(\mathbb R^N)$ if and only if
$S_F^{(k)}(f)\in L^p(\mathbb R^N)$. In this case,
$$
(Cp^{k/2})^{-1}\|f\|_{L^p(\mathbb R^N)}
\leq\|S_F^{(k)}(f)\|_{L^p(\mathbb R^N)}
\leq Cp^{k/2}\|f\|_{L^p(\mathbb R^N)}.
$$
The factor $p^{k/2}$ is sharp in both directions.
\end{itemize}
The constants $c,C>0$ depend only on $k,d_1,\ldots,d_k$.
\end{mainthm}

For the full Poisson area, the upper bound in \textup{(ii)} fails if both norms are taken on $R$.
Indeed, if $f$ is the indicator of a cube and $R$ is an adjacent cube, then both $f$ and $Q_Rf$
vanish almost everywhere on $R$, yet $S_F^{(k)}(f)$ has positive $L^2(R)$ norm. For functions supported in $R$, the upper bound holds on $R$
with the sharp factor $p^{k/2}$ (Remark~\ref{r:local-upper}).

If $0<\|S_F^{(k)}(f)\|_{L^\infty(R)}<\infty$, we apply Markov's inequality to the exponential
in \textup{(i)} to obtain the concentration bound
$$
|\{x\in R:|Q_Rf(x)|>\tau\|S_F^{(k)}(f)\|_{L^\infty(R)}\}|
\leq 2|R|e^{-c\tau^{2/k}},\qquad \tau>0.
$$
For $k=1$ this is a Gaussian tail bound for $f-f_R$. For $k=2$ the doubly centred oscillation has an exponential tail.

Earlier local results for product areas estimate a restricted Calder\'on sum.
Pipher~\cite[Theorem~2]{Pipher86} obtained an exponential estimate on an open set
$\Omega\subset\mathbb R^2$ for the Calder\'on sum over dyadic rectangles contained in $\Omega$,
under a bound on the product area over $\mathbb R^2$.
Li, Pipher and Wu~\cite[Theorem~5.2]{LiPipherWu} proved the corresponding estimate on a product of two Shilov boundaries.
Their dyadic theorem also allows almost orthogonal families without exact cancellation~\cite[Theorem~B]{LiPipherWu}.
In Theorem~\ref{t:concentration}\,\textup{(i)}, we centre the full function $f$ on $R$ and measure
its area on that same $R$. The flag reconstruction also includes terms whose integral vanishes
in only one of the two blocks. For these terms the scale-ratio estimates in
Subsection~\ref{s:intro-method} control the sum in the other block.

To deduce the global reverse estimate in \textup{(ii)} from \textup{(i)}, use the rectangles
$R_L=\prod_{j=1}^k[-L,L]^{d_j}$. Since $f\in L^u$ with $u<\infty$, H\"older's inequality on almost every slice
shows that each coordinate average in the expansion of $Q_{R_L}f$ tends to zero.
Thus $\mathbf1_{R_L}Q_{R_L}f\to f$ almost everywhere. We obtain the global bound by Fatou's lemma
(Subsection~\ref{s:flag-reverse}).

For $f\in L^p(\mathbb R^N)$, the global reverse inequality also follows by iterating the
one-parameter estimate on the flag subspaces. At each step the Poisson derivatives form an
$L^2$-valued function whose pointwise norm is the intermediate square function. The next
Poisson operator acts on each component.
As in the product argument of Fefferman and Pipher~\cite[Section~5]{FP}, each application of the
Hilbert-valued estimate contributes $\sqrt p$ (Remark~\ref{r:global-intermediate}).

For the global upper bound in \textup{(ii)}, we iterate the weighted area estimate of Chanillo and Wheeden~\cite[Theorem~1(a)]{ChanilloWheeden}
on the flag subspaces (Subsection~\ref{s:flag-lp}). To see that $p^{k/2}$ is necessary in both global bounds,
we test the upper $L^p$ bound for $S_F^{(k)}$ on a cube indicator and the reverse estimate on functions
with separated coordinate frequencies
(Subsection~\ref{s:flag-sharpness}).

\subsection{Coordinate averages in the full oscillation}\label{s:intro-coordinate}

For $R=I\times J$, Theorem~\ref{t:concentration}\textup{(i)} bounds the doubly centred term $Q_Rf$.
The two remaining terms in $f-f_R$ are the oscillations of $f_J$ on $I$ and of $f_I$ on $J$.
To control them, we use the one-parameter area function defined with the second Poisson
derivative. For $g$ on $\mathbb R^q$, write
\begin{equation}\label{e:intro-one-area}
\mathcal S_P^{(q)}g(x)^2
=\int_0^\infty\int_{|x-y|<t}
|t^2(-\Delta)e^{-t\sqrt{-\Delta}}g(y)|^2\frac{dy}{t^q}\frac{dt}{t}.
\end{equation}
For the coordinate averages, set
$$
\mathcal E_{P,R}(f)(x,y)^2
=\mathbf1_J(y)|\mathcal S_P^{(n)}f_J(x)|^2
+\mathbf1_I(x)|\mathcal S_P^{(m)}f_I(y)|^2.
$$
Combining the bound for $Q_Rf$ on $R$ with the one-parameter bounds for $f_J-f_R$ on $I$
and $f_I-f_R$ on $J$ as in Theorem~\ref{t:two-coordinate}, we obtain, for
$f\in L^u(\mathbb R^{n+m})$ with $1\leq u<\infty$ and $2\leq p<\infty$,
\begin{align*}
\|f-f_R\|_{L^p(R)}
\leq &Cp\|S_F(f)\|_{L^p(R)} 
 +C\sqrt p\,\bigg\|\Big(|\mathcal S_P^{(n)}f_J(x)|^2
                +|\mathcal S_P^{(m)}f_I(y)|^2\Big)^{1/2}\bigg\|_{L^p(R)}.
\end{align*}

Both coordinate areas are necessary even for $p=2$. For each integer $M\geq1$,
Proposition~\ref{p:coord-necessity} constructs a real Schwartz function $f$ with
$\|S_F(f)\|_\infty\leq C\sqrt M$. On the unit cube $R=I\times J$ centred at the origin,
$f_I=0$ and
$$
\inf_{\zeta\in\mathbb C}\bigg(\frac1{|R|}\int_R|f(x,y)-\zeta|^2\,dx\,dy\bigg)^{1/2}\geq cM.
$$
Here $Q_Rf$ stays uniformly bounded on $R$, so the large oscillation comes from $f_J$.
The second example satisfies the same bounds with $f_J=0$ on a rectangle whose
horizontal side length tends to zero as $M$ increases. Here the term involving $f_I$ is necessary.

For $k\geq2$ and $R=E_1\times\cdots\times E_k$, we organize the terms in $f-f_R$ according to
the variables that remain after averaging. For $\emptyset\ne A\subseteq\{1,\ldots,k\}$, write
$x_A=(x_j)_{j\in A}$, $X_A=\bigoplus_{j\in A}X_j$ and $R_A=\prod_{j\in A}E_j$.
If $A^c=\{1,\ldots,k\}\setminus A$ is nonempty, define
$$
f_A(x_A)=\frac1{|R_{A^c}|}\int_{R_{A^c}}f(x_A,u_{A^c})\,du_{A^c}.
$$
Thus $A$ indexes the variables left after averaging. For $R=I\times J$, we have
$f_{\{1\}}=f_J$ and $f_{\{2\}}=f_I$. For $A=\{1,\ldots,k\}$, we have $f_A=f$. Expanding the identity
$\prod_{j=1}^k[P_j+(I-P_j)]f=f$ on $R$, we obtain one term for each subset $A$.
The empty subset contributes $f_R$. The terms with $\ell$ coordinates left after averaging sum to
$$
F_R^{(\ell)}=\sum_{|A|=\ell}\prod_{j\in A}(I-P_j)f_A,
\qquad 1\leq\ell\leq k.
$$
Thus $f-f_R=\sum_{\ell=1}^k F_R^{(\ell)}$ on $R$.
For $k=2$, $F_R^{(1)}=(f_J-f_R)+(f_I-f_R)$ and $F_R^{(2)}=Q_Rf$.
For $A=\{i_1<\cdots<i_\ell\}$ with $\ell\geq2$, let $S_A$ denote the full-gradient
Poisson area \eqref{e:flag-area} on the flag
$X_{i_1}\oplus\cdots\oplus X_{i_\ell}\supset\cdots\supset X_{i_\ell}$.
For a single coordinate, we use $S_{\{i\}}=\mathcal S_P^{(d_i)}$ from
\eqref{e:intro-one-area}.
In Theorem~\ref{t:coordinate} one may instead use $S_{\mathrm{cl}}^{(d_i)}$ for every singleton,
by Theorem~\ref{t:concentration}\textup{(i)} with $k=1$. For each
$1\leq\ell\leq k$, write
$$
\mathcal E_{\ell,R}(f)(x)^2=\sum_{|A|=\ell}\mathbf1_{R_{A^c}}(x_{A^c})\,|(S_A f_A)(x_A)|^2,
$$
where $\mathbf1_{R_\emptyset}=1$. Thus $\mathcal E_{k,R}(f)=S_F^{(k)}(f)$, and for $k=2$,
$\mathcal E_{1,R}(f)=\mathcal E_{P,R}(f)$. Applying the local reverse estimate on the induced
flags and summing over $|A|=\ell$, we obtain
a moment bound of order $p^{\ell/2}$ for $F_R^{(\ell)}$ and the corresponding
Chang--Wilson--Wolff estimate with power $2/\ell$ (Subsection~\ref{s:coord-expansion}).

\begin{mainthm}\label{t:coordinate}
Let $k\geq2$ and $f\in L^u(\mathbb R^N)$ for some $1\leq u<\infty$.
Fix $2\leq p<\infty$. For every product $R=E_1\times\cdots\times E_k$
of axis-parallel cubes and each $1\leq\ell\leq k$, we have
$$
\bigg\|\sum_{|A|=\ell}\prod_{j\in A}(I-P_j)f_A\bigg\|_{L^p(R)}
\leq Cp^{\ell/2}
\bigg\|\bigg(\sum_{|A|=\ell}|S_Af_A|^2\bigg)^{1/2}\bigg\|_{L^p(R)}.
$$
Consequently,
$$
\|f-f_R\|_{L^p(R)}
\leq C\sum_{\ell=1}^kp^{\ell/2}
\bigg\|\bigg(\sum_{|A|=\ell}|S_A f_A|^2\bigg)^{1/2}\bigg\|_{L^p(R)}.
$$
If the estimate for $f-f_R$ holds for all such $f,R$ with finite coefficients $a_{p,\ell}\geq0$ in place of
$Cp^{\ell/2}$, then $a_{p,\ell}\geq c p^{\ell/2}$ for every $1\leq\ell\leq k$.

If $0<\|\mathcal E_{\ell,R}(f)\|_{L^\infty(R)}<\infty$, then
$$
\frac1{|R|}\int_R\exp\bigg[c\bigg(\frac{|F_R^{(\ell)}|}{\|\mathcal E_{\ell,R}(f)\|_{L^\infty(R)}}\bigg)^{2/\ell}\bigg]\,dx\leq2.
$$
The power $2/\ell$ is optimal. If a layer area vanishes almost everywhere on $R$,
so does the corresponding oscillation. The constants depend only on $k$ and the dimensions.
\end{mainthm}

For $k=2$, Proposition~\ref{p:two-orders} tests the term of order $p$ and the term of order
$\sqrt p$ on different families of functions. For general $k$, the examples in
Proposition~\ref{p:layer-sharp} have zero averages $f_A$ for $1\leq |A|<\ell$,
and the areas of all higher layers tend uniformly to zero on the testing rectangle.
The resulting distribution estimate for $f-f_R$ is stated in \eqref{e:coord-tails}.

\subsection{Bounded flag area and approximation in flag BMO}\label{s:intro-approx}

For the two-step flag, we next ask how closely a bounded function can approximate a function
with bounded flag area. In one parameter, Garnett and Jones~\cite{GJ} proved that the $\mathrm{BMO}$ distance
from $f$ to $L^\infty$ is comparable to the infimum of those $\varepsilon>0$ for which
$$
\sup_Q\frac1{|Q|}\int_Q\exp\big(|f-f_Q|/\varepsilon\big)\,dx<\infty.
$$
When the Poisson area is bounded, the exponential-square estimate makes this supremum finite for
every $\varepsilon>0$. Chang, Wilson and Wolff deduced that $f$ lies in the closure of $L^\infty$ in
$\mathrm{BMO}$ and noted the stronger inclusion in the real interpolation space
$(L^\infty,\mathrm{BMO})_{1/2,\infty}$~\cite[Corollary~3.3 and the discussion on p.~243]{CWW}.
In the flag setting, we measure the distance in $\mathrm{BMO}_F=(H_F^1)^*$, where $H_F^1$ is
the flag Hardy space of~\cite{HL,HLLW} (see also~\cite{HLW} for weighted flag Hardy spaces
and~\cite{ChangFefferman} for product Hardy spaces). We use the heat maximal function norm
on $H_F^1$ defined at the beginning of Section~\ref{s:bounded-area}.
For $0<\|S_F(f)\|_\infty<\infty$, Theorem~\ref{t:concentration}\textup{(i)} bounds the average of
$\exp(c|Q_Rf|/\|S_F(f)\|_\infty)$ on every rectangle $R$ for a fixed $c>0$.
Approximation in $\mathrm{BMO}_F$ requires a single bounded function on the whole space $\mathbb R^{n+m}$.
To construct such an approximant, we first select terms from the reconstruction in $x$ while
controlling the $y$-derivatives of each resulting sum in $x$, which permits a second selection
in $y$.
The selected sum is bounded, and pairing the discarded terms with $H_F^1$ controls the error
in $\mathrm{BMO}_F$.

\begin{mainthm}\label{t:approximation}
Suppose $f\in L^2(\mathbb R^{n+m})$ and $S_F(f)\in L^\infty(\mathbb R^{n+m})$.
For every $0<\varepsilon\leq\|S_F(f)\|_\infty$, there is a bounded function $b_\varepsilon$
with $\|f-b_\varepsilon\|_{\mathrm{BMO}_F}\leq\varepsilon$ and
$$
\|b_\varepsilon\|_\infty
\leq C\frac{\|S_F(f)\|_\infty^4}{\varepsilon^3}
\bigg[\log\bigg(e+\frac{\|S_F(f)\|_\infty}{\varepsilon}\bigg)\bigg]^{3/2}.
$$
If $\|S_F(f)\|_\infty=0$, then $f=0$ almost everywhere. The constant depends only on $n,m$.
\end{mainthm}

In the proof, the first selection increases the row bounds used in the second.
The two error estimates in Subsection~\ref{s:approximation-proof} determine the choice of parameters
and hence the upper bound above. For every sufficiently small $\varepsilon>0$, there is a real Schwartz
function $f$ with $\|S_F(f)\|_\infty=1$ such that every bounded $b$ with
$\|f-b\|_{\mathrm{BMO}_F}\leq\varepsilon$ satisfies $\|b\|_\infty\geq c/\varepsilon$
(Remark~\ref{r:approx-lower}).
For the classical Poisson area, the best uniform bound on $\|b_\varepsilon\|_\infty$ under the same
normalization has order $\varepsilon^{-1}$ as $\varepsilon\to0+$ (Remark~\ref{r:scalar-approx}).
For real $f\in L^2(\mathbb R^d)$ with $\|S_{\mathrm{cl}}^{(d)}(f)\|_\infty\leq1$, the classical estimate is
$$
[e^{\lambda f}]_{A_2}
=\sup_Q\bigg(\frac1{|Q|}\int_Q e^{\lambda f}\bigg)
       \bigg(\frac1{|Q|}\int_Q e^{-\lambda f}\bigg)
\leq C_de^{C_d\lambda^2},\qquad \lambda\geq1,
$$
where the supremum is over cubes in $\mathbb R^d$. In Remark~\ref{r:scalar-approx}
we use this estimate with $\lambda$ comparable to $\varepsilon^{-1}$.
In the two-step example of Proposition~\ref{p:coord-necessity}, the product of the two exponential averages
on one rectangle tends to infinity for every fixed $\lambda\ne0$ (see \eqref{e:two-weight-growth}).
The optimal dependence on $\varepsilon$ for flag $\mathrm{BMO}$ remains open.

Each bounded $b$ represents the functional $\langle\cdot,b\rangle$ on $H_F^1$,
where $\langle g,b\rangle=\int_{\mathbb R^{n+m}}g\,\overline b$. Define
$$
\mathcal X_F=
\overline{\{\langle\cdot,b\rangle:b\in L^\infty(\mathbb R^{n+m})\}}^{\, (H_F^1)^*}.
$$
By Theorem~\ref{t:approximation}, the class
$\mathscr S_F=\{f\in L^2(\mathbb R^{n+m}):S_F(f)\in L^\infty(\mathbb R^{n+m})\}$
lies in $\mathcal X_F\cap L^2(\mathbb R^{n+m})$, where
$$
\mathcal X_F\cap L^2(\mathbb R^{n+m})
=\big\{f\in L^2(\mathbb R^{n+m}):\inf_{b\in L^\infty}\|f-b\|_{\mathrm{BMO}_F}=0\big\}.
$$
Does a uniform exponential bound for $Q_Rf$ on rectangles characterize $\mathcal X_F\cap L^2(\mathbb R^{n+m})$? Let
$$
\mathscr E_F=\bigg\{f\in\mathrm{BMO}_F\cap L^2(\mathbb R^{n+m}):
\sup_{R=I\times J}\frac1{|R|}\int_Re^{c|Q_Rf|}\,dz<\infty
\text{ for some }c>0\bigg\},
$$
where the supremum runs over all products of axis-parallel cubes $I\subset\mathbb R^n$ and
$J\subset\mathbb R^m$.

\begin{mainthm}\label{t:classes}
Let $n,m\geq1$.
\begin{itemize}
\item[(i)] The classes $\mathcal X_F\cap L^2(\mathbb R^{n+m})$ and $\mathscr E_F$ are incomparable.
Both strictly contain $\mathscr S_F$, and
$$
\mathcal X_F\cap L^2(\mathbb R^{n+m})\subsetneq\mathrm{BMO}_F\cap L^2(\mathbb R^{n+m}).
$$
\item[(ii)] For every $f\in\mathcal X_F\cap L^2(\mathbb R^{n+m})$ and every $\lambda>0$,
$$
\sup_R\frac1{|R|}\int_R\exp\big(\lambda |Q_Rf|^{1/2}\big)\,dz<\infty.
$$
There exist a real function $F\in\mathcal X_F\cap L^2(\mathbb R^{n+m})$ and a fixed flag rectangle $R_0$
such that
$$
\int_{R_0}\exp\big(\lambda|Q_{R_0}F|^\gamma\big)\,dz=\infty
\qquad\text{for all }\lambda>0\text{ and }\gamma>1/2.
$$
\end{itemize}
\end{mainthm}

Chang and Fefferman proved exponential integrability with power $1/2$ for the sum of the
product $\mathrm{BMO}$ terms over rectangles contained in an open set~\cite[Section~III]{ChangFefferman}.
Proposition~\ref{p:rect-moments} bounds the normalized $L^p$ norm of $Q_Rf$ by
$Cp^2\|f\|_{\mathrm{BMO}_F}$ for $f\in\mathrm{BMO}_F\cap L^2(\mathbb R^{n+m})$.
One fixed function attains order $p^2$ on one flag rectangle (Remark~\ref{r:bmo-moment-sharp}).
For $f\in\mathcal X_F\cap L^2(\mathbb R^{n+m})$, approximation by bounded functions improves this bound to
$$
\lim_{p\to\infty}\frac1{p^2}\sup_R
\bigg(\frac1{|R|}\int_R |Q_Rf|^p\,dz\bigg)^{1/p}=0.
$$
Indeed, $Q_Rb$ contributes at most $4\|b\|_\infty$ for a bounded approximant $b$, and
$\|f-b\|_{\mathrm{BMO}_F}$ can be arbitrarily small.
Consequently, the exponential average with exponent $1/2$ is finite for every $\lambda>0$.
To show that this exponent is sharp, we apply two Riesz transforms along the flag to a
continuous function with compact support.
A different example, with a vertical logarithmic singularity, belongs to $\mathscr E_F$ while
remaining at positive $\mathrm{BMO}_F$ distance from $L^\infty$.
It also satisfies the exponential estimate with exponent $1/2$ for every $\lambda>0$, so that
estimate does not characterize $\mathcal X_F\cap L^2(\mathbb R^{n+m})$.

\subsection{Reconstruction and intermediate square functions}\label{s:intro-method}

We construct the smooth atoms directly from the Calder\'on reproducing formula. In one parameter,
Chen, Duong, Wu and Yan~\cite[(3.1)--(3.4) and Lemma~3.1]{CDWY} integrated this formula
over the region above each dyadic cube to obtain a localized term. For the two-step flag,
we reconstruct $f$ from the mixed time derivative of its Poisson extension using smooth kernels
as in~\cite{Wilson89}. Finite propagation ensures that the kernels have compact support~\cite[Lemma~5.1.1]{HLLW}.
After splitting the kernel according to its cancellation, we partition the centres and scales
into regions indexed by rectangles from adjacent dyadic systems in the two coordinate spaces
(Subsection~\ref{s:dyadic-supports}).
Integration over each region
produces smooth atoms supported strictly inside its rectangle $I\times J$ (Proposition~\ref{p:two-atoms}).

Multiplying both reconstruction scales by the same sufficiently large constant $\sigma$ makes the Poisson
averaging region contain $I\times J$. Fix a pair $\kappa$ of dyadic shifts and a finite selection
of rectangles from the region $t<s$. For each selected $I\times J$, define the nonnegative coefficient $c_{I,J}$ by
$$
c_{I,J}^2=\frac1{|I||J|}\int_{\mathcal T^{\mathrm{str}}_{\kappa,I,J}}
 |H_f^\sigma(u,v,t,s)|^2\,du\,dv\,\frac{dt}{t}\frac{ds}{s},
$$
where $\mathcal T^{\mathrm{str}}_{\kappa,I,J}$ is the parameter region defined in
Subsection~\ref{s:dyadic-supports}, and \eqref{e:poisson-coeff} expresses $H_f^\sigma$ in terms of the mixed time derivative of $U_f$. With $c_{I,J}=0$ on the other rectangles, we have
$$
\sum_{I,J}c_{I,J}^2\mathbf1_{I\times J}(x,y)\leq C S_F(f)(x,y)^2
$$
with a constant independent of the selection (Proposition~\ref{p:coeff-control}).
In fact, the mixed time derivative alone suffices for the local and global reverse estimates at
any fixed positive aperture (Remark~\ref{r:time-area}).

\noindent\textbf{Cancellation in the reconstruction kernels.}
For product areas the reconstruction kernel is a tensor product of two kernels of integral zero.
The two-step reconstruction kernel $K_{t,s}$ in \eqref{e:compact-kernel} is a partial convolution in $y$. It has integral zero in $y$ for each $x$, while
its integral in $x$, $b_{t,s}(y)=\int_{\mathbb R^n}K_{t,s}(x,y)\,dx$, need not vanish.
For $t<s$, write $\theta_t(x)=t^{-n}\theta(x/t)$ with $\theta$ a smooth bump of integral one
as in Subsection~\ref{s:scale-regions}. Then
$$
K_{t,s}(x,y)
=\big[K_{t,s}(x,y)-\theta_t(x)b_{t,s}(y)\big]+\theta_t(x)b_{t,s}(y).
$$
The first term has integral zero in each variable separately. In the second, the horizontal average need not vanish, but
$$
|\theta_t(x)b_{t,s}(y)|\leq Ct^{-n}s^{-m}(t/s)^2,\qquad t<s.
$$
The factor $(t/s)^2$ makes the sum over smaller horizontal scales converge at fixed $s$, leaving a one-parameter
estimate in $y$. This quantifies the weak cancellation of flag kernels in~\cite[Section~1.2 and Theorem~6.8]{NRSW1}.
In the other scale range, the kernels for a fixed $t$ fit in one cube of side comparable to $t$, and
$$
|K_{t,s}(x,y)|\leq Ct^{-(n+m)}(s/t)^2,\qquad s\leq t.
$$
Here $(s/t)^2$ makes the sum over the smaller vertical scales converge.

\noindent\textbf{Iteration in the two variables.}
We estimate sums of these smooth atoms using the ordinary dyadic martingale differences in
each coordinate block. For a finite sum $G$ of atoms with cancellation in both variables,
write $d_{U,V}=\Delta_U^x\Delta_V^yG$, where $\Delta_U^x$ and $\Delta_V^y$ are the
martingale differences on the dyadic cubes $U$ and $V$ (Subsection~\ref{s:both-cancellations}).
At fixed relative generations, $U$ and $V$ each have
a unique dyadic ancestor of the required length. The derivative bounds on the atoms then yield
the pointwise square-function estimate \eqref{e:product-square} by summation over generations,
without using Journ\'e's covering lemma for rectangles.

In the intermediate square function of Li, Pipher and Wu~\cite{LiPipherWu}, the sum in $V$
stays inside the $\ell^2$ norm in $U$. We can thus apply \eqref{e:scalar-square} in $x$ and,
for fixed $x$, apply \eqref{e:hilbert-square} in $y$ to the $\ell^2$-valued differences
$(d_{U,V})_U$. For finite martingale sums, by \eqref{e:product-square} and the coefficient
estimate above, we obtain
\begin{align*}
\bigg\|\sum_{U,V}d_{U,V}\bigg\|_p
\leq &C\sqrt p\,
   \bigg\|\bigg(\sum_U\Big|\sum_Vd_{U,V}\Big|^2\bigg)^{1/2}\bigg\|_p
\leq Cp\bigg\|\bigg(\sum_{U,V}|d_{U,V}|^2\bigg)^{1/2}\bigg\|_p\\
\leq &Cp\bigg\|\bigg(\sum_{I,J}c_{I,J}^2\mathbf1_{I\times J}\bigg)^{1/2}\bigg\|_p
 \leq Cp\|S_F(f)\|_p.
\end{align*}
For the one-parameter estimates see~\cite{IT,Pipher86}.

For the local estimate, an atom that meets the testing rectangle $R$ can extend beyond its boundary.
Consider one coordinate block, with testing cube $E$. For atoms whose supports meet $R$,
the supporting cubes $I$ with $\ell(I)\leq\ell(E)$ lie in the concentric triple $3E$.
The margin between each atom and the boundary of its cube permits a single map from
$3E$ to $E$ that sends every such $I$ into $I\cap E$. Changing variables on each branch,
we bound the $L^p$ norm of the coefficient square function on $3E$ by its norm on $E$,
uniformly in $p$ (Lemma~\ref{l:cube-fold}).
If $\ell(I)>\ell(E)$, the derivative bound and subtraction of the mean over $E$ produce
a factor $\ell(E)/\ell(I)$, summable over the larger dyadic scales. Applying these arguments in each coordinate block keeps the coefficient
norm on $R$ for every ordering of its side lengths.

For a $k$-step flag, each term in the reconstruction has cancellation in a collection of coordinate blocks.
We first sum over the side lengths in the other blocks using decay in the scale ratios.
For atoms with cancellation in the same $j$ blocks, the local Hoeffding-type estimate
bounds the $L^p(R)$ norm of their centred sum by $Cp^{j/2}$ times the $L^p(R)$ norm of
their coefficient square function (Proposition~\ref{p:flag-local-sums}).
By \eqref{e:flag-coeff}, the flag area controls this square function pointwise on $R$.
Summing over the choices of blocks, we obtain the factor $p^{k/2}$ in Theorem~\ref{t:concentration}\,\textup{(i)}.
For $k=3$, the table following Proposition~\ref{p:flag-kernel-bounds} displays all four terms and their scale factors.

The intermediate square function also enters the reverse estimates for Zygmund
dilations~\cite{Zygmund}. There the sharp power of $p$ can change with the aperture even for a fixed kernel pair.
When the radial dual kernel has sufficiently small support for the given aperture,
the ratio of the last two scales in Table~\ref{tab:two-scales} permits a single Haar estimate
in those two variables. At fixed relative generations, the first dyadic scale and one Haar interval
determine the length of the intervals in the other variable.
Distinct intervals of this length are disjoint, and coefficients from different first-scale
intervals have disjoint supports in the first scale parameter. For a flag, it is the cancellation of
each reconstruction term that determines the coordinates in which to iterate the
one-parameter estimate. Decay in the scale ratios accounts for the other coordinates.

\subsection{Weights, commutators, maximal functions and the heat area}\label{s:intro-consequences}

Corollary~\ref{c:weighted-local} estimates the moments of $f-f_R$ with respect to a weight on $R$.
For exponential weights, suppose $\|S_F^{(k)}(f)\|_\infty\leq1$ and $f$ has mean zero in each coordinate block of $R$.
Then $Q_Rf=f$ on $R$. For $k=2$, Theorem~\ref{t:concentration}\textup{(i)}
bounds the product of the averages of $e^{\lambda f}$ and $e^{-\lambda f}$ for real $f$ and small $|\lambda|$.
This product can grow without bound for large fixed $|\lambda|$ when $k=2$, and for every fixed
$\lambda\ne0$ when $k\geq3$ (Corollary~\ref{c:flag-linear-exp}).

For a two-step flag, a bound on the areas of the coordinate averages allows every real $\lambda$.
Let $\mathcal D$ be the standard dyadic system in $\mathbb R^{n+m}$. For each $Q\in\mathcal D$,
the coordinate area $\mathcal E_{P,Q}$ uses the two block projections of $Q$, which have the same side length.
If $f\in L^2(\mathbb R^{n+m})$ is real and $\sup_{Q\in\mathcal D}\|\mathcal E_{P,Q}(f)\|_{L^\infty(Q)}\leq1$,
Proposition~\ref{p:coordinate-weights} bounds the same product by $C e^{C\lambda^2}$ on every product of cubes.
Even on one fixed flag rectangle, the supremum of this product over such $f$ is at least $e^{c\lambda^2}$ for large $|\lambda|$.
We use this weight estimate for the flag Riesz commutators in Subsection~\ref{s:commutators},
following the conjugation argument of Duong, Li, Ou, Pipher and Wick~\cite{DLOPW} for little flag $\mathrm{BMO}$.

The global reverse inequality also bounds the Poisson maximal functions by the area function,
with the same sharp factor $p^{k/2}$ (Corollary~\ref{c:poisson-maximal}). In the two-step
reconstruction we may keep the compact kernels and change the spectral multipliers in the coefficients
(Proposition~\ref{p:spectral-multipliers}). For the heat area this yields the same-rectangle reverse estimate
for $f\in L^u$, $1\leq u<\infty$. On the whole space, its $L^p$ norm is equivalent to that of $f$
for $2\leq p<\infty$ (Theorem~\ref{t:heat-area}).

This paper is organized as follows. Sections~\ref{s:reconstruction} and~\ref{s:local-hoeffding}
establish the two-step estimates. We use the same decomposition for bounded approximation and for the heat area. Section~\ref{s:finite-flags} proves Theorem~\ref{t:concentration} for finite
flags, and Section~\ref{s:coord-oscillation} proves Theorem~\ref{t:coordinate} for the
coordinate averages. Section~\ref{s:bounded-area} proves Theorems~\ref{t:approximation}
and~\ref{t:classes}. Section~\ref{s:consequences} gives weighted moment estimates and bounds for
Riesz commutators and Poisson maximal functions. In Section~\ref{s:multipliers-heat} we change the multipliers in the reconstruction
and prove the two-step heat estimates.

\section{Compact reconstruction for the two-step flag}\label{s:reconstruction}


We reconstruct $f$ from the mixed time derivative of its Poisson extension. At reconstruction scales
$r_1,r_2$, the kernel has integral zero in $y$. Its support lies in a rectangle with side lengths
comparable to $r_1$ in $x$ and $\max\{r_1,r_2\}$ in $y$.
For $r_1<r_2$, we separate the term with nonzero integral in $x$ and bound it with a factor $(r_1/r_2)^2$.
For $r_2\leq r_1$, the whole kernel gains $(r_2/r_1)^2$. We bound the square sum of the dyadic coefficients
pointwise by $C S_F(f)^2$.

Write $d=n+m$ and recall that
$$
\nabla_1=(\partial_t,\nabla_x,\nabla_y),\qquad
\nabla_2=(\partial_s,\nabla_y).
$$
In the formula for $S_F(f)$ in Subsection~\ref{s:intro-background}, write $sa_2=(0,z)$ and
$(u,v)=(x,y)-ta_1-(0,z)$. Then $da_1\,da_2=t^{-(n+m)}s^{-m}\,du\,dv\,dz$.
For a displacement $(a,b)\in\mathbb R^n\times\mathbb R^m$, integration in $z\in\mathbb R^m$ leaves the weight
$$
\chi_{t,s}(a,b)=\int_{|z|<s}\mathbf1_{\{|a|^2+|b-z|^2<t^2\}}\,dz.
$$
Hence
\begin{equation}\label{e:poisson-area-weight}
S_F(f)(x,y)^2=\int_0^\infty\int_0^\infty\int_{\mathbb R^{n+m}}
\chi_{t,s}(x-u,y-v)|t\nabla_1s\nabla_2U_f(u,v,t,s)|^2
\frac{du\,dv}{t^{n+m}s^m}\frac{dt}{t}\frac{ds}{s}.
\end{equation}

\subsection{A compact kernel dual to the Poisson derivative}\label{s:compact-kernels}
Let $b_\sigma(r)=re^{-\sigma r}$ for $r\geq0$ and $\sigma>0$. For $f\in L^2(\mathbb R^{n+m})$, define its Poisson coefficients by
$$
H_f^\sigma(u,v,r_1,r_2)
=\Big [b_\sigma(r_1\sqrt{-\Delta_{x,y}})
       b_\sigma(r_2\sqrt{-\Delta_y})f\Big ]\, (u,v).
$$
At the Poisson scales $t=\sigma r_1$ and $s=\sigma r_2$, we have
\begin{equation}\label{e:poisson-coeff}
H_f^\sigma(u,v,r_1,r_2)
= r_1r_2\sqrt{-\Delta_{x,y}}\sqrt{-\Delta_y}
e^{-t\sqrt{-\Delta_{x,y}}}e^{-s\sqrt{-\Delta_y}}f(u,v)
= \sigma^{-2}t\partial_t s\partial_sU_f(u,v,t,s).
\end{equation}
We choose $\sigma$ in Subsection~\ref{s:dyadic-supports} so that the dyadic covering rectangles lie in the Poisson cones.

We write $L^2_{\mathrm{par}}$ for the space of jointly measurable, square-integrable functions of the spatial variables
$(u,v)\in\mathbb R^{n+m}$ and the scale parameters $(r_1,r_2)\in(0,\infty)^2$, with norm
$$
\|H\|_{L^2_{\mathrm{par}}}^2
=\int_0^\infty\int_0^\infty
\|H(\cdot,r_1,r_2)\|_{L^2(\mathbb R^{n+m})}^2\,
\frac{dr_1}{r_1}\frac{dr_2}{r_2}.
$$
For the Fourier transform, we use the normalization
$$
\widehat h(\xi)=\int_{\mathbb R^D}h(x)e^{-ix\cdot\xi}\,dx,
\qquad h\in\mathcal S(\mathbb R^D).
$$
With this choice, $b_\sigma(r\sqrt{-\Delta})$ has symbol
$b_\sigma(r|\xi|)$, and Plancherel's identity is
$$
\|h\|_{L^2(\mathbb R^D)}^2=(2\pi)^{-D}\|\widehat h\|_{L^2(\mathbb R^D)}^2.
$$

Following the finite-propagation construction in \cite[Lemma~5.1.1]{HLLW},
choose a real even $\theta_0\in C_c^\infty((-1,1))$ of integral one,
independently of the dimension and of $\sigma$. Write $a_0(r)=r^2\widehat\theta_0(r)$ and define
$$
c_\sigma=\int_0^\infty a_0(r)b_\sigma(r)\,\frac{dr}{r}.
$$
Since $\widehat\theta_0(0)=1$, the substitution $u=\sigma r$ and dominated convergence imply
$$
\sigma^3 c_\sigma
=\int_0^\infty u^2\widehat\theta_0(u/\sigma)e^{-u}\,du
\to \int_0^\infty u^2e^{-u}\,du=2
\qquad(\sigma\to\infty).
$$
Thus $c_\sigma>0$ for all $\sigma\geq\sigma_0$, where $\sigma_0>0$ depends only on $\theta_0$.
For such $\sigma$, let $a_\sigma=a_0/c_\sigma$. Then
\begin{equation}\label{e:dual-normalization}
\int_0^\infty a_\sigma(r)b_\sigma(r)\,\frac{dr}{r}=1.
\end{equation}
For $0<r_-<r_+<\infty$ and $H\in L^2_{\mathrm{par}}$ supported in $r_-\leq r_1,r_2\leq r_+$, define
$$
\mathcal R_\sigma H=\int_{r_-}^{r_+}\int_{r_-}^{r_+}
a_\sigma(r_1\sqrt{-\Delta_{x,y}})  \,  a_\sigma(r_2\sqrt{-\Delta_y})H(\cdot,r_1,r_2)
\,\frac{dr_1}{r_1}\frac{dr_2}{r_2}.
$$

\begin{proposition}\label{p:reconstruction}
Let $\sigma\geq\sigma_0$.
\begin{enumerate}[(i)]
\item For every $D\geq1$, the operator $a_\sigma(\sqrt{-\Delta})$
on $\mathbb R^D$ has a real, radial convolution kernel
$\varphi_D^\sigma\in C_c^\infty(\mathbb R^D)$ supported in
the closed Euclidean unit ball $\overline{B(0,1)}$. Write $\mathbb N_0=\{0,1,2,\ldots\}$.
For each multi-index $\gamma\in\mathbb N_0^D$ of order $|\gamma|=\sum_{j=1}^D\gamma_j\leq1$,
$$
\int_{\mathbb R^D} z^\gamma \varphi_D^\sigma(z)\,dz=0.
$$
For $r>0$, the kernel of $a_\sigma(r\sqrt{-\Delta})$ is
$r^{-D}\varphi_D^\sigma(\cdot/r)$.

\item The reconstruction operator $\mathcal R_\sigma$ extends to a bounded
operator from $L^2_{\mathrm{par}}$ to $L^2(\mathbb R^{n+m})$.
For every $H\in L^2_{\mathrm{par}}$ and every measurable function
$\mu=\mu(u,v,r_1,r_2)$ with $|\mu|\leq1$ almost everywhere, we have
\begin{equation}\label{e:recon-l2}
\|\mathcal R_\sigma(\mu H)\|_{L^2(\mathbb R^{n+m})}
\leq \bigg(\int_0^\infty |a_\sigma(r)|^2\,\frac{dr}{r}\bigg) \,
\|H\|_{L^2_{\mathrm{par}}}.
\end{equation}

\item For every $f\in L^2(\mathbb R^{n+m})$, the Poisson coefficients
$H_f^\sigma$ satisfy
$$
\|H_f^\sigma\|_{L^2_{\mathrm{par}}}^2
=\frac{1}{16\sigma^4}\|f\|_{L^2(\mathbb R^{n+m})}^2.
$$
The reconstruction formula is
\begin{equation}\label{e:reproduce}
\mathcal R_\sigma H_f^\sigma=f
\qquad\text{in }L^2(\mathbb R^{n+m}).
\end{equation}
\end{enumerate}
\end{proposition}

\begin{proof}
\noindent{\it Step 1.}
For $D\geq1$, define $\zeta_D\in\mathcal S(\mathbb R^D)$ by
$\widehat\zeta_D(\xi)=\widehat\theta_0(|\xi|)$.
By the evenness of $\theta_0$ and Fourier inversion,
$$
\zeta_D
=\int_{-1}^1\theta_0(\tau)\cos(\tau\sqrt{-\Delta})\delta_0\,d\tau
\qquad\text{in }\mathcal S'(\mathbb R^D),
$$
where $\delta_0$ is the Dirac mass at the origin.
By finite propagation speed, the wave kernel
$\cos(\tau\sqrt{-\Delta})\delta_0$ is supported in
$\overline{B(0,|\tau|)}$. Since $|\tau|\leq1$ in the integral above,
$\operatorname{supp}\zeta_D\subseteq\overline{B(0,1)}$.

Let $\varphi_D^\sigma=-c_\sigma^{-1}\Delta\zeta_D$.
Then $\varphi_D^\sigma$ is real, radial and smooth, with
$\operatorname{supp}\varphi_D^\sigma\subseteq\overline{B(0,1)}$.
For every $r>0$,
$$
\widehat{r^{-D}\varphi_D^\sigma(\cdot/r)}(\xi)
=\widehat{\varphi_D^\sigma}(r\xi)
=c_\sigma^{-1}r^2|\xi|^2\widehat\theta_0(r|\xi|)
=a_\sigma(r|\xi|).
$$
For every multi-index $\gamma$ with $|\gamma|\leq1$, we have
$\Delta(z^\gamma)=0$. Integrating by parts, we obtain
$$
\int_{\mathbb R^D}z^\gamma\varphi_D^\sigma(z)\,dz
=-c_\sigma^{-1}\int_{\mathbb R^D}\zeta_D(z)\Delta(z^\gamma)\,dz=0.
$$

\smallskip
\noindent{\it Step 2.}
We begin with $H\in L^2_{\mathrm{par}}$ supported where
$r_-\leq r_1,r_2\leq r_+$, with $0<r_-<r_+<\infty$. In the integral defining
$\mathcal R_\sigma H$, write
$T_{r_1,r_2}=a_\sigma(r_1\sqrt{-\Delta_{x,y}})
a_\sigma(r_2\sqrt{-\Delta_y})$.
For each fixed $g\in L^2(\mathbb R^{n+m})$, the map $(r_1,r_2)\mapsto T_{r_1,r_2}g$
is continuous in the norm of this space. Since $H$ is jointly measurable, it follows that
the integrand $(r_1,r_2)\mapsto T_{r_1,r_2}H(\cdot,r_1,r_2)$
is strongly measurable. By Plancherel in the spatial variables and
Cauchy--Schwarz in the scale variables, we have
\begin{align*}
\int_{[r_-,r_+]^2}\|T_{r_1,r_2}H(\cdot,r_1,r_2)\|_{L^2(\mathbb R^{n+m})}\,
\frac{dr_1}{r_1}\frac{dr_2}{r_2}
\leq &\|a_\sigma\|_\infty^2
\int_{[r_-,r_+]^2}\|H(\cdot,r_1,r_2)\|_{L^2(\mathbb R^{n+m})}\,
\frac{dr_1}{r_1}\frac{dr_2}{r_2}\\
\leq &\|a_\sigma\|_\infty^2\log(r_+/r_-)\|H\|_{L^2_{\mathrm{par}}}.
\end{align*}
Thus the integral defining $\mathcal R_\sigma H$ converges in $L^2(\mathbb R^{n+m})$.

To extend $\mathcal R_\sigma$ to all of $L^2_{\mathrm{par}}$, we need
a bound independent of $r_-,r_+$. For $g\in L^2(\mathbb R^{n+m})$, by Plancherel's identity and Tonelli's theorem,
\begin{align*}
 \int_0^\infty\int_0^\infty
\|T_{r_1,r_2}g\|_{L^2(\mathbb R^{n+m})}^2\,\frac{dr_1}{r_1}\frac{dr_2}{r_2} 
= &(2\pi)^{-(n+m)}\int_{\mathbb R^{n+m}}|\widehat g(\xi,\eta)|^2
\bigg(\int_0^\infty |a_\sigma(r_1|(\xi,\eta)|)|^2\,\frac{dr_1}{r_1}\bigg)\\
&\qquad \qquad \qquad \times 
\bigg(\int_0^\infty |a_\sigma(r_2|\eta|)|^2\,\frac{dr_2}{r_2}\bigg)
\,d\xi\,d\eta\\
= &\|a_\sigma\|_{L^2(du/u)}^4\|g\|_{L^2(\mathbb R^{n+m})}^2,
\end{align*}
where, outside the null set $\mathbb R^n\times\{0\}$, the substitutions
$u=r_1|(\xi,\eta)|$ and $u=r_2|\eta|$ in the respective inner integrals
show that both equal
$$
\|a_\sigma\|_{L^2(du/u)}^2
=\int_0^\infty |a_\sigma(u)|^2\,\frac{du}{u}
=c_\sigma^{-2}\int_0^\infty u^3|\widehat\theta_0(u)|^2\,du<\infty.
$$

Since $T_{r_1,r_2}$ is self-adjoint, we obtain by Cauchy--Schwarz that
\begin{align*}
|\langle\mathcal R_\sigma H,g\rangle|
= &\bigg|\int_{[r_-,r_+]^2}
\langle H(\cdot,r_1,r_2),T_{r_1,r_2}g\rangle\,
\frac{dr_1}{r_1}\frac{dr_2}{r_2}\bigg|\\
\leq &\|H\|_{L^2_{\mathrm{par}}}
\bigg(\int_0^\infty\int_0^\infty
\|T_{r_1,r_2}g\|_{L^2(\mathbb R^{n+m})}^2\,\frac{dr_1}{r_1}\frac{dr_2}{r_2}\bigg)^{1/2}
=\|a_\sigma\|_{L^2(du/u)}^2\|H\|_{L^2_{\mathrm{par}}}\|g\|_{L^2(\mathbb R^{n+m})}.
\end{align*}
By duality, $\|\mathcal R_\sigma H\|_{L^2(\mathbb R^{n+m})}\leq \|a_\sigma\|_{L^2(du/u)}^2\|H\|_{L^2_{\mathrm{par}}}$
for every such $H$, with the same constant for all $r_-,r_+$.

Now, for arbitrary $H\in L^2_{\mathrm{par}}$, let $H_N=m_NH$, where
$m_N=\mathbf1_{[N^{-1},N]^2}$ for $N\geq1$.
By dominated convergence, $H_N\to H$ in $L^2_{\mathrm{par}}$.
For $N\geq M$, the difference $H_N-H_M$ is supported where
$N^{-1}\leq r_1,r_2\leq N$. By the uniform bound just proved,
$$
\|\mathcal R_\sigma H_N-\mathcal R_\sigma H_M\|_{L^2(\mathbb R^{n+m})}
\leq \|a_\sigma\|_{L^2(du/u)}^2\|H_N-H_M\|_{L^2_{\mathrm{par}}}
\to 0
\qquad(M,N\to\infty).
$$
We therefore define $\mathcal R_\sigma H$ as the limit in $L^2(\mathbb R^{n+m})$ of
$\mathcal R_\sigma H_N$. Letting $N\to\infty$ in the estimate for $H_N$,
we obtain
$$
\|\mathcal R_\sigma H\|_{L^2(\mathbb R^{n+m})}\leq \|a_\sigma\|_{L^2(du/u)}^2\|H\|_{L^2_{\mathrm{par}}}.
$$
Applying this estimate to $\mu H$ and using $\|\mu H\|_{L^2_{\mathrm{par}}}\leq\|H\|_{L^2_{\mathrm{par}}}$, we obtain \eqref{e:recon-l2}.

\smallskip
\noindent{\it Step 3.}
Repeating the Plancherel calculation in Step~2 with $b_\sigma$ in place
of $a_\sigma$, we obtain
$$
\|H_f^\sigma\|_{L^2_{\mathrm{par}}}^2
=\bigg(\int_0^\infty |b_\sigma(r)|^2\,\frac{dr}{r}\bigg)^2\|f\|_{L^2(\mathbb R^{n+m})}^2
=\bigg(\int_0^\infty r e^{-2\sigma r}\,dr\bigg)^2\|f\|_{L^2(\mathbb R^{n+m})}^2
=\frac{1}{16\sigma^4}\|f\|_{L^2(\mathbb R^{n+m})}^2.
$$
Since $H_f^\sigma\in L^2_{\mathrm{par}}$, Step~2 implies that
$\mathcal R_\sigma(m_NH_f^\sigma)\to\mathcal R_\sigma H_f^\sigma$ in $L^2(\mathbb R^{n+m})$.

For almost every $(\xi,\eta)$ with $\eta\ne0$, we have 
$$
\widehat{\mathcal R_\sigma(m_NH_f^\sigma)}(\xi,\eta)
=q_N(|(\xi,\eta)|)q_N(|\eta|)\widehat f(\xi,\eta),
$$
where
$$
q_N(\rho)
=\int_{N^{-1}}^N a_\sigma(r\rho)b_\sigma(r\rho)\,\frac{dr}{r}
=\int_{\rho/N}^{N\rho}a_\sigma(u)b_\sigma(u)\,\frac{du}{u},\qquad \rho>0.
$$
By Cauchy--Schwarz again,
\begin{equation}\label{e:truncated-dual-bound}
|q_N(\rho)|
\leq\int_0^\infty |a_\sigma(u)b_\sigma(u)|\,\frac{du}{u}
\leq \|a_\sigma\|_{L^2(du/u)}
\bigg(\int_0^\infty u e^{-2\sigma u}\,du\bigg)^{1/2}
=\frac{\|a_\sigma\|_{L^2(du/u)}}{2\sigma}.
\end{equation}
By \eqref{e:dual-normalization}, $q_N(\rho)\to1$ for every $\rho>0$.
Using the uniform bound \eqref{e:truncated-dual-bound} and dominated
convergence in Plancherel's formula, we obtain
\begin{align*}
\|\mathcal R_\sigma H_f^\sigma-f\|_{L^2(\mathbb R^{n+m})}^2
= &\lim_{N\to\infty}\|\mathcal R_\sigma(m_NH_f^\sigma)-f\|_{L^2(\mathbb R^{n+m})}^2\\
= &(2\pi)^{-(n+m)}\lim_{N\to\infty}\int_{\mathbb R^{n+m}}
|q_N(|(\xi,\eta)|)q_N(|\eta|)-1|^2
|\widehat f(\xi,\eta)|^2\,d\xi\,d\eta 
=  0.
\end{align*}
That is,  \eqref{e:reproduce} holds. We complete the proof of Proposition \ref{p:reconstruction}.
\end{proof}

Write $\Phi_1^\sigma=\varphi_{n+m}^\sigma$ and $\Phi_2^\sigma=\varphi_m^\sigma$.
Their dilates are
$$
\Phi_{1,r}^\sigma(x,y)=r^{-(n+m)}\Phi_1^\sigma(x/r,y/r),
\qquad
\Phi_{2,r}^\sigma(y)=r^{-m}\Phi_2^\sigma(y/r).
$$
For $r_1,r_2>0$, the reconstruction kernel is
\begin{equation}\label{e:compact-kernel}
K_{r_1,r_2}^\sigma(x,y)
=\int_{\mathbb R^m}
\Phi_{1,r_1}^\sigma(x,y-z)\Phi_{2,r_2}^\sigma(z)\,dz.
\end{equation}
Note that the integrand vanishes unless $|(x,y-z)|\leq r_1$ and $|z|\leq r_2$, thus 
$$
\operatorname{supp} K_{r_1,r_2}^\sigma
\subseteq\big\{(x,y):\, |x|\leq r_1,\ |y|\leq r_1+r_2\big\}.
$$
The corresponding scales are $r_1$ in $x$ and $r_1+r_2$ in $y$.
In the area function below, we use the comparable $y$-scale $\max\{r_1,r_2\}$. Moreover, 
\eqref{e:reproduce} with this kernel becomes
$$
f=\lim_{N\to\infty}\int_{1/N}^N\int_{1/N}^N
K_{r_1,r_2}^\sigma*H_f^\sigma(\cdot,r_1,r_2)\,\frac{dr_1}{r_1}\frac{dr_2}{r_2}
\quad\text{in }  \     L^2(\mathbb R^{n+m}).
$$

\subsection{Area estimates for the Poisson coefficients}\label{s:coefficient-area}

For a jointly measurable function $H$ and $\lambda>0$, define
\begin{equation}\label{e:aux-area}
\mathcal A_{\lambda,F}(H)(x,y)^2=\int_0^\infty\int_0^\infty
\int_{\substack{|x-u|<\lambda r_1\\|y-v|<\lambda\max\{r_1,r_2\}}}
\frac{|H(u,v,r_1,r_2)|^2}{r_1^n\max\{r_1,r_2\}^m}\,du\,dv\,\frac{dr_1}{r_1}\frac{dr_2}{r_2}.
\end{equation}
Recall that the weight in the Poisson area formula \eqref{e:poisson-area-weight} is
$$
\chi_{t,s}(a,b)=\big|\{z\in\mathbb R^m:|z|<s,\ |a|^2+|b-z|^2<t^2\}\big|.
$$
For $t,s>0$ and every $(a,b)\in\mathbb R^n\times\mathbb R^m$, the two weights satisfy
\begin{equation}\label{e:cone-lower}
\frac{\mathbf1_{B(0,t/8)\times B(0,\max\{t,s\}/8)}(a,b)}{t^n\max\{t,s\}^m}
\leq C_m\frac{\chi_{t,s}(a,b)}{t^{n+m}s^m}.
\end{equation}
For $n=m=1$, the indicator in \eqref{e:cone-lower} equals one precisely inside the dashed rectangle
in each panel of Figure~\ref{fig:flag-cones}.

The left-hand side vanishes unless $|a|<t/8$ and $|b|<\max\{t,s\}/8$.
For such $(a,b)$, the inequality $|(a,b-z)|\leq |a|+|b-z|$ implies
$$
\chi_{t,s}(a,b)\geq |B(0,s)\cap B(b,t-|a|)|.
$$
If $s\leq t$, this intersection contains $B(0,s/2)$, since $|b-z|<t/8+s/2<t-|a|$ on that ball.
If $t<s$, it contains $B(b,t/2)$, since $t/2<t-|a|$ and $|z|<s/8+t/2<s$ there.
In either case,
$$
\frac{\chi_{t,s}(a,b)}{t^{n+m}s^m}
\geq c_m\frac{\min\{t,s\}^m}{t^{n+m}s^m}
=\frac{c_m}{t^n\max\{t,s\}^m},
$$
as required in \eqref{e:cone-lower}.

\begin{proposition}\label{p:coeff-area}
Let $\sigma\geq\sigma_0$ and $0<\lambda\leq\sigma/8$.
For every $f\in L^2(\mathbb R^{n+m})$,
\begin{equation}\label{e:coeff-area}
\mathcal A_{\lambda,F}(H_f^\sigma)(x,y)
\leq C\sigma^{(n+m-4)/2}S_F(f)(x,y),
\end{equation}
where $C$ depends only on $n,m$.
\end{proposition}

\begin{proof}
In \eqref{e:aux-area}, change variables to $t=\sigma r_1$ and $s=\sigma r_2$.
Since $\lambda/\sigma\leq1/8$, we may use \eqref{e:cone-lower} on the resulting
region of integration. By \eqref{e:poisson-coeff} and
$|\partial_t\partial_sU_f|\leq|\nabla_1\nabla_2U_f|$, we obtain
\begin{align*}
\mathcal A_{\lambda,F}(H_f^\sigma)(x,y)^2
= &\sigma^{n+m-4}\int_0^\infty\int_0^\infty
\int_{\substack{|x-u|<(\lambda/\sigma)t\\
                 |y-v|<(\lambda/\sigma)\max\{t,s\}}}
|t\partial_t s\partial_sU_f(u,v,t,s)|^2\,
\frac{du\,dv}{t^n\max\{t,s\}^m}\frac{dt}{t}\frac{ds}{s}\\
\leq &C\sigma^{n+m-4}\int_0^\infty\int_0^\infty\int_{\mathbb R^{n+m}}
\chi_{t,s}(x-u,y-v)|t\nabla_1s\nabla_2U_f(u,v,t,s)|^2\,
\frac{du\,dv}{t^{n+m}s^m}\frac{dt}{t}\frac{ds}{s}\\
= &C\sigma^{n+m-4}S_F(f)(x,y)^2,
\end{align*}
where the last equality follows from \eqref{e:poisson-area-weight}.
Hence we obtain \eqref{e:coeff-area}.
\end{proof}

\subsection{Kernel bounds in the two scale regions}\label{s:scale-regions}

For the rest of this section, let $(t,s)$ and $(\sigma t,\sigma s)$ denote
the reconstruction scales and the corresponding Poisson scales, respectively.
Fix $\sigma\geq\sigma_0$. We omit the superscript $\sigma$ on $\Phi_i^\sigma$, its dilates and $K_{t,s}^\sigma$.
By Proposition~\ref{p:reconstruction}\,\textup{(i)},
$$
\widehat\Phi_1(0,\eta)=a_\sigma(|\eta|)=\widehat\Phi_2(\eta),
$$
so Fourier inversion in $y$ yields $\int_{\mathbb R^n}\Phi_1(x,y)\,dx=\Phi_2(y)$.
Integrating \eqref{e:compact-kernel} in $x$, we obtain
\begin{equation}\label{e:kernel-integral}
b_{t,s}(y):=\int_{\mathbb R^n}K_{t,s}(x,y)\,dx
=(\Phi_{2,t}*\Phi_{2,s})(y).
\end{equation}
For fixed $t,s>0$, the expansion $a_\sigma(r)=c_\sigma^{-1}r^2+O(r^4)$ as $r\to0+$
shows that, for every sufficiently small $|\eta|>0$,
$$
\widehat b_{t,s}(\eta)=a_\sigma(t|\eta|)a_\sigma(s|\eta|)\ne0.
$$

For $t<s$, we subtract a kernel with the same $x$-integral as $K_{t,s}$.
Let $\theta(x)=\prod_{\ell=1}^n\theta_0(x_\ell)$ and
$\theta_t(x)=t^{-n}\theta(x/t)$. Then $\theta\in C_c^\infty((-1,1)^n)$
and $\int_{\mathbb R^n}\theta_t(x)\,dx=1$.
With $K_{t,s}^{\mathrm{rem}}(x,y)=\theta_t(x)b_{t,s}(y)$, we write
\begin{equation}\label{e:kernel-split}
K_{t,s}=K_{t,s}^{\mathrm{bi}}+K_{t,s}^{\mathrm{rem}}.
\end{equation}
By \eqref{e:kernel-integral}, the two kernels $K_{t,s}$ and
$K_{t,s}^{\mathrm{rem}}$ have $x$-integral $b_{t,s}(y)$, so
$\int_{\mathbb R^n}K_{t,s}^{\mathrm{bi}}(x,y)\,dx=0$ for every $y$.

For $s\leq t$, we estimate $K_{t,s}$ directly. By Proposition~\ref{p:reconstruction}\,\textup{(i)},
the mean and first moments of $\Phi_{2,s}$ vanish, so only the
second-order Taylor remainder contributes, with a factor $(s/t)^2$.

For the support bounds in each coordinate block we use the maximum norm
$$
|z|_\infty=\max_{1\leq j\leq q}|z_j|,\qquad z=(z_1,\ldots,z_q)\in\mathbb R^q.
$$

\begin{proposition}\label{p:kernel-bounds}
Let $t,s>0$ and let $\alpha\in\mathbb N_0^n$ and $\beta\in\mathbb N_0^m$ be multi-indices.
\begin{enumerate}[(i)]
\item Suppose $t<s$. Both kernels in \eqref{e:kernel-split} are supported
where $|x|_\infty\leq t$ and $|y|_\infty\leq t+s$, and both have
integral zero in $y$ for each fixed $x$.
The kernel $K_{t,s}^{\mathrm{bi}}$ also has integral zero in $x$ for
each fixed $y$ and satisfies
\begin{equation}\label{e:cancel-kernel}
|\partial_x^\alpha\partial_y^\beta K_{t,s}^{\mathrm{bi}}(x,y)|
\leq C_{\alpha,\beta,\sigma}\,t^{-n-|\alpha|}s^{-m-|\beta|}.
\end{equation}
For the remainder, the bound has an additional factor $(t/s)^2$,
\begin{equation}\label{e:remainder-kernel}
|\partial_x^\alpha\partial_y^\beta K_{t,s}^{\mathrm{rem}}(x,y)|
\leq C_{\alpha,\beta,\sigma}(t/s)^2\,t^{-n-|\alpha|}s^{-m-|\beta|}.
\end{equation}

\item Suppose $s\leq t$. The kernel $K_{t,s}$ is supported where
$|x|_\infty\leq t$ and $|y|_\infty\leq2t$, and its integral in $y$
vanishes for each fixed $x$. In this range,
\begin{equation}\label{e:cubic-kernel}
|\partial_x^\alpha\partial_y^\beta K_{t,s}(x,y)|
\leq C_{\alpha,\beta,\sigma}(s/t)^2\,t^{-n-m-|\alpha|-|\beta|}.
\end{equation}
\end{enumerate}
The constants depend only on $n,m,\alpha,\beta,\sigma$ and $\theta_0$.
\end{proposition}

\begin{proof}
We first verify the support and cancellation properties in \textup{(i)} and \textup{(ii)}.

By \eqref{e:compact-kernel}, $K_{t,s}(x,y)$ vanishes unless
$|x|_\infty\leq t$ and $|y|_\infty\leq t+s$.
For $t<s$, the same support bounds hold for $\theta_t(x)b_{t,s}(y)$
by \eqref{e:kernel-integral}, and hence for $K_{t,s}^{\mathrm{bi}}$
by \eqref{e:kernel-split}. For $s\leq t$, we have $t+s\leq2t$.

Since $\int_{\mathbb R^m}\Phi_{2,s}=0$, we have for every $x$ that
$$
\int_{\mathbb R^m}K_{t,s}(x,y)\,dy
=\int_{\mathbb R^m}\Phi_{1,t}(x,u)
  \bigg(\int_{\mathbb R^m}\Phi_{2,s}(y-u)\,dy\bigg)\,du=0.
$$
Integrating in $x$, we obtain $\int_{\mathbb R^m}b_{t,s}(y)\,dy=0$.
Thus, for $t<s$, both terms in \eqref{e:kernel-split} have integral zero
in $y$. For every $y$, we also have
$$
\int_{\mathbb R^n}K_{t,s}^{\mathrm{bi}}(x,y)\,dx
=b_{t,s}(y)-b_{t,s}(y)\int_{\mathbb R^n}\theta_t(x)\,dx=0.
$$

\smallskip
\noindent\textup{(i)} Suppose $t<s$. In \eqref{e:compact-kernel},
we differentiate $\Phi_{1,t}$ in $x$ and $\Phi_{2,s}$ in $y$. Thus
\begin{align}
|\partial_x^\alpha\partial_y^\beta K_{t,s}(x,y)|
= &\bigg|\int_{\mathbb R^m}\partial_x^\alpha\Phi_{1,t}(x,u)
       \partial^\beta\Phi_{2,s}(y-u)\,du\bigg|\notag\\
\leq &\|\partial_x^\alpha\Phi_{1,t}(x,\cdot)\|_{L^1(\mathbb R^m)}
       \|\partial^\beta\Phi_{2,s}\|_\infty
\leq C_{\alpha,\beta,\sigma}t^{-n-|\alpha|}s^{-m-|\beta|}.
\label{e:kernel-size}
\end{align}
To estimate $b_{t,s}$, we expand $\partial^\beta\Phi_{2,s}(y-u)$
at $u=0$ in \eqref{e:kernel-integral}. Since the mean and first moments
of $\Phi_{2,t}$ vanish, it follows from Taylor's formula that
\begin{align}
|\partial^\beta b_{t,s}(y)|
= &\bigg|\int_{\mathbb R^m}\Phi_{2,t}(u)
   \big[\partial^\beta\Phi_{2,s}(y-u)-\partial^\beta\Phi_{2,s}(y)
       +u\cdot\nabla\partial^\beta\Phi_{2,s}(y)\big]\,du\bigg|\notag\\
\leq &C_{\beta,\sigma}s^{-m-|\beta|-2}
      \int_{\mathbb R^m}|u|^2|\Phi_{2,t}(u)|\,du
\leq C_{\beta,\sigma}(t/s)^2s^{-m-|\beta|}.
\label{e:marginal-remainder}
\end{align}
Multiplying by $|\partial^\alpha\theta_t(x)|\leq C_\alpha t^{-n-|\alpha|}$,
we obtain \eqref{e:remainder-kernel}. Since $t/s<1$, we obtain
\eqref{e:cancel-kernel} by applying \eqref{e:kernel-size} and
\eqref{e:remainder-kernel} to
$K_{t,s}^{\mathrm{bi}}=K_{t,s}-K_{t,s}^{\mathrm{rem}}$.

\smallskip
\noindent\textup{(ii)} Suppose $s\leq t$. We expand
$\partial_x^\alpha\partial_y^\beta\Phi_{1,t}(x,y-v)$ at $v=0$.
Since the mean and first moments of $\Phi_{2,s}$ vanish, only the
Taylor remainder contributes, and we obtain
\begin{align*}
|\partial_x^\alpha\partial_y^\beta K_{t,s}(x,y)|
= &\bigg|\int_{\mathbb R^m}\Phi_{2,s}(v)\int_0^1(1-\tau)
   (v\cdot\nabla_y)^2\partial_x^\alpha\partial_y^\beta
   \Phi_{1,t}(x,y-\tau v)\,d\tau\,dv\bigg|\\
\leq &C_{\alpha,\beta,\sigma}t^{-n-m-|\alpha|-|\beta|-2}
      \int_{\mathbb R^m}|v|^2|\Phi_{2,s}(v)|\,dv\\
\leq &C_{\alpha,\beta,\sigma}(s/t)^2t^{-n-m-|\alpha|-|\beta|}.
\end{align*}
This proves \eqref{e:cubic-kernel} and completes the proof.
\end{proof}

\subsection{Dyadic rectangles and parameter regions}\label{s:dyadic-supports}

Following \cite[Proposition~2.2 and Lemmas~3.1--3.2]{CDWY}, we choose a dyadic
rectangle that contains the translated kernel support in its interior and
lies within the cone in \eqref{e:aux-area}.

For $u=(u_1,\ldots,u_q)\in\mathbb R^q$ and $r>0$, the open ball in the maximum norm is the cube
$$
B_\infty(u,r)=\{v\in\mathbb R^q:|v-u|_\infty<r\}
=\prod_{j=1}^q(u_j-r,u_j+r).
$$
For a nonempty set $E\subseteq\mathbb R^q$, write
$\operatorname{dist}_\infty(u,E)=\inf_{v\in E}|u-v|_\infty$.
We use the adjacent systems of half-open dyadic cubes
$$
\mathcal D_{\mathbb R^q}^{\omega}
=\big\{2^k([0,1)^q+j+(-1)^k\omega):k\in\mathbb Z,\ j\in\mathbb Z^q\big\},
\qquad \omega\in\{0,1/3,2/3\}^q.
$$
The partitions in each system are nested since $3\omega\in\mathbb Z^q$.
For $r>0$, set
$$
L(r)=\min\{2^j:j\in\mathbb Z,\ 2^j>24r\},
\qquad 24r<L(r)\leq48r.
$$
In each coordinate the endpoints of the three grids at this scale are spaced
$L(r)/3$ apart. Since $8r<L(r)/3$, at most one grid has an endpoint in
$[u_i-4r,u_i+4r]$. Choosing a shift in each coordinate thus places
$\overline{B_\infty(u,4r)}$ in the interior of a dyadic cube of side $L(r)$.

Fix an ordering of the shifts in $\{0,1/3,2/3\}^q$. For $u\in\mathbb R^q$ and
$r>0$, let $\kappa_q(u,r)$ be the first shift for which a cube of side $L(r)$
contains $\overline{B_\infty(u,4r)}$ in its interior.
For any dyadic cube $Q$, define $Q_q(u,r)$ by
\begin{equation}\label{e:dyadic-choice}
Q_q(u,r)=Q\quad  \text{if and only if }\quad
Q\in\mathcal D_{\mathbb R^q}^{\kappa_q(u,r)},\quad
\ell(Q)=L(r),\quad \overline{B_\infty(u,4r)}\subset Q^\circ.
\end{equation}
At a fixed scale only one cube of a given system contains $u$. For $u\in Q^\circ$,
$\operatorname{dist}_\infty(u,\partial Q)$ is the smallest distance from $u$ to a face of $Q$. Thus $\overline{B_\infty(u,4r)}\subset Q^\circ$
holds exactly when $u\in Q^\circ$ and $\operatorname{dist}_\infty(u,\partial Q)>4r$.
These conditions and $\ell(Q)=L(r)$ are Borel in $(u,r)$, so both selections are Borel.

For $(u,v,t,s)\in\mathbb R^{n+m}\times(0,\infty)^2$, the selected rectangle and
its pair of shifts are
$$
I\times J=Q_n(u,t)\times Q_m(v,\max\{t,s\}),\qquad
\kappa=(\kappa_n(u,t),\kappa_m(v,\max\{t,s\})).
$$
For fixed $\kappa,I,J$, let $\mathcal T^{\mathrm{str}}_{\kappa,I,J}$ consist of
the parameters selecting this pair of shifts and this rectangle with $t<s$.
Here $\ell(I)=L(t)\leq L(s)=\ell(J)$, with equality possible, and
\begin{equation}\label{e:rect-scales}
\frac{\ell(I)}{48}\leq t<\frac{\ell(I)}{24},\qquad
\frac{\ell(J)}{48}\leq s<\frac{\ell(J)}{24}.
\end{equation}
For $h=0,1,\ldots$, let $\mathcal T^{\mathrm{cub}}_{\kappa,I,J,h}$ consist of
the parameters selecting the same shifts and rectangle with
$2^{-h-1}t<s\leq2^{-h}t$. In this region $I\times J$ is a cube in
$\mathcal D_{\mathbb R^{n+m}}^\kappa$ and
\begin{equation}\label{e:cube-scales}
2^{-h-1}<s/t\leq2^{-h},\qquad
\ell(I)=\ell(J),\qquad
\frac{\ell(I)}{48}\leq t<\frac{\ell(I)}{24}.
\end{equation}
The nonempty regions form a countable disjoint Borel partition of
$\mathbb R^{n+m}\times(0,\infty)^2$.

On the translated kernel support, $|x-u|_\infty\leq t$ and
$|y-v|_\infty\leq t+s$ by Proposition~\ref{p:kernel-bounds}.
The interior margins in the two blocks are therefore
$$
\operatorname{dist}_\infty(x,\partial I)\geq4t-t\geq\ell(I)/24,
\qquad
\operatorname{dist}_\infty(y,\partial J)
\geq4\max\{t,s\}-(t+s)\geq\ell(J)/24.
$$
For the whole rectangle, take $\lambda_0=49\max\{\sqrt n,\sqrt m\}$.
Using $|z|\leq\sqrt q\,|z|_\infty$ in each block and $L(r)\leq48r$, we obtain
\begin{align}
I\times J
\subset &B(u,\sqrt n\,\ell(I))\times B(v,\sqrt m\,\ell(J)) 
\subset  B(u,\lambda_0t)\times B(v,\lambda_0\max\{t,s\}).
\label{e:rectangle-cone}
\end{align}
We now fix $\sigma=\max\{\sigma_0,8\lambda_0\}$. By Proposition~\ref{p:coeff-area},
$\mathcal A_{\lambda_0,F}(H_f^\sigma)\leq C S_F(f)$.

\subsection{Finite sums in the reconstruction}\label{s:finite-sums}

We first restrict the reconstruction integral to finitely many of the
parameter regions in Subsection~\ref{s:dyadic-supports}. The nonempty
regions are indexed by
$$
\mathcal I_{\rm par}
=\{(\kappa,I,J):\mathcal T^{\mathrm{str}}_{\kappa,I,J}\ne\emptyset\}
\cup
\{(\kappa,I,J,h):\mathcal T^{\mathrm{cub}}_{\kappa,I,J,h}\ne\emptyset\}.
$$
For $\gamma\in\mathcal I_{\rm par}$, let $\mathcal T_\gamma$ denote
the corresponding region and $R_\gamma=I\times J$ its rectangle.
We integrate with respect to $d\nu=du\,dv\,(dt/t)(ds/s)$.
Since the regions are disjoint, the indicator of a finite selection
$\mathcal E\subset\mathcal I_{\rm par}$ is
$$
m_{\mathcal E}
=\mathbf1_{\bigcup_{\gamma\in\mathcal E}\mathcal T_\gamma}
=\sum_{\gamma\in\mathcal E}\mathbf1_{\mathcal T_\gamma}.
$$

Fix $H\in L^2_{\mathrm{par}}$ and write $z=(x,y)$.
When $t<s$, the two terms in \eqref{e:kernel-split} are integrated
over the same region and have the common coefficient
$$
\lambda_{\kappa,I,J}^{\mathrm{str}}(\mathcal E)
=\bigg(\int_{\mathcal T^{\mathrm{str}}_{\kappa,I,J}}
       |m_{\mathcal E}H|^2\,d\nu\bigg)^{1/2}.
$$
The term with cancellation in both variables is
\begin{equation}\label{e:bi-atom}
A_{\kappa,I,J}^{\mathrm{bi}}(\mathcal E)(z)
=\int_{\mathcal T^{\mathrm{str}}_{\kappa,I,J}}
K_{t,s}^{\mathrm{bi}}(x-u,y-v)(m_{\mathcal E}H)(u,v,t,s)\,d\nu.
\end{equation}
The remainder is
\begin{equation}\label{e:rem-atom}
A_{\kappa,I,J}^{\mathrm{rem}}(\mathcal E)(z)
=\int_{\mathcal T^{\mathrm{str}}_{\kappa,I,J}}
K_{t,s}^{\mathrm{rem}}(x-u,y-v)(m_{\mathcal E}H)(u,v,t,s)\,d\nu.
\end{equation}

When $s\leq t$, we use $K_{t,s}$ directly. Each region
$\mathcal T^{\mathrm{cub}}_{\kappa,I,J,h}$ has coefficient
$$
\lambda_{\kappa,I,J,h}^{\mathrm{cub}}(\mathcal E)
=\bigg(\int_{\mathcal T^{\mathrm{cub}}_{\kappa,I,J,h}}
       |m_{\mathcal E}H|^2\,d\nu\bigg)^{1/2}
$$
and contributes one term,
\begin{equation}\label{e:cube-sums}
A_{\kappa,I,J,h}^{\mathrm{cub}}(\mathcal E)(z)
=\int_{\mathcal T^{\mathrm{cub}}_{\kappa,I,J,h}}
K_{t,s}(x-u,y-v)(m_{\mathcal E}H)(u,v,t,s)\,d\nu.
\end{equation}

For $f\in L^2(\mathbb R^{n+m})$, take $H=H_f^\sigma$ and define
the finite reconstruction by
$$
f_{\mathcal E}
=\sum_{(\kappa,I,J)\in\mathcal E}
 \big(A_{\kappa,I,J}^{\mathrm{bi}}(\mathcal E)
     +A_{\kappa,I,J}^{\mathrm{rem}}(\mathcal E)\big)
+\sum_{(\kappa,I,J,h)\in\mathcal E}
 A_{\kappa,I,J,h}^{\mathrm{cub}}(\mathcal E).
$$
To identify this sum with the reconstruction integral, note that
$(u,v)\in R_\gamma$ on $\mathcal T_\gamma$. By \eqref{e:rect-scales}
and \eqref{e:cube-scales}, using $ds/s=d(s/t)/(s/t)$ at fixed $t$
in the cubic regions, we have
$$
\nu(\mathcal T_\gamma)\leq(\log 2)^2|R_\gamma|.
$$
Thus $m_{\mathcal E}H_f^\sigma\in L^1(d\nu)$ by Cauchy--Schwarz.
Moreover, $t$ is bounded away from zero on this finite union, and
\eqref{e:compact-kernel} implies
$$
\|K_{t,s}\|_{L^2(\mathbb R^{n+m})}
\leq \|\Phi_{1,t}\|_{L^2(\mathbb R^{n+m})}
     \|\Phi_{2,s}\|_{L^1(\mathbb R^m)}
=t^{-(n+m)/2}\|\Phi_1\|_{L^2(\mathbb R^{n+m})}
              \|\Phi_2\|_{L^1(\mathbb R^m)}.
$$
By Fubini's theorem and \eqref{e:kernel-split}, we obtain
\begin{align}
f_{\mathcal E}
&=\int_0^\infty\int_0^\infty
  K_{t,s}*(m_{\mathcal E}H_f^\sigma)(\cdot,t,s)\,
  \frac{dt}{t}\frac{ds}{s} 
 =\mathcal R_\sigma(m_{\mathcal E}H_f^\sigma)
  \qquad\text{in }L^2(\mathbb R^{n+m}).
\label{e:finite-recon}
\end{align}

For general $H\in L^2_{\mathrm{par}}$, each selected region contributes
once to the square sum of the coefficients. Since $m_{\mathcal E}=1$
on these regions and vanishes on the others, we have
\begin{align}
&\sum_{\kappa,I,J}
 \frac{\lambda_{\kappa,I,J}^{\mathrm{str}}(\mathcal E)^2}{|I||J|}
 \mathbf1_{I\times J}(z)
+\sum_{\kappa,I,J,h}
 \frac{\lambda_{\kappa,I,J,h}^{\mathrm{cub}}(\mathcal E)^2}{|I||J|}
 \mathbf1_{I\times J}(z)
=\sum_{\gamma\in\mathcal E}
 \frac{\mathbf1_{R_\gamma}(z)}{|R_\gamma|}
 \int_{\mathcal T_\gamma}|H|^2\,d\nu.
\label{e:coefficient-square}
\end{align}
We use the indicator of the whole rectangle $R_\gamma$ because a
dyadic average can be nonzero throughout a subrectangle even when
the summand is supported on only part of it.

\begin{proposition}\label{p:two-atoms}
Fix $H\in L^2_{\mathrm{par}}$ and a finite $\mathcal E\subset\mathcal I_{\rm par}$.
The three summands in \eqref{e:bi-atom}--\eqref{e:cube-sums} are smooth and supported
inside $I\times J$ with margins of at least $\ell(I)/24$ in $x$ and $\ell(J)/24$ in $y$.
Each has integral zero in $y$ for every $x$, and the summand in \eqref{e:bi-atom}
also has integral zero in $x$ for every $y$.

For $\alpha\in\mathbb N_0^n$ and $\beta\in\mathbb N_0^m$, the derivative bounds are
\begin{equation}\label{e:bi-deriv}
|\partial_x^\alpha\partial_y^\beta A_{\kappa,I,J}^{\mathrm{bi}}(\mathcal E)(z)|
\leq C_{\alpha,\beta}
\frac{\lambda_{\kappa,I,J}^{\mathrm{str}}(\mathcal E)\,\mathbf1_{I\times J}(z)}
{|I|^{1/2}|J|^{1/2}\ell(I)^{|\alpha|}\ell(J)^{|\beta|}},
\end{equation}
\begin{equation}\label{e:rem-deriv}
|\partial_x^\alpha\partial_y^\beta A_{\kappa,I,J}^{\mathrm{rem}}(\mathcal E)(z)|
\leq C_{\alpha,\beta}
\frac{\lambda_{\kappa,I,J}^{\mathrm{str}}(\mathcal E)\,\mathbf1_{I\times J}(z)}
{|I|^{1/2}|J|^{1/2}\ell(I)^{|\alpha|}\ell(J)^{|\beta|}}
\bigg(\frac{\ell(I)}{\ell(J)}\bigg)^2,
\end{equation}
and
\begin{equation}\label{e:cube-derivatives}
|\partial_x^\alpha\partial_y^\beta A_{\kappa,I,J,h}^{\mathrm{cub}}(\mathcal E)(z)|
\leq C_{\alpha,\beta}2^{-2h}
\frac{\lambda_{\kappa,I,J,h}^{\mathrm{cub}}(\mathcal E)\,\mathbf1_{I\times J}(z)}
{|I|^{1/2}|J|^{1/2}\ell(I)^{|\alpha|}\ell(J)^{|\beta|}}.
\end{equation}
The constants are independent of $I,J,h$ and $\mathcal E$.
\end{proposition}

\begin{proof}
By Proposition~\ref{p:kernel-bounds}, the kernel derivatives are bounded on each fixed parameter region.
Since $m_{\mathcal E}H\in L^1(d\nu)$, we may differentiate under the integral.
Each summand has the support margins from Subsection~\ref{s:dyadic-supports},
so it suffices to estimate its derivatives for $z=(x,y)\in I\times J$.

\smallskip
\noindent{\it Case 1. $t<s$.}
For fixed $t,s$, the support of either translated kernel in $(u,v)$ has measure
at most $Ct^ns^m$. Applying Cauchy--Schwarz in \eqref{e:bi-atom} and using
\eqref{e:cancel-kernel} and \eqref{e:rect-scales}, we obtain
\begin{align*}
|\partial_x^\alpha\partial_y^\beta A_{\kappa,I,J}^{\mathrm{bi}}(\mathcal E)(z)|^2
\leq &\lambda_{\kappa,I,J}^{\mathrm{str}}(\mathcal E)^2
\int_{\mathcal T^{\mathrm{str}}_{\kappa,I,J}}
|\partial_x^\alpha\partial_y^\beta K_{t,s}^{\mathrm{bi}}(x-u,y-v)|^2\,d\nu\\
\leq &C_{\alpha,\beta}\lambda_{\kappa,I,J}^{\mathrm{str}}(\mathcal E)^2
\int_{\ell(I)/48}^{\ell(I)/24}t^{-n-2|\alpha|}\,\frac{dt}{t}
\int_{\ell(J)/48}^{\ell(J)/24}s^{-m-2|\beta|}\,\frac{ds}{s}\\
\leq &\frac{C_{\alpha,\beta}\lambda_{\kappa,I,J}^{\mathrm{str}}(\mathcal E)^2}
{|I||J|\ell(I)^{2|\alpha|}\ell(J)^{2|\beta|}}.
\end{align*}
For the remainder, squaring \eqref{e:remainder-kernel} introduces the additional
factor $(t/s)^4$. Since $t/s\leq2\ell(I)/\ell(J)$ on this region, the same calculation
in \eqref{e:rem-atom} gives
\begin{align*}
|\partial_x^\alpha\partial_y^\beta A_{\kappa,I,J}^{\mathrm{rem}}(\mathcal E)(z)|^2
\leq &\lambda_{\kappa,I,J}^{\mathrm{str}}(\mathcal E)^2
\int_{\mathcal T^{\mathrm{str}}_{\kappa,I,J}}
|\partial_x^\alpha\partial_y^\beta K_{t,s}^{\mathrm{rem}}(x-u,y-v)|^2\,d\nu\\
\leq &\frac{C_{\alpha,\beta}\lambda_{\kappa,I,J}^{\mathrm{str}}(\mathcal E)^2}
{|I||J|\ell(I)^{2|\alpha|}\ell(J)^{2|\beta|}}
\bigg(\frac{\ell(I)}{\ell(J)}\bigg)^4.
\end{align*}
Taking square roots, we obtain \eqref{e:bi-deriv} and \eqref{e:rem-deriv}.

\smallskip
\noindent{\it Case 2. $s\leq t$.}
Here $\ell(I)=\ell(J)$, and the support of the translated kernel in $(u,v)$
has measure at most $Ct^{n+m}$. At fixed $t$, the substitution $r=s/t$ in
\eqref{e:cube-scales} restricts $r$ to $(2^{-h-1},2^{-h}]$.
By Cauchy--Schwarz in \eqref{e:cube-sums} and the kernel bound \eqref{e:cubic-kernel}, we have
\begin{align*}
|\partial_x^\alpha\partial_y^\beta A_{\kappa,I,J,h}^{\mathrm{cub}}(\mathcal E)(z)|^2
\leq &\lambda_{\kappa,I,J,h}^{\mathrm{cub}}(\mathcal E)^2
\int_{\mathcal T^{\mathrm{cub}}_{\kappa,I,J,h}}
|\partial_x^\alpha\partial_y^\beta K_{t,s}(x-u,y-v)|^2\,d\nu\\
\leq &C_{\alpha,\beta}\lambda_{\kappa,I,J,h}^{\mathrm{cub}}(\mathcal E)^2
\int_{\ell(I)/48}^{\ell(I)/24}t^{-(n+m)-2|\alpha|-2|\beta|}\,\frac{dt}{t}
\int_{2^{-h-1}}^{2^{-h}}r^4\,\frac{dr}{r}\\
\leq &\frac{C_{\alpha,\beta}2^{-4h}\lambda_{\kappa,I,J,h}^{\mathrm{cub}}(\mathcal E)^2}
{|I||J|\ell(I)^{2|\alpha|}\ell(J)^{2|\beta|}}.
\end{align*}
Taking square roots proves \eqref{e:cube-derivatives}.

Finally, for each $y$, Fubini's theorem and the cancellation in $x$ from
Proposition~\ref{p:kernel-bounds} imply
$$
\int_{\mathbb R^n}A_{\kappa,I,J}^{\mathrm{bi}}(\mathcal E)(x,y)\,dx
=\int_{\mathcal T^{\mathrm{str}}_{\kappa,I,J}}(m_{\mathcal E}H)(u,v,t,s)
\bigg(\int_{\mathbb R^n}K_{t,s}^{\mathrm{bi}}(x-u,y-v)\,dx\bigg)\,d\nu=0.
$$
Integrating in $y$ instead, we obtain the corresponding cancellation for all three summands.
\end{proof}

For fixed $\kappa,I,J$, the cubic summands are all supported on $I\times J$.
Summing \eqref{e:cube-derivatives} over $h$ and applying Cauchy--Schwarz, we obtain
\begin{align}
\bigg|\partial_x^\alpha\partial_y^\beta
\sum_{h\geq0}A_{\kappa,I,J,h}^{\mathrm{cub}}(\mathcal E)(z)\bigg|
\leq &\frac{C_{\alpha,\beta}\mathbf1_{I\times J}(z)}
{|I|^{1/2}|J|^{1/2}\ell(I)^{|\alpha|}\ell(J)^{|\beta|}}
\sum_{h\geq0}2^{-2h}\lambda_{\kappa,I,J,h}^{\mathrm{cub}}(\mathcal E)\notag\\
\leq &\frac{C_{\alpha,\beta}\mathbf1_{I\times J}(z)}
{|I|^{1/2}|J|^{1/2}\ell(I)^{|\alpha|}\ell(J)^{|\beta|}}
\bigg(\sum_{h\geq0}\lambda_{\kappa,I,J,h}^{\mathrm{cub}}(\mathcal E)^2\bigg)^{1/2}.
\label{e:cube-sum}
\end{align}

\subsection{Coefficient bounds and convergence of the finite sums}\label{s:coefficient-control}

\begin{proposition}\label{p:coeff-control}
Let $f\in L^2(\mathbb R^{n+m})$ and set $H=H_f^\sigma$.
\begin{enumerate}[(i)]
\item For every finite $\mathcal E\subset\mathcal I_{\rm par}$ and every $z\in\mathbb R^{n+m}$,
\begin{equation}\label{e:coefficient-control}
\bigg(\sum_{\gamma\in\mathcal E}\frac{\mathbf1_{R_\gamma}(z)}{|R_\gamma|}\int_{\mathcal T_\gamma}|H|^2\,d\nu\bigg)^{1/2}
\leq C_1\mathcal A_{\lambda_0,F}(m_{\mathcal E}H)(z)\leq C_2 S_F(f)(z),
\end{equation}
where $C_1,C_2$ depend only on $n,m$.
For every $H\in L^2_{\mathrm{par}}$, the first inequality holds and the left-hand side is nondecreasing in $\mathcal E$.

\item If the finite sets $\mathcal E_N$ increase to $\mathcal I_{\rm par}$, then
$f_{\mathcal E_N}\to f$ in $L^2(\mathbb R^{n+m})$.
\end{enumerate}
\end{proposition}

\begin{proof}
\noindent\textup{(i)}
Fix $H\in L^2_{\mathrm{par}}$ and $z=(x,y)$. On $\mathcal T_\gamma$,
we have $|R_\gamma|\simeq t^n\max\{t,s\}^m$ by \eqref{e:rect-scales} and \eqref{e:cube-scales}.
Combining this with \eqref{e:rectangle-cone}, we obtain
$$
\frac{\mathbf1_{R_\gamma}(z)}{|R_\gamma|}
\leq C\,\frac{\mathbf1_{B(u,\lambda_0t)\times B(v,\lambda_0\max\{t,s\})}(z)}{t^n\max\{t,s\}^m}.
$$
Integrating against $|H|^2\,d\nu$ on $\mathcal T_\gamma$ and summing over $\gamma\in\mathcal E$, we have
\begin{align*}
\sum_{\gamma\in\mathcal E}\frac{\mathbf1_{R_\gamma}(z)}{|R_\gamma|}\int_{\mathcal T_\gamma}|H|^2\,d\nu
\leq &C\int_0^\infty\int_0^\infty
\int_{\substack{|x-u|<\lambda_0t\\|y-v|<\lambda_0\max\{t,s\}}}
\frac{|(m_{\mathcal E}H)(u,v,t,s)|^2}{t^n\max\{t,s\}^m}\,du\,dv\,\frac{dt}{t}\frac{ds}{s}\\
= &C\mathcal A_{\lambda_0,F}(m_{\mathcal E}H)(z)^2,
\end{align*}
since $m_{\mathcal E}$ is the indicator of the disjoint union of these regions.
Taking square roots, we obtain the first inequality in \eqref{e:coefficient-control}.
The sum on the left is nondecreasing in $\mathcal E$ because its terms are nonnegative.
For $H=H_f^\sigma$, we use $|m_{\mathcal E}|\leq1$ and \eqref{e:coeff-area} to obtain
$$
\mathcal A_{\lambda_0,F}(m_{\mathcal E}H_f^\sigma)(z)
\leq\mathcal A_{\lambda_0,F}(H_f^\sigma)(z)\leq C S_F(f)(z).
$$

\smallskip
\noindent\textup{(ii)}
Let $\mathcal E_N\uparrow\mathcal I_{\rm par}$ through finite sets.
Since the regions $\mathcal T_\gamma$ partition the parameter space, $m_{\mathcal E_N}\to1$ pointwise.
For every $H\in L^2_{\mathrm{par}}$, it follows from \eqref{e:recon-l2} that
\begin{align}
\|\mathcal R_\sigma(m_{\mathcal E_N}H)-\mathcal R_\sigma H\|_{L^2(\mathbb R^{n+m})}
= &\|\mathcal R_\sigma((m_{\mathcal E_N}-1)H)\|_{L^2(\mathbb R^{n+m})}\notag\\
\leq &\|a_\sigma\|_{L^2(dr/r)}^2
\bigg(\int_{\{m_{\mathcal E_N}=0\}}|H|^2\,d\nu\bigg)^{1/2}\to0,\label{e:finite-recon-limit}
\end{align}
where the limit follows by dominated convergence.

Taking $H=H_f^\sigma$ and using \eqref{e:finite-recon} and \eqref{e:reproduce},
we obtain $f_{\mathcal E_N}\to f$ in $L^2(\mathbb R^{n+m})$.
\end{proof}


\section{Intermediate square functions and local Hoeffding estimates}\label{s:local-hoeffding}

We begin with the atoms $A^{\mathrm{bi}}$ in Proposition~\ref{p:two-atoms},
which have cancellation in both $x$ and $y$. We first apply the dyadic
square-function estimate in $x$. The sum in $y$ remains linear inside the
intermediate square function, so we can apply the estimate once more in $y$.
Each step contributes a factor $\sqrt p$.

For $A^{\mathrm{rem}}$ and $A^{\mathrm{cub}}$, we use cancellation only in $y$.
For the remainders, the factor $(\ell(I)/\ell(J))^2$ in \eqref{e:rem-deriv}
is summable over the generations of $I$ at fixed $J$.
For the cubic terms, we use \eqref{e:cube-sum} to sum over $h$.
At fixed $J$, the remaining cubes $I$ have side length $\ell(J)$
and are pairwise disjoint within each dyadic system.
In both cases, a single square-function estimate in $y$ contributes the factor $\sqrt p$.

\subsection{Cancellation in both variables}\label{s:both-cancellations}

Fix $H\in L^2_{\mathrm{par}}$, a finite set $\mathcal E\subset\mathcal I_{\rm par}$
and a pair of shifts $\kappa=(\omega_x,\omega_y)$.
Write $\mathcal D_x=\mathcal D_{\mathbb R^n}^{\omega_x}$ and
$\mathcal D_y=\mathcal D_{\mathbb R^m}^{\omega_y}$.
For $I\in\mathcal D_x$ and $J\in\mathcal D_y$, write
$A_{I,J}=A_{\kappa,I,J}^{\mathrm{bi}}(\mathcal E)$.
We estimate the finite sum $G=\sum_{I,J}A_{I,J}$ in terms of the normalized coefficients
$$
c_{I,J}
=\frac{\lambda_{\kappa,I,J}^{\mathrm{str}}(\mathcal E)}{(|I||J|)^{1/2}}
=\bigg(\frac1{|I||J|}\int_{\mathcal T_{\kappa,I,J}^{\mathrm{str}}}
             |m_{\mathcal E}H|^2\,d\nu\bigg)^{1/2}.
$$

For a dyadic cube $Q$, let $\langle h\rangle_Q=|Q|^{-1}\int_Q h$ and denote its dyadic children by $\operatorname{ch}(Q)$.
The martingale difference on $Q$ is
$$
\Delta_Qh=\sum_{Q'\in\operatorname{ch}(Q)}
(\langle h\rangle_{Q'}-\langle h\rangle_Q)\mathbf1_{Q'}.
$$
We use $\Delta_U^x$ and $\Delta_V^y$ for the differences in $x$ and $y$,
respectively, and $\Delta_U^{(j)}$ for the difference in the $j$th coordinate block.

In the $x$ variable we use the scalar estimate
\begin{equation}\label{e:scalar-square}
\bigg\|\sum_{Q\in\mathcal F}\Delta_Qh\bigg\|_p
\leq C_q\sqrt p\,\bigg\|\bigg(\sum_{Q\in\mathcal F}|\Delta_Qh|^2\bigg)^{1/2}\bigg\|_p,
\qquad 2\leq p<\infty,
\end{equation}
where $h$ is locally integrable and $\mathcal F$ is a finite family in a single dyadic system on $\mathbb R^q$.
To see this, let $g=\sum_{Q\in\mathcal F}\Delta_Qh$ and fix a maximal cube $Q_0\in\mathcal F$.
Then $\langle g\rangle_{Q_0}=0$, and for every dyadic $R\subseteq Q_0$,
$$
\Delta_Rg=\mathbf1_{\{R\in\mathcal F\}}\Delta_Rh.
$$
Equip $Q_0$ with the normalized measure $d\mu_0=dx/|Q_0|$.
Each dyadic child has $2^{-q}$ times the measure of its parent.
For the dyadic maximal function $g^*(x)=\sup_{x\in R\subseteq Q_0}|\langle g\rangle_R|$,
Ivanisvili and Treil~\cite[Theorem~4.3 of the preprint]{IT} proved that
$$
\|g\|_{L^p(d\mu_0)}\leq\|g^*\|_{L^p(d\mu_0)}
\leq C_q\sqrt p\,\bigg\|\bigg(\sum_{R\subseteq Q_0}|\Delta_Rg|^2\bigg)^{1/2}\bigg\|_{L^p(d\mu_0)}.
$$
Taking $p$th powers, multiplying by $|Q_0|$ and summing over the disjoint maximal cubes, we obtain \eqref{e:scalar-square}.
For complex $h$, we apply the estimate to its real and imaginary parts.

After taking the square function in $x$, the martingale differences in $y$ take values in $\ell^2$.
We therefore need the Hilbert-valued form of \eqref{e:scalar-square}.
Let $\mathcal H$ be a Hilbert space and let $d_Q=\Delta_Qd_Q$ be $\mathcal H$-valued differences
indexed by a finite family in a single dyadic system on $\mathbb R^q$. Then
\begin{equation}\label{e:hilbert-square}
\bigg\|\sum_Qd_Q\bigg\|_{L^p(\mathcal H)}
\leq C_q\sqrt p\,\bigg\|\bigg(\sum_Q\|d_Q\|_{\mathcal H}^2\bigg)^{1/2}\bigg\|_p,
\qquad 2\leq p<\infty,
\end{equation}
with $C_q$ independent of $\mathcal H$.

Each $d_Q$ is supported on $Q$ and constant on its children.
We may thus identify the real span of their values with a subspace of some $\mathbb R^L$,
using the real part of the inner product when $\mathcal H$ is complex.
By rotational invariance of the normalized surface measure $\nu_L$ on $S^{L-1}$,
$$
\bigg(\int_{S^{L-1}}|\langle\omega,v\rangle|^p\,d\nu_L(\omega)\bigg)^{1/p}
=\|\omega_1\|_{L^p(d\nu_L)}\|v\|_{\mathcal H}.
$$
For each $\omega$, the $Q$th difference of $\sum_U\langle\omega,d_U\rangle$ is
$\langle\omega,d_Q\rangle$, since $\Delta_Q\Delta_U=0$ for $Q\ne U$.
Applying \eqref{e:scalar-square} for each $\omega$, then using Fubini's theorem
and the triangle inequality in $L^{p/2}(d\nu_L)$, we obtain
\begin{align*}
\|\omega_1\|_{L^p(d\nu_L)}\,\bigg\|\sum_Qd_Q\bigg\|_{L^p(\mathcal H)}
= &\bigg\|\sum_Q\langle\omega,d_Q\rangle\bigg\|_{L^p(dx\,d\nu_L)}\\
\leq &C_q\sqrt p\,\bigg\|\bigg(\sum_Q|\langle\omega,d_Q\rangle|^2\bigg)^{1/2}\bigg\|_{L^p(dx\,d\nu_L)}\\
\leq &C_q\sqrt p\,\bigg\|\bigg(\sum_Q
\|\langle\omega,d_Q(x)\rangle\|_{L^p(d\nu_L)}^2\bigg)^{1/2}\bigg\|_{L^p(dx)}\\
= &C_q\sqrt p\,\|\omega_1\|_{L^p(d\nu_L)}\,\bigg\|\bigg(\sum_Q\|d_Q\|_{\mathcal H}^2\bigg)^{1/2}\bigg\|_p.
\end{align*}
That is, we obtain \eqref{e:hilbert-square}.

\begin{proposition}\label{p:product-square}
Let $G=\sum_{I,J}A_{I,J}$ be the finite sum above. Its product dyadic square function satisfies
\begin{equation}\label{e:product-square}
S_{xy}G(x,y):=\bigg(\sum_{U\in\mathcal D_x}\sum_{V\in\mathcal D_y}
|\Delta_U^x\Delta_V^yG(x,y)|^2\bigg)^{1/2}
\leq C\,\bigg(\sum_{I,J}c_{I,J}^2\mathbf1_{I\times J}(x,y)\bigg)^{1/2}.
\end{equation}
Moreover, for $2\leq p<\infty$,
\begin{equation}\label{e:product-lp}
\|G\|_p\leq Cp\,\bigg\|\bigg(\sum_{I,J}c_{I,J}^2\mathbf1_{I\times J}\bigg)^{1/2}\bigg\|_p.
\end{equation}
The constants depend only on $n,m$.
\end{proposition}

\begin{proof}
Since $A_{I,J}$ is supported in $I\times J$ and has integral zero in each variable,
$$
\Delta_U^x\Delta_V^yA_{I,J}=0\qquad\text{unless }U\subseteq I,\ V\subseteq J.
$$
For $a\geq0$, let $U^{(a)}$ denote the dyadic ancestor of $U$ with side length $2^a\ell(U)$.
If $I=U^{(a)}$ and $J=V^{(b)}$, then by the definition of the two differences and \eqref{e:bi-deriv},
\begin{align*}
|\Delta_U^x\Delta_V^yA_{I,J}|
\leq &\sup_{\substack{u,u'\in U\\v,v'\in V}}
|A_{I,J}(u,v)-A_{I,J}(u',v)-A_{I,J}(u,v')+A_{I,J}(u',v')|\mathbf1_{U\times V}\\
\leq &C\ell(U)\ell(V)\|\nabla_x\nabla_yA_{I,J}\|_\infty\mathbf1_{U\times V}\\
\leq &C2^{-a-b}c_{I,J}\mathbf1_{U\times V}.
\end{align*}
Set $c_{I,J}=0$ for pairs outside the finite sum. Summing over the ancestors of $U$ and $V$, we obtain
$$
|\Delta_U^x\Delta_V^yG(x,y)|
\leq C\sum_{a,b\geq0}2^{-a-b}c_{U^{(a)},V^{(b)}}\mathbf1_{U\times V}(x,y).
$$
For fixed $a,b$, the rectangles $U\times V$ with $U^{(a)}=I$ and $V^{(b)}=J$ partition $I\times J$.
By Minkowski's inequality in $\ell^2(U,V)$,
\begin{align*}
S_{xy}G(x,y)
\leq &C\sum_{a,b\geq0}2^{-a-b}
\bigg(\sum_{U,V}c_{U^{(a)},V^{(b)}}^2\mathbf1_{U\times V}(x,y)\bigg)^{1/2}\\
= &C\sum_{a,b\geq0}2^{-a-b}
\bigg(\sum_{I,J}c_{I,J}^2
\sum_{\substack{U^{(a)}=I\\V^{(b)}=J}}\mathbf1_{U\times V}(x,y)\bigg)^{1/2}\\
\leq &C\,\bigg(\sum_{I,J}c_{I,J}^2\mathbf1_{I\times J}(x,y)\bigg)^{1/2},
\end{align*}
namely,  \eqref{e:product-square} holds.

It suffices to  prove \eqref{e:product-lp}, we first truncate the martingale expansion of $G$.
Let $\mathbb E_j^x$ denote dyadic averaging in $x$ at side length $2^{-j}$,
$$
(\mathbb E_j^xG)(x,y)=\sum_{\substack{U\in\mathcal D_x\\\ell(U)=2^{-j}}}
\frac{\mathbf1_U(x)}{|U|}\int_U G(u,y)\,du,
$$
and define $\mathbb E_l^y$ similarly. Choose $j_0,l_0\in\mathbb Z$ so that
$\ell(I)\leq2^{-j_0}$ and $\ell(J)\leq2^{-l_0}$ for every pair in the sum.
By the separate cancellation of each atom, $\mathbb E_{j_0}^xG=\mathbb E_{l_0}^yG=0$.
For $M\geq\max\{j_0,l_0\}$, let
$$
\mathcal U_M=\{U\in\mathcal D_x:2^{-M}\leq\ell(U)\leq2^{-j_0}\},\qquad
\mathcal V_M=\{V\in\mathcal D_y:2^{-M}\leq\ell(V)\leq2^{-l_0}\}.
$$
Write $d_{U,V}=\Delta_U^x\Delta_V^yG$. Note that
\begin{align}
G_M:=\sum_{U\in\mathcal U_M}\sum_{V\in\mathcal V_M}d_{U,V}
= &(\mathbb E_{M+1}^x-\mathbb E_{j_0}^x)(\mathbb E_{M+1}^y-\mathbb E_{l_0}^y)G
= \mathbb E_{M+1}^x\mathbb E_{M+1}^yG\to G\quad\text{in }L^p(\mathbb R^{n+m}),\label{e:martingale-sums}
\end{align}
where only finitely many terms are nonzero and the convergence follows from $G\in C_c^\infty(\mathbb R^{n+m})$.

We now apply Pipher's vector-valued argument~\cite[pp.~73--76]{Pipher86}
using the pointwise intermediate square function of Li, Pipher and Wu~\cite{LiPipherWu}.
In the following sums, $U\in\mathcal U_M$ and $V\in\mathcal V_M$. Taking the square function in $x$ leaves the sum in $V$ inside the $\ell^2$ norm,
$$
S_xG_M(x,y)=\bigg(\sum_U\bigg|\sum_Vd_{U,V}(x,y)\bigg|^2\bigg)^{1/2}
=\bigg\|\sum_V\big(d_{U,V}(x,y)\big)_U\bigg\|_{\ell^2}.
$$
For fixed $y$, we apply \eqref{e:scalar-square} in $x$ to obtain
\begin{equation}\label{e:horizontal-square}
\|G_M(\cdot,y)\|_{L^p(\mathbb R^n)}
\leq C_n\sqrt p\,\|S_xG_M(\cdot,y)\|_{L^p(\mathbb R^n)}.
\end{equation}
For fixed $x$, each vector $(d_{U,V}(x,\cdot))_U$ is an $\ell^2$-valued martingale difference on $V$.
The square function in $y$ is therefore
$$
S_{xy}G_M(x,y)=\bigg(\sum_V\|(d_{U,V}(x,y))_U\|_{\ell^2}^2\bigg)^{1/2}
=\bigg(\sum_{U,V}|d_{U,V}(x,y)|^2\bigg)^{1/2}.
$$
By \eqref{e:hilbert-square} in $y$,
\begin{equation}\label{e:vertical-hilbert}
\|S_xG_M(x,\cdot)\|_{L^p(\mathbb R^m)}
\leq C_m\sqrt p\,\|S_{xy}G_M(x,\cdot)\|_{L^p(\mathbb R^m)}.
\end{equation}
Taking $p$th powers in \eqref{e:horizontal-square} and \eqref{e:vertical-hilbert}
and integrating in $y$ and $x$, respectively, we obtain
\begin{align*}
\|G_M\|_p
\leq &C_n\sqrt p\,\|S_xG_M\|_p
\leq C_nC_m p\,\|S_{xy}G_M\|_p
\leq C_nC_m p\,\|S_{xy}G\|_p
\leq Cp\bigg\|\bigg(\sum_{I,J}c_{I,J}^2\mathbf1_{I\times J}\bigg)^{1/2}\bigg\|_p,
\end{align*}
where the last inequality follows from \eqref{e:product-square}.
Letting $M\to\infty$ in \eqref{e:martingale-sums}, we obtain \eqref{e:product-lp}
with the same constant for every finite subfamily.
\end{proof}

\subsection{Vertical cancellation}\label{s:vertical-cancellation}

Keep $\kappa$ and $\mathcal E$ fixed and omit them from the notation for $A$ and $\lambda$.
Both the remainder and the cubic terms have integral zero in $y$ for each $x$.
To apply the martingale estimate in $y$, we first sum the terms with the same vertical cube $J$.

For the remainders, let $B_J^{\mathrm{rem}}=\sum_I A_{I,J}^{\mathrm{rem}}$ and define its coefficient by
$$
c_J^{\mathrm{rem}}(x)^2=\sum_I\frac{(\lambda_{I,J}^{\mathrm{str}})^2}{|I||J|}\mathbf1_I(x).
$$
For the cubic terms, the sum also runs over $h$. Let $B_J^{\mathrm{cub}}=\sum_{I,h}A_{I,J,h}^{\mathrm{cub}}$, with coefficient
$$
c_J^{\mathrm{cub}}(x)^2=\sum_{I,h}\frac{(\lambda_{I,J,h}^{\mathrm{cub}})^2}{|I||J|}\mathbf1_I(x).
$$
Summing over $J$, we write $G^{\mathrm{rem}}=\sum_J B_J^{\mathrm{rem}}$ and $G^{\mathrm{cub}}=\sum_J B_J^{\mathrm{cub}}$.

\begin{proposition}\label{p:vertical-square}
Let $G=\sum_J B_J$, where $(B_J,c_J)$ denotes either the remainder family
$(B_J^{\mathrm{rem}},c_J^{\mathrm{rem}})$ or the cubic family $(B_J^{\mathrm{cub}},c_J^{\mathrm{cub}})$.
Then
\begin{equation}\label{e:vertical-square}
S_yG(x,y):=\bigg(\sum_{V\in\mathcal D_y}|\Delta_V^yG(x,y)|^2\bigg)^{1/2}
\leq C\bigg(\sum_J c_J(x)^2\mathbf1_J(y)\bigg)^{1/2},
\end{equation}
and, for $2\leq p<\infty$,
\begin{equation}\label{e:vertical-lp}
\|G\|_p\leq C\sqrt p\,\bigg\|\bigg(\sum_J c_J(x)^2\mathbf1_J(y)\bigg)^{1/2}\bigg\|_{L^p(dx\,dy)}.
\end{equation}
The constants depend only on $n,m$.
\end{proposition}

\begin{proof}
Fix $x$ and $J$. For the remainders, at most one cube $I$ of each side length $2^{-k}\ell(J)$ contains $x$.
By \eqref{e:rem-deriv} and Cauchy--Schwarz,
\begin{align*}
|\partial_y^\beta B_J^{\mathrm{rem}}(x,y)|
\leq &\frac{C_\beta\mathbf1_J(y)}{\ell(J)^{|\beta|}}
\sum_{I\ni x}\Big(\frac{\ell(I)}{\ell(J)}\Big)^2
\frac{\lambda_{I,J}^{\mathrm{str}}}{(|I||J|)^{1/2}}\\
\leq &\frac{C_\beta\mathbf1_J(y)}{\ell(J)^{|\beta|}}
\bigg(\sum_{k\geq0}2^{-4k}\bigg)^{1/2}
\bigg(\sum_{I\ni x}\frac{(\lambda_{I,J}^{\mathrm{str}})^2}{|I||J|}\bigg)^{1/2}\\
\leq &C_\beta\ell(J)^{-|\beta|}c_J^{\mathrm{rem}}(x)\mathbf1_J(y).
\end{align*}
For the cubic terms, the cubes $I$ have side length $\ell(J)$ and are pairwise disjoint. By \eqref{e:cube-sum},
\begin{align*}
|\partial_y^\beta B_J^{\mathrm{cub}}(x,y)|
\leq &\frac{C_\beta\mathbf1_J(y)}{\ell(J)^{|\beta|}}
\sum_{I\ni x}\bigg(\sum_h\frac{(\lambda_{I,J,h}^{\mathrm{cub}})^2}{|I||J|}\bigg)^{1/2}\\
= &C_\beta\ell(J)^{-|\beta|}c_J^{\mathrm{cub}}(x)\mathbf1_J(y).
\end{align*}
Thus both families satisfy
\begin{equation}\label{e:vertical-size}
|\partial_y^\beta B_J(x,y)|
\leq C_\beta\ell(J)^{-|\beta|}c_J(x)\mathbf1_J(y),
\qquad \beta\in\mathbb N_0^m.
\end{equation}

By support and cancellation in $y$, $\Delta_V^yB_J=0$ unless $V\subseteq J$.
For $J=V^{(k)}$, the first derivative bound in \eqref{e:vertical-size} implies
\begin{align*}
|\Delta_V^yB_J(x,y)|
\leq &\sup_{v,v'\in V}|B_J(x,v)-B_J(x,v')|\mathbf1_V(y)\\
\leq &C\frac{\ell(V)}{\ell(J)}c_J(x)\mathbf1_V(y)
=C2^{-k}c_{V^{(k)}}(x)\mathbf1_V(y).
\end{align*}
Set $c_J=0$ outside the finite family. By Minkowski's inequality in $\ell^2(V)$,
\begin{align*}
S_yG(x,y)
\leq &C\sum_{k\geq0}2^{-k}
\bigg(\sum_{V\in\mathcal D_y}c_{V^{(k)}}(x)^2\mathbf1_V(y)\bigg)^{1/2}\\
= &C\sum_{k\geq0}2^{-k}
\bigg(\sum_J c_J(x)^2\sum_{V^{(k)}=J}\mathbf1_V(y)\bigg)^{1/2}\\
\leq &C\bigg(\sum_J c_J(x)^2\mathbf1_J(y)\bigg)^{1/2},
\end{align*}
since the cubes $V$ with $V^{(k)}=J$ partition $J$. We obtain \eqref{e:vertical-square}.

Choose $l_0$ so that $\ell(J)\leq2^{-l_0}$ for every $J$ in the sum. By cancellation,
$\mathbb E_{l_0}^yG=0$. As in \eqref{e:martingale-sums},
$$
\sum_{\substack{V\in\mathcal D_y\\2^{-M}\leq\ell(V)\leq2^{-l_0}}}\Delta_V^yG(x,\cdot)
=\mathbb E_{M+1}^yG(x,\cdot)\to G(x,\cdot)\quad\text{in }L^p(\mathbb R^m).
$$
Applying \eqref{e:scalar-square} to these sums and using \eqref{e:vertical-square}, we obtain
\begin{equation}\label{e:vertical-fixed}
\|G(x,\cdot)\|_{L^p(dy)}
\leq C_m\sqrt p\,\|S_yG(x,\cdot)\|_{L^p(dy)}
\leq C\sqrt p\,\bigg\|\bigg(\sum_J c_J(x)^2\mathbf1_J\bigg)^{1/2}\bigg\|_{L^p(dy)}.
\end{equation}
Taking $p$th powers and integrating in $x$, we obtain \eqref{e:vertical-lp}.

For fixed $x$, the proof of \eqref{e:vertical-fixed} uses only support in $J$, mean zero on $J$
and \eqref{e:vertical-size} for $|\beta|\leq1$. The estimate thus holds for any finite family
satisfying these conditions in a dyadic system on $\mathbb R^q$, with a constant depending
only on $q$ and the constants in those bounds.
\end{proof}

Let $2\leq p<\infty$ and let $f\in L^2(\mathbb R^{n+m})$ with $S_F(f)\in L^p(\mathbb R^{n+m})$.
For each pair of shifts, we apply \eqref{e:product-lp} to $\sum_{I,J}A_{I,J}^{\mathrm{bi}}$
and \eqref{e:vertical-lp} to $G^{\mathrm{rem}}$ and $G^{\mathrm{cub}}$.
Summing over the $3^{n+m}$ pairs of shifts and using \eqref{e:coefficient-control}, we obtain
$$
\|f_{\mathcal E}\|_p\leq C(p+\sqrt p)\,\bigg\|\bigg(\sum_{\gamma\in\mathcal E}\frac{\mathbf1_{R_\gamma}}{|R_\gamma|}\int_{\mathcal T_\gamma}|H_f^\sigma|^2\,d\nu\bigg)^{1/2}\bigg\|_p\leq Cp\|S_F(f)\|_p.
$$
Choose finite sets $\mathcal E_N\uparrow\mathcal I_{\rm par}$.
By \eqref{e:finite-recon-limit}, $f_{\mathcal E_N}\to f$ in $L^2(\mathbb R^{n+m})$ and hence almost everywhere along a subsequence.
Along this subsequence,
$$
\|f\|_p\leq\liminf_{j\to\infty}\|f_{\mathcal E_{N_j}}\|_p\leq Cp\|S_F(f)\|_p,
$$
where the first inequality follows from Fatou's lemma.

\subsection{Local estimates with cancellation in both variables}\label{s:local-centring}

Fix a product $R_0=I_0\times J_0$ of axis-parallel cubes. To estimate the oscillation on $R_0$,
we centre in each variable separately. In the $x$ variable, let
$$
P_xg(x,y)=\frac{\mathbf1_{I_0}(x)}{|I_0|}\int_{I_0}g(u,y)\,du
$$
and set $Q_xg=\mathbf1_{I_0}g-P_xg$. Define $P_y,Q_y$ similarly on $J_0$.
Since $P_x,P_y$ are contractions on $L^r(R_0)$ for $1\leq r\leq\infty$,
both $Q_x$ and $Q_y$ have norm at most $2$ on this space.
For the sum $G$ in Proposition~\ref{p:product-square} and $(x,y)\in R_0$, we have
\begin{align*}
Q_{R_0}G(x,y)=Q_xQ_yG(x,y)
= &G(x,y)-\frac1{|I_0|}\int_{I_0}G(u,y)\,du\\
&\quad-\frac1{|J_0|}\int_{J_0}G(x,v)\,dv
+\frac1{|R_0|}\int_{R_0}G(u,v)\,du\,dv.
\end{align*}

For a cube $I$, write $I_*=\{z\in I:\operatorname{dist}_\infty(z,\partial I)\geq\ell(I)/48\}$.
By Proposition~\ref{p:two-atoms}, each atom is supported in $I_*\times J_*$.
Its coefficient in \eqref{e:product-square}, however, is counted on all of $I\times J$,
which may extend beyond $R_0$ even when $\ell(I)\leq\ell(I_0)$ and $\ell(J)\leq\ell(J_0)$.
For rectangles of these sizes with $(I_*\times J_*)\cap R_0\ne\varnothing$,
the next lemma bounds the $L^p(\mathbb R^{n+m})$ norm of the coefficient square function
by its $L^p(R_0)$ norm.

\begin{lemma}\label{l:cube-fold}
Let $E=\prod_{j=1}^h E_j\subset\mathbb R^d$ be a product of axis-parallel cubes.
Let $\mathcal I$ be a finite family of products $\mathbf I=\prod_{j=1}^h I_j$ of axis-parallel cubes satisfying for $1\leq j\leq h$
$$
\ell(I_j)\leq\ell(E_j),\qquad (I_j)_*\cap E_j\ne\emptyset .
$$
For $a_{\mathbf I}\geq0$ and $2\leq p<\infty$, we have
$$
\int_{\mathbb R^d}\bigg(\sum_{\mathbf I\in\mathcal I}a_{\mathbf I}^2\mathbf1_{\mathbf I}(x)\bigg)^{p/2}dx
\leq C_d\int_E\bigg(\sum_{\mathbf I\in\mathcal I}a_{\mathbf I}^2\mathbf1_{\mathbf I}(x)\bigg)^{p/2}dx,
$$
where $C_d$ depends only on $d$.
\end{lemma}

\begin{proof}
Fix a coordinate interval $[a,b]$ of $E$. On its concentric triple, define
$$
T(x)=
\begin{cases}
a+(a-x)/64,&x<a,\\
x,&a\leq x\leq b,\\
b-(x-b)/64,&x>b.
\end{cases}
$$
Let $I=[u,v]$ be the corresponding interval of $\mathbf I$, with $r=v-u\leq b-a$.
By hypothesis, $[u+r/48,v-r/48]$ meets $[a,b]$. If $x\in I$ and $x<a$, then $a-x\leq r$ and
$$
u\leq x<a<T(x)\leq a+r/64<a+r/48\leq\min\{v,b\}.
$$
Thus $T(x)\in I\cap[a,b]$. The case $x>b$ follows by reflection, and $T(x)=x$ on $[a,b]$.

Applying these maps in the $d$ scalar coordinates, we obtain a single map $T$ on $3E=\prod_j3E_j$.
Every $\mathbf I\in\mathcal I$ lies in $3E$, and $T(\mathbf I)\subseteq\mathbf I\cap E$ up to boundaries. Hence
$$
\mathbf1_{\mathbf I}(x)\leq\mathbf1_{\mathbf I}(Tx),
\qquad\text{a.e. }\   x\in3E.
$$
In each coordinate the three inverse branches have absolute slopes $64,1,64$.
Changing variables on each piece of $3E$, we obtain
\begin{align*}
\int_{\mathbb R^d}\bigg(\sum_{\mathbf I\in\mathcal I}a_{\mathbf I}^2\mathbf1_{\mathbf I}(x)\bigg)^{p/2}dx
\leq &\int_{3E}\bigg(\sum_{\mathbf I\in\mathcal I}a_{\mathbf I}^2\mathbf1_{\mathbf I}(Tx)\bigg)^{p/2}dx\\
\leq &(64+1+64)^d\int_E\bigg(\sum_{\mathbf I\in\mathcal I}a_{\mathbf I}^2\mathbf1_{\mathbf I}(y)\bigg)^{p/2}dy.
\end{align*}
This proves the estimate with $C_d=129^d$.
\end{proof}

\begin{proposition}\label{p:local-product}
Fix $\kappa,\mathcal E$ and let $G$ be any finite sum of the atoms $A_{I,J}=A_{\kappa,I,J}^{\mathrm{bi}}(\mathcal E)$.
For every product $R_0=I_0\times J_0$ of axis-parallel cubes,
\begin{equation}\label{e:local-product}
\|Q_xQ_yG\|_{L^p(R_0)}\leq Cp\,\bigg\|\bigg(\sum_{I,J}c_{I,J}^2\mathbf1_{I\times J}\bigg)^{1/2}\bigg\|_{L^p(R_0)},\qquad 2\leq p<\infty.
\end{equation}
For the fixed reconstruction kernels, $C$ depends only on $n,m$.
\end{proposition}

\begin{proof}
Only pairs with $I_*\cap I_0\ne\emptyset$ and $J_*\cap J_0\ne\emptyset$ contribute to $Q_xQ_yG$ on $R_0$.

For $i,j\in\{0,1\}$, let $\mathcal F^{ij}$ consist of the contributing pairs with
$i=\mathbf1_{\{\ell(I)>\ell(I_0)\}}$ and $j=\mathbf1_{\{\ell(J)>\ell(J_0)\}}$.
Write $G_{ij}=\sum_{(I,J)\in\mathcal F^{ij}}A_{I,J}$.

\smallskip
\noindent{\it Case 1. $\ell(I)\leq\ell(I_0)$ and $\ell(J)\leq\ell(J_0)$.}
By \eqref{e:product-lp} and Lemma~\ref{l:cube-fold},
$$
\|Q_xQ_yG_{00}\|_{L^p(R_0)}\leq4\|G_{00}\|_p
\leq Cp\,\bigg\|\bigg(\sum_{(I,J)\in\mathcal F^{00}}c_{I,J}^2\mathbf1_{I\times J}\bigg)^{1/2}\bigg\|_p
\leq Cp\,\bigg\|\bigg(\sum_{(I,J)\in\mathcal F^{00}}c_{I,J}^2\mathbf1_{I\times J}\bigg)^{1/2}\bigg\|_{L^p(R_0)}.
$$

\smallskip
 
\noindent{\it Case 2. $\ell(I)>\ell(I_0)$ and $\ell(J)\leq\ell(J_0)$,
or $\ell(I)\leq\ell(I_0)$ and $\ell(J)>\ell(J_0)$.}
We begin with $G_{10}$. For $x,u\in I_0$ and $|\beta|\leq1$,
the mean value theorem and \eqref{e:bi-deriv} imply
$$
|\partial_y^\beta A_{I,J}(x,y)-\partial_y^\beta A_{I,J}(u,y)|
\leq C\frac{\ell(I_0)}{\ell(I)}c_{I,J}\ell(J)^{-|\beta|}
\big(\mathbf1_I(x)+\mathbf1_I(u)\big)\mathbf1_J(y),
$$
since $|x-u|\leq\sqrt n\,\ell(I_0)$ and the difference vanishes when $x,u\notin I$.
Averaging over $u\in I_0$, we obtain
\begin{equation}\label{e:horizontal-mean}
|\partial_y^\beta Q_xA_{I,J}(x,y)|
\leq C\frac{\ell(I_0)}{\ell(I)}c_{I,J}\ell(J)^{-|\beta|}
\bigg(\mathbf1_I(x)+\frac{|I\cap I_0|}{|I_0|}\bigg)\mathbf1_J(y),\qquad x\in I_0.
\end{equation}

For fixed $J$, let $W_J=\sum_{I:(I,J)\in\mathcal F^{10}}Q_xA_{I,J}$ and, on $I_0$, set
$$
d_J=(\operatorname{Id}+P_x)\sum_{I:(I,J)\in\mathcal F^{10}}
\frac{\ell(I_0)}{\ell(I)}c_{I,J}\mathbf1_I.
$$
For each $x\in I_0$, $W_J(x,\cdot)$ is supported in $J$ and has integral zero. By \eqref{e:horizontal-mean},
\begin{equation}\label{e:centred-row-deriv}
|\partial_y^\beta W_J(x,y)|\leq C\ell(J)^{-|\beta|}d_J(x)\mathbf1_J(y),\qquad |\beta|\leq1.
\end{equation}
At most one $I$ of each generation contains $x$, so
$$
\sum_{\ell(I)>\ell(I_0)}\bigg(\frac{\ell(I_0)}{\ell(I)}\bigg)^2\mathbf1_I(x)
\leq\sum_{k\geq0}2^{-2k}.
$$
By Minkowski's inequality in $\ell^2(J)$ and Cauchy--Schwarz in $I$, we have on $I_0\times\mathbb R^m$
\begin{align}
\bigg(\sum_Jd_J^2\mathbf1_J\bigg)^{1/2}
\leq &(\operatorname{Id}+P_x)
\bigg(\sum_J\bigg(\sum_{I:(I,J)\in\mathcal F^{10}}
\frac{\ell(I_0)}{\ell(I)}c_{I,J}\mathbf1_I\bigg)^2\mathbf1_J\bigg)^{1/2}\notag\\
\leq &C(\operatorname{Id}+P_x)
\bigg(\sum_{(I,J)\in\mathcal F^{10}}c_{I,J}^2\mathbf1_{I\times J}\bigg)^{1/2}.
\label{e:centred-coefficients}
\end{align}
For each $x\in I_0$, we may now apply \eqref{e:vertical-fixed} to $W_J(x,\cdot)$
and Lemma~\ref{l:cube-fold} to the cubes $J$ meeting $J_0$ in $J_*$.
Integrating in $x$ and using \eqref{e:centred-coefficients}, we obtain
\begin{align*}
\|Q_xQ_yG_{10}\|_{L^p(R_0)}
\leq &2\bigg\|\sum_JW_J\bigg\|_{L^p(I_0\times\mathbb R^m)}\\
\leq &C\sqrt p\,\bigg\|\bigg(\sum_Jd_J^2\mathbf1_J\bigg)^{1/2}\bigg\|_{L^p(I_0\times\mathbb R^m)}
\leq C\sqrt p\,\bigg\|\bigg(\sum_Jd_J^2\mathbf1_J\bigg)^{1/2}\bigg\|_{L^p(R_0)}\\
\leq &C\sqrt p\,\bigg\|\bigg(\sum_{(I,J)\in\mathcal F^{10}}c_{I,J}^2\mathbf1_{I\times J}\bigg)^{1/2}\bigg\|_{L^p(R_0)}.
\end{align*}

For $G_{01}$, we centre first in $y$, with
\begin{equation}\label{e:vertical-weights}
(\operatorname{Id}+P_y)\mathbf1_J(y)=\mathbf1_J(y)+\frac{|J\cap J_0|}{|J_0|},\qquad y\in J_0.
\end{equation}
Since $Q_yA_{I,J}$ has integral zero in $x$, the same argument with $x$ and $y$ interchanged yields
$$
\|Q_xQ_yG_{01}\|_{L^p(R_0)}\leq C\sqrt p\,\bigg\|\bigg(\sum_{(I,J)\in\mathcal F^{01}}c_{I,J}^2\mathbf1_{I\times J}\bigg)^{1/2}\bigg\|_{L^p(R_0)}.
$$

\smallskip
\noindent{\it Case 3. $\ell(I)>\ell(I_0)$ and $\ell(J)>\ell(J_0)$.}
For $x,u\in I_0$ and $y,v\in J_0$, the mixed difference vanishes if $x,u\notin I$ or $y,v\notin J$.
Using \eqref{e:bi-deriv} for one derivative in each variable, we obtain on $R_0$
\begin{align*}
|Q_xQ_yA_{I,J}(x,y)|
\leq &\frac1{|R_0|}\int_{I_0}\int_{J_0}
|A_{I,J}(x,y)-A_{I,J}(u,y)-A_{I,J}(x,v)+A_{I,J}(u,v)|\,dv\,du\\
\leq &C\frac{\ell(I_0)}{\ell(I)}\frac{\ell(J_0)}{\ell(J)}c_{I,J}
\bigg(\mathbf1_I(x)+\frac{|I\cap I_0|}{|I_0|}\bigg)
\bigg(\mathbf1_J(y)+\frac{|J\cap J_0|}{|J_0|}\bigg).
\end{align*}
At each pair of generations, at most one $I\times J$ contains $(x,y)$. Hence
$$
\sum_{(I,J)\in\mathcal F^{11}}
\bigg(\frac{\ell(I_0)\ell(J_0)}{\ell(I)\ell(J)}\bigg)^2\mathbf1_{I\times J}(x,y)
\leq\bigg(\sum_{k\geq0}2^{-2k}\bigg)^2.
$$
Summing the bounds for $Q_xQ_yA_{I,J}$ and applying Cauchy--Schwarz before taking the averages, we obtain
\begin{align*}
|Q_xQ_yG_{11}|
\leq &C(\operatorname{Id}+P_x)(\operatorname{Id}+P_y)
\sum_{(I,J)\in\mathcal F^{11}}\frac{\ell(I_0)\ell(J_0)}{\ell(I)\ell(J)}c_{I,J}\mathbf1_{I\times J}\\
\leq &C(\operatorname{Id}+P_x)(\operatorname{Id}+P_y)
\bigg(\sum_{(I,J)\in\mathcal F^{11}}c_{I,J}^2\mathbf1_{I\times J}\bigg)^{1/2}.
\end{align*}
Since $P_x,P_y$ are contractions on $L^p(R_0)$,
$$
\|Q_xQ_yG_{11}\|_{L^p(R_0)}
\leq C\bigg\|\bigg(\sum_{(I,J)\in\mathcal F^{11}}c_{I,J}^2\mathbf1_{I\times J}\bigg)^{1/2}\bigg\|_{L^p(R_0)}.
$$
Adding the four estimates, we conclude that
$$
\|Q_xQ_yG\|_{L^p(R_0)}
\leq C(p+2\sqrt p+1)\bigg\|\bigg(\sum_{I,J}c_{I,J}^2\mathbf1_{I\times J}\bigg)^{1/2}\bigg\|_{L^p(R_0)},
$$
which proves \eqref{e:local-product} for $p\geq2$.
\end{proof}

\subsection{Local estimates with vertical cancellation}\label{s:local-vertical}

\begin{proposition}\label{p:local-vertical}
Let $G=\sum_J B_J$ be either the remainder sum or the cubic sum in Proposition~\ref{p:vertical-square},
with coefficients $c_J$. For every product $R_0=I_0\times J_0$ of axis-parallel cubes,
\begin{equation}\label{e:local-vertical}
\|Q_xQ_yG\|_{L^p(R_0)}
\leq C\sqrt p\,\bigg\|\bigg(\sum_J c_J(x)^2\mathbf1_J(y)\bigg)^{1/2}\bigg\|_{L^p(R_0)},\qquad 2\leq p<\infty.
\end{equation}
For the fixed reconstruction kernels, $C$ depends only on $n,m$.
\end{proposition}

\begin{proof}
Fix $x\in I_0$. Since $B_J(x,\cdot)$ is supported in $J_*$, only cubes $J$ with
$J_*\cap J_0\ne\emptyset$ contribute to $Q_yG$ on $J_0$.

\smallskip
\noindent{\it Case 1. $\ell(J)\leq\ell(J_0)$.}
By \eqref{e:vertical-fixed} and Lemma~\ref{l:cube-fold},
\begin{align*}
\bigg\|Q_y\sum_{\ell(J)\leq\ell(J_0)}B_J(x,\cdot)\bigg\|_{L^p(J_0)}
\leq &2\bigg\|\sum_{\ell(J)\leq\ell(J_0)}B_J(x,\cdot)\bigg\|_{L^p(\mathbb R^m)}\\
\leq &C\sqrt p\,\bigg\|\bigg(\sum_{\ell(J)\leq\ell(J_0)}c_J(x)^2\mathbf1_J\bigg)^{1/2}\bigg\|_{L^p(\mathbb R^m)}\\
\leq &C\sqrt p\,\bigg\|\bigg(\sum_J c_J(x)^2\mathbf1_J\bigg)^{1/2}\bigg\|_{L^p(J_0)}.
\end{align*}

\smallskip
\noindent{\it Case 2. $\ell(J)>\ell(J_0)$.}
Using the first derivative bound in \eqref{e:vertical-size} and averaging over $v\in J_0$ as in \eqref{e:horizontal-mean}, we obtain for $y\in J_0$
$$
|Q_yB_J(x,y)|\leq C\frac{\ell(J_0)}{\ell(J)}c_J(x)
\bigg(\mathbf1_J(y)+\frac{|J\cap J_0|}{|J_0|}\bigg).
$$
At each generation at most one $J$ contains $y$, so
$$
\sum_{\ell(J)>\ell(J_0)}\bigg(\frac{\ell(J_0)}{\ell(J)}\bigg)^2\mathbf1_J(y)
\leq\sum_{k\geq0}2^{-2k}.
$$
By \eqref{e:vertical-weights}, the contraction of $P_y$ on $L^p(J_0)$ and Cauchy--Schwarz in $J$,
\begin{align*}
\bigg\|Q_y\sum_{\ell(J)>\ell(J_0)}B_J(x,\cdot)\bigg\|_{L^p(J_0)}
\leq &2C\,\bigg\|\sum_{\ell(J)>\ell(J_0)}\frac{\ell(J_0)}{\ell(J)}c_J(x)\mathbf1_J\bigg\|_{L^p(J_0)} 
\leq  C\,\bigg\|\bigg(\sum_J c_J(x)^2\mathbf1_J\bigg)^{1/2}\bigg\|_{L^p(J_0)}.
\end{align*}

Adding the two bounds and taking the $L^p(I_0)$ norm in $x$, we obtain
$$
\|Q_xQ_yG\|_{L^p(R_0)}\leq2\|Q_yG\|_{L^p(R_0)}
\leq C\sqrt p\,\bigg\|\bigg(\sum_J c_J(x)^2\mathbf1_J(y)\bigg)^{1/2}\bigg\|_{L^p(R_0)},
$$
as desired.
\end{proof}

\subsection{Local \texorpdfstring{$L^p$}{Lp} estimates and exponential integrability}\label{s:local-exponential}

Suppose the coefficient square function is bounded by one on $R_0$.
For the centred sum of the $A^{\mathrm{bi}}$ terms, the factor $p$ in \eqref{e:local-product}
implies an exponential tail. For the centred remainder and cubic sums, the factor $\sqrt p$
in \eqref{e:local-vertical} implies Gaussian tails.
We now combine the three estimates to control the doubly centred oscillation $Q_xQ_yf$
by the flag area function on the same rectangle.

\begin{theorem}\label{t:two-local}
Let $f\in L^2(\mathbb R^{n+m})$ and let $R_0=I_0\times J_0$ be a product of axis-parallel cubes.
Then
\begin{equation}\label{e:two-step-moments}
\|Q_xQ_yf\|_{L^p(R_0)}\leq Cp\|S_F(f)\|_{L^p(R_0)},\qquad 2\leq p<\infty.
\end{equation}
If $0<\|S_F(f)\|_{L^\infty(R_0)}<\infty$, then
\begin{equation}\label{e:two-step-exp}
\frac1{|R_0|}\int_{R_0}\exp\bigg(c\frac{|Q_xQ_yf(z)|}{\|S_F(f)\|_{L^\infty(R_0)}}\bigg)\,dz\leq2.
\end{equation}
In particular,
$$
\big|\{z\in R_0:|Q_xQ_yf(z)|>\tau\}\big|
\leq2|R_0|\exp\bigg(-\frac{c\tau}{\|S_F(f)\|_{L^\infty(R_0)}}\bigg),\qquad \tau>0.
$$
The constants $C,c>0$ depend only on $n,m$.
\end{theorem}

\begin{proof}
For each $\kappa$, we apply \eqref{e:local-product} to the sum of the $A^{\mathrm{bi}}$ terms
and \eqref{e:local-vertical} to the remainder and cubic sums. Summing over $\kappa$ and using
\eqref{e:coefficient-control}, we obtain
$$
\|Q_xQ_yf_{\mathcal E}\|_{L^p(R_0)}
\leq Cp\,\bigg\|\bigg(\sum_{\gamma\in\mathcal E}\frac{\mathbf1_{R_\gamma}}{|R_\gamma|}\int_{\mathcal T_\gamma}|H_f^\sigma|^2\,d\nu\bigg)^{1/2}\bigg\|_{L^p(R_0)}
\leq Cp\|S_F(f)\|_{L^p(R_0)}.
$$
Choose finite sets $\mathcal E_N\uparrow\mathcal I_{\rm par}$.
Since $Q_x,Q_y$ are orthogonal projections on $L^2(R_0)$, by \eqref{e:finite-recon-limit},
$$
\|Q_xQ_y(f_{\mathcal E_N}-f)\|_{L^2(R_0)}
\leq\|f_{\mathcal E_N}-f\|_{L^2(\mathbb R^{n+m})}\to0.
$$
Passing to an almost everywhere convergent subsequence on $R_0$ and applying Fatou's lemma, we obtain \eqref{e:two-step-moments}.

Suppose $0<\|S_F(f)\|_{L^\infty(R_0)}<\infty$. For every integer $j\geq2$, \eqref{e:two-step-moments} implies
$$
\frac1{|R_0|}\int_{R_0}\bigg(\frac{|Q_xQ_yf(z)|}{\|S_F(f)\|_{L^\infty(R_0)}}\bigg)^j\,dz
\leq(C_{n,m}j)^j.
$$
By Cauchy--Schwarz, the bound for $j=2$ also yields the case $j=1$ after increasing $C_{n,m}$.
Choose $c>0$ so that $ceC_{n,m}\leq1/2$. By Tonelli's theorem and $j!\geq(j/e)^j$,
\begin{align}
\frac1{|R_0|}\int_{R_0}\exp\bigg(c\frac{|Q_xQ_yf(z)|}{\|S_F(f)\|_{L^\infty(R_0)}}\bigg)\,dz
= &1+\sum_{j=1}^\infty\frac{c^j}{j!|R_0|}\int_{R_0}\bigg(\frac{|Q_xQ_yf(z)|}{\|S_F(f)\|_{L^\infty(R_0)}}\bigg)^j\,dz\notag\\
\leq &1+\sum_{j=1}^\infty\frac{(cC_{n,m}j)^j}{j!}
\leq1+\sum_{j=1}^\infty2^{-j}=2.\label{e:exponential-series}
\end{align}
This proves \eqref{e:two-step-exp}. The tail estimate follows by Markov's inequality.
\end{proof}

\subsection{Global area bounds and local estimates for \texorpdfstring{$L^u$}{Lu} data}\label{s:lp-extension}

To prove the global $L^p$ upper bound for $S_F$, $2\leq p<\infty$, we use a weighted
estimate for the classical Poisson area function in $\mathbb R^{n+m}$ and in the $y$ variable.
We then extend the local reverse estimate to $f\in L^u(\mathbb R^{n+m})$, $1\leq u<\infty$,
and obtain the global reverse estimate by letting the rectangles increase to $\mathbb R^{n+m}$.

Let $M_q$ be the Hardy--Littlewood maximal operator on $\mathbb R^q$,
$$
M_qw(z)=\sup_{Q\ni z}\frac1{|Q|}\int_Q|w(u)|\,du,
$$
where the supremum is over axis-parallel cubes. In the $y$ variable, write
$M_yw(x,y)=M_m(w(x,\cdot))(y)$.

Let $P_r^{[q]}g=P_r^{(q)}*g$ and write
$$
D_{r,\nu}^{[q]}g=
\begin{cases}
r\partial_r P_r^{[q]}g, & \nu=0,\\[2pt]
r\partial_{x_\nu}P_r^{[q]}g, & 1\leq\nu\leq q.
\end{cases}
$$
The classical Poisson area function uses all $q+1$ derivatives. With the change of variables $y=x-ra$, we have
\begin{align*}
S_{\mathrm{cl}}^{(q)}g(x)^2
= &\sum_{\nu=0}^q\int_0^\infty\int_{|x-y|<r}|D_{r,\nu}^{[q]}g(y)|^2\,\frac{dy\,dr}{r^{q+1}} 
=  \sum_{\nu=0}^q\int_0^\infty\int_{|a|<1}|D_{r,\nu}^{[q]}g(x-ra)|^2\,da\,\frac{dr}{r}.
\end{align*}
For the coordinate averages, the area function $\mathcal S_P^{(q)}$ in \eqref{e:intro-one-area}
uses the second derivative $r^2(-\Delta)P_r^{[q]}$.

For $g\in\mathcal S(\mathbb R^q)$ and bounded $w\geq0$, the weighted estimate of
Chanillo and Wheeden~\cite[Theorem~1(a)]{ChanilloWheeden} at exponent $2$ is
\begin{equation}\label{e:one-weighted}
\int_{\mathbb R^q}|S_{\mathrm{cl}}^{(q)}g|^2w\,dx
\leq C_q\int_{\mathbb R^q}|g|^2M_qw\,dx.
\end{equation}
For $g\in L^2(\mathbb R^q)$, choose $g_k\in\mathcal S(\mathbb R^q)$ with $g_k\to g$ in $L^2(\mathbb R^q)$.
At each $r>0$, Cauchy--Schwarz and the scaling of the Poisson kernel imply
$$
\|D_{r,\nu}^{[q]}(g_k-g)\|_\infty\leq C_qr^{-q/2}\|g_k-g\|_2\to0.
$$
By Fatou's lemma and \eqref{e:one-weighted},
$$
\int_{\mathbb R^q}|S_{\mathrm{cl}}^{(q)}g|^2w\,dx
\leq\liminf_{k\to\infty}\int_{\mathbb R^q}|S_{\mathrm{cl}}^{(q)}g_k|^2w\,dx
\leq C_q\lim_{k\to\infty}\int_{\mathbb R^q}|g_k|^2M_qw\,dx
=C_q\int_{\mathbb R^q}|g|^2M_qw\,dx,
$$
where the last equality follows from $M_qw\leq\|w\|_\infty$ and $g_k\to g$ in $L^2(\mathbb R^q)$.

\begin{proposition}\label{p:two-step-upper}
Let $d=n+m$. For $f\in L^2(\mathbb R^d)$ and every nonnegative locally integrable weight $w$,
\begin{equation}\label{e:two-step-weighted}
\int_{\mathbb R^d}S_F(f)^2w\,dz\leq C_{n,m}\int_{\mathbb R^d}|f|^2M_yM_dw\,dz.
\end{equation}
The estimate also holds with $M_dM_yw$ on the right. In particular,
\begin{equation}\label{e:two-step-upper}
\|S_F(f)\|_p\leq C_{n,m}p\|f\|_p,\qquad f\in L^p(\mathbb R^d),\quad 2\leq p<\infty.
\end{equation}
\end{proposition}

\begin{proof}
Write $z=(x,y)$ and let $D_{s,\nu_2}^y$ denote $D_{s,\nu_2}^{[m]}$ acting in $y$.
In \eqref{e:poisson-area-weight}, write the two displacements as $ta_1$ and $(0,sa_2)$,
with $a_1\in B_1$ and $a_2\in B_2$. The Jacobians cancel the factors $t^{-d}$ and $s^{-m}$.
Integrating first in $a_1,t$, we obtain
\begin{align*}
S_F(f)(z)^2
= &\sum_{\nu_1=0}^d\sum_{\nu_2=0}^m\int_0^\infty\int_0^\infty\int_{B_1}\int_{B_2}
|D_{t,\nu_1}^{[d]}D_{s,\nu_2}^yf(z-ta_1-(0,sa_2))|^2\,da_2\,da_1\,\frac{dt}{t}\frac{ds}{s}\\
= &\sum_{\nu_2=0}^m\int_0^\infty\int_{B_2}
|S_{\mathrm{cl}}^{(d)}(D_{s,\nu_2}^yf)(x,y-sa_2)|^2\,da_2\,\frac{ds}{s}.
\end{align*}
Let $S_{\mathrm{cl},y}f(x,y)=S_{\mathrm{cl}}^{(m)}(f(x,\cdot))(y)$.
Since the Poisson derivatives commute, integration in $a_2,s$ first yields
$$
S_F(f)(z)^2
=\sum_{\nu_1=0}^d\int_0^\infty\int_{B_1}
|S_{\mathrm{cl},y}(D_{t,\nu_1}^{[d]}f)(z-ta_1)|^2\,da_1\,\frac{dt}{t}.
$$

Suppose first that $w$ is bounded. Applying \eqref{e:one-weighted} in $\mathbb R^d$ and then in $y$,
and using translation invariance of $M_d$, we obtain
\begin{align*}
\int_{\mathbb R^d}S_F(f)^2w\,dz
= &\sum_{\nu_2=0}^m\int_0^\infty\int_{B_2}\int_{\mathbb R^d}
|S_{\mathrm{cl}}^{(d)}(D_{s,\nu_2}^yf)(z)|^2w(z+(0,sa_2))\,dz\,da_2\,\frac{ds}{s}\\
\leq &C_d\sum_{\nu_2=0}^m\int_0^\infty\int_{B_2}\int_{\mathbb R^d}
|D_{s,\nu_2}^yf(z)|^2(M_dw)(z+(0,sa_2))\,dz\,da_2\,\frac{ds}{s}\\
= &C_d\int_{\mathbb R^d}|S_{\mathrm{cl},y}f(z)|^2M_dw(z)\,dz
\leq C_{n,m}\int_{\mathbb R^d}|f|^2M_yM_dw\,dz.
\end{align*}
Using the second formula for $S_F(f)^2$ and applying \eqref{e:one-weighted} in the opposite order, we have
$$
\int_{\mathbb R^d}S_F(f)^2w\,dz
\leq C_m\int_{\mathbb R^d}|S_{\mathrm{cl}}^{(d)}f(z)|^2M_yw(z)\,dz
\leq C_{n,m}\int_{\mathbb R^d}|f|^2M_dM_yw\,dz.
$$
For general $w$, apply both estimates to $\min\{w,N\}$ and let $N\to\infty$,
using monotone convergence on the left and monotonicity of the maximal operators on the right.

For $p=2$, take $w=1$ in \eqref{e:two-step-weighted}. For $p>2$, let $q=(p/2)'$ and
$f\in L^2(\mathbb R^d)\cap L^p(\mathbb R^d)$. Both maximal operators have $L^q$ norm at most
$C_{n,m}q'$, where $q'=p/2$. By \eqref{e:two-step-weighted} and H\"older's inequality,
\begin{align*}
\int_{\mathbb R^d}S_F(f)^2w\,dz
\leq &C_{n,m}\int_{\mathbb R^d}|f|^2M_yM_dw\,dz
\leq C_{n,m}\|f\|_p^2\|M_yM_dw\|_q\\
\leq &C_{n,m}(q')^2\|f\|_p^2\|w\|_q
\leq C_{n,m}p^2\|f\|_p^2\|w\|_q.
\end{align*}
Duality in $L^{p/2}$ now gives $\|S_F(f)\|_p^2\leq C_{n,m}p^2\|f\|_p^2$.

For $f\in L^p(\mathbb R^d)$, choose $f_k\in C_c^\infty(\mathbb R^d)$ with $f_k\to f$ in $L^p(\mathbb R^d)$.
At fixed $t,s>0$, the $L^{p'}$ bound for the kernel of $D_{t,\nu_1}^{[d]}$ and the $L^1$ bound
for the kernel of $D_{s,\nu_2}^y$ imply
$$
\|D_{t,\nu_1}^{[d]}D_{s,\nu_2}^y(f_k-f)\|_\infty
\leq C_{d,p}t^{-d/p}\|D_{s,\nu_2}^y(f_k-f)\|_p
\leq C_{n,m,p}t^{-d/p}\|f_k-f\|_p\to0.
$$
Applying Fatou's lemma first in the area integral and then in $z$, we obtain
$$
\|S_F(f)\|_p\leq\liminf_{k\to\infty}\|S_F(f_k)\|_p\leq C_{n,m}p\|f\|_p.
$$
For any sequence $f_k\to f$ in $L^p(\mathbb R^d)$, the triangle inequality in the area integral now implies
\begin{equation}\label{e:two-step-cont}
\|S_F(f_k)-S_F(f)\|_p\leq\|S_F(f_k-f)\|_p\leq C_{n,m}p\|f_k-f\|_p\to0.
\end{equation}
This completes the proof of Proposition~\ref{p:two-step-upper}.
\end{proof}

\begin{corollary}\label{c:two-extension}
Let $f\in L^u(\mathbb R^d)$ for some $1\leq u<\infty$.
The local moment, exponential and tail estimates of Theorem~\ref{t:two-local} hold on every product $R=I\times J$ of axis-parallel cubes.
For $2\leq p<\infty$, $f\in L^p(\mathbb R^d)$ if and only if $S_F(f)\in L^p(\mathbb R^d)$, and in this case
$$
(Cp)^{-1}\|f\|_p\leq\|S_F(f)\|_p\leq Cp\|f\|_p.
$$
The constants depend only on $n,m$.
\end{corollary}

\begin{proof}
For $2\leq p\leq u$, choose $f_i\in C_c^\infty(\mathbb R^d)$ with $f_i\to f$ in $L^u(\mathbb R^d)$.
By H\"older's inequality on $R$ and \eqref{e:two-step-cont},
$$
\|Q_R(f_i-f)\|_{L^p(R)}+\|S_F(f_i)-S_F(f)\|_{L^p(R)}
\leq C_{n,m}u|R|^{1/p-1/u}\|f_i-f\|_u\to0.
$$
Passing to the limit in \eqref{e:two-step-moments} proves the local estimate in this range.

For $p>u$, suppose $S_F(f)\in L^p(R)$. Choose concentric products $R_i\uparrow R^\circ$
and nonnegative smooth mollifiers $\rho_i$ with supports shrinking to the origin and $R_i-\operatorname{supp}\rho_i\subset R$.
By Young's inequality, $g_i=\rho_i*f\in L^u\cap L^\infty\subset L^p$, so the case $u=p$ applies to $g_i$ on $R_i$.
Since convolution commutes with the Poisson derivatives, Minkowski's inequality in the area integral yields
$$
\|Q_{R_i}g_i\|_{L^p(R_i)}
\leq Cp\|S_F(g_i)\|_{L^p(R_i)}
\leq Cp\|\rho_i*(\mathbf1_R S_F(f))\|_{L^p(R_i)}
\leq Cp\|S_F(f)\|_{L^p(R)}.
$$
Extend the centred functions by zero to $R$. By Fubini's theorem and $f\in L^1(R)$, the coordinate averages over $R_i$
converge almost everywhere to those over $R$, so $Q_{R_i}f\to Q_Rf$ almost everywhere on $R$. Since
$$
\|Q_{R_i}(g_i-f)\|_{L^u(R)}\leq4\|g_i-f\|_u\to0,
$$
we may apply Fatou's lemma along a subsequence. Thus in both ranges,
\begin{equation}\label{e:local-lp-extension}
\|Q_Rf\|_{L^p(R)}\leq Cp\|S_F(f)\|_{L^p(R)},\qquad 2\leq p<\infty.
\end{equation}
For $0<\|S_F(f)\|_{L^\infty(R)}<\infty$, summing the moments as in \eqref{e:exponential-series}
yields the exponential estimate. The tail bound follows by Markov's inequality.

For the global estimate, take $R_L=I_L\times J_L$, where $I_L=[-L,L]^n$ and $J_L=[-L,L]^m$.
For almost every $(x,y)$ and all sufficiently large $L$, H\"older's inequality gives
$$
|Q_{R_L}f(x,y)-f(x,y)|
\leq\frac{\|f(\cdot,y)\|_{L^u(\mathbb R^n)}}{|I_L|^{1/u}}
+\frac{\|f(x,\cdot)\|_{L^u(\mathbb R^m)}}{|J_L|^{1/u}}
+\frac{\|f\|_u}{|R_L|^{1/u}}\to0.
$$
By Fatou's lemma and \eqref{e:local-lp-extension},
$$
\|f\|_p\leq\liminf_{L\to\infty}\|Q_{R_L}f\|_{L^p(R_L)}
\leq Cp\|S_F(f)\|_p.
$$
Combining this with \eqref{e:two-step-upper}, we complete the proof.
\end{proof}

In the next section we prove the corresponding estimates for $k$-step flags.
For $k=2$ and $(d_1,d_2)=(n,m)$, Proposition~\ref{p:flag-area-sharp} shows that
the factor $p$ in \eqref{e:two-step-upper} is sharp and that
$S_F(\mathbf1_{[0,1]^d})\notin L^\infty(\mathbb R^d)$.


\section{Local oscillation and sharp square-function estimates for Euclidean flags}\label{s:finite-flags}

We extend the local reverse estimate of Section~\ref{s:local-hoeffding} to $k$-step flags,
with $p^{k/2}$ in place of $p$, and prove the sharp global $L^p$ bounds.
Recall that $\mathbb R^N=X_1\oplus\cdots\oplus X_k$, where $X_j=\mathbb R^{d_j}$.
The $i$th reconstruction kernel acts on $V_i=X_i\oplus\cdots\oplus X_k$,
so only the factors with $i\leq j$ contribute to the displacement in $X_j$.
At reconstruction scales $t_1,\ldots,t_k$, this displacement is at most
$$
\sum_{i=1}^j t_i\leq kT_j,\qquad T_j=\max_{1\leq i\leq j}t_i.
$$
Thus the compactly supported reconstruction kernel lies in a product of cubes of side lengths $2kT_j$.
Several dyadic scale choices can lead to the same enclosing rectangle.
We sum their contributions before estimating the local oscillation.

Let $D_j=\dim V_j$ and let $B_j$ be the unit ball in $V_j$. All cubes in $X_j$ are axis-parallel.
For the Poisson semigroup $P_{j,t}=e^{-t\sqrt{-\Delta_{V_j}}}$ on $V_j$, write
$$
D_{j,t,\nu}=
\begin{cases}
t\partial_tP_{j,t}, &\nu=0,\\[2pt]
t\partial_{V_j,\nu}P_{j,t}, &1\leq\nu\leq D_j,
\end{cases}
$$
where $\partial_{V_j,\nu}$ denotes differentiation in the $\nu$th coordinate of $V_j$.
For $\mathbf t=(t_1,\ldots,t_k)$, the Poisson extension
$U_f(\cdot,\mathbf t)=P_{1,t_1}\cdots P_{k,t_k}f$ has mixed derivatives
$$
\mathcal G_{\mathbf t}f
=\big(D_{1,t_1,\nu_1}\cdots D_{k,t_k,\nu_k}f\big)_{0\leq\nu_j\leq D_j}.
$$
We use the Euclidean norm on this array.

Let $K_{\mathbf t,\boldsymbol\nu}$ be the convolution kernel of the component indexed by
$\boldsymbol\nu=(\nu_1,\ldots,\nu_k)$. The first factor acts on $V_1=\mathbb R^N$.
Each remaining factor has uniformly bounded $L^1(V_j)$ norm and defines a finite measure on $\mathbb R^N$.
By scaling in the first factor and Young's inequality,
$$
\|K_{\mathbf t,\boldsymbol\nu}\|_{L^q(\mathbb R^N)}
\leq C_qt_1^{-N(1-1/q)},\qquad 1\leq q\leq\infty,
$$
where $C_q$ depends only on $q$ and the dimensions.
Taking $q=u'$ shows that every component of $\mathcal G_{\mathbf t}f$ is defined by an absolutely
convergent integral for $f\in L^u(\mathbb R^N)$, $1\leq u<\infty$.

Regard each $V_j$ as a coordinate subspace of $\mathbb R^N$ and set
$$
\chi_{\mathbf t}(x)
=\int_{t_2B_2}\cdots\int_{t_kB_k}
\mathbf1_{\{|x-v_2-\cdots-v_k|<t_1\}}\,dv_k\cdots dv_2.
$$
For $k=1$, set $\chi_{t_1}(x)=\mathbf1_{\{|x|<t_1\}}$. By Tonelli's theorem,
$$
\int_{\mathbb R^N}\chi_{\mathbf t}(x)\,dx=\prod_{j=1}^k|B_j|t_j^{D_j}.
$$
We use Lebesgue measure on each ball, as in \cite[Definitions~1.2.2 and~1.3.1]{HLLW} for $k=2$.
The flag area function is
$$
S_F^{(k)}(f)(x)^2
=\int_{(0,\infty)^k}\int_{\mathbb R^N}
\chi_{\mathbf t}(x-y)|\mathcal G_{\mathbf t}f(y)|^2\,dy
\prod_{j=1}^k\frac{dt_j}{t_j^{D_j+1}}.
$$
Writing $x-y=v_1+\cdots+v_k$ and substituting $v_j=t_ja_j$, we obtain
\begin{equation}\label{e:flag-area}
S_F^{(k)}(f)(x)^2
=\sum_{\nu_1=0}^{D_1}\cdots\sum_{\nu_k=0}^{D_k}
\int_{(0,\infty)^k}\int_{B_1\times\cdots\times B_k}
\bigg|\big(D_{1,t_1,\nu_1}\cdots D_{k,t_k,\nu_k}f\big)
\Big(x-\sum_{j=1}^kt_ja_j\Big)\bigg|^2
\prod_{j=1}^k da_j\frac{dt_j}{t_j}.
\end{equation}
For $k=1$, $S_F^{(1)}=S_{\mathrm{cl}}^{(N)}$. For $k=2$ and $(d_1,d_2)=(n,m)$, $S_F^{(2)}=S_F$.

Convolutions along the nested subspaces $V_1,\ldots,V_k$ also appear in \cite[Section~8.6]{NRSW1}
for compactly supported or Schwartz kernels. Here the Poisson derivatives have polynomial decay.
To apply the one-parameter estimates successively in \eqref{e:flag-area}, write
$$
T_{j,\omega_j}g(x)=(D_{j,t_j,\nu_j}g)(x-t_ja_j),\qquad \omega_j=(t_j,a_j,\nu_j),
$$
and define $\mu_j$ by
\begin{equation}\label{e:flag-parameter}
\int F(\omega_j)\,d\mu_j(\omega_j)
=\sum_{\nu_j=0}^{D_j}\int_0^\infty\int_{B_j}F(t_j,a_j,\nu_j)\,da_j\frac{dt_j}{t_j}.
\end{equation}
For fixed parameters, the $T_{j,\omega_j}$ are commuting convolution operators on $L^2(\mathbb R^N)$, and
$$
S_F^{(k)}(f)(x)^2
=\int\big|T_{1,\omega_1}\cdots T_{k,\omega_k}f(x)\big|^2\,d\mu_1(\omega_1)\cdots d\mu_k(\omega_k).
$$

For $\xi=(\xi_1,\ldots,\xi_k)$, write $\xi_{V_j}=(\xi_j,\ldots,\xi_k)$.
On $V_j$, the squared multipliers of the time and spatial derivatives sum to
$$
(t_j|\xi_{V_j}|)^2e^{-2t_j|\xi_{V_j}|}
+t_j^2\sum_{\nu=1}^{D_j}|\xi_{V_j,\nu}|^2e^{-2t_j|\xi_{V_j}|}
=2(t_j|\xi_{V_j}|)^2e^{-2t_j|\xi_{V_j}|}.
$$
For $\xi_k\ne0$, each $|\xi_{V_j}|$ is positive, so
$$
\prod_{j=1}^k\int_0^\infty 2(t_j|\xi_{V_j}|)^2e^{-2t_j|\xi_{V_j}|}\,\frac{dt_j}{t_j}
=\bigg(2\int_0^\infty ue^{-2u}\,du\bigg)^k=2^{-k}.
$$
The translations in \eqref{e:flag-area} preserve the $L^2$ norm, and integration over each $B_j$
contributes $|B_j|$. By Plancherel's identity, for every $f\in L^2(\mathbb R^N)$,
\begin{equation}\label{e:flag-plancherel}
\|S_F^{(k)}(f)\|_{L^2(\mathbb R^N)}^2
=2^{-k}\bigg(\prod_{j=1}^k|B_j|\bigg)\|f\|_{L^2(\mathbb R^N)}^2.
\end{equation}

\subsection{Weighted estimates and global \texorpdfstring{$L^p$}{Lp} bounds for the flag area function}\label{s:flag-lp}

\begin{theorem}\label{t:flag-area-bounds}
Let $M_j$ be the Hardy--Littlewood maximal operator on $V_j$, acting with $x_1,\ldots,x_{j-1}$ fixed.
For $f\in L^2(\mathbb R^N)$, every nonnegative locally integrable weight $w$ and every permutation $\pi$ of $\{1,\ldots,k\}$,
\begin{equation}\label{e:flag-weighted}
\int_{\mathbb R^N}S_F^{(k)}(f)^2w\,dx
\leq C\int_{\mathbb R^N}|f|^2M_{\pi(k)}\cdots M_{\pi(1)}w\,dx.
\end{equation}
In particular,
\begin{equation}\label{e:flag-area-lp}
\|S_F^{(k)}(f)\|_p\leq Cp^{k/2}\|f\|_p,\qquad f\in L^p(\mathbb R^N),\quad 2\leq p<\infty.
\end{equation}
The constants depend only on $k,d_1,\ldots,d_k$.
\end{theorem}

\begin{proof}
For $g\in L^2(\mathbb R^N)$ and almost every fixed $x_1,\ldots,x_{j-1}$,
the $L^2(d\mu_j)$ norm of $T_{j,\omega_j}g(x)$ is the classical Poisson area function
of $g(x_1,\ldots,x_{j-1},\cdot)$ at $(x_j,\ldots,x_k)$.
Applying \eqref{e:one-weighted} on $V_j$ and integrating in the remaining variables, we obtain for bounded $v\geq0$
$$
\int\int_{\mathbb R^N}|T_{j,\omega_j}g(x)|^2v(x)\,dx\,d\mu_j(\omega_j)
\leq C_{D_j}\int_{\mathbb R^N}|g|^2M_jv\,dx.
$$

Suppose first that $w$ is bounded. The operators $T_{j,\omega_j}$ commute and are bounded on $L^2(\mathbb R^N)$
for fixed $\omega_j$. We may therefore apply the preceding estimate first on $V_{\pi(1)}$,
then on $V_{\pi(2)}$, and so on. By Tonelli's theorem,
\begin{align*}
\int_{\mathbb R^N}S_F^{(k)}(f)^2w\,dx
= &\int\int_{\mathbb R^N}
\big|T_{\pi(1),\omega_{\pi(1)}}\cdots T_{\pi(k),\omega_{\pi(k)}}f(x)\big|^2
w(x)\,dx\prod_{i=1}^k d\mu_{\pi(i)}(\omega_{\pi(i)})\\
\leq &C_{D_{\pi(1)}}\int\int_{\mathbb R^N}
\big|T_{\pi(2),\omega_{\pi(2)}}\cdots T_{\pi(k),\omega_{\pi(k)}}f(x)\big|^2
(M_{\pi(1)}w)(x)\,dx\prod_{i=2}^k d\mu_{\pi(i)}(\omega_{\pi(i)})\\
\leq &\bigg(\prod_{j=1}^k C_{D_j}\bigg)
\int_{\mathbb R^N}|f|^2M_{\pi(k)}\cdots M_{\pi(1)}w\,dx,
\end{align*}
where each intermediate weight is bounded by $\|w\|_\infty$.
For general $w$, apply this bound to $\min\{w,L\}$ and let $L\to\infty$ as in Proposition~\ref{p:two-step-upper}.

For $p=2$, use \eqref{e:flag-plancherel}. For $p>2$, let $f\in L^2(\mathbb R^N)\cap L^p(\mathbb R^N)$
and $q=(p/2)'$. Since $\|M_j\|_{L^q\to L^q}\leq C_{D_j}q'$ and $q'=p/2$, we have for bounded, compactly supported $w\geq0$
$$
\int_{\mathbb R^N}S_F^{(k)}(f)^2w\,dx
\leq C\|f\|_p^2\|M_{\pi(k)}\cdots M_{\pi(1)}w\|_q
\leq C(q')^k\|f\|_p^2\|w\|_q
\leq Cp^k\|f\|_p^2\|w\|_q.
$$
Taking the supremum over such $w$ with $\|w\|_q\leq1$, we obtain
$\|S_F^{(k)}(f)\|_p^2\leq Cp^k\|f\|_p^2$.

At fixed scales, $K_{\mathbf t,\boldsymbol\nu}\in L^{p'}(\mathbb R^N)$, so the density argument
in Proposition~\ref{p:two-step-upper} applies to every $f\in L^p(\mathbb R^N)$.
For $f_i\to f$ in $L^p(\mathbb R^N)$, the triangle inequality in \eqref{e:flag-area} then yields
$$
\|S_F^{(k)}(f_i)-S_F^{(k)}(f)\|_p
\leq\|S_F^{(k)}(f_i-f)\|_p\leq Cp^{k/2}\|f_i-f\|_p\to0.
$$
This completes the proof.
\end{proof}

\subsection{Local moment bounds and exponential integrability on products of cubes}\label{s:flag-local}

Fix $R_0=E_1\times\cdots\times E_k$, where each $E_j\subset X_j$ is an axis-parallel cube.
Let $P_j$ denote averaging over $E_j$ in the $j$th coordinate block,
$$
P_jg(x)=\frac{\mathbf1_{E_j}(x_j)}{|E_j|}
\int_{E_j}g(x_1,\ldots,x_{j-1},u_j,x_{j+1},\ldots,x_k)\,du_j.
$$
With $Q_jg(x)=\mathbf1_{E_j}(x_j)g(x)-P_jg(x)$, we have $Q_{R_0}f=Q_1\cdots Q_kf$ on $R_0$.

\begin{theorem}\label{t:flag-local}
Let $k\geq1$ and let $f\in L^u(\mathbb R^N)$ for some $1\leq u<\infty$.
For every product $R_0$ of axis-parallel cubes,
\begin{equation}\label{e:flag-local-lp}
\|Q_{R_0}f\|_{L^p(R_0)}
\leq Cp^{k/2}\|S_F^{(k)}(f)\|_{L^p(R_0)},\qquad 2\leq p<\infty.
\end{equation}
If $0<\|S_F^{(k)}(f)\|_{L^\infty(R_0)}<\infty$, then
\begin{equation}\label{e:flag-exponential}
\frac1{|R_0|}\int_{R_0}
\exp\bigg[c\bigg(\frac{|Q_{R_0}f(x)|}{\|S_F^{(k)}(f)\|_{L^\infty(R_0)}}\bigg)^{2/k}\bigg]\,dx\leq2.
\end{equation}
The constants depend only on $k,d_1,\ldots,d_k$.
The order $p^{k/2}$ and the exponential power $2/k$ are sharp even on the unit cube.
\end{theorem}

We prove Theorem~\ref{t:flag-local} at the end of this subsection, after establishing
the local estimates for finite sums in Proposition~\ref{p:flag-local-sums}.
The reconstruction kernel has integral zero in $x_k$. As in the two-step case,
we decompose it according to cancellation in the remaining blocks.
For a term with cancellation in $\ell$ blocks, the one-parameter square-function
estimate contributes a factor $p^{\ell/2}$. In the other blocks, we sum over
the dyadic generations using decay in the ratios of successive side lengths.
The examples proving sharpness are constructed in Subsection~\ref{s:flag-sharpness}.

\subsubsection{Compact reconstruction from Poisson coefficients}

Let $t_1,\ldots,t_k$ be the reconstruction scales and $T_j=\max_{i\leq j}t_i$.
Fix $\sigma_k=\max\{\sigma_0,98k\sqrt{N\max_i d_i}\}$.
At the corresponding Poisson scales $\sigma_kt_1,\ldots,\sigma_kt_k$, write
$$
H_f(y,\mathbf t)=\bigg[\prod_{j=1}^k b_{\sigma_k}\big(t_j\sqrt{-\Delta_{V_j}}\big)f\bigg](y),
\qquad f\in L^2(\mathbb R^N).
$$
On the parameter space we use $d\nu_k(y,\mathbf t)=dy\prod_{j=1}^kdt_j/t_j$.

Let $\Phi_j=\varphi_{D_j}^{\sigma_k}$ be the kernel from Proposition~\ref{p:reconstruction}
and set $\Phi_{j,t}(v)=t^{-D_j}\Phi_j(v/t)$. For a kernel $g$ on $V_j$, write
$$
(F*_{V_j}g)(x)=\int_{V_j}F(x-v)g(v)\,dv.
$$
The reconstruction kernel and its Fourier transform are
$$
K_{\mathbf t}=\Phi_{1,t_1}*_{V_2}\Phi_{2,t_2}*_{V_3}\cdots*_{V_k}\Phi_{k,t_k},
\qquad
\widehat K_{\mathbf t}(\xi)=\prod_{j=1}^ka_{\sigma_k}(t_j|\xi_{V_j}|).
$$
For $k=1$, $K_{\mathbf t}=\Phi_{1,t_1}$. Each $\Phi_j$ is smooth, is supported in the closed unit ball of $V_j$
and has vanishing moments through order one. Since only the first $j$ factors act in $X_j$,
$$
\operatorname{supp}K_{\mathbf t}
\subseteq\prod_{j=1}^k\bigg\{x_j\in X_j:|x_j|_\infty\leq\sum_{i=1}^jt_i\bigg\}
\subseteq\prod_{j=1}^k\{x_j\in X_j:|x_j|_\infty\leq kT_j\}.
$$

Let $\widetilde K_{\mathbf t}(x)=\overline{K_{\mathbf t}(-x)}$.
For $\xi_k\ne0$, each $|\xi_{V_j}|$ is positive. By Plancherel and the change of variables
$u=t_j|\xi_{V_j}|$ in each scale integral, we have for $g\in L^2(\mathbb R^N)$
\begin{equation}\label{e:flag-recon-energy}
\int_{(0,\infty)^k}\|\widetilde K_{\mathbf t}*g\|_{L^2(\mathbb R^N)}^2\prod_{j=1}^k\frac{dt_j}{t_j}
=\|a_{\sigma_k}\|_{L^2(du/u)}^{2k}\|g\|_{L^2(\mathbb R^N)}^2.
\end{equation}
With $b_{\sigma_k}$ in place of $a_{\sigma_k}$, the same calculation yields
$$
\|H_f\|_{L^2(d\nu_k)}^2
=\bigg(\int_0^\infty u e^{-2\sigma_k u}\,du\bigg)^k\|f\|_{L^2(\mathbb R^N)}^2
=(4\sigma_k^2)^{-k}\|f\|_{L^2(\mathbb R^N)}^2.
$$

Let $\mu$ be a Borel function with $|\mu|\leq1$.
For every measurable $E\subset[L^{-1},L]^k$, $L>1$, duality, Cauchy--Schwarz and \eqref{e:flag-recon-energy} imply
$$
\bigg\|\int_E K_{\mathbf t}*(\mu H_f)(\cdot,\mathbf t)\prod_{j=1}^k\frac{dt_j}{t_j}\bigg\|_{L^2(\mathbb R^N)}
\leq\|a_{\sigma_k}\|_{L^2(du/u)}^k\|\mathbf1_E(\mathbf t)\mu H_f\|_{L^2(d\nu_k)}.
$$
Since $\mu H_f\in L^2(d\nu_k)$, applying this bound to differences of truncations shows that the limit
$$
f_\mu=\lim_{L\to\infty}\int_{[L^{-1},L]^k}K_{\mathbf t}*(\mu H_f)(\cdot,\mathbf t)\prod_{j=1}^k\frac{dt_j}{t_j}
\qquad\text{in }L^2(\mathbb R^N)
$$
exists and satisfies
\begin{equation}\label{e:flag-recon-l2}
\|f_\mu\|_{L^2(\mathbb R^N)}
\leq\|a_{\sigma_k}\|_{L^2(du/u)}^k\|\mu H_f\|_{L^2(d\nu_k)}.
\end{equation}

For $\mu=1$ and $\xi_k\ne0$, the multiplier of the truncated integral satisfies
$$
\prod_{j=1}^k\int_{L^{-1}}^L
a_{\sigma_k}(t_j|\xi_{V_j}|)b_{\sigma_k}(t_j|\xi_{V_j}|)\frac{dt_j}{t_j}
=\prod_{j=1}^k\int_{|\xi_{V_j}|/L}^{L|\xi_{V_j}|}
a_{\sigma_k}(u)b_{\sigma_k}(u)\,\frac{du}{u}\to1
$$
by \eqref{e:dual-normalization}. Each factor is uniformly bounded by \eqref{e:truncated-dual-bound}.
By Plancherel and dominated convergence, we obtain $f_1=f$ in $L^2(\mathbb R^N)$.

\subsubsection{Coordinate cancellation and decay in the scale ratios}

For $1\leq j<k$, let $\theta_j(u)=\prod_{\ell=1}^{d_j}\theta_0(u_\ell)$ and
$\theta_{j,T}(u)=T^{-d_j}\theta_j(u/T)$. Since $\int_{X_j}\theta_j=1$, the operator
$$
\mathsf E_jF(x)=\theta_{j,T_j}(x_j)
\int_{X_j}F(x_1,\ldots,x_{j-1},u_j,x_{j+1},\ldots,x_k)\,du_j
$$
preserves the integral in $x_j$, and $(I-\mathsf E_j)F$ has integral zero in $x_j$.
The projections $\mathsf E_j$ commute. In $X_k$, the kernel $K_{\mathbf t}$ already has
integral zero, since its Fourier transform contains the factor $a_{\sigma_k}(t_k|\xi_k|)$.

For $k\in\Gamma\subseteq\{1,\ldots,k\}$, define
\begin{equation}\label{e:flag-kernel-split}
K_{\mathbf t}^\Gamma
=\prod_{j\notin\Gamma}\mathsf E_j
\prod_{\substack{j\in\Gamma\\j<k}}(I-\mathsf E_j)K_{\mathbf t}.
\end{equation}
The two-step estimates illustrate the two scale ratios below. With $t_1=t$ and $t_2=s$,
the factor $(t_2/T_2)^2=(s/t)^2$ in \eqref{e:cubic-kernel} is summable over the dyadic scales $s\leq t$
at fixed $t$. For $t<s$, the remainder in \eqref{e:remainder-kernel} has $\Gamma=\{2\}$:
integration in $x_1$ supplies the factor $(T_1/T_2)^2=(t/s)^2$, summable over $t$ at fixed $s$.

\begin{proposition}\label{p:flag-kernel-bounds}
We have $K_{\mathbf t}=\sum_{\Gamma\ni k}K_{\mathbf t}^\Gamma$.
Each $K_{\mathbf t}^\Gamma$ is supported where $|x_j|_\infty\leq kT_j$ for $1\leq j\leq k$
and has integral zero in $x_j$ for every $j\in\Gamma$.
For $\alpha=(\alpha_1,\ldots,\alpha_k)$ with $\alpha_j\in\mathbb N_0^{d_j}$,
\begin{equation}\label{e:flag-kernel-deriv}
|\partial^\alpha K_{\mathbf t}^\Gamma(x)|
\leq C_{\alpha,k}\prod_{j=2}^k\Big(\frac{t_j}{T_j}\Big)^2
\prod_{j\notin\Gamma}\Big(\frac{T_j}{T_{j+1}}\Big)^2
\prod_{j=1}^kT_j^{-d_j-|\alpha_j|}.
\end{equation}
The constant also depends on the block dimensions. Empty products equal one.
\end{proposition}

\begin{proof}
Summing \eqref{e:flag-kernel-split} over $\Gamma\ni k$, we obtain
$$
\sum_{\Gamma\ni k}K_{\mathbf t}^\Gamma
=\prod_{j=1}^{k-1}\big[(I-\mathsf E_j)+\mathsf E_j\big]K_{\mathbf t}
=K_{\mathbf t}.
$$
Each $\mathsf E_j$ preserves the support bound, since $\theta_{j,T_j}$ is supported where
$|x_j|_\infty\leq T_j$. All these projections preserve cancellation in $x_k$.
For $j\in\Gamma$ with $j<k$, cancellation in $x_j$ follows from the factor $I-\mathsf E_j$.

To prove \eqref{e:flag-kernel-deriv}, first integrate out the coordinates outside $G$.
For $k\in G\subseteq\Gamma$, write $X_G=\prod_{\ell\in G}X_\ell$ and
$$
\mathsf I_GK_{\mathbf t}(x_G)=\int_{X_{G^c}}K_{\mathbf t}(x_G,u_{G^c})\,du_{G^c}.
$$
Expanding the factors $I-\mathsf E_j$ in \eqref{e:flag-kernel-split}, we have
\begin{equation}\label{e:flag-projection-expansion}
K_{\mathbf t}^\Gamma(x)
=\sum_{\substack{G\subseteq\Gamma\\k\in G}}(-1)^{|\Gamma\setminus G|}
\bigg(\prod_{j\notin G}\theta_{j,T_j}(x_j)\bigg)\mathsf I_GK_{\mathbf t}(x_G).
\end{equation}
Let $\nu_G(j)=\min\{\ell\in G:\ell\geq j\}$, which exists since $k\in G$,
and write $\partial_G^\alpha=\prod_{\ell\in G}\partial_{x_\ell}^{\alpha_\ell}$.
We first show that
\begin{equation}\label{e:flag-integral-bound}
|\partial_G^\alpha\mathsf I_GK_{\mathbf t}(x_G)|
\leq C_{\alpha,k}\prod_{j=1}^k\Big(\frac{t_j}{T_{\nu_G(j)}}\Big)^2
\prod_{\ell\in G}T_\ell^{-d_\ell-|\alpha_\ell|}.
\end{equation}

Integration in $X_{G^c}$ sets $\xi_\ell=0$ for $\ell\notin G$. Thus
$$
\widehat{\mathsf I_GK_{\mathbf t}}(\xi_G)
=\prod_{j=1}^k a_{\sigma_k}\bigg(t_j
\Big(\sum_{\substack{\ell\in G\\\ell\geq j}}|\xi_\ell|^2\Big)^{1/2}\bigg).
$$
Set $\eta_\ell=T_\ell\xi_\ell$ for $\ell\in G$. Since $T_\ell\geq T_{\nu_G(j)}$ for $\ell\in G$, $\ell\geq j$,
the argument of the $j$th multiplier satisfies
$$
u_j:=t_j\bigg(\sum_{\substack{\ell\in G\\\ell\geq j}}
\frac{|\eta_\ell|^2}{T_\ell^2}\bigg)^{1/2}
\leq\frac{t_j}{T_{\nu_G(j)}}|\eta|.
$$
Conversely, for each $\ell\in G$, choose $i\leq\ell$ with $t_i=T_\ell$.
Then $u_i\geq|\eta_\ell|$, so $\max_j u_j\geq |G|^{-1/2}|\eta|$.
By the definition of $a_{\sigma_k}$ in Subsection~\ref{s:compact-kernels},
$$
|a_{\sigma_k}(u)|=c_{\sigma_k}^{-1}u^2|\widehat\theta_0(u)|
\leq C_Mu^2(1+u)^{-M},\qquad M>0.
$$
Using the upper bounds for the $u_j$ in the quadratic factors and the lower bound for $\max_j u_j$ in the decay factor, we obtain
\begin{align*}
|\widehat{\mathsf I_GK_{\mathbf t}}(\xi_G)|
\leq &C_{M,k}\bigg(\prod_{j=1}^ku_j^2\bigg)(1+\max_j u_j)^{-M}\\
\leq &C_{M,k}\frac{|\eta|^{2k}}{(1+|\eta|)^M}
\prod_{j=1}^k\Big(\frac{t_j}{T_{\nu_G(j)}}\Big)^2.
\end{align*}
By Fourier inversion and the change of variables $\eta_\ell=T_\ell\xi_\ell$,
\begin{align*}
|\partial_G^\alpha\mathsf I_GK_{\mathbf t}(x_G)|
\leq &C\int_{X_G}\prod_{\ell\in G}|\xi_\ell|^{|\alpha_\ell|}
|\widehat{\mathsf I_GK_{\mathbf t}}(\xi_G)|\,d\xi_G\\
\leq &C_{M,k}\prod_{j=1}^k\Big(\frac{t_j}{T_{\nu_G(j)}}\Big)^2
\prod_{\ell\in G}T_\ell^{-d_\ell-|\alpha_\ell|}
\int_{X_G}\frac{|\eta|^{2k+\sum_{\ell\in G}|\alpha_\ell|}}{(1+|\eta|)^M}\,d\eta.
\end{align*}
Taking $M>N+2k+|\alpha|$, we obtain \eqref{e:flag-integral-bound}.

For $j\in G$, $\nu_G(j)=j$. For $j\notin G$, $T_{\nu_G(j)}\geq T_{j+1}$.
Since $t_1=T_1$ and $G\subseteq\Gamma$,
\begin{equation}\label{e:flag-scale-gain}
\prod_{j=1}^k\Big(\frac{t_j}{T_{\nu_G(j)}}\Big)^2
=\prod_{j=2}^k\Big(\frac{t_j}{T_j}\Big)^2
\prod_{j\notin G}\Big(\frac{T_j}{T_{\nu_G(j)}}\Big)^2
\leq\prod_{j=2}^k\Big(\frac{t_j}{T_j}\Big)^2
\prod_{j\notin\Gamma}\Big(\frac{T_j}{T_{j+1}}\Big)^2.
\end{equation}
In \eqref{e:flag-projection-expansion}, each derivative in a block $j\notin G$ falls on
$\theta_{j,T_j}$, with $|\partial^{\alpha_j}\theta_{j,T_j}|\leq C_{\alpha_j}T_j^{-d_j-|\alpha_j|}$.
Combining this bound with \eqref{e:flag-integral-bound} and \eqref{e:flag-scale-gain}, we obtain \eqref{e:flag-kernel-deriv}.
\end{proof}

For $k=3$, all four bounds in \eqref{e:flag-kernel-deriv} contain the factor $(t_2/T_2)^2(t_3/T_3)^2$.
The middle column shows the additional decay from integration in the blocks outside $\Gamma$.
Each block in $\Gamma$ contributes a factor $\sqrt p$ in Proposition~\ref{p:flag-local-sums}, as recorded in the last column.
$$
\begin{array}{c|c|c}
\Gamma & \prod_{j\notin\Gamma}(T_j/T_{j+1})^2 & p^{|\Gamma|/2}\\ \hline
\{1,2,3\} & 1             & p^{3/2}\\[1mm]
\{1,3\}   & (T_2/T_3)^2   & p\\[1mm]
\{2,3\}   & (T_1/T_2)^2   & p\\[1mm]
\{3\}     & (T_1/T_3)^2   & \sqrt p
\end{array}
$$


\subsubsection{Parameter regions and smooth sums on dyadic rectangles}

For each $(y,\mathbf t)$, choose the cube and shift in $X_j$ by \eqref{e:dyadic-choice} at scale $kT_j$,
$$
I_j=Q_{d_j}(y_j,kT_j),\qquad \kappa_j=\kappa_{d_j}(y_j,kT_j),\qquad 1\leq j\leq k.
$$
Then $\overline{B_\infty(y_j,4kT_j)}\subset I_j^\circ$ and $24kT_j<\ell(I_j)\leq48kT_j$.
Since $T_1\leq\cdots\leq T_k$ and the same side-length rule applies in every block,
$\ell(I_1)\leq\cdots\leq\ell(I_k)$.

Write $\kappa=(\kappa_1,\ldots,\kappa_k)$ and $R=I_1\times\cdots\times I_k$.
For $\mathbf h=(h_2,\ldots,h_k)\in\mathbb N_0^{k-1}$ and $\gamma=(\kappa,R,\mathbf h)$, define
$$
\mathcal T_\gamma=
\big\{(y,\mathbf t):
Q_{d_j}(y_j,kT_j)=I_j,\ \kappa_{d_j}(y_j,kT_j)=\kappa_j\ (1\leq j\leq k),\
2^{-h_j-1}<t_j/T_j\leq2^{-h_j}\ (2\leq j\leq k)\big\}.
$$
The nonempty regions form a countable disjoint Borel partition of $\mathbb R^N\times(0,\infty)^k$.
On $\mathcal T_\gamma$, we have $y\in R_\gamma:=R$ and
$$
\frac{\ell(I_1)}{48k}\leq t_1<\frac{\ell(I_1)}{24k},\qquad
\frac{2^{-h_j}\ell(I_j)}{96k}<t_j<\frac{2^{-h_j}\ell(I_j)}{24k},\quad 2\leq j\leq k.
$$
Integrating in $y$ and the scale variables, we obtain
$$
\nu_k(\mathcal T_\gamma)\leq |R_\gamma|(\log 2)(\log 4)^{k-1}.
$$

For a finite set $\mathcal E$ of indices, let $m_{\mathcal E}=\mathbf1_{\bigcup_{\gamma\in\mathcal E}\mathcal T_\gamma}$.
The contribution from $\mathcal T_\gamma$ with cancellation in the blocks $j\in\Gamma$ is
$$
A_\gamma^\Gamma(x)=\int_{\mathcal T_\gamma}K_{\mathbf t}^\Gamma(x-y)H_f(y,\mathbf t)\,d\nu_k.
$$
For all $\Gamma\ni k$, this term has the same coefficient
$$
c_\gamma=\bigg(\frac1{|R_\gamma|}\int_{\mathcal T_\gamma}|H_f(y,\mathbf t)|^2\,d\nu_k\bigg)^{1/2}.
$$
On the translated kernel support, $|x_j-y_j|_\infty\leq kT_j$ by Proposition~\ref{p:flag-kernel-bounds}.
Since $\overline{B_\infty(y_j,4kT_j)}\subset I_j^\circ$, the distance from $x_j$ to $\partial I_j$
is at least $3kT_j\geq\ell(I_j)/16$. In particular,
\begin{equation}\label{e:flag-atom-support}
\operatorname{supp}A_\gamma^\Gamma\subseteq
\prod_{j=1}^k\big\{x_j\in I_j:\operatorname{dist}_\infty(x_j,\partial I_j)\geq\ell(I_j)/48\big\}.
\end{equation}
By Fubini's theorem and the cancellation of $K_{\mathbf t}^\Gamma$,
\begin{equation}\label{e:flag-atom-cancel}
\int_{X_j}A_\gamma^\Gamma(x)\,dx_j=0,\qquad j\in\Gamma.
\end{equation}

On each fixed $\mathcal T_\gamma$, the scales are bounded away from zero and
$\mathbf1_{\mathcal T_\gamma}H_f\in L^1(d\nu_k)$, so we may differentiate under the integral.
To apply \eqref{e:flag-kernel-deriv}, note that
$$
\frac{t_j}{T_j}\leq2^{-h_j}\quad(2\leq j\leq k),\qquad
\frac{T_j}{T_{j+1}}\leq2\frac{\ell(I_j)}{\ell(I_{j+1})}\quad(1\leq j<k).
$$
By Cauchy--Schwarz, \eqref{e:flag-kernel-deriv} and $T_j\simeq_k\ell(I_j)$,
\begin{align*}
|\partial^\alpha A_\gamma^\Gamma(x)|
\leq &|R_\gamma|^{1/2}c_\gamma
\bigg(\int_{\mathcal T_\gamma}|\partial^\alpha K_{\mathbf t}^\Gamma(x-y)|^2\,d\nu_k\bigg)^{1/2}\\
\leq &C_{\alpha,k}
\bigg(\frac{\nu_k(\mathcal T_\gamma)}{|R_\gamma|}\bigg)^{1/2}
2^{-2\sum_{j=2}^kh_j}c_\gamma
\prod_{j\notin\Gamma}\Big(\frac{\ell(I_j)}{\ell(I_{j+1})}\Big)^2
\prod_{j=1}^k\ell(I_j)^{-|\alpha_j|}\mathbf1_{R_\gamma}(x)\\
\leq &C_{\alpha,k}2^{-2\sum_{j=2}^kh_j}c_\gamma
\prod_{j\notin\Gamma}\Big(\frac{\ell(I_j)}{\ell(I_{j+1})}\Big)^2
\prod_{j=1}^k\ell(I_j)^{-|\alpha_j|}\mathbf1_{R_\gamma}(x).
\end{align*}

For fixed $\kappa,R$, sum over the scale indices $\mathbf h$ for which $(\kappa,R,\mathbf h)\in\mathcal E$.
Write
$$
A_{\kappa,R}^\Gamma=\sum_{\mathbf h:(\kappa,R,\mathbf h)\in\mathcal E}A_{\kappa,R,\mathbf h}^\Gamma,
\qquad
c_{\kappa,R}=\bigg(\sum_{\mathbf h:(\kappa,R,\mathbf h)\in\mathcal E}c_{\kappa,R,\mathbf h}^2\bigg)^{1/2}.
$$
By Cauchy--Schwarz in $\mathbf h$,
\begin{align}
|\partial^\alpha A_{\kappa,R}^\Gamma(x)|
\leq &C_{\alpha,k}
\prod_{j\notin\Gamma}\Big(\frac{\ell(I_j)}{\ell(I_{j+1})}\Big)^2
\prod_{j=1}^k\ell(I_j)^{-|\alpha_j|}\mathbf1_R(x)
\sum_{\mathbf h:(\kappa,R,\mathbf h)\in\mathcal E}2^{-2\sum_{j=2}^kh_j}c_{\kappa,R,\mathbf h}\notag\\
\leq &C_{\alpha,k}\,c_{\kappa,R}
\prod_{j\notin\Gamma}\Big(\frac{\ell(I_j)}{\ell(I_{j+1})}\Big)^2
\prod_{j=1}^k\ell(I_j)^{-|\alpha_j|}\mathbf1_R(x),\label{e:flag-sum-deriv}
\end{align}
since $\sum_{\mathbf h\in\mathbb N_0^{k-1}}2^{-4\sum_{j=2}^kh_j}=(1-2^{-4})^{1-k}$.
All terms in this sum have the same rectangle $R$ and cancellation in every block $j\in\Gamma$.
Thus \eqref{e:flag-atom-support} and \eqref{e:flag-atom-cancel} also hold for $A_{\kappa,R}^\Gamma$.

\subsubsection{Pointwise coefficient bounds and convergence of the finite sums}

For the finite sum $f_{\mathcal E}=\sum_{\kappa,R}\sum_{\Gamma\ni k}A_{\kappa,R}^\Gamma$, we have
$$
\sum_{\kappa,R}c_{\kappa,R}^2\mathbf1_R(x)
=\sum_{\gamma\in\mathcal E}\frac{\mathbf1_{R_\gamma}(x)}{|R_\gamma|}
\int_{\mathcal T_\gamma}|H_f(y,\mathbf t)|^2\,d\nu_k.
$$
To compare the right-hand side with $S_F^{(k)}(f)(x)^2$, we first show that
\begin{equation}\label{e:flag-cone-weight}
\frac{\prod_{j=1}^k\mathbf1_{\{|x_j|<T_j/(2\sqrt N)\}}}{\prod_{j=1}^k T_j^{d_j}}
\leq C\frac{\chi_{\mathbf t}(x)}{\prod_{j=1}^k t_j^{D_j}}.
\end{equation}

Let $r_j=t_j/\sqrt N$ and write $q_{r,d}(u)=(2r)^{-d}\mathbf1_{B_\infty(0,r)}(u)$ on $\mathbb R^d$.
For $v_j=(v_{j,j},\ldots,v_{j,k})\in V_j$, the cube $B_\infty(0,r_j)$ in $V_j$ lies in the Euclidean ball of radius $t_j$. Hence
$$
t_j^{-D_j}\mathbf1_{\{|v_j|<t_j\}}
\geq c\prod_{\ell=j}^k q_{r_j,d_\ell}(v_{j,\ell}).
$$
Convolving in $V_1,\ldots,V_k$ and separating the integrals over the blocks $X_\ell$, we obtain
$$
\frac{\chi_{\mathbf t}(x)}{\prod_{j=1}^k t_j^{D_j}}
\geq c\prod_{\ell=1}^k(q_{r_1,d_\ell}*\cdots*q_{r_\ell,d_\ell})(x_\ell).
$$
Fix $\ell$ and choose $i\leq\ell$ with $r_i=\max_{j\leq\ell}r_j=T_\ell/\sqrt N$.
For $|x_\ell|_\infty<r_i/2$, restrict the other variables to $|u_j|_\infty<r_j/(2\ell)$.
On this region the $i$th density is constant, since
$$
\bigg|x_\ell-\sum_{\substack{j\leq\ell\\j\ne i}}u_j\bigg|_\infty
<\frac{r_i}{2}+\frac{\ell-1}{2\ell}r_i<r_i.
$$
It follows that
$$
(q_{r_1,d_\ell}*\cdots*q_{r_\ell,d_\ell})(x_\ell)
\geq (2r_i)^{-d_\ell}\prod_{\substack{j\leq\ell\\j\ne i}}
\int_{|u_j|_\infty<r_j/(2\ell)}q_{r_j,d_\ell}(u_j)\,du_j
=(2r_i)^{-d_\ell}(2\ell)^{-d_\ell(\ell-1)}
\geq cT_\ell^{-d_\ell}.
$$
Since $|x_\ell|<T_\ell/(2\sqrt N)$ implies $|x_\ell|_\infty<r_i/2$, multiplying over $\ell$ proves \eqref{e:flag-cone-weight}.

Now fix $(y,\mathbf t)\in\mathcal T_\gamma$ and pass to the Poisson scales $s_j=\sigma_kt_j$.
For $x\in R_\gamma$, the choice of $\sigma_k$ ensures that
$$
|x_j-y_j|\leq\sqrt{d_j}\,\ell(I_j)\leq48k\sqrt{d_j}\,T_j
<\frac{\sigma_kT_j}{2\sqrt N}
=\frac{\max_{i\leq j}s_i}{2\sqrt N}.
$$
Since $\ell(I_j)>24kT_j$, we obtain from \eqref{e:flag-cone-weight}, for every $x\in\mathbb R^N$,
$$
\frac{\mathbf1_{R_\gamma}(x)}{|R_\gamma|}
\leq C\sigma_k^N
\prod_{j=1}^k\frac{\mathbf1_{\{|x_j-y_j|<(\max_{i\leq j}s_i)/(2\sqrt N)\}}}
{(\max_{i\leq j}s_i)^{d_j}}
\leq C\sigma_k^N\frac{\chi_{\mathbf s}(x-y)}{\prod_{j=1}^k s_j^{D_j}}.
$$
The component of $\mathcal G_{\mathbf s}f$ with a time derivative in every parameter satisfies
$$
|H_f(y,\mathbf s/\sigma_k)|^2
=\sigma_k^{-2k}\big|[\mathcal G_{\mathbf s}f(y)]_{(0,\ldots,0)}\big|^2
\leq\sigma_k^{-2k}|\mathcal G_{\mathbf s}f(y)|^2.
$$
Integrating over the disjoint regions $\mathcal T_\gamma$ and using $dt_j/t_j=ds_j/s_j$, we obtain
\begin{align}
\sum_{\kappa,R}c_{\kappa,R}^2\mathbf1_R(x)
\leq &C\sigma_k^N\int_{(0,\infty)^k}\int_{\mathbb R^N}
\chi_{\mathbf s}(x-y)|H_f(y,\mathbf s/\sigma_k)|^2\,dy
\prod_{j=1}^k\frac{ds_j}{s_j^{D_j+1}}\notag\\
\leq &C\sigma_k^{N-2k}\int_{(0,\infty)^k}\int_{\mathbb R^N}
\chi_{\mathbf s}(x-y)|\mathcal G_{\mathbf s}f(y)|^2\,dy
\prod_{j=1}^k\frac{ds_j}{s_j^{D_j+1}}
\leq CS_F^{(k)}(f)(x)^2.\label{e:flag-coeff}
\end{align}
The constant is independent of $\mathcal E$ and of the side lengths of the rectangles.

On a finite union of parameter regions, $y$ is bounded and all scales lie in a compact subinterval of $(0,\infty)$.
By Proposition~\ref{p:flag-kernel-bounds} and the definition of $m_{\mathcal E}$,
$$
f_{\mathcal E}
=\sum_{\gamma\in\mathcal E}\int_{\mathcal T_\gamma}
\bigg(\sum_{\Gamma\ni k}K_{\mathbf t}^\Gamma(\cdot-y)\bigg)H_f(y,\mathbf t)\,d\nu_k
=f_{m_{\mathcal E}}\qquad\text{in }L^2(\mathbb R^N).
$$
Choose finite sets $\mathcal E_\nu\uparrow\{\gamma:\mathcal T_\gamma\ne\emptyset\}$.
Since $m_{\mathcal E_\nu}\to1$ pointwise and $f_1=f$, it follows from \eqref{e:flag-recon-l2} that
\begin{equation}\label{e:flag-finite-limit}
\|f_{\mathcal E_\nu}-f\|_{L^2(\mathbb R^N)}
\leq\|a_{\sigma_k}\|_{L^2(du/u)}^k
\|(m_{\mathcal E_\nu}-1)H_f\|_{L^2(d\nu_k)}\to0,
\end{equation}
where the limit follows by dominated convergence.

\subsubsection{Local estimates for each cancellation pattern}

For each cancellation pattern, we first use the scale ratios in \eqref{e:flag-sum-deriv}
to sum over the blocks without cancellation. In the remaining blocks we use the following product estimate.

\begin{lemma}\label{l:active-cubes}
Let $\emptyset\ne J\subseteq\{1,\ldots,k\}$ and fix a dyadic system in each $X_j$, $j\in J$.
Write $X_J=\prod_{j\in J}X_j$ and $dx_J=\prod_{j\in J}dx_j$.
For each tuple $\mathbf I=(I_j)_{j\in J}$ in a finite family of distinct dyadic rectangles
$R_{\mathbf I}=\prod_{j\in J}I_j$, let $W_{\mathbf I}$ be smooth, supported in $R_{\mathbf I}$
and of integral zero in each $x_j$. Suppose that, for some constants $h_{\mathbf I}\geq0$,
$$
|\partial^\alpha W_{\mathbf I}(x_J)|
\leq h_{\mathbf I}\mathbf1_{R_{\mathbf I}}(x_J)\prod_{j\in J}\ell(I_j)^{-|\alpha_j|},
\qquad |\alpha_j|\leq1\quad(j\in J).
$$
Then $G=\sum_{\mathbf I}W_{\mathbf I}$ satisfies
$$
\bigg(\sum_{\mathbf U}\Big|\prod_{j\in J}\Delta_{U_j}^{(j)}G(x_J)\Big|^2\bigg)^{1/2}
\leq C\bigg(\sum_{\mathbf I}h_{\mathbf I}^2\mathbf1_{R_{\mathbf I}}(x_J)\bigg)^{1/2},
$$
where $\mathbf U=(U_j)_{j\in J}$ ranges over all tuples of cubes in the chosen systems. Moreover,
$$
\|G\|_{L^p(dx_J)}
\leq Cp^{|J|/2}\bigg\|\bigg(\sum_{\mathbf I}h_{\mathbf I}^2\mathbf1_{R_{\mathbf I}}\bigg)^{1/2}\bigg\|_{L^p(dx_J)},
\qquad 2\leq p<\infty.
$$
The constants depend only on the dimensions $(d_j)_{j\in J}$.
\end{lemma}

\begin{proof}
Write $\Delta_{\mathbf U}=\prod_{j\in J}\Delta_{U_j}^{(j)}$. By support and separate cancellation,
$\Delta_{\mathbf U}W_{\mathbf I}=0$ unless $U_j\subseteq I_j$ for every $j\in J$.
For $I_j=U_j^{(n_j)}$, applying the mean value theorem once in each block as in Proposition~\ref{p:product-square}, we obtain
$$
|\Delta_{\mathbf U}W_{\mathbf I}|
\leq Ch_{\mathbf I}\mathbf1_{R_{\mathbf U}}\prod_{j\in J}\frac{\ell(U_j)}{\ell(I_j)}
=C2^{-\sum_{j\in J}n_j}h_{\mathbf I}\mathbf1_{R_{\mathbf U}}.
$$
For fixed $\mathbf n=(n_j)_{j\in J}$, the rectangles $R_{\mathbf U}$ with $U_j^{(n_j)}=I_j$ for all $j\in J$
partition $R_{\mathbf I}$. By Minkowski's inequality in $\ell^2(\mathbf U)$,
\begin{align*}
\bigg(\sum_{\mathbf U}|\Delta_{\mathbf U}G|^2\bigg)^{1/2}
\leq &C\sum_{\mathbf n\in\mathbb N_0^J}2^{-\sum_{j\in J}n_j}
\bigg(\sum_{\mathbf I}h_{\mathbf I}^2
\sum_{\mathbf U:\,U_j^{(n_j)}=I_j\ (j\in J)}\mathbf1_{R_{\mathbf U}}\bigg)^{1/2}\\
= &C\sum_{\mathbf n\in\mathbb N_0^J}2^{-\sum_{j\in J}n_j}
\bigg(\sum_{\mathbf I}h_{\mathbf I}^2\mathbf1_{R_{\mathbf I}}\bigg)^{1/2}
\leq C\bigg(\sum_{\mathbf I}h_{\mathbf I}^2\mathbf1_{R_{\mathbf I}}\bigg)^{1/2}.
\end{align*}

Let $\mathbb E_{j,r}$ denote dyadic averaging in $X_j$ at side length $2^{-r}$.
Choose $m_0$ with $2^{-m_0}\geq\ell(I_j)$ for every cube in the family.
By cancellation, $\mathbb E_{j,m_0}G=0$ for every $j\in J$. Thus, as in \eqref{e:martingale-sums},
$$
G_M:=\sum_{\substack{\mathbf U:\,2^{-M}\leq\ell(U_j)\leq2^{-m_0}\\j\in J}}\Delta_{\mathbf U}G
=\prod_{j\in J}(\mathbb E_{j,M+1}-\mathbb E_{j,m_0})G
=\prod_{j\in J}\mathbb E_{j,M+1}G\to G\qquad\text{in }L^p(dx_J).
$$
For $G_M$, we apply \eqref{e:scalar-square} in the first block of $J$.
In each subsequent block, the preceding cube indices are the $\ell^2$ coordinates.
Since differences in distinct blocks commute, we may apply \eqref{e:hilbert-square} in turn, obtaining
\begin{align*}
\|G\|_{L^p(dx_J)}
= &\lim_{M\to\infty}\|G_M\|_{L^p(dx_J)}
\leq Cp^{|J|/2}\bigg\|\bigg(\sum_{\mathbf U}|\Delta_{\mathbf U}G|^2\bigg)^{1/2}\bigg\|_{L^p(dx_J)}\\
\leq &Cp^{|J|/2}\bigg\|\bigg(\sum_{\mathbf I}h_{\mathbf I}^2\mathbf1_{R_{\mathbf I}}\bigg)^{1/2}\bigg\|_{L^p(dx_J)},
\end{align*}
as desired.
\end{proof}

Fix $\kappa=(\kappa_1,\ldots,\kappa_k)$ and a finite set $\mathcal E$ of parameter indices.
For $\mathbf I=(I_1,\ldots,I_k)$, write $R_{\mathbf I}=\prod_{j=1}^k I_j$, $r_j=\ell(I_j)$,
$A_{\mathbf I}^\Gamma=A_{\kappa,R_{\mathbf I}}^\Gamma$ and $c_{\mathbf I}=c_{\kappa,R_{\mathbf I}}$.
The following local Hoeffding-type inequality estimates the centred sum and its coefficient square function on the same rectangle.

\begin{proposition}\label{p:flag-local-sums}
For every $k\in\Gamma\subseteq\{1,\ldots,k\}$ and every product $R_0=\prod_{j=1}^k E_j$ of axis-parallel cubes,
\begin{equation}\label{e:flag-local-sums}
\bigg\|\prod_{j=1}^kQ_j\sum_{\mathbf I}A_{\mathbf I}^\Gamma\bigg\|_{L^p(R_0)}
\leq Cp^{|\Gamma|/2}\bigg\|\bigg(\sum_{\mathbf I}c_{\mathbf I}^2\mathbf1_{R_{\mathbf I}}\bigg)^{1/2}\bigg\|_{L^p(R_0)},
\qquad 2\leq p<\infty.
\end{equation}
The constant depends only on $k,d_1,\ldots,d_k$.
\end{proposition}

\begin{proof}
\noindent{\it Step 1.}
Fix $\mathbf I_\Gamma$ and sum over the coordinates without cancellation. Write
$$
B_{\mathbf I_\Gamma}(x)=\sum_{\mathbf I_{\Gamma^c}}A_{\mathbf I}^\Gamma(x),
\qquad
b_{\mathbf I_\Gamma}(x_{\Gamma^c})^2
=\sum_{\mathbf I_{\Gamma^c}}c_{\mathbf I}^2\prod_{j\notin\Gamma}\mathbf1_{I_j}(x_j).
$$
For $j\notin\Gamma$, let $g_j=\log_2(r_{j+1}/r_j)\in\mathbb N_0$.
Since $k\in\Gamma$, the lengths in $\mathbf I_\Gamma$ and the gaps $g_j$ determine all remaining lengths:
starting at $j=k-1$ and proceeding downwards, use $r_j=2^{-g_j}r_{j+1}$ whenever $j\notin\Gamma$.
At each such length there is at most one cube $I_j$ containing $x_j$. Hence
$$
\sum_{\mathbf I_{\Gamma^c}}\prod_{j\notin\Gamma}
\Big(\frac{r_j}{r_{j+1}}\Big)^4\mathbf1_{I_j}(x_j)
\leq\sum_{(g_j)_{j\notin\Gamma}\in\mathbb N_0^{\Gamma^c}}2^{-4\sum_{j\notin\Gamma}g_j}
=(1-2^{-4})^{-|\Gamma^c|}.
$$
By \eqref{e:flag-sum-deriv} and Cauchy--Schwarz in $\mathbf I_{\Gamma^c}$, for $|\alpha_j|\leq1$ in each block $j\in\Gamma$,
\begin{align}
|\partial_\Gamma^\alpha B_{\mathbf I_\Gamma}(x)|
\leq &C\bigg(\sum_{\mathbf I_{\Gamma^c}}\prod_{j\notin\Gamma}
\Big(\frac{r_j}{r_{j+1}}\Big)^4\mathbf1_{I_j}(x_j)\bigg)^{1/2}
b_{\mathbf I_\Gamma}(x_{\Gamma^c})\prod_{j\in\Gamma}r_j^{-|\alpha_j|}\mathbf1_{I_j}(x_j)\notag\\
\leq &Cb_{\mathbf I_\Gamma}(x_{\Gamma^c})\prod_{j\in\Gamma}r_j^{-|\alpha_j|}\mathbf1_{I_j}(x_j).
\label{e:flag-block-deriv}
\end{align}
For each fixed $x_{\Gamma^c}$, these functions have integral zero in every block $j\in\Gamma$, and
\begin{equation}\label{e:flag-coeff-id}
\sum_{\mathbf I_\Gamma}b_{\mathbf I_\Gamma}(x_{\Gamma^c})^2\prod_{j\in\Gamma}\mathbf1_{I_j}(x_j)
=\sum_{\mathbf I}c_{\mathbf I}^2\mathbf1_{R_{\mathbf I}}(x).
\end{equation}
Thus, for $G=\sum_{\mathbf I_\Gamma}B_{\mathbf I_\Gamma}$, Lemma~\ref{l:active-cubes} implies
\begin{equation}\label{e:flag-active-lp}
\|G(\cdot,x_{\Gamma^c})\|_{L^p(dx_\Gamma)}
\leq Cp^{|\Gamma|/2}
\bigg\|\bigg(\sum_{\mathbf I}c_{\mathbf I}^2\mathbf1_{R_{\mathbf I}}(\cdot,x_{\Gamma^c})\bigg)^{1/2}\bigg\|_{L^p(dx_\Gamma)}.
\end{equation}

\smallskip
\noindent{\it Step 2.}
To restrict the estimate to $R_0$, retain only those $\mathbf I_\Gamma$ for which
$(I_j)_*\cap E_j\ne\emptyset$ for every $j\in\Gamma$.
By \eqref{e:flag-atom-support}, all other terms vanish on $R_0$ after centring.
Let $G_A$ be the sum over the retained tuples with
$$
A=\{j\in\Gamma:r_j>\ell(E_j)\},\qquad S=\Gamma\setminus A.
$$
We first centre in the blocks $j\in A$. Write $E_A=\prod_{j\in A}E_j$ and
let $\tau_v^{(j)}$ denote translation by $v$ in $x_j$.
For $x_A\in E_A$,
$$
\prod_{j\in A}Q_jB_{\mathbf I_\Gamma}(x)
=\int_{E_A}\bigg[\prod_{j\in A}(I-\tau_{u_j-x_j}^{(j)})B_{\mathbf I_\Gamma}\bigg](x)
\prod_{j\in A}\frac{du_j}{|E_j|}.
$$
Applying the fundamental theorem of calculus along $[x_j,u_j]$ in each block $j\in A$
and using \eqref{e:flag-block-deriv}, we obtain for $|\alpha_j|\leq1$, $j\in S$,
\begin{align*}
\bigg|\prod_{j\in A}(I-\tau_{u_j-x_j}^{(j)})\partial_S^\alpha B_{\mathbf I_\Gamma}(x)\bigg|
\leq &Cb_{\mathbf I_\Gamma}(x_{\Gamma^c})
\prod_{j\in S}r_j^{-|\alpha_j|}\mathbf1_{I_j}(x_j)
 \prod_{j\in A}\frac{|x_j-u_j|}{r_j}
\big(\mathbf1_{I_j}(x_j)+\mathbf1_{I_j}(u_j)\big),
\end{align*}
since the mixed difference vanishes if $x_j,u_j\notin I_j$ for some $j\in A$.
Averaging in $u_A$ and using $|x_j-u_j|\leq\sqrt{d_j}\,\ell(E_j)$, we have
$$
\bigg|\partial_S^\alpha\prod_{j\in A}Q_jB_{\mathbf I_\Gamma}(x)\bigg|
\leq Cb_{\mathbf I_\Gamma}(x_{\Gamma^c})
\prod_{j\in A}\frac{\ell(E_j)}{r_j}
\bigg(\mathbf1_{I_j}(x_j)+\frac{|I_j\cap E_j|}{|E_j|}\bigg)
\prod_{j\in S}r_j^{-|\alpha_j|}\mathbf1_{I_j}(x_j).
$$

For fixed $x_A\in E_A$ and $x_{\Gamma^c}$, put
$$
W_{\mathbf I_S}=\sum_{\mathbf I_A}\prod_{j\in A}Q_jB_{\mathbf I_\Gamma},
$$
where the sum runs over the tuples contributing to $G_A$. On $E_A$, the corresponding coefficient is
$$
h_{\mathbf I_S}
=\prod_{j\in A}(\operatorname{Id}+P_j)
\bigg[\sum_{\mathbf I_A}b_{\mathbf I_\Gamma}(x_{\Gamma^c})
\prod_{j\in A}\frac{\ell(E_j)}{r_j}\mathbf1_{I_j}(x_j)\bigg].
$$
Indeed, $(\operatorname{Id}+P_j)\mathbf1_{I_j}(x_j)
=\mathbf1_{I_j}(x_j)+|I_j\cap E_j|/|E_j|$ on $E_j$, and $b_{\mathbf I_\Gamma}$ is independent of $x_A$. Thus
\begin{equation}\label{e:flag-local-block}
|\partial_S^\alpha W_{\mathbf I_S}(x)|
\leq Ch_{\mathbf I_S}\prod_{j\in S}r_j^{-|\alpha_j|}\mathbf1_{I_j}(x_j),
\qquad |\alpha_j|\leq1\quad(j\in S).
\end{equation}
Each $W_{\mathbf I_S}$ still has integral zero in $x_j$ for every $j\in S$.

\smallskip
\noindent{\it Step 3.}
For $j\in S$, we have $r_j\leq\ell(E_j)$ and $(I_j)_*\cap E_j\ne\emptyset$.
If $S\ne\emptyset$, apply Lemma~\ref{l:active-cubes} on $X_S$ and then Lemma~\ref{l:cube-fold} on $E_S=\prod_{j\in S}E_j$.
Since $\|Q_j\|_{p\to p}\leq2$,
\begin{align*}
\bigg\|\prod_{j\in S}Q_j\sum_{\mathbf I_S}W_{\mathbf I_S}\bigg\|_{L^p(E_S)}
\leq &Cp^{|S|/2}\bigg\|\bigg(\sum_{\mathbf I_S}h_{\mathbf I_S}^2
\prod_{j\in S}\mathbf1_{I_j}\bigg)^{1/2}\bigg\|_{L^p(X_S)}\\
\leq &Cp^{|S|/2}\bigg\|\bigg(\sum_{\mathbf I_S}h_{\mathbf I_S}^2
\prod_{j\in S}\mathbf1_{I_j}\bigg)^{1/2}\bigg\|_{L^p(E_S)}.
\end{align*}
For $S=\emptyset$, the corresponding pointwise bound follows directly from \eqref{e:flag-local-block}.

It remains to estimate the coefficients $h_{\mathbf I_S}$. For each $j\in A$, at most one $I_j$
of each dyadic side length contains $x_j$, so
$$
\sum_{\ell(I_j)>\ell(E_j)}
\Big(\frac{\ell(E_j)}{\ell(I_j)}\Big)^2\mathbf1_{I_j}(x_j)
\leq\sum_{r\geq0}2^{-2r}=\frac43.
$$
By Minkowski's inequality for the averages $P_j$, Cauchy--Schwarz in $\mathbf I_A$ and \eqref{e:flag-coeff-id},
\begin{align*}
\bigg(\sum_{\mathbf I_S}h_{\mathbf I_S}^2\prod_{j\in S}\mathbf1_{I_j}\bigg)^{1/2}
\leq &\prod_{j\in A}(\operatorname{Id}+P_j)
\bigg[\sum_{\mathbf I_S}
\bigg(\sum_{\mathbf I_A}b_{\mathbf I_\Gamma}
\prod_{j\in A}\frac{\ell(E_j)}{r_j}\mathbf1_{I_j}\bigg)^2
\prod_{j\in S}\mathbf1_{I_j}\bigg]^{1/2}\\
\leq &(4/3)^{|A|/2}\prod_{j\in A}(\operatorname{Id}+P_j)
\bigg(\sum_{\mathbf I_\Gamma}b_{\mathbf I_\Gamma}^2
\prod_{j\in\Gamma}\mathbf1_{I_j}\bigg)^{1/2}\\
\leq &C\prod_{j\in A}(\operatorname{Id}+P_j)
\bigg(\sum_{\mathbf I}c_{\mathbf I}^2\mathbf1_{R_{\mathbf I}}\bigg)^{1/2}.
\end{align*}
Taking the $L^p$ norm in the fixed variables and using the contraction of each $P_j$ on $R_0$, we obtain
$$
\bigg\|\prod_{j\in\Gamma}Q_jG_A\bigg\|_{L^p(R_0)}
\leq Cp^{|S|/2}\bigg\|\bigg(\sum_{\mathbf I_S}h_{\mathbf I_S}^2\prod_{j\in S}\mathbf1_{I_j}\bigg)^{1/2}\bigg\|_{L^p(R_0)}
\leq Cp^{|S|/2}\bigg\|\bigg(\sum_{\mathbf I}c_{\mathbf I}^2\mathbf1_{R_{\mathbf I}}\bigg)^{1/2}\bigg\|_{L^p(R_0)}.
$$
Finally, $\|\prod_{j\notin\Gamma}Q_j\|_{p\to p}\leq2^{|\Gamma^c|}$ and $S=\Gamma\setminus A$. Summing over $A\subseteq\Gamma$, we have
\begin{align*}
\bigg\|\prod_{j=1}^kQ_jG\bigg\|_{L^p(R_0)}
\leq &C\sum_{A\subseteq\Gamma}p^{(|\Gamma|-|A|)/2}
\bigg\|\bigg(\sum_{\mathbf I}c_{\mathbf I}^2\mathbf1_{R_{\mathbf I}}\bigg)^{1/2}\bigg\|_{L^p(R_0)}\\
= &C(1+\sqrt p)^{|\Gamma|}
\bigg\|\bigg(\sum_{\mathbf I}c_{\mathbf I}^2\mathbf1_{R_{\mathbf I}}\bigg)^{1/2}\bigg\|_{L^p(R_0)},
\end{align*}
which proves \eqref{e:flag-local-sums} for $p\geq2$.
\end{proof}


For the heat area function we shall use the one-parameter case of the preceding argument.
Let $\mathcal F$ be a finite family of cubes in a single dyadic system on $\mathbb R^q$.
For each $Q\in\mathcal F$, let $A_Q$ be smooth, supported in $Q_*$ and of integral zero, with
$\|\partial^\alpha A_Q\|_\infty\leq c_Q\ell(Q)^{-|\alpha|}$ for $|\alpha|\leq1$.
Then, for every axis-parallel cube $E$,
\begin{equation}\label{e:one-coeff-local}
\bigg\|\sum_{Q\in\mathcal F}\big(A_Q-(A_Q)_E\big)\bigg\|_{L^p(E)}
\leq C_q\sqrt p\,
\bigg\|\Big(\sum_{Q\in\mathcal F}c_Q^2\mathbf1_Q\Big)^{1/2}\bigg\|_{L^p(E)},
\qquad 2\leq p<\infty.
\end{equation}
For $\ell(Q)\leq\ell(E)$, only cubes with $Q_*\cap E\ne\emptyset$ contribute,
so Lemma~\ref{l:cube-fold} applies after Lemma~\ref{l:active-cubes}.
For $\ell(Q)>\ell(E)$, the first derivative bound gives, as in \eqref{e:horizontal-mean},
$$
|A_Q(x)-(A_Q)_E|
\leq C_q\frac{\ell(E)}{\ell(Q)}c_Q
\bigg(\mathbf1_Q(x)+\frac{|Q\cap E|}{|E|}\bigg),\qquad x\in E.
$$
Summing over the larger cubes by Cauchy--Schwarz completes the estimate.

\smallskip

\begin{proof}[Proof of Theorem~\ref{t:flag-local}]
First let $f\in L^2(\mathbb R^N)$. For each finite $\mathcal E$, Proposition~\ref{p:flag-local-sums}
and \eqref{e:flag-coeff} imply
\begin{align*}
\bigg\|\prod_{j=1}^kQ_jf_{\mathcal E}\bigg\|_{L^p(R_0)}
\leq &C\sum_{\Gamma\ni k}p^{|\Gamma|/2}\sum_\kappa
\bigg\|\bigg(\sum_Rc_{\kappa,R}^2\mathbf1_R\bigg)^{1/2}\bigg\|_{L^p(R_0)}\\
\leq &C\sqrt p(1+\sqrt p)^{k-1}
\bigg\|\bigg(\sum_{\kappa,R}c_{\kappa,R}^2\mathbf1_R\bigg)^{1/2}\bigg\|_{L^p(R_0)}\\
\leq &Cp^{k/2}\|S_F^{(k)}(f)\|_{L^p(R_0)},
\end{align*}
where there are $3^N$ choices of $\kappa$ and $\sum_{\Gamma\ni k}p^{|\Gamma|/2}=\sqrt p(1+\sqrt p)^{k-1}$.
By \eqref{e:flag-finite-limit} and the $L^2(R_0)$ contraction of each $Q_j$,
$$
\|Q_{R_0}f_{\mathcal E_\nu}-Q_{R_0}f\|_{L^2(R_0)}
\leq\|f_{\mathcal E_\nu}-f\|_{L^2(\mathbb R^N)}\to0.
$$
Fatou's lemma along an almost everywhere convergent subsequence proves \eqref{e:flag-local-lp} for $f\in L^2(\mathbb R^N)$.

For $2\leq p\leq u$, choose $f_i\in C_c^\infty(\mathbb R^N)$ with $f_i\to f$ in $L^u(\mathbb R^N)$.
By H\"older's inequality on $R_0$ and \eqref{e:flag-area-lp},
$$
\|Q_{R_0}(f_i-f)\|_{L^p(R_0)}
+\|S_F^{(k)}(f_i)-S_F^{(k)}(f)\|_{L^p(R_0)}
\leq Cu^{k/2}|R_0|^{1/p-1/u}\|f_i-f\|_u\to0.
$$
Passing to the limit proves \eqref{e:flag-local-lp} in this range.

For $1\leq u<p$, suppose $S_F^{(k)}(f)\in L^p(R_0)$. Take $R_i\uparrow R_0^\circ$ and mollifiers $\rho_i$
as in Corollary~\ref{c:two-extension}, with $R_i-\operatorname{supp}\rho_i\subset R_0$.
Since $g_i=\rho_i*f\in L^u\cap L^\infty\subset L^p$, the case $u=p$ and Minkowski's inequality in the flag area integral imply
\begin{align*}
\|Q_{R_i}g_i\|_{L^p(R_i)}
\leq &Cp^{k/2}\|S_F^{(k)}(g_i)\|_{L^p(R_i)}
\leq Cp^{k/2}\|\rho_i*(\mathbf1_{R_0}S_F^{(k)}(f))\|_{L^p(R_i)}
\leq Cp^{k/2}\|S_F^{(k)}(f)\|_{L^p(R_0)}.
\end{align*}
As in Corollary~\ref{c:two-extension}, the zero extensions of $Q_{R_i}g_i$ converge along a subsequence
to $Q_{R_0}f$ almost everywhere on $R_0$. Fatou's lemma proves \eqref{e:flag-local-lp}.

For \eqref{e:flag-exponential}, normalize $\|S_F^{(k)}(f)\|_{L^\infty(R_0)}=1$.
For each integer $m\geq1$, take $q=\max\{2,2m/k\}$. By H\"older's inequality on $(R_0,dx/|R_0|)$ and \eqref{e:flag-local-lp},
$$
\frac1{|R_0|}\int_{R_0}|Q_{R_0}f|^{2m/k}\,dx
\leq\bigg(\frac{\|Q_{R_0}f\|_{L^q(R_0)}}{|R_0|^{1/q}}\bigg)^{2m/k}
\leq(Cq^{k/2})^{2m/k}\leq(C_1m)^m.
$$
Choose $c>0$ with $cC_1e\leq1/2$. Since $m!\geq(m/e)^m$,
\begin{equation}\label{e:flag-exp-series}
\frac1{|R_0|}\int_{R_0}\exp(c|Q_{R_0}f|^{2/k})\,dx
=1+\sum_{m=1}^\infty\frac{c^m}{m!|R_0|}\int_{R_0}|Q_{R_0}f|^{2m/k}\,dx
\leq1+\sum_{m=1}^\infty(cC_1e)^m\leq2.
\end{equation}
Rescaling $f$, we obtain \eqref{e:flag-exponential}.
When $0<\|S_F^{(k)}(f)\|_\infty<\infty$, the same calculation with the global norm in the denominator is uniform over $R_0$.

For $f_p$ and the flag rectangle $R_p$ in Proposition~\ref{p:flag-reverse-sharp},
$$
\|Q_{R_p}f_p\|_{L^p(R_p)}
=\|f_p\|_{L^p(R_p)}
\geq cp^{k/2}\|S_F^{(k)}(f_p)\|_{L^p(\mathbb R^N)}
\geq cp^{k/2}\|S_F^{(k)}(f_p)\|_{L^p(R_p)}.
$$
The same proposition proves sharpness of $2/k$ even with the global area norm in the denominator.
Remark~\ref{r:fixed-unit} gives both lower bounds on the fixed unit cube.
\end{proof}


\subsection{The global reverse inequality}\label{s:flag-reverse}

The local estimate \eqref{e:flag-local-lp} implies the global reverse inequality
once the coordinate averages in $Q_Rf$ tend to zero.

\begin{theorem}\label{t:flag-reverse}
Let $f\in L^u(\mathbb R^N)$ for some $1\leq u<\infty$, and let $2\leq p<\infty$.
If $S_F^{(k)}(f)\in L^p(\mathbb R^N)$, then $f\in L^p(\mathbb R^N)$ and
$$
\|f\|_p\leq Cp^{k/2}\|S_F^{(k)}(f)\|_p,
$$
where $C$ depends only on $k,d_1,\ldots,d_k$.
\end{theorem}

\begin{proof}
Let $E_{j,L}=[-L,L]^{d_j}$ and $R_L=\prod_{j=1}^k E_{j,L}$, $L\geq1$.
For $\emptyset\ne A\subseteq\{1,\ldots,k\}$, write $R_{A,L}=\prod_{j\in A}E_{j,L}$.
Extending $Q_{R_L}f$ by zero outside $R_L$ and expanding the product of the centring operators, we obtain
$$
Q_{R_L}f(x)
=\mathbf1_{R_L}(x)\bigg[f(x)+\sum_{\emptyset\ne A\subseteq\{1,\ldots,k\}}
\frac{(-1)^{|A|}}{|R_{A,L}|}\int_{R_{A,L}}f(y_A,x_{A^c})\,dy_A\bigg].
$$
By Fubini's theorem and H\"older's inequality, for each such $A$ and almost every $x_{A^c}$,
$$
\bigg|\frac1{|R_{A,L}|}\int_{R_{A,L}}f(y_A,x_{A^c})\,dy_A\bigg|
\leq |R_{A,L}|^{-1/u}\|f(\cdot,x_{A^c})\|_{L^u(X_A)}\to0.
$$
Thus $Q_{R_L}f\to f$ almost everywhere. By Fatou's lemma and \eqref{e:flag-local-lp},
$$
\|f\|_p
\leq\liminf_{L\to\infty}\|Q_{R_L}f\|_{L^p(R_L)}
\leq Cp^{k/2}\lim_{L\to\infty}\|S_F^{(k)}(f)\|_{L^p(R_L)}
=Cp^{k/2}\|S_F^{(k)}(f)\|_p,
$$
as desired.
\end{proof}

\begin{remark}\label{r:time-area}
The local and global reverse estimates also hold for the area function formed from the mixed time derivative.
For $\eta>0$, define
$$
S_{F,\eta}^{\mathrm{time}}(f)(x)^2
=\int_{(0,\infty)^k}\int_{\eta B_1\times\cdots\times\eta B_k}
\bigg|\bigg(\prod_{j=1}^kt_j\partial_{t_j}\bigg)U_f
\bigg(x-\sum_{j=1}^kt_ja_j,\mathbf t\bigg)\bigg|^2
\,da_1\cdots da_k\prod_{j=1}^k\frac{dt_j}{t_j}.
$$
In Theorems~\ref{t:flag-local} and~\ref{t:flag-reverse}, we may replace $S_F^{(k)}$ on the right
by $S_{F,\eta}^{\mathrm{time}}$, with constants depending also on $\eta$.

By monotonicity in $\eta$, it suffices to consider $0<\eta\leq1$. Choose
$$
\sigma_k=\max\Big\{\sigma_0,98k\sqrt{N\max_i d_i}/\eta\Big\}.
$$
At the Poisson scales $s_j=\sigma_kt_j$, we have
$$
49k\sqrt{d_j}\,T_j\leq\frac{\eta}{2\sqrt N}\max_{i\leq j}s_i,\qquad 1\leq j\leq k.
$$
Since $H_f(y,\mathbf t)=(-1)^k\sigma_k^{-k}[\mathcal G_{\mathbf s}f(y)]_{(0,\ldots,0)}$,
applying \eqref{e:flag-cone-weight} at scales $\eta\mathbf s$ in the proof of \eqref{e:flag-coeff}, we obtain
\begin{align*}
\sum_{\kappa,R}c_{\kappa,R}^2\mathbf1_R(x)
\leq &C_\eta\int_{(0,\infty)^k}\int_{\mathbb R^N}
\chi_{\eta\mathbf s}(x-y)\big|[\mathcal G_{\mathbf s}f(y)]_{(0,\ldots,0)}\big|^2
\,dy\prod_{j=1}^k\frac{ds_j}{s_j^{D_j+1}}
= C_\eta S_{F,\eta}^{\mathrm{time}}(f)(x)^2.
\end{align*}
With this coefficient bound, the proof of Theorem~\ref{t:flag-local} applies to $f\in L^2(\mathbb R^N)$.
For $2\leq p\leq u$, the same density argument applies since
$S_{F,\eta}^{\mathrm{time}}(h)\leq S_F^{(k)}(h)$.
For $1\leq u<p$, Minkowski's inequality on the inner rectangles gives
$$
S_{F,\eta}^{\mathrm{time}}(\rho_i*f)
\leq \rho_i*S_{F,\eta}^{\mathrm{time}}(f).
$$
The local moment estimate therefore extends to $f\in L^u(\mathbb R^N)$ as before.
Summing the moments as in \eqref{e:flag-exp-series} proves the exponential estimate.
Letting $R_L$ increase to $\mathbb R^N$ proves the global reverse inequality.
\end{remark}


\begin{remark}\label{r:global-intermediate}
For $f\in L^p(\mathbb R^N)$, the global reverse estimate also follows by iterating
the one-parameter estimate on the nested subspaces $V_1\supset\cdots\supset V_k$.
As in the product-space argument of R.~Fefferman and Pipher~\cite[Section~5]{FP},
we keep the parameters from each preceding step in a Hilbert space.

Fix $2\leq p<\infty$. For real functions $g_1,\ldots,g_M\in L^p(\mathbb R^N)$,
let $Z_1,\ldots,Z_M$ be independent standard Gaussian variables.
Applying the case $k=1$ of Theorem~\ref{t:flag-reverse} to $\sum_i Z_i g_i$ on each $V_j$-slice,
then using Minkowski's inequality in $L^{p/2}$ of the Gaussian variables, we obtain
\begin{align*}
(\mathbb E|Z_1|^p)^{1/p}\bigg\|\bigg(\sum_i|g_i|^2\bigg)^{1/2}\bigg\|_p
= &\bigg(\int_{\mathbb R^N}\mathbb E\bigg|\sum_i Z_i g_i(x)\bigg|^p\,dx\bigg)^{1/p}\\
\leq &C_{D_j}\sqrt p\,\bigg(\int_{\mathbb R^N}\mathbb E
\bigg(\int\bigg|\sum_i Z_i T_{j,\omega_j}g_i(x)\bigg|^2\,d\mu_j\bigg)^{p/2}dx\bigg)^{1/p}\\
\leq &C_{D_j}\sqrt p\,(\mathbb E|Z_1|^p)^{1/p}
\bigg\|\bigg(\int\sum_i|T_{j,\omega_j}g_i|^2\,d\mu_j\bigg)^{1/2}\bigg\|_p.
\end{align*}
Cancel the common Gaussian moment and apply the estimate to the first $M$ coordinates
of $g\in L^p(\mathbb R^N;\mathcal H)$ in an orthonormal basis of a separable Hilbert space $\mathcal H$.
Since $T_{j,\omega_j}$ acts coordinatewise, monotone convergence yields
$$
\|g\|_{L^p(\mathbb R^N;\mathcal H)}
\leq C_{D_j}\sqrt p\,\bigg\|\bigg(\int\|T_{j,\omega_j}g\|_{\mathcal H}^2\,d\mu_j\bigg)^{1/2}\bigg\|_p.
$$
For complex $\mathcal H$, we use its underlying real Hilbert space.

For $k\geq2$, the first two steps are
$$
\|f\|_p
\leq C\sqrt p\,\bigg\|\bigg(\int|T_{1,\omega_1}f|^2\,d\mu_1\bigg)^{1/2}\bigg\|_p
\leq Cp\,\bigg\|\bigg(\iint|T_{2,\omega_2}T_{1,\omega_1}f|^2\,d\mu_1\,d\mu_2\bigg)^{1/2}\bigg\|_p.
$$
The first intermediate square function is the pointwise $L^2(d\mu_1)$ norm of
$(T_{1,\omega_1}f)_{\omega_1}$. To pass to the second estimate, we apply $T_{2,\omega_2}$
to each member of this family, just as the martingale difference in \eqref{e:vertical-hilbert}
acts on each $\ell^2$ coordinate.

At step $j$, the preceding parameters form the Hilbert space $L^2(d\mu_1\cdots d\mu_{j-1})$.
The pointwise norm of the family belongs to $L^p$ by \eqref{e:flag-area-lp}
for the shorter flag ending at $V_{j-1}$. Iterating and using \eqref{e:flag-area}, we obtain
$$
\|f\|_p
\leq Cp^{k/2}\bigg\|\bigg(\int
|T_{k,\omega_k}\cdots T_{1,\omega_1}f|^2\,d\mu_1\cdots d\mu_k\bigg)^{1/2}\bigg\|_p
=Cp^{k/2}\|S_F^{(k)}(f)\|_p.
$$
The exhaustion argument in Theorem~\ref{t:flag-reverse} also proves $f\in L^p$
when one assumes only $f\in L^u$, $1\leq u<\infty$, and $S_F^{(k)}(f)\in L^p$.
\end{remark}

\subsection{Sharpness of the upper and reverse square-function estimates}\label{s:flag-sharpness}

The local estimates in Theorem~\ref{t:concentration}\textup{(i)} are \eqref{e:flag-local-lp}
and \eqref{e:flag-exponential}. The global estimates in \textup{(ii)} follow from
\eqref{e:flag-area-lp} and Theorem~\ref{t:flag-reverse}.
For sharpness of the upper bound, we estimate the area function of the unit cube indicator
near the intersection of $k$ faces, one from each coordinate block.
For the reverse bound, we separate the coordinate frequencies so that on each $V_j$
the frequency in $X_j$ dominates those in the subsequent blocks.
Remark~\ref{r:fixed-unit} places the local examples on the fixed unit cube.

\begin{proposition}\label{p:flag-area-sharp}
Let $f=\mathbf1_{[0,1]^N}$. Then $S_F^{(k)}(f)\notin L^\infty(\mathbb R^N)$ and
$$
\|S_F^{(k)}(f)\|_p\geq cp^{k/2},\qquad 2\leq p<\infty.
$$
For each such $p$, there is also a nonzero real $f_p\in C_c^\infty(\mathbb R^N)$ such that
$$
\|S_F^{(k)}(f_p)\|_p\geq cp^{k/2}\|f_p\|_p.
$$
The constant depends only on $k,d_1,\ldots,d_k$.
\end{proposition}

\begin{proof}
\noindent{\it Step 1.}
We estimate the component of $\mathcal G_{\mathbf t}f$ obtained by differentiating in the first
coordinate of each block. By subordination,
$$
e^{-t\sqrt{-\Delta}}
=\frac1{2\sqrt\pi}\int_0^\infty e^{t^2s\Delta}e^{-1/(4s)}\,\frac{ds}{s^{3/2}}.
$$
Let $G_a(v)=(4\pi a)^{-1/2}e^{-v^2/(4a)}$ and $q_a=G_a*\mathbf1_{[0,1]}$.
Since $\Delta_{V_j}=\sum_{i=j}^k\Delta_{X_i}$, the heat time in $X_i$ is
$A_i=\sum_{j=1}^it_j^2s_j$. Thus
$$
\bigg(\prod_{j=1}^ke^{t_j^2s_j\Delta_{V_j}}\bigg)f(y)
=\prod_{i=1}^k\prod_{\ell=1}^{d_i}q_{A_i}(y_{i,\ell}),
$$
and
$$
\Big(\prod_{i=1}^kt_i\partial_{y_{i,1}}\Big)U_f(y,\mathbf t)
=\int_{(0,\infty)^k}
\prod_{i=1}^k\bigg[t_iq_{A_i}'(y_{i,1})\prod_{\ell=2}^{d_i}q_{A_i}(y_{i,\ell})\bigg]
\prod_{j=1}^k\frac{e^{-1/(4s_j)}}{2\sqrt\pi}\frac{ds_j}{s_j^{3/2}}.
$$
Indeed, $0\leq q_a\leq1$ and $t_i\|q_{A_i}'\|_\infty\leq Ct_iA_i^{-1/2}\leq Cs_i^{-1/2}$.
Since $\int_0^\infty e^{-1/(4s)}s^{-2}\,ds<\infty$, we may differentiate under the integral.

For $v<1/2$ and every $a>0$,
$$
q_a'(v)=G_a(v)-G_a(v-1)=G_a(v)\big(1-e^{-(1-2v)/(4a)}\big)>0.
$$
We may therefore restrict the integral to $1\leq s_j\leq2$ whenever $y_{i,1}<1/2$ for every $i$.
Let $0<t_1\leq\cdots\leq t_k\leq t_0=1/(8\sqrt k)$ and suppose that
$$
|y_{i,1}|\leq t_i,\qquad y_i'\in[1/4,3/4]^{d_i-1},\qquad 1\leq i\leq k.
$$
On $[1,2]^k$, we have $t_i^2\leq A_i\leq2kt_i^2\leq1/32$ and $1-2y_{i,1}\geq3/4$. Hence
$$
t_iq_{A_i}'(y_{i,1})
=\frac{t_i e^{-y_{i,1}^2/(4A_i)}}{\sqrt{4\pi A_i}}
\big(1-e^{-(1-2y_{i,1})/(4A_i)}\big)
\geq\frac{e^{-1/4}(1-e^{-6})}{\sqrt{8\pi k}}.
$$
For $2\leq\ell\leq d_i$, we also have
$$
q_{A_i}(y_{i,\ell})\geq\int_{-1/4}^{1/4}G_{A_i}(v)\,dv
\geq\frac1{\sqrt\pi}\int_{-1/\sqrt2}^{1/\sqrt2}e^{-v^2}\,dv.
$$
Restricting the subordination integral to $[1,2]^k$, we conclude that
$$
\Big(\prod_{i=1}^kt_i\partial_{y_{i,1}}\Big)U_f(y,\mathbf t)
\geq c\bigg(\frac1{2\sqrt\pi}\int_1^2e^{-1/(4s)}\,\frac{ds}{s^{3/2}}\bigg)^k
\geq c'>0.
$$

\smallskip
\noindent{\it Step 2.}
For $0<\delta<t_0/2$, let
$$
\Omega_\delta=\big\{x:0<x_{i,1}<\delta,\ x_i'\in[3/8,5/8]^{d_i-1},\ 1\leq i\leq k\big\}.
$$
Fix $x\in\Omega_\delta$ and $2\delta\leq t_1\leq\cdots\leq t_k\leq t_0$.
If $a_j\in V_j$ and $|a_j|<1/(4k)$, then $y=x-\sum_jt_ja_j$ satisfies
$$
|y_i-x_i|\leq\sum_{j\leq i}t_j|a_j|\leq t_i/4.
$$
Thus $|y_{i,1}|\leq\delta+t_i/4\leq3t_i/4$ and $y_i'\in[1/4,3/4]^{d_i-1}$.
Using the component from Step~1 in \eqref{e:flag-area}, we obtain
\begin{align*}
S_F^{(k)}(f)(x)^2
\geq &\int_{2\delta\leq t_1\leq\cdots\leq t_k\leq t_0}
\int_{\prod_{j=1}^k B_{V_j}(0,1/(4k))}
\bigg|\Big(\prod_{i=1}^kt_i\partial_{y_{i,1}}\Big)U_f
\Big(x-\sum_{j=1}^kt_ja_j,\mathbf t\Big)\bigg|^2
\prod_{j=1}^k da_j\frac{dt_j}{t_j}\\
\geq &c\int_{2\delta\leq t_1\leq\cdots\leq t_k\leq t_0}\prod_{j=1}^k\frac{dt_j}{t_j}
=\frac{c}{k!}\log^k\frac{t_0}{2\delta}.
\end{align*}
Since $|\Omega_\delta|=4^{-(N-k)}\delta^k>0$ for every $\delta$, the area function is unbounded.
Taking $\delta=(t_0/2)e^{-p}$, we further obtain
$$
\|S_F^{(k)}(f)\|_p
\geq cp^{k/2}|\Omega_\delta|^{1/p}
=cp^{k/2}\big(4^{-(N-k)}(t_0/2)^k\big)^{1/p}e^{-k}
\geq c'p^{k/2}.
$$

Finally, let $f_\varepsilon=\rho_\varepsilon*f$ for a nonnegative $C_c^\infty$ approximate identity.
Then $\|f_\varepsilon\|_p\leq\|f\|_p=1$, and by \eqref{e:flag-area-lp},
$$
\|S_F^{(k)}(f_\varepsilon)-S_F^{(k)}(f)\|_p
\leq Cp^{k/2}\|f_\varepsilon-f\|_p\to0
\qquad(\varepsilon\to0+).
$$
For each $p$, choose $\varepsilon_p>0$ small enough that $f_p=f_{\varepsilon_p}$ satisfies
$$
\|S_F^{(k)}(f_p)\|_p
\geq\tfrac12\|S_F^{(k)}(f)\|_p
\geq cp^{k/2}\geq cp^{k/2}\|f_p\|_p,
$$
as desired.
\end{proof}

\begin{remark}\label{r:local-upper}
Even for a unit cube $R$, the full Poisson area on $R$ cannot be bounded in terms of $\|f\|_{L^p(R)}$.
Indeed, let $f=\mathbf1_{[0,1]^N}$ and take the adjacent cube $R=[-1,0]\times[0,1]^{N-1}$.
Reflect $\Omega_\delta$ across $x_{1,1}=0$. The reflected set lies in $R$ and has the same positive measure.
The lower bound in the proof of Proposition~\ref{p:flag-area-sharp} still holds there, since the estimate
in the first coordinate uses only $|x_{1,1}|<\delta$. Thus, for every $2\leq p<\infty$,
$$
\|S_F^{(k)}(f)\|_{L^p(R)}>0=\|f\|_{L^p(R)}.
$$

For the area of the zero extension of $f\in L^p(R)$, Theorem~\ref{t:flag-area-bounds} implies
$$
\|S_F^{(k)}(f\mathbf1_R)\|_{L^p(R)}
\leq\|S_F^{(k)}(f\mathbf1_R)\|_{L^p(\mathbb R^N)}
\leq Cp^{k/2}\|f\|_{L^p(R)}
$$
on every product $R$ of axis-parallel cubes. The order $p^{k/2}$ is sharp on $R=[0,1]^N$ itself.
Indeed, $\Omega_\delta\subset R$, and the same choice $\delta=(t_0/2)e^{-p}$ in the preceding proof gives
$\|S_F^{(k)}(\mathbf1_R)\|_{L^p(R)}\geq cp^{k/2}$.
\end{remark}

For the reverse inequality, we use $M$ dyadic dilates of a single function.
Their sum has size $M$ near the origin, while its Poisson area has size at most $C\sqrt M$.
A correction with bounded coefficient makes the sum have mean zero on a fixed cube.
We also use this construction in Subsection~\ref{s:coord-sharp}.

\begin{lemma}\label{l:one-lacunary}
For each $q\geq1$, there are a cube $E\subset\mathbb R^q$ centred at zero with $\ell(E)\leq1$,
constants $c,C>0$ and an integer $M_0\geq1$ such that the following holds.
For every integer $M\geq M_0$, there is a real-valued $h_M\in\mathcal S(\mathbb R^q)$
with mean zero on $E$ and compact Fourier support away from the origin such that
$$
\|h_M\|_{L^2(\mathbb R^q)}+\|\mathcal S_P^{(q)}h_M\|_{L^2(\mathbb R^q)}\leq C.
$$
Moreover, $\|\mathcal S_P^{(q)}h_M\|_\infty\leq C\sqrt M$, and $h_M\geq M/2$
on a ball $E_M\subset E$ of measure at least $c2^{-qM}$.
\end{lemma}

\begin{proof}
Choose $b=b_q$ with nonnegative radial Fourier transform
$\widehat b\in C_c^\infty(\{1<|\zeta|<2\})$, normalized by $b(0)=1$.
Then $b$ is real and Schwartz. Choose a cube $E$ centred at zero with $\ell(E)\leq1$
and $b\geq1/2$ on $E$. For $G_M(x)=\sum_{n=1}^M b(2^nx)$, let
\begin{equation}\label{e:one-lacunary-sums}
h_M=G_M-\frac{(G_M)_E}{b_E}\,b.
\end{equation}
Then $(h_M)_E=0$ and
$$
\operatorname{supp}\widehat h_M
\subseteq\bigcup_{n=0}^M2^n\operatorname{supp}\widehat b
\subseteq\{\zeta\in\mathbb R^q:1<|\zeta|<2^{M+1}\}.
$$
Since $b_E\geq1/2$, the coefficient of $b$ satisfies
$$
\frac{|(G_M)_E|}{b_E}
\leq\frac{2}{|E|}\sum_{n=1}^M\int_E|b(2^nx)|\,dx
\leq\frac{2\|b\|_1}{|E|}\sum_{n=1}^M2^{-nq}\leq C_q.
$$

Fix $\varepsilon>0$ with $B(0,\varepsilon)\subset E$ and let $E_M=B(0,\varepsilon2^{-M})$.
For $x\in E_M$, we have
$$
|G_M(x)-M|\leq\sum_{n=1}^M|b(2^nx)-b(0)|
\leq\|\nabla b\|_\infty|x|\sum_{n=1}^M2^n\leq2\varepsilon\|\nabla b\|_\infty.
$$
Thus $h_M\geq M-C_q\geq M/2$ on $E_M$ for all sufficiently large $M$,
and $|E_M|=2^{-qM}|B(0,\varepsilon)|$.

Write $T_t=t^2(-\Delta)e^{-t\sqrt{-\Delta}}$. By Fourier inversion,
$$
\|T_t[b(2^n\cdot)]\|_\infty
\leq(2\pi)^{-q}\int_{1<|\xi|<2}(t2^n|\xi|)^2e^{-t2^n|\xi|}|\widehat b(\xi)|\,d\xi
\leq C_q(t2^n)^2e^{-t2^n}.
$$
Including the correction in \eqref{e:one-lacunary-sums} at $n=0$, we obtain from \eqref{e:intro-one-area}
\begin{align*}
\mathcal S_P^{(q)}h_M(x)^2
&=\int_{|a|<1}\int_0^\infty|T_th_M(x-ta)|^2\,\frac{dt}{t}\,da
\leq C_q\int_0^\infty\bigg(\sum_{n=0}^M(t2^n)^2e^{-t2^n}\bigg)^2\frac{dt}{t}\\
&=C_q\sum_{j,n=0}^M\frac{6\,2^{2j+2n}}{(2^j+2^n)^4}
\leq C_q\sum_{j,n=0}^M2^{-2|j-n|}\leq C_qM.
\end{align*}

Finally, by dilation and the bound for the mean correction,
$$
\|h_M\|_{L^2(\mathbb R^q)}
\leq\|b\|_{L^2(\mathbb R^q)}
\bigg(\sum_{n=1}^M2^{-nq/2}+\frac{|(G_M)_E|}{b_E}\bigg)\leq C_q.
$$
By Plancherel's identity and Tonelli's theorem,
$$
\|\mathcal S_P^{(q)}h_M\|_{L^2(\mathbb R^q)}^2
=|B(0,1)|\int_0^\infty r^3e^{-2r}\,dr\,\|h_M\|_{L^2(\mathbb R^q)}^2
=\frac{3|B(0,1)|}{8}\|h_M\|_{L^2(\mathbb R^q)}^2\leq C_q.
$$
This completes the proof of Lemma~\ref{l:one-lacunary}.
\end{proof}


For the sharpness examples we use flag rectangles $R=\prod_{j=1}^k I_j$ with
$\ell(I_1)\leq\cdots\leq\ell(I_k)$. Write
$$
P_{I_j}^{(j)}f(x)=\frac1{|I_j|}\int_{I_j}
f(x_1,\ldots,x_{j-1},u,x_{j+1},\ldots,x_k)\,du.
$$
If $P_{I_j}^{(j)}f\equiv0$ for every $j$, then $Q_Rf=f$ on $R$.

\begin{proposition}\label{p:flag-reverse-sharp}
For every $2\leq p<\infty$, there are a nonzero real $f_p\in\mathcal S(\mathbb R^N)$
and a flag rectangle $R_p=\prod_{j=1}^k I_{j,p}$ such that
$P_{I_{j,p}}^{(j)}f_p\equiv0$ for every $j$ and
$$
\|f_p\|_{L^p(R_p)}\geq cp^{k/2}\|S_F^{(k)}(f_p)\|_{L^p(\mathbb R^N)},
$$
where $c>0$ depends only on $k,d_1,\ldots,d_k$.
The pairs $(f_p,R_p)$ may be chosen so that, for every $\alpha>2/k$ and $c_1>0$,
$$
\frac1{|R_p|}\int_{R_p}
\exp\bigg[c_1\bigg(\frac{|f_p(x)|}{\|S_F^{(k)}(f_p)\|_\infty}\bigg)^\alpha\bigg]\,dx
\to\infty\qquad(p\to\infty).
$$
\end{proposition}

\begin{proof}
In each block $X_j$, take the function and cube from Lemma~\ref{l:one-lacunary},
$$
h_{j,M}(x_j)=\sum_{n=1}^M b_{d_j}(2^nx_j)-c_{j,M}b_{d_j}(x_j),
\qquad |c_{j,M}|\leq C.
$$
Thus $(h_{j,M})_{E_j}=0$, $\|h_{j,M}\|_{L^2(X_j)}\leq C$, and
$h_{j,M}\geq M/2$ on a ball $E_{j,M}\subset E_j$ of measure at least $c2^{-d_jM}$.

Let $L_M=2^{M+2}$ and $\lambda_j=L_M^{k-j}$. On $R_M=\prod_{j=1}^k\lambda_j^{-1}E_j$, consider
$$
f_M(x)=\prod_{j=1}^kh_{j,M}(\lambda_jx_j).
$$
Then $P_{\lambda_j^{-1}E_j}^{(j)}f_M\equiv0$ for every $j$. Since
$$
\frac{\ell(\lambda_j^{-1}E_j)}{\ell(\lambda_{j+1}^{-1}E_{j+1})}
=\frac{\ell(E_j)}{L_M\ell(E_{j+1})},\qquad 1\leq j<k,
$$
$R_M$ is a flag rectangle for all $M\geq M_0$, after increasing $M_0$ if necessary.

The choice of $\lambda_j$ ensures that $|\xi_{V_j}|$ is comparable to $|\xi_j|$
on the Fourier support of each term in the expansion of $f_M$.
Indeed, each such term is a bounded multiple of
$\prod_{j=1}^k b_{d_j}(\lambda_j2^{n_j}x_j)$ with $0\leq n_j\leq M$.
On its Fourier support, $\lambda_j2^{n_j}\leq|\xi_j|\leq2\lambda_j2^{n_j}$, so
\begin{align*}
(\lambda_j2^{n_j})^2
&\leq|\xi_{V_j}|^2=\sum_{i=j}^k|\xi_i|^2
\leq4(\lambda_j2^{n_j})^2+4\lambda_j^22^{2M}\sum_{h=1}^{k-j}L_M^{-2h}\\
&\leq(\lambda_j2^{n_j})^2\bigg(4+\frac{4\cdot2^{2M}}{L_M^2-1}\bigg)
\leq5(\lambda_j2^{n_j})^2.
\end{align*}
Write $\beta_c(u)=ue^{-cu}$ for $c>0$ and $u\geq0$. Each derivative multiplier on $V_j$ has modulus at most
$$
t_j|\xi_{V_j}|e^{-t_j|\xi_{V_j}|}\leq\sqrt5\,\beta_1(t_j\lambda_j2^{n_j}).
$$
Since dilation preserves the $L^1$ norm of the Fourier transform, Fourier inversion gives, uniformly in $y$,
$$
|\mathcal G_{\mathbf t}f_M(y)|
\leq C\sum_{n_1,\ldots,n_k=0}^M\prod_{j=1}^k\beta_1(t_j\lambda_j2^{n_j})
=C\prod_{j=1}^k\bigg(\sum_{n=0}^M\beta_1(t_j\lambda_j2^n)\bigg).
$$
For $c,\lambda,\mu>0$, we have
\begin{equation}\label{e:one-scale-integral}
\int_0^\infty\beta_c(t\lambda)\beta_c(t\mu)\frac{dt}{t}
=\frac{\lambda\mu}{c^2(\lambda+\mu)^2}
\leq c^{-2}\min\{\lambda/\mu,\mu/\lambda\}.
\end{equation}
By \eqref{e:flag-area} and \eqref{e:one-scale-integral},
\begin{align}
\|S_F^{(k)}(f_M)\|_\infty^2
&\leq C\prod_{j=1}^k\int_0^\infty
\bigg(\sum_{n=0}^M\beta_1(t\lambda_j2^n)\bigg)^2\frac{dt}{t}\notag\\
&=C\bigg(\sum_{n,l=0}^M\frac{2^{n+l}}{(2^n+2^l)^2}\bigg)^k
\leq C\bigg(\sum_{n,l=0}^M2^{-|n-l|}\bigg)^k\leq CM^k.\label{e:flag-sharp-sup}
\end{align}
By \eqref{e:flag-plancherel} and the $L^2$ bounds from Lemma~\ref{l:one-lacunary},
$$
\|S_F^{(k)}(f_M)\|_{L^2(\mathbb R^N)}
=C\|f_M\|_{L^2(\mathbb R^N)}
=C|R_M|^{1/2}\prod_{j=1}^k\frac{\|h_{j,M}\|_{L^2(X_j)}}{|E_j|^{1/2}}
\leq C'|R_M|^{1/2}.
$$
Hence, for $2\leq p<\infty$,
\begin{equation}\label{e:flag-sharp-lp}
\|S_F^{(k)}(f_M)\|_p
\leq\|S_F^{(k)}(f_M)\|_\infty^{1-2/p}\|S_F^{(k)}(f_M)\|_2^{2/p}
\leq CM^{k/2-k/p}|R_M|^{1/p}.
\end{equation}

On $\Omega_M=\prod_{j=1}^k\lambda_j^{-1}E_{j,M}$, each factor of $f_M$ is at least $M/2$. Thus
\begin{equation}\label{e:flag-large-values}
\frac{|\{x\in R_M:f_M(x)\geq(M/2)^k\}|}{|R_M|}
\geq\frac{|\Omega_M|}{|R_M|}
=\prod_{j=1}^k\frac{|E_{j,M}|}{|E_j|}\geq c2^{-NM}.
\end{equation}
Choose $M=\max\{M_0,\lceil p\rceil\}$, so $p\leq M\leq Cp$.
Combining \eqref{e:flag-sharp-lp} and \eqref{e:flag-large-values}, we obtain
$$
\frac{\|f_M\|_{L^p(R_M)}}{\|S_F^{(k)}(f_M)\|_p}
\geq\frac{(M/2)^k(c|R_M|2^{-NM})^{1/p}}{CM^{k/2-k/p}|R_M|^{1/p}}
\geq cM^{k/2+k/p}2^{-NM/p}\geq c'p^{k/2}.
$$
The denominator is positive by \eqref{e:flag-plancherel}.

For the same pairs, \eqref{e:flag-sharp-sup} and \eqref{e:flag-large-values} imply
\begin{align*}
\frac1{|R_M|}\int_{R_M}
\exp\bigg[c_1\bigg(\frac{|f_M(x)|}{\|S_F^{(k)}(f_M)\|_\infty}\bigg)^\alpha\bigg]\,dx
&\geq\frac{|\Omega_M|}{|R_M|}
\exp\bigg[c_1\bigg(\frac{(M/2)^k}{CM^{k/2}}\bigg)^\alpha\bigg]\\
&\geq c\exp\big(c_1(c'M^{k/2})^\alpha-NM\log2\big).
\end{align*}
Since $k\alpha/2>1$, the last expression tends to infinity as $p\to\infty$ for every $c_1>0$.
This completes the proof of Proposition~\ref{p:flag-reverse-sharp}.
\end{proof}

\begin{remark}\label{r:fixed-unit}
Both the moment order $p^{k/2}$ and the exponential power $2/k$ are sharp on the fixed unit cube
$R_*=\prod_{j=1}^k C_j$, where $C_j=[-1/2,1/2]^{d_j}$.
For the functions $f_M$ in Proposition~\ref{p:flag-reverse-sharp}, write
$g_{j,M}(x_j)=h_{j,M}(\lambda_jx_j)$. Since $|C_j|=1$ and $\lambda_j\geq1$,
$$
|(g_{j,M})_{C_j}|\leq\|g_{j,M}\|_{L^2(X_j)}
=\lambda_j^{-d_j/2}\|h_{j,M}\|_{L^2(X_j)}\leq C_0.
$$
The cubes $E_j$ are centred at zero and have side length at most one, so
$\Omega_M\subset R_M\subset R_*$. On $\Omega_M$, each $g_{j,M}\geq M/2$. Hence, for all sufficiently large $M$,
$$
Q_{R_*}f_M(x)=\prod_{j=1}^k\big(g_{j,M}(x_j)-(g_{j,M})_{C_j}\big)
\geq(M/2-C_0)^k\geq(M/4)^k.
$$
By \eqref{e:flag-large-values} and the choice of $\lambda_j$,
$$
|\Omega_M|\geq c2^{-NM}|R_M|
=c2^{-NM}\prod_{j=1}^k2^{-(M+2)(k-j)d_j}|E_j|\geq e^{-c_0M}.
$$
Combining these bounds with \eqref{e:flag-sharp-sup}, we obtain
$$
\frac{\|Q_{R_*}f_M\|_{L^p(R_*)}}{\|S_F^{(k)}(f_M)\|_\infty}
\geq\frac{(M/4)^k|\Omega_M|^{1/p}}{CM^{k/2}}
\geq cM^{k/2}e^{-c_0M/p}.
$$
Take $M=\lceil p\rceil$ for large $p$. Since $|R_*|=1$, it follows that
$$
\|Q_{R_*}f_M\|_{L^p(R_*)}\geq cp^{k/2}\|S_F^{(k)}(f_M)\|_\infty
\geq cp^{k/2}\|S_F^{(k)}(f_M)\|_{L^p(R_*)}.
$$
For every $\alpha>2/k$ and $c_1>0$, integration over $\Omega_M$ gives
$$
\frac1{|R_*|}\int_{R_*}\exp\bigg[c_1\bigg(\frac{|Q_{R_*}f_M(x)|}{\|S_F^{(k)}(f_M)\|_\infty}\bigg)^\alpha\bigg]\,dx
\geq\exp\big(c_1(cM^{k/2})^\alpha-c_0M\big)\to\infty,
$$
since $k\alpha/2>1$.
\end{remark}

We next consider the $A_2$ expression for the weight $e^{\lambda f}$ on a flag rectangle $R$.
If the coordinate averages of $f$ vanish on $R$, then $Q_Rf=f$ there.
For $k=2$ and $\|S_F(f)\|_\infty\leq1$, Theorem~\ref{t:concentration}\textup{(i)}
bounds the averages of $e^{\lambda f}$ and $e^{-\lambda f}$ by $2$ for all sufficiently small $|\lambda|$.

\begin{corollary}\label{c:flag-linear-exp}
For each $k\geq2$, there are real Schwartz functions $f_\nu$ with
$\|S_F^{(k)}(f_\nu)\|_{L^\infty(\mathbb R^N)}=1$ and flag rectangles
$R_\nu=\prod_{j=1}^k I_{j,\nu}$ such that $P_{I_{j,\nu}}^{(j)}f_\nu\equiv0$ for every $j$.
If $k=2$, let $|\lambda|>\lambda_0$ for a constant $\lambda_0$ depending only on $d_1,d_2$.
If $k\geq3$, let $\lambda$ be any nonzero real number. Then
$$
\bigg(\frac1{|R_\nu|}\int_{R_\nu}e^{\lambda f_\nu}\,dx\bigg)
\bigg(\frac1{|R_\nu|}\int_{R_\nu}e^{-\lambda f_\nu}\,dx\bigg)\to\infty
\qquad(\nu\to\infty).
$$
For the same sequence,
$$
\frac1{|R_\nu|}\int_{R_\nu}e^{|\lambda f_\nu|}\,dx\to\infty.
$$
\end{corollary}

\begin{proof}
Let $h_\nu,R_\nu$ be the examples from the proof of Proposition~\ref{p:flag-reverse-sharp}
with $M=\nu$. By \eqref{e:flag-plancherel} and \eqref{e:flag-sharp-sup},
$0<\|S_F^{(k)}(h_\nu)\|_\infty\leq C\nu^{k/2}$. For
$$
f_\nu=\frac{h_\nu}{\|S_F^{(k)}(h_\nu)\|_\infty},
$$
we have $\|S_F^{(k)}(f_\nu)\|_\infty=1$, and each coordinate average on $R_\nu$ still vanishes.
By \eqref{e:flag-large-values}, $f_\nu\geq c'\nu^{k/2}$ on $\Omega_\nu$, whence
$$
\frac1{|R_\nu|}\int_{R_\nu}e^{|\lambda|f_\nu}\,dx
\geq\frac{|\Omega_\nu|}{|R_\nu|}\exp(c'|\lambda|\nu^{k/2})
\geq c\exp\big(c'|\lambda|\nu^{k/2}-N\nu\log2\big).
$$
For $k=2$, the last expression tends to infinity when $|\lambda|>N\log2/c'$.
For $k\geq3$, it tends to infinity for every $\lambda\ne0$.
Since $(f_\nu)_{R_\nu}=0$, Jensen's inequality gives
$$
\frac1{|R_\nu|}\int_{R_\nu}e^{-|\lambda|f_\nu}\,dx
\geq e^{-|\lambda|(f_\nu)_{R_\nu}}=1.
$$
The product therefore diverges for either sign of $\lambda$.
The last conclusion follows from $e^{|\lambda f_\nu|}\geq e^{|\lambda|f_\nu}$.
\end{proof}


\section{Coordinate averages and the full oscillation}\label{s:coord-oscillation}

On $R=I\times J$, the full oscillation has the decomposition
$$
f(x,y)-f_R=Q_xQ_yf(x,y)+\big(f_J(x)-f_R\big)+\big(f_I(y)-f_R\big).
$$
Theorem~\ref{t:two-local} controls the first term with a factor $p$.
Each of the other two terms depends on one variable and is controlled by the Poisson
area of the corresponding average, with a factor $\sqrt p$. More generally, for a
$k$-step flag, the term with $\ell$ coordinate blocks left after averaging has a
bound of order $p^{\ell/2}$.

\subsection{Reverse estimates for the one-parameter Poisson area function}\label{s:one-poisson-area}

For the terms depending on one variable, we use the area function $\mathcal S_P^{(q)}$
defined in \eqref{e:intro-one-area} by means of $t^2(-\Delta)e^{-t\sqrt{-\Delta}}$.
For $L=-\Delta$ and aperture one, \cite[Theorem~5.2]{CDWY} gives the global reverse
estimate with factor $C_q\sqrt p$ as $p\to\infty$. Here we need the oscillation and
the area on the same cube. For a cube $E\subset\mathbb R^q$, write
$P_Eg=\mathbf1_Eg_E$ and $Q_Eg=\mathbf1_E(g-g_E)$, where $g_E=|E|^{-1}\int_Eg$.

\begin{proposition}\label{p:one-area}
Let $g\in L^u(\mathbb R^q)$ for some $1\leq u<\infty$, and let $2\leq p<\infty$.
For every axis-parallel cube $E\subset\mathbb R^q$,
\begin{equation}\label{e:one-local}
\|Q_Eg\|_{L^p(E)}\leq C_q\sqrt p\,\|\mathcal S_P^{(q)}g\|_{L^p(E)}.
\end{equation}
On $\mathbb R^q$, we have
\begin{equation}\label{e:one-reverse}
\|g\|_p\leq C_q\sqrt p\,\|\mathcal S_P^{(q)}g\|_p.
\end{equation}
The constant depends only on $q$.
\end{proposition}

\begin{proof}
By Remark~\ref{r:time-area} with $k=1$, it suffices to show that
$S_{F,1/4}^{\mathrm{time}}(g)(x)\leq C_q\mathcal S_P^{(q)}g(x)$.
Let $v(y,t)=(P_t^{(q)}*g)(y)$. By scaling,
$\|\partial_sv(\cdot,s)\|_\infty\leq C_{q,u}s^{-1-q/u}\|g\|_u\to0$ as $s\to\infty$.
Integrating from $t$ to infinity and applying Cauchy--Schwarz, we obtain
$$
|t\partial_tv(y,t)|^2
=\bigg|\int_t^\infty\frac{t}{s}\,s^2\partial_s^2v(y,s)\,\frac{ds}{s}\bigg|^2
\leq\int_t^\infty\frac{t}{s}|s^2\partial_s^2v(y,s)|^2\,\frac{ds}{s}.
$$

For $|y-x|<t/4$ and $s\geq t$, the ball of radius $s/8$ about $(y,s)$ lies in
$\{(w,r):7s/8<r<9s/8,\ |w-x|<r\}$. Since $\partial_s^2v$ is harmonic,
the mean-value inequality implies
$$
|s^2\partial_s^2v(y,s)|^2
\leq C_q\int_{7s/8}^{9s/8}\int_{|w-x|<r}|r^2\partial_r^2v(w,r)|^2
\,\frac{dw\,dr}{r^{q+1}}.
$$
Integrating over the cone $|y-x|<t/4$ and using $\int_0^s(t/s)\,dt/t=1$, we obtain
\begin{align*}
S_{F,1/4}^{\mathrm{time}}(g)(x)^2
&\leq C_q\int_0^\infty\int_{7s/8}^{9s/8}\int_{|w-x|<r}
|r^2\partial_r^2v(w,r)|^2\,\frac{dw\,dr}{r^{q+1}}\frac{ds}{s}\\
&=C_q\log(9/7)\int_0^\infty\int_{|w-x|<r}
|r^2(-\Delta)e^{-r\sqrt{-\Delta}}g(w)|^2\,\frac{dw\,dr}{r^{q+1}}
=C_q'\mathcal S_P^{(q)}g(x)^2,
\end{align*}
where $8r/9<s<8r/7$ for each fixed $r$. This proves \eqref{e:one-local} and \eqref{e:one-reverse}.
\end{proof}


\subsection{The full oscillation on a product of two cubes}\label{s:two-coord-averages}

For $R_0=I_0\times J_0$, write
$$
f_{J_0}(x)=\frac1{|J_0|}\int_{J_0}f(x,v)\,dv,\qquad
f_{I_0}(y)=\frac1{|I_0|}\int_{I_0}f(u,y)\,du.
$$
For $f\in L^u(\mathbb R^{n+m})$, $1\leq u<\infty$, Jensen's inequality gives
$$
\|f_{J_0}\|_{L^u(\mathbb R^n)}^u
\leq\frac1{|J_0|}\int_{\mathbb R^n}\int_{J_0}|f(x,v)|^u\,dv\,dx
\leq\frac1{|J_0|}\|f\|_{L^u(\mathbb R^{n+m})}^u.
$$
Similarly, $\|f_{I_0}\|_{L^u(\mathbb R^m)}\leq |I_0|^{-1/u}\|f\|_{L^u(\mathbb R^{n+m})}$.
By Fubini's theorem, $(f_{J_0})_{I_0}=(f_{I_0})_{J_0}=f_{R_0}$. Thus on $R_0$,
\begin{equation}\label{e:coord-splitting}
f(x,y)-f_{R_0}=Q_xQ_yf(x,y)+Q_{I_0}f_{J_0}(x)+Q_{J_0}f_{I_0}(y).
\end{equation}
For the last two terms we use the Poisson areas of the coordinate averages,
$$
\mathcal E_{P,R_0}(f)(x,y)^2
=\mathbf1_{J_0}(y)\big(\mathcal S_P^{(n)}f_{J_0}(x)\big)^2
+\mathbf1_{I_0}(x)\big(\mathcal S_P^{(m)}f_{I_0}(y)\big)^2.
$$

\begin{theorem}\label{t:two-coordinate}
Let $f\in L^u(\mathbb R^{n+m})$ for some $1\leq u<\infty$, and let $R_0=I_0\times J_0$
be a product of axis-parallel cubes. For $2\leq p<\infty$,
\begin{equation}\label{e:two-scalar}
\|f-f_{R_0}\|_{L^p(R_0)}
\leq Cp\|S_F(f)\|_{L^p(R_0)}+C\sqrt p\,\|\mathcal E_{P,R_0}(f)\|_{L^p(R_0)}.
\end{equation}
If $0<\|\mathcal E_{P,R_0}(f)\|_{L^\infty(R_0)}<\infty$, then
\begin{equation}\label{e:coord-exponential}
\frac1{|R_0|}\int_{R_0}
\exp\bigg(c\frac{|f_{J_0}(x)+f_{I_0}(y)-2f_{R_0}|^2}
{\|\mathcal E_{P,R_0}(f)\|_{L^\infty(R_0)}^2}\bigg)\,dx\,dy\leq2.
\end{equation}
If both $S_F(f)$ and $\mathcal E_{P,R_0}(f)$ are bounded on $R_0$, then for every $\tau>0$,
\begin{equation}\label{e:two-tails}
|\{z\in R_0:|f(z)-f_{R_0}|>\tau\}|
\leq 2|R_0|\exp\bigg(-\frac{c\tau}{\|S_F(f)\|_{L^\infty(R_0)}}\bigg)
+2|R_0|\exp\bigg(-\frac{c\tau^2}{\|\mathcal E_{P,R_0}(f)\|_{L^\infty(R_0)}^2}\bigg).
\end{equation}
A term with zero denominator is omitted. The constants depend only on $n,m$.
\end{theorem}

\begin{proof}
Applying \eqref{e:one-local} to $f_{J_0}$ on $I_0$ and to $f_{I_0}$ on $J_0$, we obtain
\begin{align}
\|f_{J_0}+f_{I_0}-2f_{R_0}\|_{L^p(R_0)}
&\leq |J_0|^{1/p}\|f_{J_0}-f_{R_0}\|_{L^p(I_0)}
     +|I_0|^{1/p}\|f_{I_0}-f_{R_0}\|_{L^p(J_0)}\notag\\
&\leq C\sqrt p\,\Big(|J_0|^{1/p}\|\mathcal S_P^{(n)}f_{J_0}\|_{L^p(I_0)}
     +|I_0|^{1/p}\|\mathcal S_P^{(m)}f_{I_0}\|_{L^p(J_0)}\Big)\notag\\
&\leq C\sqrt p\,\|\mathcal E_{P,R_0}(f)\|_{L^p(R_0)}.\label{e:coord-mean-lp}
\end{align}
Combining \eqref{e:coord-splitting} and \eqref{e:coord-mean-lp} with
\eqref{e:flag-local-lp} for $k=2$, we obtain \eqref{e:two-scalar}.

Suppose $0<\|\mathcal E_{P,R_0}(f)\|_{L^\infty(R_0)}<\infty$.
Taking $p=2j$ in \eqref{e:coord-mean-lp}, we have for every integer $j\geq1$,
$$
\frac1{|R_0|}\int_{R_0}
\bigg(\frac{|f_{J_0}(x)+f_{I_0}(y)-2f_{R_0}|}
{\|\mathcal E_{P,R_0}(f)\|_{L^\infty(R_0)}}\bigg)^{2j}\,dx\,dy
\leq (C\sqrt{2j})^{2j}\leq(C'j)^j.
$$
Summing the exponential series as in \eqref{e:flag-exp-series}, we obtain \eqref{e:coord-exponential}.

Finally, by \eqref{e:coord-splitting},
\begin{align*}
|\{z\in R_0:|f(z)-f_{R_0}|>\tau\}|
&\leq |\{z\in R_0:|Q_xQ_yf(z)|>\tau/2\}|\\
&\quad+|\{(x,y)\in R_0:|f_{J_0}(x)+f_{I_0}(y)-2f_{R_0}|>\tau/2\}|.
\end{align*}
For the first term we apply Markov's inequality to \eqref{e:flag-exponential} with $k=2$,
and for the second we use \eqref{e:coord-exponential}. If either area norm is zero,
the corresponding term vanishes by \eqref{e:flag-local-lp} or \eqref{e:coord-mean-lp}
at $p=2$. This proves \eqref{e:two-tails}.
\end{proof}


\subsection{Neither coordinate average can be omitted}\label{s:coord-necessity}

We construct two products. In the first, the average in $x$ vanishes and the average in $y$
accounts for an oscillation of order $M$. In the second, the roles of the averages are reversed.
In both examples the flag area is at most $C\sqrt M$.

Let $E_q=[-1/2,1/2]^q$, $I_0=E_n$, $J_0=E_m$ and $R_0=I_0\times J_0$.
For the function $b_q$ in the proof of Lemma~\ref{l:one-lacunary}, write $a_q=\partial_1b_q$.
Since $b_q$ is radial, $a_q$ is odd in its first coordinate and $(a_q)_{E_q}=0$. Moreover,
$$
\partial_1a_q(0)=-(2\pi)^{-q}\int_{\mathbb R^q}\zeta_1^2\widehat b_q(\zeta)\,d\zeta<0,
$$
so $\|a_q\|_{L^2(E_q)}>0$.

Fix $\rho\geq2$ and $0<\varepsilon\leq1$, to be chosen below. For $\rho_r=\rho 2^{r-1}$, let
$$
u_M(x,y)=a_n(x)H_M(y),\qquad H_M(y)=\sum_{r=1}^M b_m(y/\rho_r).
$$
We use $u_M$ on $R_0$. For the second example, let
$$
v_M(x,y)=G_M(x)a_m(y),\qquad G_M(x)=\sum_{r=1}^M b_n(2^rx).
$$
We use $v_M$ on $R_M=I_M\times J_0$, where $I_M$ is the cube centred at zero
of side length $\varepsilon2^{-M}$.

\begin{proposition}\label{p:coord-necessity}
For suitable fixed $\rho,\varepsilon$ and every $M\geq1$, we have
$(u_M)_{I_0}\equiv0$, $(v_M)_{J_0}\equiv0$ and
$$
\|S_F(u_M)\|_\infty+\|S_F(v_M)\|_\infty\leq C\sqrt M.
$$
For either pair $(g,R)=(u_M,R_0)$ or $(v_M,R_M)$, the doubly centred part satisfies
$\|Q_Rg\|_{L^\infty(R)}\leq C$, while
$$
\inf_{\zeta\in\mathbb C}\bigg(\frac1{|R|}\int_R|g(z)-\zeta|^2\,dz\bigg)^{1/2}\geq c_0M.
$$
In each case the area of the nonzero coordinate average satisfies
$$
\|\mathcal E_{P,R}(g)\|_{L^\infty(R)}\simeq M.
$$
The constants are independent of $M$.
\end{proposition}

\begin{proof}
Since $b_n(0)=b_m(0)=1$, the mean value theorem implies
\begin{equation}\label{e:coord-flat}
\sup_{y\in J_0}|H_M(y)-M|
\leq\sum_{r=1}^M\sup_{y\in J_0}|b_m(y/\rho_r)-1|
\leq C\sum_{r=1}^M\rho_r^{-1}\leq C/\rho.
\end{equation}
Similarly,
$$
\sup_{x\in I_M}|G_M(x)-M|
\leq\sum_{r=1}^M\sup_{x\in I_M}|b_n(2^rx)-1|
\leq C\varepsilon2^{-M}\sum_{r=1}^M2^r\leq C\varepsilon.
$$
Choose $\rho$ large and $\varepsilon$ small so that both errors are at most $1/4$.

Since $(a_n)_{I_0}=(a_m)_{J_0}=0$, the only nonzero coordinate averages are
$(u_M)_{J_0}=(H_M)_{J_0}a_n$ and $(v_M)_{I_M}=(G_M)_{I_M}a_m$. Hence
\begin{align*}
\|Q_{R_0}u_M\|_{L^\infty(R_0)}+\|Q_{R_M}v_M\|_{L^\infty(R_M)}
&=\|a_n\|_{L^\infty(I_0)}\|H_M-(H_M)_{J_0}\|_{L^\infty(J_0)}\\
&\quad+\|G_M-(G_M)_{I_M}\|_{L^\infty(I_M)}\|a_m\|_{L^\infty(J_0)}\\
&\leq\tfrac12\big(\|a_n\|_\infty+\|a_m\|_\infty\big).
\end{align*}
For either pair $(g,R)$, we have $g_R=0$. Since $|I_0|=|J_0|=1$, for every $\zeta\in\mathbb C$,
$$
\frac1{|R|}\int_R|g-\zeta|^2
=|\zeta|^2+\frac1{|R|}\int_R|g|^2
\geq(M-1/4)^2\min\big\{\|a_n\|_{L^2(I_0)}^2,\|a_m\|_{L^2(J_0)}^2\big\}
\geq c_0^2M^2.
$$

On the Fourier support of $a_n(x)b_m(y/\rho_r)$, we have $|(\xi,\eta)|\simeq1$ and $|\eta|\simeq\rho_r^{-1}$.
For $b_n(2^rx)a_m(y)$, the corresponding sizes are $2^r$ and $1$.
Each derivative multiplier is bounded by
$t|(\xi,\eta)|e^{-t|(\xi,\eta)|}\,s|\eta|e^{-s|\eta|}$,
and the Fourier $L^1$ norms of both products are independent of $r$.
By Fourier inversion, \eqref{e:flag-area} and \eqref{e:one-scale-integral}, we obtain
\begin{align*}
\|S_F(u_M)\|_\infty^2+\|S_F(v_M)\|_\infty^2
&\leq C\int_0^\infty\int_0^\infty
\bigg[\beta_c(t)^2\bigg(\sum_{r=1}^M\beta_c(s/\rho_r)\bigg)^2
+\beta_c(s)^2\bigg(\sum_{r=1}^M\beta_c(2^rt)\bigg)^2\bigg]\frac{dt}{t}\frac{ds}{s}\\
&\leq C\sum_{j,r=1}^M\frac{2^{j+r}}{(2^j+2^r)^2}
\leq C\sum_{j,r=1}^M2^{-|j-r|}\leq CM,
\end{align*}
where $\beta_c(t)=te^{-ct}$ and $\rho_r^{-1}=(2/\rho)2^{-r}$.

Finally, $(a_q)_{E_q}=0$ and $\widehat a_q$ is supported in $\{1<|\xi|<2\}$.
By \eqref{e:one-local} at $p=2$ and Fourier inversion,
\begin{align*}
0<\|a_q\|_{L^2(E_q)}
&\leq C_q\|\mathcal S_P^{(q)}a_q\|_{L^2(E_q)}
\leq C_q\|\mathcal S_P^{(q)}a_q\|_{L^\infty(E_q)}\\
&\leq C_q\bigg(\int_0^\infty\|t^2(-\Delta)e^{-t\sqrt{-\Delta}}a_q\|_\infty^2\,\frac{dt}{t}\bigg)^{1/2}
\leq C_q\bigg(\int_0^\infty t^4e^{-2t}\,\frac{dt}{t}\bigg)^{1/2}<\infty.
\end{align*}
Substituting the nonzero coordinate average in each case, we obtain
\begin{align*}
\|\mathcal E_{P,R}(g)\|_{L^\infty(R)}
&=\begin{cases}\displaystyle
|(H_M)_{J_0}|\,\|\mathcal S_P^{(n)}a_n\|_{L^\infty(I_0)}, &(g,R)=(u_M,R_0),\\[4pt]
|(G_M)_{I_M}|\,\|\mathcal S_P^{(m)}a_m\|_{L^\infty(J_0)}, &(g,R)=(v_M,R_M)
\end{cases}\\
&\simeq M.
\end{align*}
Omitting either coordinate average in \eqref{e:two-scalar} at $p=2$ would therefore imply
$c_0M\leq C\sqrt M$ for every $M$. This completes the proof.
\end{proof}


\smallskip

Enlarging the rectangle on the area side cannot compensate for the missing coordinate averages.
Fix $\Lambda\geq1$ and write $\Lambda R=(\Lambda I)\times(\Lambda J)$, with each cube dilated about its centre.
If
$$
\|f-f_R\|_{L^2(R)}\leq C_\Lambda\|S_F(f)\|_{L^2(\Lambda R)}
$$
held for every real Schwartz function $f$ and every flag rectangle $R$, then
Proposition~\ref{p:coord-necessity}, applied to either $(u_M,R_0)$ or $(v_M,R_M)$, would imply
$$
c_0M\leq\frac{\|f-f_R\|_{L^2(R)}}{|R|^{1/2}}
\leq C_\Lambda\frac{\|S_F(f)\|_{L^2(\Lambda R)}}{|R|^{1/2}}
\leq C_\Lambda\Lambda^{(n+m)/2}\|S_F(f)\|_\infty
\leq C'_\Lambda\sqrt M,
$$
a contradiction as $M\to\infty$.

For $k=2$, Theorem~\ref{t:concentration}\textup{(i)} controls the exponential averages
of the doubly centred oscillation for small $|\lambda|$. The first example shows what happens
for $f-f_{R_0}$ when only one coordinate average vanishes.
Let $f_M=u_M/\|S_F(u_M)\|_\infty$, where the denominator is positive by \eqref{e:flag-plancherel}.
Then $(f_M)_{R_0}=0$, and Proposition~\ref{p:coord-necessity} implies
$$
\frac1{|R_0|}\int_{R_0}|f_M|^2
=\frac{\|u_M\|_{L^2(R_0)}^2}{|R_0|\,\|S_F(u_M)\|_\infty^2}
\geq cM.
$$
For every real $\lambda\ne0$, by symmetry in $z,z'$ and $\cosh t\geq1+t^2/2$, we obtain
\begin{align}
\bigg(\frac1{|R_0|}\int_{R_0}e^{\lambda f_M}\bigg)
\bigg(\frac1{|R_0|}\int_{R_0}e^{-\lambda f_M}\bigg)
&=\frac1{|R_0|^2}\int_{R_0}\int_{R_0}
\cosh\big(\lambda(f_M(z)-f_M(z'))\big)\,dz\,dz'\notag\\
&\geq 1+\frac{\lambda^2}{2|R_0|^2}\int_{R_0}\int_{R_0}|f_M(z)-f_M(z')|^2\,dz\,dz'\notag\\
&=1+\frac{\lambda^2}{|R_0|}\int_{R_0}|f_M|^2\geq1+c\lambda^2M.\label{e:two-weight-growth}
\end{align}
Thus the $A_2$ expression diverges on the fixed rectangle $R_0$ for every $\lambda\ne0$.

Moreover, $|a_n|\geq c_n>0$ on a fixed subset of $I_0$ of positive measure.
Since $H_M\geq M-1/4$ on $J_0$ and $\|S_F(u_M)\|_\infty\leq C\sqrt M$, we have
$|f_M|\geq c\sqrt M$ on a fixed subset of $R_0$ of positive measure. Hence
$$
\frac1{|R_0|}\int_{R_0}e^{|\lambda f_M|}\,dz\geq c\exp(c'|\lambda|\sqrt M).
$$



\subsection{Sharpness of the flag and coordinate terms}\label{s:coord-sharp}

When both coordinate averages vanish, \eqref{e:two-scalar} reduces to the reverse estimate
for the flag area. To test the coordinate term, we keep the rectangle fixed and choose
functions whose flag areas tend uniformly to zero while one coordinate average converges
to the function $h_L$ of Lemma~\ref{l:one-lacunary}.

\begin{proposition}\label{p:two-orders}
Fix $\Lambda\geq1$. Let $f$ range over real Schwartz functions on $\mathbb R^{n+m}$
and $R$ over flag rectangles.

\begin{enumerate}[(i)]
\item Fix $2\leq p<\infty$. If
$$
\|f-f_R\|_{L^p(R)}
\leq a_p\|S_F(f)\|_{L^p(\Lambda R)}
+b_p\|\mathcal E_{P,R}(f)\|_{L^p(\mathbb R^{n+m})}
$$
for all $f,R$, then $a_p\geq c_0p$ and $b_p\geq c_0\sqrt p$, with $c_0=c_0(n,m)>0$.
For the bound on $b_p$, it suffices to test on a single flag rectangle $R_*=R_*(n,m)$.

\smallskip

\item Fix $\alpha,\beta>0$, with the convention $\tau/0=+\infty$. If
\begin{align*}
|\{z\in R:|f(z)-f_R|>\tau\}|
&\leq C|R|\exp\bigg[-c\bigg(\frac{\tau}{\|S_F(f)\|_{L^\infty(\Lambda R)}}\bigg)^\alpha\bigg]\\
&\quad+C|R|\exp\bigg[-c\bigg(\frac{\tau}{\|\mathcal E_{P,R}(f)\|_{L^\infty(\mathbb R^{n+m})}}\bigg)^\beta\bigg]
\end{align*}
for all $f,R$ and $\tau>0$ with uniform constants $c,C>0$, then $\alpha\leq1$ and $\beta\leq2$.
For any $\alpha>0$, the bound on $\beta$ follows by testing on $R_*$ with both area norms positive there.
The power $2$ in \eqref{e:coord-exponential} is also sharp on $R_*$.
\end{enumerate}
\end{proposition}


\begin{proof}
\noindent\textit{Step 1.}
Use the functions $f_L(x,y)=h_{1,L}(2^{L+2}x)h_{2,L}(y)$ and flag rectangles
$R_L=(2^{-L-2}E_1)\times E_2$ from the proof of Proposition~\ref{p:flag-reverse-sharp}.
Both coordinate averages vanish identically, so $(f_L)_{R_L}=0$ and
$\mathcal E_{P,R_L}(f_L)=0$. For the choice of $L$ in that proposition, \textup{(i)} implies
$$
cp\|S_F(f_L)\|_p\leq\|f_L\|_{L^p(R_L)}
\leq a_p\|S_F(f_L)\|_{L^p(\Lambda R_L)}\leq a_p\|S_F(f_L)\|_p.
$$
Hence $a_p\geq cp$.

For the tail estimate, \eqref{e:flag-sharp-sup} with $k=2$ reads
\begin{equation}\label{e:trace-free-sup}
0<\|S_F(f_L)\|_\infty\leq CL.
\end{equation}
Let $\Omega_L=\{z\in R_L:f_L(z)\geq L^2/4\}$. By \eqref{e:flag-large-values},
\begin{equation}\label{e:zero-trace-values}
|\Omega_L|\geq c|R_L|2^{-(n+m)L}.
\end{equation}
Applying \textup{(ii)} at $\tau=L^2/8$, we obtain
$$
c|R_L|2^{-(n+m)L}\leq|\Omega_L|
\leq C|R_L|\exp\bigg[-c_1\bigg(\frac{L^2}{8\|S_F(f_L)\|_\infty}\bigg)^\alpha\bigg]
\leq C|R_L|e^{-c_2L^\alpha}.
$$
This requires $\alpha\leq1$.

\smallskip

\noindent\textit{Step 2.}
Take $h_L,E,E_L$ from Lemma~\ref{l:one-lacunary} in dimension $n$.
The rectangle $R_*=E\times J_0$ is fixed and $\ell(E)\leq1=\ell(J_0)$.
For $L\geq M_0$ and $M\geq1$, let
$$
F_{L,M}(x,y)=h_L(x)\frac{H_M(y)}M.
$$
By \eqref{e:coord-flat}, $|H_M-M|\leq1/4$ on $J_0$. Since $(h_L)_E=0$, the rectangular
mean vanishes and the only nonzero coordinate average is $(F_{L,M})_{J_0}=(H_M)_{J_0}h_L/M$. Thus
\begin{equation}\label{e:trace-area}
\mathcal E_{P,R_*}(F_{L,M})(x,y)
=\frac{(H_M)_{J_0}}M\,\mathbf1_{J_0}(y)\mathcal S_P^{(n)}h_L(x).
\end{equation}

On the Fourier support of $b_n(2^jx)b_m(y/\rho_r)$, we have
$|(\xi,\eta)|\simeq2^j$ and $|\eta|\simeq\rho_r^{-1}$, including $j=0$ for the correction in $h_L$.
Each product has uniformly bounded Fourier $L^1$ norm, and the correction coefficient is bounded independently of $L$.
As in the proof of Proposition~\ref{p:coord-necessity}, Fourier inversion and \eqref{e:one-scale-integral} imply
\begin{equation}\label{e:trace-flag-decay}
\|S_F(F_{L,M})\|_\infty^2
\leq\frac C{M^2}\int_0^\infty\bigg(\sum_{j=0}^L\beta_c(2^jt)\bigg)^2\frac{dt}{t}
\int_0^\infty\bigg(\sum_{r=1}^M\beta_c(s/\rho_r)\bigg)^2\frac{ds}{s}
\leq\frac{C(L+1)M}{M^2}\leq\frac{CL}{M}.
\end{equation}
For the coordinate area, \eqref{e:trace-area} and Lemma~\ref{l:one-lacunary} imply
\begin{equation}\label{e:trace-area-sup}
\|\mathcal E_{P,R_*}(F_{L,M})\|_\infty
=\frac{(H_M)_{J_0}}M\|\mathcal S_P^{(n)}h_L\|_\infty\leq C\sqrt L.
\end{equation}

Fix $p$ and $L$. Since $|J_0|=1$ and $(H_M)_{J_0}/M\geq3/4$, Jensen's inequality in $y$ and \textup{(i)} imply
\begin{align*}
\|h_L\|_{L^p(E)}
&=\frac{M}{(H_M)_{J_0}}\|(F_{L,M})_{J_0}\|_{L^p(E)}
\leq\frac{M}{(H_M)_{J_0}}\|F_{L,M}\|_{L^p(R_*)}\\
&\leq Ca_p|\Lambda R_*|^{1/p}\sqrt{L/M}
+b_p\|\mathcal S_P^{(n)}h_L\|_{L^p(\mathbb R^n)}.
\end{align*}
Letting $M\to\infty$ and using $h_L\geq L/2$ on $E_L$, we obtain
\begin{align*}
\tfrac L2(c2^{-nL})^{1/p}
&\leq\|h_L\|_{L^p(E)}
\leq b_p\|\mathcal S_P^{(n)}h_L\|_{L^p(\mathbb R^n)}\\
&\leq b_p\|\mathcal S_P^{(n)}h_L\|_2^{2/p}\|\mathcal S_P^{(n)}h_L\|_\infty^{1-2/p}
\leq Cb_pL^{1/2-1/p}.
\end{align*}
Take $L=\lceil Kp\rceil$ with $K\geq M_0$ fixed. Then $cp\leq Cb_p\sqrt p$, so $b_p\geq c_0\sqrt p$.
The same lower bound for $\|h_L\|_{L^p(E)}$ and upper bound for $\|\mathcal S_P^{(n)}h_L\|_p$
also prove sharpness of $\sqrt p$ in \eqref{e:one-local} and \eqref{e:one-reverse}, since $(h_L)_E=0$.

For the tail estimate, $H_M/M\geq3/4$ on $J_0$, so $F_{L,M}\geq3L/8$ on $E_L\times J_0$. Hence
\begin{equation}\label{e:trace-values}
|\{z\in R_*:|F_{L,M}(z)|>L/4\}|\geq|E_L||J_0|\geq c|R_*|2^{-nL}.
\end{equation}
Both area norms are positive on $R_*$. For the coordinate area this follows from
\eqref{e:trace-area} and \eqref{e:one-local}, since $(h_L)_E=0$ and $h_L\geq L/2$ on $E_L$.
For the flag area, $\partial_1^2H_M(0)=\partial_1^2b_m(0)\sum_{r=1}^M\rho_r^{-2}<0$,
so $H_M$ is nonconstant on $J_0$. By \eqref{e:two-step-moments} at $p=2$,
$$
0<\frac{\|h_L\|_{L^2(E)}}M\|H_M-(H_M)_{J_0}\|_{L^2(J_0)}
=\|Q_xQ_yF_{L,M}\|_{L^2(R_*)}\leq C\|S_F(F_{L,M})\|_{L^2(R_*)}.
$$

Taking $\tau=L/4$ in \textup{(ii)} and using \eqref{e:trace-values},
\eqref{e:trace-flag-decay} and \eqref{e:trace-area-sup}, we obtain
$$
c2^{-nL}\leq C\exp\big(-c_1(LM)^{\alpha/2}\big)+C\exp\big(-c_2L^{\beta/2}\big).
$$
For every fixed $\alpha>0$, letting $M\to\infty$ gives
$c2^{-nL}\leq Ce^{-c_2L^{\beta/2}}$ for all $L\geq M_0$. Hence $\beta\leq2$.

Finally, take $M=1$. The factor $(H_1)_{J_0}$ in the coordinate average cancels
with the same factor in \eqref{e:trace-area}. For every $\beta>2$ and $c_*>0$,
\begin{align*}
\frac1{|R_*|}\int_{R_*}
\exp\bigg[c_*\bigg(\frac{|(F_{L,1})_{J_0}(x)|}{\|\mathcal E_{P,R_*}(F_{L,1})\|_\infty}\bigg)^\beta\bigg]\,dx\,dy
&=\frac1{|E|}\int_E\exp\bigg[c_*\bigg(\frac{|h_L(x)|}{\|\mathcal S_P^{(n)}h_L\|_\infty}\bigg)^\beta\bigg]\,dx\\
&\geq c2^{-nL}\exp\big(c_*(c'\sqrt L)^\beta\big)\to\infty.
\end{align*}
This proves sharpness of the power $2$ in \eqref{e:coord-exponential} on $R_*$.
\end{proof}


\subsection{Decomposition of the full oscillation for a finite flag}\label{s:coord-expansion}

Fix $k\geq2$ and let $R=\prod_{j=1}^kE_j$, where the $E_j\subset X_j$ are axis-parallel cubes
of arbitrary side lengths. For $A\subseteq\{1,\ldots,k\}$, write
$X_A=\bigoplus_{j\in A}X_j$, $R_A=\prod_{j\in A}E_j$ and $x_A=(x_j)_{j\in A}$.
The index $A$ records the coordinates left after averaging. Thus, for each nonempty proper subset $A$,
$$
f_A(x_A)=\frac1{|R_{A^c}|}\int_{R_{A^c}}f(x_A,u_{A^c})\,du_{A^c},\qquad x_A\in X_A,
$$
where $A^c=\{1,\ldots,k\}\setminus A$. At the two extremes, $f_\emptyset=f_R$ and
$f_{\{1,\ldots,k\}}=f$.

On $R$, each $P_j+Q_j$ is the identity and $\prod_{j\in A^c}P_jf=f_A$.
Expanding $\prod_{j=1}^k(P_j+Q_j)f$ and subtracting $f_R$, we obtain
\begin{equation}\label{e:coord-expansion}
f-f_R=\sum_{\emptyset\ne A\subseteq\{1,\ldots,k\}}\prod_{j\in A}Q_jf_A
=\sum_{\ell=1}^k F_R^{(\ell)},
\end{equation}
where the terms centred in exactly $\ell$ coordinates sum to
$$
F_R^{(\ell)}=\sum_{|A|=\ell}\prod_{j\in A}Q_jf_A.
$$
For $k=2$, $f_{\{1\}}=f_J$ and $f_{\{2\}}=f_I$, so \eqref{e:coord-expansion} is precisely
\eqref{e:coord-splitting}. For $k=3$, the term centred in all three coordinates is
$$
Q_1Q_2Q_3f=f-f_{\{1,2\}}-f_{\{1,3\}}-f_{\{2,3\}}+f_{\{1\}}+f_{\{2\}}+f_{\{3\}}-f_R.
$$

For $A=\{i_1<\cdots<i_\ell\}$ with $\ell\geq2$, let $S_A$ be the full-gradient Poisson
area function \eqref{e:flag-area} on the induced flag
$$
X_{i_1}\oplus\cdots\oplus X_{i_\ell}
\supset X_{i_2}\oplus\cdots\oplus X_{i_\ell}\supset\cdots\supset X_{i_\ell}.
$$
For a singleton, use $S_{\{i\}}=\mathcal S_P^{(d_i)}$ from \eqref{e:intro-one-area}.
The area corresponding to $F_R^{(\ell)}$ is
$$
\mathcal E_{\ell,R}(f)(x)^2
=\sum_{|A|=\ell}\mathbf1_{R_{A^c}}(x_{A^c})\,|(S_Af_A)(x_A)|^2,
$$
with $|R_\emptyset|=1$ and $\mathbf1_{R_\emptyset}=1$.
In particular, $F_R^{(k)}=Q_Rf$ on $R$ and $\mathcal E_{k,R}(f)=S_F^{(k)}(f)$.

\begin{proof}[Proof of Theorem~\ref{t:coordinate}]
For each nonempty $A$, Jensen's inequality implies
$$
\|f_A\|_{L^u(X_A)}^u
\leq\frac1{|R_{A^c}|}\int_{X_A}\int_{R_{A^c}}|f(x_A,u_{A^c})|^u\,du_{A^c}\,dx_A
\leq\frac{\|f\|_{L^u(\mathbb R^N)}^u}{|R_{A^c}|}.
$$
Apply \eqref{e:one-local} on $R_A$ when $|A|=1$ and \eqref{e:flag-local-lp} on the induced
flag when $|A|\geq2$. Since $\prod_{j\in A}Q_jf_A$ depends only on $x_A$, we obtain
\begin{align*}
\|F_R^{(\ell)}\|_{L^p(R)}
&\leq\sum_{|A|=\ell}|R_{A^c}|^{1/p}
\bigg\|\prod_{j\in A}Q_jf_A\bigg\|_{L^p(R_A)}\\
&\leq Cp^{\ell/2}\sum_{|A|=\ell}|R_{A^c}|^{1/p}\|S_Af_A\|_{L^p(R_A)}
\leq C\binom{k}{\ell}p^{\ell/2}\|\mathcal E_{\ell,R}(f)\|_{L^p(R)},
\end{align*}
where the last inequality follows from $|(S_Af_A)(x_A)|\leq\mathcal E_{\ell,R}(f)(x)$ on $R$.
Summing in \eqref{e:coord-expansion}, we have
\begin{equation}\label{e:coord-moments}
\|f-f_R\|_{L^p(R)}
\leq\sum_{\ell=1}^k\|F_R^{(\ell)}\|_{L^p(R)}
\leq C\sum_{\ell=1}^kp^{\ell/2}\|\mathcal E_{\ell,R}(f)\|_{L^p(R)}.
\end{equation}

If $\mathcal E_{\ell,R}(f)=0$ almost everywhere on $R$, the estimate for $F_R^{(\ell)}$
at $p=2$ implies $F_R^{(\ell)}=0$ there almost everywhere.
For $0<\|\mathcal E_{\ell,R}(f)\|_{L^\infty(R)}<\infty$, the same moment bound gives
$$
\bigg(\frac1{|R|}\int_R
\bigg|\frac{F_R^{(\ell)}(x)}{\|\mathcal E_{\ell,R}(f)\|_{L^\infty(R)}}\bigg|^p\,dx\bigg)^{1/p}
\leq Cp^{\ell/2},\qquad 2\leq p<\infty.
$$
Summing the exponential series as in \eqref{e:flag-exp-series}, with $\ell$ in place of $k$, we obtain
\begin{equation}\label{e:coord-order-exp}
\frac1{|R|}\int_R\exp\bigg[c\bigg(\frac{|F_R^{(\ell)}(x)|}
{\|\mathcal E_{\ell,R}(f)\|_{L^\infty(R)}}\bigg)^{2/\ell}\bigg]\,dx\leq2.
\end{equation}
For each $\ell$, Proposition~\ref{p:layer-sharp} proves sharpness of $p^{\ell/2}$ in
\eqref{e:coord-moments} and of $2/\ell$ in \eqref{e:coord-order-exp}.
\end{proof}


Suppose $\mathcal E_{\ell,R}(f)$ is bounded on $R$ for every $1\leq\ell\leq k$.
By \eqref{e:coord-expansion}, \eqref{e:coord-order-exp} and Markov's inequality, for $\tau>0$ we have
\begin{align}
|\{x\in R:|f(x)-f_R|>\tau\}|
&\leq\sum_{\ell=1}^k|\{x\in R:|F_R^{(\ell)}(x)|>\tau/k\}|\notag\\
&\leq 2|R|\sum_{\ell=1}^k
\exp\bigg[-c\bigg(\frac{\tau}{\|\mathcal E_{\ell,R}(f)\|_{L^\infty(R)}}\bigg)^{2/\ell}\bigg],
\label{e:coord-tails}
\end{align}
where $c>0$ depends only on $k$ and the block dimensions. A term with zero area norm is omitted,
since the corresponding $F_R^{(\ell)}$ vanishes almost everywhere on $R$.

\begin{proposition}\label{p:layer-sharp}
Let $f$ range over real Schwartz functions and $R$ over flag rectangles.
\begin{enumerate}[(i)]
\item For every $1\leq\ell\leq k$, the moment estimate for $F_R^{(\ell)}$ has optimal order $p^{\ell/2}$.
If \eqref{e:coord-moments} holds with finite coefficients $a_{p,\ell}\geq0$ in place of $Cp^{\ell/2}$, then
$$
a_{p,\ell}\geq cp^{\ell/2},\qquad 1\leq\ell\leq k,\quad 2\leq p<\infty.
$$
\item The power $2/\ell$ in \eqref{e:coord-order-exp} is optimal for every $\ell$.
If \eqref{e:coord-tails} holds with powers $\gamma_\ell>0$ in place of $2/\ell$
and constants $C_*,c_*>0$ independent of $f,R,\tau$, then $\gamma_\ell\leq2/\ell$ for every $\ell$.
\end{enumerate}
\end{proposition}

\begin{proof}
\noindent\textit{Step 1.}
Fix $1\leq\ell\leq k$ and let $A=\{1,\ldots,\ell\}$.
We choose $f$ so that all terms of order less than $\ell$ vanish. A second parameter $M$
will allow us to make the areas of order greater than $\ell$ arbitrarily small.

For positive integers $L,M$, with $L$ sufficiently large, let
$$
\lambda_i=
\begin{cases}
2^{(L+2)(\ell-i)}, & i\leq\ell,\\[2pt]
2^{-(M+3)(i-\ell)}, & i>\ell.
\end{cases}
$$
In the first $\ell$ blocks take the function from the proof of Proposition~\ref{p:flag-reverse-sharp},
$$
g_L(x_A)=\prod_{i=1}^\ell h_{i,L}(\lambda_i x_i).
$$
For $i\leq\ell$, let $E_i$ be the $\lambda_i^{-1}$-dilate of the cube in Lemma~\ref{l:one-lacunary}
on which $h_{i,L}$ has mean zero. By \eqref{e:flag-large-values}, there is a set $U_L\subset R_A$ such that
$$
g_L(x_A)\geq(L/2)^\ell\quad(x_A\in U_L),\qquad |U_L|\geq c|R_A|2^{-L\dim X_A}.
$$
For $i>\ell$, take $E_i=[-1/2,1/2]^{d_i}$ and let
$$
v_{i,M}(z)=\frac1M\sum_{r=0}^{M-1}b_{d_i}(2^r\lambda_i z).
$$
Since $b_{d_i}(0)=1$ and $2^M\lambda_i\leq2^{-3}$,
$$
\sup_{z\in E_i}|v_{i,M}(z)-1|
\leq\frac{C}{M}\sum_{r=0}^{M-1}2^r\lambda_i
\leq\frac{C}{M}.
$$

On $R=\prod_{i=1}^kE_i$, take
$$
f(x)=g_L(x_A)\prod_{i>\ell}v_{i,M}(x_i).
$$
For large $L$, $\ell(E_1)\leq\cdots\leq\ell(E_\ell)\leq1$, so $R$ is a flag rectangle.
When $\ell=k$, this construction reduces to $f=g_L$.

Averaging out any coordinate of $A$ gives zero. Hence $f_R=0$ and $f_B=0$ unless $A\subseteq B$.
For $B=A\cup D$, $D\subseteq A^c$, we have
$$
f_B(x_B)=\bigg(\prod_{i\in A^c\setminus D}(v_{i,M})_{E_i}\bigg)
g_L(x_A)\prod_{i\in D}v_{i,M}(x_i).
$$
In particular, $F_R^{(\ell)}=f_A=(\prod_{i>\ell}(v_{i,M})_{E_i})g_L$ on $R$.
For all sufficiently large $M$, independently of $L$, it follows that
$$
\min\{f(x),F_R^{(\ell)}(x)\}\geq cL^\ell,\qquad x\in U_L\times R_{A^c}.
$$

\smallskip
\noindent\textit{Step 2.}
We estimate $S_Bf_B$ for $B=A\cup D$. In the Fourier expansion of $f_B$, a term with indices
$(r_i)_{i\in B}$ has frequencies $2^{r_i}\lambda_i\leq|\xi_i|\leq2^{r_i+1}\lambda_i$, where
$0\leq r_i\leq L$ for $i\in A$ and $0\leq r_i<M$ for $i\in D$.
The index $r_i=0$ in $A$ includes the mean correction in $h_{i,L}$, whose coefficient is uniformly bounded.
By the choice of $\lambda_i$, $2^{r_j+1}\lambda_j\leq2^{r_i}\lambda_i$ whenever $i<j$. Thus
$$
2^{r_i}\lambda_i
\leq\bigg(\sum_{\substack{j\in B\\j\geq i}}|\xi_j|^2\bigg)^{1/2}
\leq 2\sqrt{k}\,2^{r_i}\lambda_i.
$$
For $|B|\geq2$, the derivative multiplier in the $i$th parameter is therefore bounded by
$C\beta_1(2^{r_i}\lambda_i t_i)$, where $\beta_1(u)=ue^{-u}$.
By \eqref{e:one-scale-integral}, for every $\lambda>0$ and $m\geq1$,
$$
\int_0^\infty\bigg(\sum_{r=0}^{m-1}\beta_1(2^r\lambda t)\bigg)^2\frac{dt}{t}
=\sum_{r,s=0}^{m-1}\frac{2^{r+s}}{(2^r+2^s)^2}
\leq\sum_{r,s=0}^{m-1}2^{-|r-s|}\leq3m.
$$
Dilation preserves the Fourier $L^1$ norm of each bump. By Fourier inversion and
integration in the cones, as in \eqref{e:flag-sharp-sup}, we obtain
\begin{align*}
\|S_Bf_B\|_\infty^2
&\leq\frac{C}{M^{2|D|}}
\prod_{i\in A}\int_0^\infty\bigg(\sum_{r=0}^L\beta_1(2^r\lambda_i t)\bigg)^2\frac{dt}{t}
\prod_{i\in D}\int_0^\infty\bigg(\sum_{r=0}^{M-1}\beta_1(2^r\lambda_i t)\bigg)^2\frac{dt}{t}\\
&\leq\frac{C(L+1)^\ell M^{|D|}}{M^{2|D|}}\leq CL^\ell M^{-|D|}.
\end{align*}
If $|B|=1$, then $B=A=\{1\}$ and the same bound follows from
$\|\mathcal S_P^{(d_1)}h_{1,L}\|_\infty\leq C\sqrt L$ in Lemma~\ref{l:one-lacunary}.
Summing over $|B|=s$, we conclude that
$$
\|\mathcal E_{s,R}(f)\|_{L^\infty(R)}
\leq
\begin{cases}
0, & s<\ell,\\
CL^{\ell/2}M^{-(s-\ell)/2}, & \ell\leq s\leq k.
\end{cases}
$$

\smallskip
\noindent\textit{Step 3.}
Apply the moment estimate for $f-f_R$ with coefficients $a_{p,s}$. Since $f_R=0$, the bounds
on $U_L\times R_{A^c}$ and on the areas imply
\begin{align*}
cL^\ell 2^{-L\dim X_A/p}
&\leq |R|^{-1/p}\|f-f_R\|_{L^p(R)}
\leq\sum_{s=\ell}^k a_{p,s}\|\mathcal E_{s,R}(f)\|_{L^\infty(R)}
\leq CL^{\ell/2}\sum_{s=\ell}^k a_{p,s}M^{-(s-\ell)/2}.
\end{align*}
Let $M\to\infty$ with $p,L$ fixed, and then take $L=\lceil Kp\rceil$ with $K$ fixed and sufficiently large. We obtain
$$
a_{p,\ell}\geq cL^{\ell/2}2^{-L\dim X_A/p}\geq c'p^{\ell/2}.
$$
For the individual term, the same choice of $L$ gives
$$
\|F_R^{(\ell)}\|_{L^p(R)}
\geq cL^\ell |R|^{1/p}2^{-L\dim X_A/p}
\geq c'p^{\ell/2}\|\mathcal E_{\ell,R}(f)\|_{L^p(R)}.
$$

For the tail estimate, choose $\tau=c_0L^\ell$ so that $f>\tau$ on $U_L\times R_{A^c}$.
Applying \eqref{e:coord-tails} with powers $\gamma_s>0$ and cancelling $|R|$, we obtain
$$
c2^{-L\dim X_A}
\leq C_*\sum_{s=\ell}^k
\exp\bigg[-c_*\Big(cL^{\ell/2}M^{(s-\ell)/2}\Big)^{\gamma_s}\bigg].
$$
Letting $M\to\infty$ with $L$ fixed leaves
$$
c2^{-L\dim X_A}\leq C_*\exp\big[-c_*(cL^{\ell/2})^{\gamma_\ell}\big].
$$
Thus $\ell\gamma_\ell/2\leq1$.

Finally, $\|\mathcal E_{\ell,R}(f)\|_{L^\infty(R)}>0$ by the moment estimate for $F_R^{(\ell)}$.
For any $\gamma>2/\ell$ and $c_*>0$, integration over $U_L\times R_{A^c}$ gives
$$
\frac1{|R|}\int_R
\exp\bigg[c_*\bigg(\frac{|F_R^{(\ell)}(x)|}{\|\mathcal E_{\ell,R}(f)\|_{L^\infty(R)}}\bigg)^\gamma\bigg]\,dx
\geq c2^{-L\dim X_A}\exp\big[c_*(c'L^{\ell/2})^\gamma\big]\to\infty
$$
as $L\to\infty$. This proves sharpness of $2/\ell$ in \eqref{e:coord-order-exp}.
\end{proof}


\section{Approximation by bounded functions in flag BMO}\label{s:bounded-area}

For the two-step flag on $\mathbb R^n\times\mathbb R^m$, Theorem~\ref{t:two-local}
controls $Q_Rf$ on each product of cubes $R$. We now approximate $f$ on the whole
space by a bounded function. By homogeneity, it suffices to prove that if
$f\in L^2(\mathbb R^{n+m})$ and $\|S_F(f)\|_\infty=1$, then for every
$0<\varepsilon\leq1$ there is a $b_\varepsilon$ with
$\|f-b_\varepsilon\|_{\mathrm{BMO}_F}\leq\varepsilon$ and
$$
\|b_\varepsilon\|_\infty\leq C\varepsilon^{-3}[\log(e/\varepsilon)]^{3/2}.
$$
We select terms from the reconstruction in Section~\ref{s:reconstruction}
and estimate the discarded terms by duality with $H_F^1$.

For $g\in L^1(\mathbb R^{n+m})$, let
$$
\mathcal M_Fg(x,y)=\sup_{t,s>0}\big|(e^{t^2\Delta_{x,y}}e^{s^2\Delta_y}g)(x,y)\big|.
$$
We use $H_F^1=\{g\in L^1(\mathbb R^{n+m}):\mathcal M_Fg\in L^1(\mathbb R^{n+m})\}$
with norm $\|g\|_{H_F^1}=\|\mathcal M_Fg\|_1$. Taking Gaussian kernels in
\cite[Definitions~1.2.3 and~1.2.5, Theorem~1.2.6]{HLLW}, we obtain an equivalent
norm on the flag Hardy space. Its dual is the flag Carleson measure space
$\mathrm{BMO}_F$ \cite[Theorem~1.17]{HL}.

By the atomic characterization, $H_F^1\cap L^2(\mathbb R^{n+m})$ is dense in $H_F^1$
\cite[Definition~1.2.10 and Theorem~1.2.11]{HLLW}. Thus for
$u\in L^2(\mathbb R^{n+m})+L^\infty(\mathbb R^{n+m})$ we use the seminorm
$$
\|u\|_{\mathrm{BMO}_F}
=\sup_{\substack{g\in H_F^1\cap L^2(\mathbb R^{n+m})\\ \|g\|_{H_F^1}\leq1}}
|\langle g,u\rangle|,
\qquad
\langle g,u\rangle=\int_{\mathbb R^{n+m}}g\,\overline u.
$$
When finite, this is the norm of the corresponding functional on $H_F^1$.
In particular, every $b\in L^\infty(\mathbb R^{n+m})$ defines such a functional, since
$$
|\langle g,b\rangle|\leq\|b\|_\infty\|g\|_1\leq\|b\|_\infty\|g\|_{H_F^1}.
$$
If $b$ depends only on $x$, this functional is zero. Indeed, by Lemma~\ref{l:hardy-bmo},
$\int_{\mathbb R^m}g(x,y)\,dy=0$ for almost every $x$.

Fix a pair of shifts $\kappa=(\omega_x,\omega_y)$ and a finite family $\mathcal E$.
For the terms with cancellation in both variables, write
$A_{I,J}=A_{\kappa,I,J}^{\mathrm{bi}}(\mathcal E)$ and
$$
c_{I,J}=\frac{\lambda_{\kappa,I,J}^{\mathrm{str}}(\mathcal E)}{(|I||J|)^{1/2}}.
$$
Only finitely many $c_{I,J}$ are nonzero, and by \eqref{e:coefficient-control},
\begin{equation}\label{e:row-coefficients}
\sum_{I,J}c_{I,J}^2\mathbf1_{I\times J}(x,y)\leq C S_F(f)(x,y)^2
\qquad\text{for almost every }(x,y).
\end{equation}

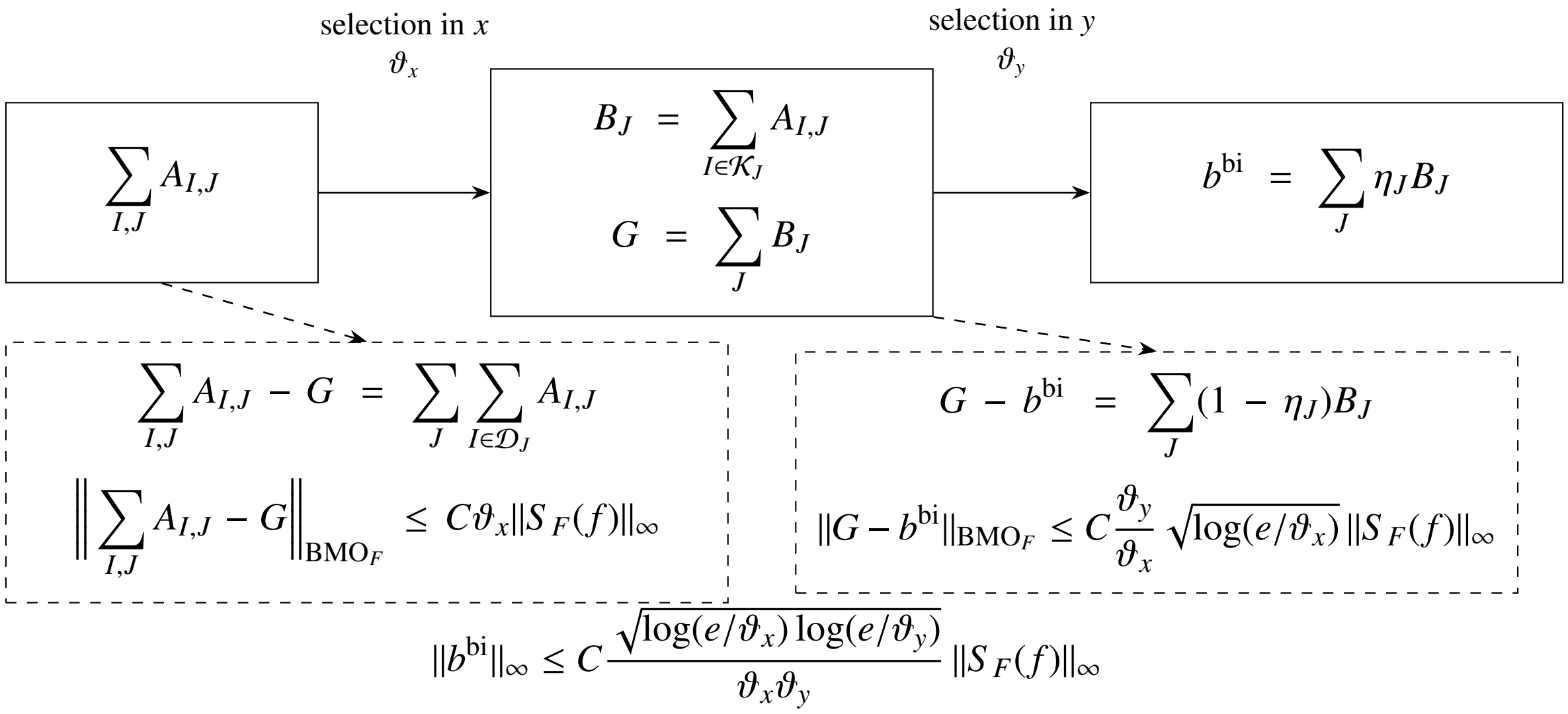
\begin{figure}[tbp]
\centering
\begin{tikzpicture}[
  >=Stealth,
  every node/.style={font=\small,align=center},
  sum/.style={draw=black,rectangle,inner sep=6pt,minimum height=18mm},
  deleted/.style={draw=black,dashed,rectangle,inner sep=6pt,minimum height=15mm},
  flow/.style={draw=black,->,line width=.5pt},
  discard/.style={draw=black,->,dashed,line width=.5pt}
]
\node[sum,text width=27mm] (original) at (1.65,0)
  {$\displaystyle\sum_{I,J}A_{I,J}$};
\node[sum,text width=40mm] (rows) at (7.15,0)
  {$\displaystyle B_J=\sum_{I\in\mathcal K_J}A_{I,J}$\\[2mm]
   $\displaystyle G=\sum_JB_J$};
\node[sum,text width=43mm] (bounded) at (13.3,0)
  {$\displaystyle b^{\mathrm{bi}}=\sum_J\eta_JB_J$};
\draw[flow] (original.east) -- node[above=10mm,font=\scriptsize]
  {selection in $x$\\$\vartheta_x$} (rows.west);
\draw[flow] (rows.east) -- node[above=10mm,font=\scriptsize]
  {selection in $y$\\$\vartheta_y$} (bounded.west);
\node[deleted,text width=68mm] (horizontal-error) at (3.7,-2.8)
  {$\displaystyle\sum_{I,J}A_{I,J}-G=\sum_J\sum_{I\in\mathcal D_J}A_{I,J}$\\[2mm]
   $\displaystyle\bigg\|\sum_{I,J}A_{I,J}-G\bigg\|_{\mathrm{BMO}_F}
   \leq C\vartheta_x\|S_F(f)\|_\infty$};
\node[deleted,text width=68mm] (vertical-error) at (11.6,-2.8)
  {$\displaystyle G-b^{\mathrm{bi}}=\sum_J(1-\eta_J)B_J$\\[2mm]
   $\displaystyle\|G-b^{\mathrm{bi}}\|_{\mathrm{BMO}_F}
   \leq C\frac{\vartheta_y}{\vartheta_x}\sqrt{\log(e/\vartheta_x)}\,\|S_F(f)\|_\infty$};
\draw[discard] (original.south) -- (horizontal-error.north);
\draw[discard] (rows.south east) -- (vertical-error.north);
\node at (7.7,-4.65)
  {$\displaystyle\|b^{\mathrm{bi}}\|_\infty\leq C
  \frac{\sqrt{\log(e/\vartheta_x)\log(e/\vartheta_y)}}{\vartheta_x\vartheta_y}\,
  \|S_F(f)\|_\infty$};
\end{tikzpicture}
\caption{The two selections for the sum with cancellation in both variables,
$0<\vartheta_x,\vartheta_y\leq1/2$. In row $J$, the horizontal selection keeps
the cubes in $\mathcal K_J$ and discards those in $\mathcal D_J$. At each $x$,
the vertical selection keeps $B_J(x,\cdot)$ when $\eta_J(x)=1$.
The dashed boxes show the two errors.}\label{fig:two-selections}
\end{figure}

For each $J$, we first select horizontal cubes $I$ so that the sum $B_J$ over the
selected cubes satisfies
$$
|B_J(x,y)|+\ell(J)|\nabla_yB_J(x,y)|
\leq C\frac{\sqrt{\log(e/\vartheta_x)}}{\vartheta_x}
\bigg(\sum_{I\ni x}c_{I,J}^2\bigg)^{1/2}\mathbf1_J(y).
$$
At fixed $x$, the bound on the right is the coefficient estimate for the selection
in $y$. Thus the factor $\vartheta_x^{-1}\sqrt{\log(e/\vartheta_x)}$ also enters
the vertical error in Figure~\ref{fig:two-selections}. For $G^{\mathrm{rem}}$
and $G^{\mathrm{cub}}$, the row estimates are already available from
\eqref{e:vertical-size}, and only the selection in $y$ is needed.

By Proposition~\ref{p:carleson-bmo}, the horizontal error is at most
$C\vartheta_x\|S_F(f)\|_\infty$ provided that, for every open
$\Omega\subset\mathbb R^{n+m}$ of finite measure,
$$
\sum_J\sum_{\substack{I\in\mathcal D_J\\I\times J\subset\Omega}}c_{I,J}^2|I||J|
\leq C\vartheta_x^2|\Omega|\,\|S_F(f)\|_\infty^2.
$$

\subsection{Horizontal selection within each row}\label{s:horizontal-selection}

For each $J$, we select the horizontal cubes using a Sobolev norm in $y$
to control both the row sum and its $y$ derivatives. For the subsequent
selection in $y$, the scalar case of the same lemma suffices.

\begin{lemma}\label{l:dyadic-selection}
Let $\mathcal F$ be a finite family of cubes in a dyadic system $\mathcal D$
on $\mathbb R^{d_0}$, and let $\mathcal H$ be a Hilbert space.
For each $I\in\mathcal F$, suppose that $u_I\in C_c^1(\mathbb R^{d_0};\mathcal H)$
is supported in $\overline I$, has integral zero and satisfies
$$
\|u_I(x)\|_{\mathcal H}
+\ell(I)\bigg(\sum_{a=1}^{d_0}\|\partial_a u_I(x)\|_{\mathcal H}^2\bigg)^{1/2}
\leq c_I\mathbf1_I(x).
$$
For $0<\vartheta\leq1/2$, there is a partition
$\mathcal F=\mathcal F_{\rm sel}\cup\mathcal F_{\rm del}$ such that, for almost every $x$,
\begin{equation}\label{e:selection-size}
\bigg\|\sum_{I\in\mathcal F_{\rm sel}}u_I(x)\bigg\|_{\mathcal H}
\leq C\frac{\sqrt{\log(e/\vartheta)}}{\vartheta}
\bigg(\sum_{I\in\mathcal F}c_I^2\mathbf1_I(x)\bigg)^{1/2},
\end{equation}
and, for every $U\in\mathcal D$,
\begin{equation}\label{e:selection-energy}
\sum_{\substack{I\in\mathcal F_{\rm del}\\I\subseteq U}}c_I^2|I|
\leq C\vartheta^2\sum_{I\in\mathcal F}c_I^2|I\cap U|.
\end{equation}
The constant depends only on $d_0$.
\end{lemma}

\begin{proof}
\emph{Step 1.} For $\mathcal P\subseteq\mathcal F$, write $h_{\mathcal P}=\sum_{I\in\mathcal P}u_I$.
By support and cancellation, $\Delta_Uu_I=0$ unless $U\subseteq I$.
For $I=U^{(k)}$, the gradient bound implies
$$
\|\Delta_Uu_I(x)\|_{\mathcal H}
\leq \mathbf1_U(x)\sup_{v,w\in U}\|u_I(v)-u_I(w)\|_{\mathcal H}
\leq C\frac{\ell(U)}{\ell(I)}c_I\mathbf1_U(x)
=C2^{-k}c_I\mathbf1_U(x).
$$
By Minkowski's inequality and the partition of $I$ by its descendants of generation $k$,
\begin{align*}
\bigg(\sum_U\|\Delta_Uh_{\mathcal P}(x)\|_{\mathcal H}^2\bigg)^{1/2}
&\leq C\sum_{k\geq0}2^{-k}
\bigg(\sum_{I\in\mathcal P}c_I^2\sum_{U^{(k)}=I}\mathbf1_U(x)\bigg)^{1/2}
\leq C\bigg(\sum_{I\in\mathcal P}c_I^2\mathbf1_I(x)\bigg)^{1/2}.
\end{align*}
The dyadic averages of $h_{\mathcal P}$ vanish at all sufficiently coarse scales.
Applying \eqref{e:hilbert-square} to the finite generation sums and passing to the limit
as in \eqref{e:martingale-sums}, we obtain, for $2\leq p<\infty$,
\begin{equation}\label{e:row-moments}
\|h_{\mathcal P}\|_{L^p(\mathcal H)}
\leq C\sqrt p\,\bigg\|\bigg(\sum_{I\in\mathcal P}c_I^2\mathbf1_I\bigg)^{1/2}\bigg\|_p.
\end{equation}
In particular, at $p=2$,
\begin{equation}\label{e:row-orthogonality}
\int_{\mathbb R^{d_0}}\|h_{\mathcal P}(x)\|_{\mathcal H}^2\,dx
\leq C\sum_{I\in\mathcal P}c_I^2|I|.
\end{equation}

\smallskip
\noindent\emph{Step 2.} Let $\delta=\vartheta^2$, $\rho=1+\delta$, and for $U\in\mathcal D$ write
$
E(U)=\sum_{\substack{I\in\mathcal F\\I\supsetneq U}}c_I^2.
$

First select every $I$ with $c_I=0$ or $c_I^2>\delta E(I)$.
Along the chain of cubes containing $x$, the intervals $[E(I),E(I)+c_I^2]$
have disjoint interiors. Hence
\begin{align*}
\bigg\|\sum_{c_I^2>\delta E(I)}u_I(x)\bigg\|_{\mathcal H}
&\leq\sum_{\substack{I\ni x\\c_I^2>\delta E(I)}}c_I
\leq\big(\sqrt{1+\delta^{-1}}+\delta^{-1/2}\big)
\sum_{\substack{I\ni x\\c_I^2>\delta E(I)}}
\big(\sqrt{E(I)+c_I^2}-\sqrt{E(I)}\big)\\
&\leq C\delta^{-1/2}\bigg(\sum_{I\ni x}c_I^2\bigg)^{1/2}.
\end{align*}

For $j\in\mathbb Z$, let $R$ range over the maximal cubes among the remaining $I$
with $\rho^j\leq E(I)<\rho^{j+1}$. Denote by $\mathcal P_{j,R}$ all such $I\subseteq R$.
For $x\in R$, let $I_x$ be the smallest cube in $\mathcal P_{j,R}$ containing $x$.
Since $c_{I_x}^2\leq\delta E(I_x)$,
\begin{equation}\label{e:row-energy}
\sum_{I\in\mathcal P_{j,R}}c_I^2\mathbf1_I(x)
\leq E(I_x)+c_{I_x}^2-E(R)
\leq (1+\delta)\rho^{j+1}-\rho^j
=\delta(1+\rho)\rho^j\leq3\delta\rho^j.
\end{equation}
The sum vanishes outside $R$. Thus for $h_{j,R}=\sum_{I\in\mathcal P_{j,R}}u_I$,
\eqref{e:row-moments} implies
$$
\bigg(\frac1{|R|}\int_R\|h_{j,R}(x)\|_{\mathcal H}^p\,dx\bigg)^{1/p}
\leq C_1\sqrt{p\delta}\,\rho^{j/2}.
$$

Take $p=2+\log(1/\delta)$ and $H_j=eC_1\sqrt{p\delta}\,\rho^{j/2}$.
Let $\mathcal V_{j,R}$ consist of the maximal dyadic cubes $V\subseteq R$
with $\|\langle h_{j,R}\rangle_V\|_{\mathcal H}>H_j$, and write
$\mathcal B_{j,R}=\bigcup_{V\in\mathcal V_{j,R}}V$.
Since $\langle h_{j,R}\rangle_R=0$, each $V$ is a proper subcube of $R$.
By Jensen's inequality and disjointness,
\begin{equation}\label{e:stopping-measure}
|\mathcal B_{j,R}|
\leq H_j^{-p}\sum_{V\in\mathcal V_{j,R}}\int_V\|h_{j,R}(x)\|_{\mathcal H}^p\,dx
\leq e^{-p}|R|\leq\delta|R|.
\end{equation}
Delete every $I\in\mathcal P_{j,R}$ contained in a stopping cube and
let $b_{j,R}$ be the sum of the selected terms in $\mathcal P_{j,R}$.

For $V\in\mathcal V_{j,R}$, let $P=V^{(1)}$. Then
$\|\langle h_{j,R}\rangle_P\|_{\mathcal H}\leq H_j$.
On $V$, only cubes $I\supseteq P$ contribute to $b_{j,R}$.
The same cubes contribute to $\langle h_{j,R}\rangle_P$, by cancellation.
Thus for $x\in V$,
\begin{align*}
\|b_{j,R}(x)\|_{\mathcal H}
&\leq \|\langle h_{j,R}\rangle_P\|_{\mathcal H}
+\sum_{\substack{I\in\mathcal P_{j,R}\\I\supseteq P}}
\|u_I(x)-\langle u_I\rangle_P\|_{\mathcal H}\\
&\leq H_j+C\sum_{\substack{I\in\mathcal P_{j,R}\\I\supseteq P}}
\frac{\ell(P)}{\ell(I)}c_I
\leq H_j+C\sqrt{\delta\rho^{j+1}}\sum_{l\geq0}2^{-l}
\leq C\sqrt{p\delta}\,\rho^{j/2}.
\end{align*}
On $R\setminus\mathcal B_{j,R}$, dyadic differentiation implies
$\|b_{j,R}(x)\|_{\mathcal H}=\|h_{j,R}(x)\|_{\mathcal H}\leq H_j$ almost everywhere.

For fixed $j$, the maximal cubes $R$ are disjoint, and
$\rho^j\leq E(R)\leq\sum_{I\ni x}c_I^2$ on $R$. Hence
\begin{align*}
\sum_{j,R}\|b_{j,R}(x)\|_{\mathcal H}
&\leq C\sqrt{p\delta}\sum_{\rho^j\leq\sum_{I\ni x}c_I^2}\rho^{j/2}\\
&\leq C\sqrt{p\delta}\,\frac{\sqrt\rho}{\sqrt\rho-1}
\bigg(\sum_{I\ni x}c_I^2\bigg)^{1/2}
\leq C\sqrt{p/\delta}\bigg(\sum_{I\ni x}c_I^2\bigg)^{1/2}.
\end{align*}
Adding the estimate for $c_I^2>\delta E(I)$ and recalling that
$p=2\log(e/\vartheta)$ and $\delta=\vartheta^2$, we obtain \eqref{e:selection-size}.

\smallskip
\noindent\emph{Step 3.} Fix $U\in\mathcal D$. By \eqref{e:row-energy} and
\eqref{e:stopping-measure}, the deleted coefficients in $\mathcal P_{j,R}$ satisfy
$$
\sum_{I\in\mathcal P_{j,R}\cap\mathcal F_{\rm del}}c_I^2|I|
\leq\int_{\mathcal B_{j,R}}\sum_{I\in\mathcal P_{j,R}}c_I^2\mathbf1_I(x)\,dx
\leq3\delta\rho^j|\mathcal B_{j,R}|\leq3\delta^2\rho^j|R|.
$$
Summing first over $R\subseteq U$, we obtain
\begin{align*}
\sum_{j,R:\,R\subseteq U}\sum_{I\in\mathcal P_{j,R}\cap\mathcal F_{\rm del}}c_I^2|I|
&\leq 3\delta^2\sum_{j,R:\,R\subseteq U}\rho^j|R|
\leq 3\delta^2\int_U\sum_{\rho^j\leq\sum_{I\ni x}c_I^2}\rho^j\,dx\\
&\leq 3\delta^2\frac{\rho}{\rho-1}\int_U\sum_{I\ni x}c_I^2\,dx
\leq C\delta\sum_{I\in\mathcal F}c_I^2|I\cap U|.
\end{align*}

For the remaining terms, $I\subseteq U\subsetneq R$, and hence
$$
\rho^j\leq E(R)\leq E(U)\leq E(I)<\rho^{j+1}.
$$
Thus $E(U)$ determines $j$, and at that level at most one maximal cube $R$ contains $U$.
For this $\mathcal P_{j,R}$, \eqref{e:row-energy} implies
$$
\sum_{\substack{I\in\mathcal P_{j,R}\cap\mathcal F_{\rm del}\\I\subseteq U}}c_I^2|I|
\leq\int_U\sum_{I\in\mathcal P_{j,R}}c_I^2\mathbf1_I(x)\,dx
\leq3\delta\rho^j|U|\leq3\delta E(U)|U|
\leq3\delta\sum_{I\in\mathcal F}c_I^2|I\cap U|.
$$
Adding the two estimates, we obtain \eqref{e:selection-energy}.
\end{proof}


For the horizontal selection, we use a Sobolev norm in $y$ that controls both
the row sum and its first derivatives. Fix an integer $M_0>m/2+1$ depending
only on $m$, and let $y_J$ be the centre of $J$. For
$v\in W^{M_0,2}(\mathbb R^m)$ supported in $\overline J$, write
$$
\|v\|_J^2
=\sum_{|\beta|\leq M_0}\frac{\ell(J)^{2|\beta|}}{|J|}\|\partial_y^\beta v\|_{L^2(\mathbb R^m)}^2
=\|v(y_J+\ell(J)\,\cdot)\|_{W^{M_0,2}(\mathbb R^m)}^2.
$$
By Sobolev embedding applied to $v(y_J+\ell(J)\,\cdot)$,
\begin{equation}\label{e:scaled-sobolev}
\|v\|_\infty+\ell(J)\|\nabla_yv\|_\infty\leq C_m\|v\|_J.
\end{equation}

\begin{proposition}\label{p:horizontal-selection}
Fix $0<\vartheta_x\leq1/2$. For each $J$, the horizontal indices admit a partition
$\mathcal K_J\cup\mathcal D_J$ such that the sum $B_J=\sum_{I\in\mathcal K_J}A_{I,J}$
satisfies, for almost every $x$,
\begin{equation}\label{e:selected-row}
\|B_J(x,\cdot)\|_J\leq C\frac{\sqrt{\log(e/\vartheta_x)}}{\vartheta_x}
\bigg(\sum_Ic_{I,J}^2\mathbf1_I(x)\bigg)^{1/2}.
\end{equation}
For every $x$, $B_J(x,\cdot)$ has integral zero and is supported in
$\{y\in J:\operatorname{dist}_\infty(y,\partial J)\geq\ell(J)/24\}$.

For every $U\in\mathcal D_x$, the deleted coefficients satisfy
\begin{equation}\label{e:deleted-row-energy}
\sum_{\substack{I\in\mathcal D_J\\I\subseteq U}}c_{I,J}^2|I|
\leq C\vartheta_x^2\sum_Ic_{I,J}^2|I\cap U|.
\end{equation}
Consequently, for every open $\Omega\subset\mathbb R^{n+m}$ of finite measure,
\begin{equation}\label{e:deleted-carleson}
\sum_J\sum_{\substack{I\in\mathcal D_J\\I\times J\subseteq\Omega}}c_{I,J}^2|I||J|
\leq C\vartheta_x^2\int_\Omega\sum_{I,J}c_{I,J}^2\mathbf1_{I\times J}(x,y)\,dx\,dy.
\end{equation}
The constants are independent of $\kappa,\mathcal E,J$ and $\vartheta_x$.
\end{proposition}

\begin{proof}
Fix $J$. By \eqref{e:bi-deriv},
\begin{align*}
\sum_{|\alpha|\leq1}\ell(I)^{2|\alpha|}\|\partial_x^\alpha A_{I,J}(x,\cdot)\|_J^2
&=\frac1{|J|}\sum_{\substack{|\alpha|\leq1\\|\beta|\leq M_0}}
\ell(I)^{2|\alpha|}\ell(J)^{2|\beta|}
\int_J|\partial_x^\alpha\partial_y^\beta A_{I,J}(x,y)|^2\,dy
\leq Cc_{I,J}^2\mathbf1_I(x).
\end{align*}
The map $x\mapsto A_{I,J}(x,\cdot)$ is smooth in this Sobolev norm, is supported in $I$,
and has integral zero, since for $|\beta|\leq M_0$,
$$
\int_I\partial_y^\beta A_{I,J}(x,y)\,dx
=\partial_y^\beta\int_I A_{I,J}(x,y)\,dx=0.
$$
Applying Lemma~\ref{l:dyadic-selection} to these functions with coefficients $Cc_{I,J}$,
we obtain \eqref{e:selected-row} and \eqref{e:deleted-row-energy}.
The support and cancellation of $B_J(x,\cdot)$ follow from Proposition~\ref{p:two-atoms}.

For each $J$, let $\mathcal U_J$ consist of the maximal $U\in\mathcal D_x$
with $U\times J\subseteq\Omega$. Every $I$ with $I\times J\subseteq\Omega$
lies in exactly one such $U$, since the admissible cubes have volume at most
$|\Omega|/|J|$. By \eqref{e:deleted-row-energy},
\begin{align*}
\sum_J\sum_{\substack{I\in\mathcal D_J\\I\times J\subseteq\Omega}}c_{I,J}^2|I||J|
&=\sum_J|J|\sum_{U\in\mathcal U_J}\sum_{\substack{I\in\mathcal D_J\\I\subseteq U}}c_{I,J}^2|I| 
\leq C\vartheta_x^2\sum_J\sum_{U\in\mathcal U_J}
\int_{U\times J}\sum_Ic_{I,J}^2\mathbf1_I(x)\,dx\,dy\\
&\leq C\vartheta_x^2\int_\Omega\sum_{I,J}c_{I,J}^2\mathbf1_{I\times J}(x,y)\,dx\,dy
\leq C\vartheta_x^2\int_\Omega S_F(f)(x,y)^2\,dx\,dy,
\end{align*}
where the last two inequalities follow from the disjointness of the cubes in each
$\mathcal U_J$ and \eqref{e:row-coefficients}, respectively.
\end{proof}


\subsection{Vertical selection of rows}\label{s:vertical-selection}

For each $x$, we apply Lemma~\ref{l:dyadic-selection} in $y$ to the rows $B_J(x,\cdot)$.

\begin{proposition}\label{p:vertical-selection}
Fix $\kappa,\mathcal E$ as in Subsection~\ref{s:vertical-cancellation}.
Let $(B_J,c_J)$ denote either the remainder rows and their coefficients or the cubic
rows and their coefficients, and write $G=\sum_JB_J$.
For $0<\vartheta\leq1/2$, there are measurable functions
$\eta_J:\mathbb R^n\to\{0,1\}$ such that $b(x,y)=\sum_J\eta_J(x)B_J(x,y)$ satisfies
$$
\|b\|_\infty\leq C\frac{\sqrt{\log(e/\vartheta)}}{\vartheta}
\bigg\|\bigg(\sum_Jc_J(x)^2\mathbf1_J(y)\bigg)^{1/2}\bigg\|_{L^\infty(\mathbb R^{n+m})}.
$$
For the discarded sum,
$$
\operatorname*{ess\,sup}_{x\in\mathbb R^n}\|G(x,\cdot)-b(x,\cdot)\|_{\mathrm{BMO}(\mathbb R^m)}
\leq C\vartheta\bigg\|\bigg(\sum_Jc_J(x)^2\mathbf1_J(y)\bigg)^{1/2}\bigg\|_{L^\infty(\mathbb R^{n+m})}.
$$
Here $C$ depends only on $n,m$, and
$\|h\|_{\mathrm{BMO}(\mathbb R^m)}=\sup_Q(|Q|^{-1}\int_Q|h-h_Q|^2)^{1/2}$,
where $Q$ ranges over cubes in $\mathbb R^m$.
\end{proposition}

\begin{proof}
Fix $x$ and write $B_J=B_J(x,\cdot)$ and $c_J=c_J(x)$. Each $B_J\in C_c^\infty(J^\circ)$
has integral zero, and by \eqref{e:vertical-size},
$$
|B_J(y)|+\ell(J)|\nabla_yB_J(y)|\leq Cc_J\mathbf1_J(y).
$$
Apply Lemma~\ref{l:dyadic-selection} with coefficients $Cc_J$. In a fixed ordering
of the subsets of the finite row family, choose the first satisfying
\eqref{e:selection-size} almost everywhere in $y$ and \eqref{e:selection-energy}
for every $U\in\mathcal D_y$. By Fubini's theorem and the countability of $\mathcal D_y$,
the selected indicators $\eta_J(x)$ are measurable in $x$. Thus
$$
|b(x,y)|\leq C\frac{\sqrt{\log(e/\vartheta)}}{\vartheta}
\bigg(\sum_Jc_J^2\mathbf1_J(y)\bigg)^{1/2}.
$$
For $\mathcal D(x)=\{J:\eta_J(x)=0\}$ and $U\in\mathcal D_y$, we have
\begin{equation}\label{e:vertical-deleted}
\sum_{\substack{J\in\mathcal D(x)\\J\subseteq U}}c_J^2|J|
\leq C\vartheta^2\sum_Jc_J^2|J\cap U|
=C\vartheta^2\int_U\sum_Jc_J^2\mathbf1_J(y)\,dy.
\end{equation}
In particular, taking $U=J$ bounds each deleted coefficient by
$C\vartheta\|(\sum_Jc_J^2\mathbf1_J)^{1/2}\|_{L^\infty(\mathbb R^m)}$.

Let $h=\sum_{J\in\mathcal D(x)}B_J$, and fix a cube $Q$ with centre $y_Q$ and side length $r$.
Every $J$ meeting $Q$ with $\ell(J)\leq r$ lies in a dyadic cube $U$ meeting $Q$
with $r\leq\ell(U)<2r$. Denote these $U$ by $\mathcal U(Q)$.
At each dyadic side length at least $r$, at most $C_m$ cubes meet $Q$.

Subtracting $\sum_{J\in\mathcal D(x),\,\ell(J)>r}B_J(y_Q)$ on $Q$, we use
\eqref{e:row-orthogonality} for $\ell(J)\leq r$ and the gradient bound for $\ell(J)>r$.
By \eqref{e:vertical-deleted},
\begin{align*}
\bigg(\frac1{|Q|}\int_Q|h-h_Q|^2\bigg)^{1/2}
&\leq |Q|^{-1/2}\bigg\|\sum_{\substack{J\in\mathcal D(x)\\\ell(J)\leq r}}B_J\bigg\|_{L^2(Q)}
+\sum_{\substack{J\in\mathcal D(x)\\\ell(J)>r,\ J\cap Q\ne\emptyset}}
\sup_{y\in Q}|B_J(y)-B_J(y_Q)|\\
&\leq C\bigg(\frac1{|Q|}\sum_{U\in\mathcal U(Q)}
\sum_{\substack{J\in\mathcal D(x)\\J\subseteq U}}c_J^2|J|\bigg)^{1/2}
+Cr\sum_{\substack{J\in\mathcal D(x)\\\ell(J)>r,\ J\cap Q\ne\emptyset}}
\frac{c_J}{\ell(J)}\\
&\leq C\vartheta\bigg[\bigg(\sum_{U\in\mathcal U(Q)}\frac{|U|}{|Q|}\bigg)^{1/2}
+\sum_{k\geq0}2^{-k}\bigg]
\bigg\|\bigg(\sum_Jc_J^2\mathbf1_J\bigg)^{1/2}\bigg\|_{L^\infty(\mathbb R^m)}\\
&\leq C\vartheta\bigg\|\bigg(\sum_Jc_J^2\mathbf1_J\bigg)^{1/2}\bigg\|_{L^\infty(\mathbb R^m)}.
\end{align*}
Taking the supremum over $Q$ and then the essential supremum in $x$, we obtain
the required $\mathrm{BMO}(\mathbb R^m)$ bound.
\end{proof}



\subsection{Flag Hardy space estimates for the discarded sums}\label{s:hardy-estimates}

For the discarded rows, we use Hardy--BMO duality on each vertical slice.
For the remainders with $\ell(J)=2^k\ell(I)$, the map $(x,y)\mapsto(x,2^{-k}y)$
takes $I\times J$ to the cube $I\times 2^{-k}J$, where $2^{-k}J=\{2^{-k}y:y\in J\}$.
We need the corresponding BMO bound uniformly in $k$.

Write $d=n+m$. In the Fourier variables $(\xi,\eta)\in\mathbb R^n\times\mathbb R^m$,
the full-space Riesz transforms $R_i^{(1)}$, $1\leq i\leq d$, have multipliers
$-\mathrm i(\xi,\eta)_i/|(\xi,\eta)|$, and the vertical Riesz transforms
$R_j^{(2)}$, $1\leq j\leq m$, have multipliers $-\mathrm i\eta_j/|\eta|$.
Let $R_0^{(1)}=R_0^{(2)}=\mathrm{Id}$ and denote the classical real Hardy space
on $\mathbb R^q$ by $H^1(\mathbb R^q)$.

\begin{lemma}\label{l:hardy-bmo}
For $g\in H_F^1\cap L^2(\mathbb R^d)$,
\begin{equation}\label{e:hardy-riesz}
\|g\|_{H^1(\mathbb R^d)}
+\int_{\mathbb R^n}\|g(x,\cdot)\|_{H^1(\mathbb R^m)}\,dx
+\sum_{i=0}^d\sum_{j=0}^m\|R_i^{(1)}R_j^{(2)}g\|_1
\leq C\|g\|_{H_F^1}.
\end{equation}
For $u\in L^2(\mathbb R^d)$, we have
\begin{equation}\label{e:vertical-bmo}
\|u\|_{\mathrm{BMO}_F}
\leq C\operatorname*{ess\,sup}_{x\in\mathbb R^n}\|u(x,\cdot)\|_{\mathrm{BMO}(\mathbb R^m)}
\end{equation}
and, for every $\rho\geq1$,
\begin{equation}\label{e:dilated-bmo}
\|u\|_{\mathrm{BMO}_F}\leq C\|u(\cdot,\rho\,\cdot)\|_{\mathrm{BMO}(\mathbb R^d)}.
\end{equation}
The constants depend only on $n,m$.
\end{lemma}

\begin{proof}
In $e^{t^2\Delta_{x,y}}e^{s^2\Delta_y}g$, fix either time at a positive rational
value and let the other tend to zero. By strong continuity in $L^2(\mathbb R^d)$,
a subsequence converges almost everywhere to the corresponding one-parameter
extension. Continuity at positive times then implies
$$
\max\bigg\{\sup_{t>0}|e^{t^2\Delta_{x,y}}g|,\,
           \sup_{s>0}|e^{s^2\Delta_y}g|\bigg\}\leq\mathcal M_Fg
\qquad\text{a.e.}
$$
By the classical heat maximal characterization \cite[Theorem~11]{FS},
\begin{align*}
\|g\|_{H^1(\mathbb R^d)}+\int_{\mathbb R^n}\|g(x,\cdot)\|_{H^1(\mathbb R^m)}\,dx
&\leq C\int_{\mathbb R^d}\bigg(\sup_{t>0}|e^{t^2\Delta_{x,y}}g|
                               +\sup_{s>0}|e^{s^2\Delta_y}g|\bigg)
\leq C\|\mathcal M_Fg\|_1=C\|g\|_{H_F^1}.
\end{align*}
The classical Riesz characterization \cite[Section~2]{FS} now bounds the terms
with $i=0$ or $j=0$ in \eqref{e:hardy-riesz}. For the mixed terms,
\cite[Corollary~1.2.13]{HLLW} implies
$$
\sum_{i=1}^d\sum_{j=1}^m\|R_i^{(1)}R_j^{(2)}g\|_1\leq C\|\mathcal M_Fg\|_1.
$$

For $h\in H^1(\mathbb R^q)\cap L^2(\mathbb R^q)$ and
$v\in\mathrm{BMO}(\mathbb R^q)\cap L^2(\mathbb R^q)$, truncate the real and
imaginary parts of $v$ at $\pm T$. Since $\|v_T\|_{\mathrm{BMO}}\leq\|v\|_{\mathrm{BMO}}$
and $v_T\to v$ in $L^2(\mathbb R^q)$, classical duality \cite[Theorem~2]{FS} implies
$$
|\langle h,v\rangle|
=\lim_{T\to\infty}|\langle h,v_T\rangle|
\leq C_q\|h\|_{H^1(\mathbb R^q)}\|v\|_{\mathrm{BMO}(\mathbb R^q)}.
$$
Applying this bound on the vertical slices and using \eqref{e:hardy-riesz}, we obtain
\begin{align*}
|\langle g,u\rangle|
&\leq C\int_{\mathbb R^n}\|g(x,\cdot)\|_{H^1(\mathbb R^m)}
                           \|u(x,\cdot)\|_{\mathrm{BMO}(\mathbb R^m)}\,dx
\leq C\|g\|_{H_F^1}\operatorname*{ess\,sup}_{x\in\mathbb R^n}
                      \|u(x,\cdot)\|_{\mathrm{BMO}(\mathbb R^m)},
\end{align*}
which proves \eqref{e:vertical-bmo}.

For $\rho\geq1$, let $g_\rho(x,y)=\rho^mg(x,\rho y)$. By scaling in $y$,
$$
|e^{t^2\Delta_{x,y}}g_\rho(x,y)|
=\rho^m\big|(e^{t^2\Delta_{x,y}}e^{(\rho^2-1)t^2\Delta_y}g)(x,\rho y)\big|
\leq\rho^m\mathcal M_Fg(x,\rho y),
$$
where $s=t\sqrt{\rho^2-1}$ for $\rho>1$ and the case $\rho=1$ follows
from the first display of the proof. Hence
$$
\|g_\rho\|_{H^1(\mathbb R^d)}
\leq C\int_{\mathbb R^d}\sup_{t>0}|e^{t^2\Delta_{x,y}}g_\rho(x,y)|\,dx\,dy
\leq C\int_{\mathbb R^d}\rho^m\mathcal M_Fg(x,\rho y)\,dx\,dy
=C\|g\|_{H_F^1}.
$$
Changing variables in the pairing and applying classical duality, we obtain
$$
|\langle g,u\rangle|
=|\langle g_\rho,u(\cdot,\rho\,\cdot)\rangle|
\leq C\|g_\rho\|_{H^1(\mathbb R^d)}\|u(\cdot,\rho\,\cdot)\|_{\mathrm{BMO}(\mathbb R^d)}
\leq C\|g\|_{H_F^1}\|u(\cdot,\rho\,\cdot)\|_{\mathrm{BMO}(\mathbb R^d)}.
$$
Taking the supremum over $\|g\|_{H_F^1}\leq1$, we obtain \eqref{e:dilated-bmo}.
\end{proof}

For the proof of Proposition~\ref{p:carleson-bmo}, we need an $L^1$ estimate
for the area function of the adjoint reconstruction kernels.
Write $k^*(z)=\overline{k(-z)}$.

\begin{lemma}\label{l:hardy-area}
For every $\lambda>0$ and $g\in H_F^1\cap L^2(\mathbb R^{n+m})$,
\begin{equation}\label{e:hardy-area}
\bigg\|\bigg(\int_0^\infty\int_0^\infty
\int_{\substack{|x-u|<\lambda t\\|y-v|<\lambda\max\{t,s\}}}
\frac{|((K_{t,s}^\sigma)^* *g)(u,v)|^2}{t^n\max\{t,s\}^m}\,du\,dv\,
\frac{dt}{t}\frac{ds}{s}\bigg)^{1/2}\bigg\|_{L^1(dx\,dy)}
\leq C\|g\|_{H_F^1},
\end{equation}
where $C$ depends only on $n,m,\lambda,\sigma$ and $\theta_0$.
\end{lemma}

\begin{proof}
\noindent{\it Step 1.}
Let $L=8\max\{1,\lambda\}$, $D_1=n+m$ and $D_2=m$. For $i=1,2$, define
$$
(\mathsf T_i h)(z)(w,r)=((\Phi_{i,r}^\sigma)^* *h)(z+rw),
\qquad w\in B_{\mathbb R^{D_i}}(0,L),\quad r>0,
$$
and let $\mathcal H_i=L^2(B_{\mathbb R^{D_i}}(0,L)\times(0,\infty),dw\,dr/r)$.
By translation invariance, Plancherel's identity and
$\widehat\Phi_i^\sigma(\xi)=a_\sigma(|\xi|)$,
\begin{align*}
\|\mathsf T_i h\|_{L^2(\mathbb R^{D_i};\mathcal H_i)}^2
&=|B_{\mathbb R^{D_i}}(0,L)|\int_0^\infty
\|(\Phi_{i,r}^\sigma)^* *h\|_{L^2(\mathbb R^{D_i})}^2\,\frac{dr}{r}\\
&=|B_{\mathbb R^{D_i}}(0,L)|\int_0^\infty|a_\sigma(r)|^2\,\frac{dr}{r}\,
\|h\|_{L^2(\mathbb R^{D_i})}^2
\leq C\|h\|_{L^2(\mathbb R^{D_i})}^2.
\end{align*}
For $i=1$, the same identity holds for $\mathcal H_2$-valued $h$ by summing
over an orthonormal basis of $\mathcal H_2$.

Write $G_g(u,v,t,s)=((K_{t,s}^\sigma)^* *g)(u,v)$ and let $\mathsf T_2$ act in $y$.
For $z=(x,y)$ and $w=(w_x,w_y)$,
$$
(\mathsf T_1\mathsf T_2g)(z)(w,t,w',s)
=G_g(x+tw_x,y+tw_y+sw',t,s).
$$
Since $\lambda\leq L/8$, applying \eqref{e:cone-lower} at scales $Lt,Ls$
and changing variables in the two ball integrals, we obtain
\begin{align*}
\mathcal A_{\lambda,F}(G_g)(z)^2
&\leq CL^{n+m}\int_0^\infty\int_0^\infty\int_{\mathbb R^{n+m}}
\frac{\chi_{Lt,Ls}(x-u,y-v)}{(Lt)^{n+m}(Ls)^m}
|G_g(u,v,t,s)|^2\,du\,dv\,\frac{dt}{t}\frac{ds}{s}\\
&=CL^{-m}\int_0^\infty\int_0^\infty\int_{|w|<L}\int_{|w'|<L}
|G_g(x+tw_x,y+tw_y+sw',t,s)|^2\,dw'\,dw\,\frac{dt}{t}\frac{ds}{s}\\
&=CL^{-m}\|\mathsf T_1\mathsf T_2g(z)\|_{\mathcal H_1\otimes\mathcal H_2}^2.
\end{align*}

\smallskip
\noindent{\it Step 2.}
Let $a_Q$ be a $(1,2)$-atom in $\mathbb R^q$, $q=D_i$, supported in a cube $Q$
of side $\rho$ and centre $z_Q$. Then $\int_Qa_Q=0$ and
$\|a_Q\|_1\leq |Q|^{1/2}\|a_Q\|_2\leq1$.

Let $B_Q=B(z_Q,2\sqrt q\,\rho)$. For $z\notin B_Q$, write $R=|z-z_Q|$.
If $u\in Q$, $|w|<L$ and $r<R/[2(L+1)]$, then
$$
|z+rw-u|\geq R-r|w|-|u-z_Q|\geq 3R/4-Lr>r.
$$
By the support of $\Phi_{i,r}^\sigma$, the cancellation of $a_Q$ and the mean value theorem,
\begin{align*}
|(\mathsf T_i a_Q)(z)(w,r)|
&=\bigg|\int_Q\big[(\Phi_{i,r}^\sigma)^*(z+rw-u)
                  -(\Phi_{i,r}^\sigma)^*(z+rw-z_Q)\big]a_Q(u)\,du\bigg|\\
&\leq C\rho r^{-q-1}\|a_Q\|_1\,\mathbf1_{\{r\geq R/[2(L+1)]\}}.
\end{align*}
Using the $L^2$ estimate of Step~1 on $B_Q$ and polar coordinates outside $B_Q$, we obtain
\begin{align*}
\|\mathsf T_i a_Q\|_{L^1(\mathcal H_i)}
&\leq |B_Q|^{1/2}\|\mathsf T_i a_Q\|_{L^2(\mathcal H_i)}
+C\rho\int_{2\sqrt q\,\rho}^\infty
\bigg(\int_{R/[2(L+1)]}^\infty r^{-2q-2}\,\frac{dr}{r}\bigg)^{1/2}R^{q-1}\,dR\\
&\leq C|Q|^{1/2}\|a_Q\|_2+C\rho\int_{2\sqrt q\,\rho}^\infty\frac{dR}{R^2}
\leq C.
\end{align*}

For $h\in H^1(\mathbb R^q)\cap L^2(\mathbb R^q)$, take an atomic decomposition
$h=\sum_\nu c_\nu a_\nu$ as in \cite{Latter}. Since $\|a_\nu\|_1\leq1$ and
$\sum_\nu|c_\nu|\leq C\|h\|_{H^1}$, convolution with $(\Phi_{i,r}^\sigma)^*$
may be taken term by term at each $r>0$. By Minkowski's inequality and the
Riesz characterization \cite[Section~2]{FS},
\begin{equation}\label{e:atomic-area}
\|\mathsf T_i h\|_{L^1(\mathcal H_i)}
\leq\sum_\nu|c_\nu|\,\|\mathsf T_i a_\nu\|_{L^1(\mathcal H_i)}
\leq C\sum_\nu|c_\nu|
\leq C\|h\|_{H^1(\mathbb R^q)}
\leq C\sum_{j=0}^{D_i}\|R_j^{(i)}h\|_1.
\end{equation}

\smallskip
\noindent{\it Step 3.}
By \eqref{e:hardy-riesz} and Fubini's theorem, the vertical slices of $R_i^{(1)}g$
and their Riesz transforms belong to $L^1(\mathbb R^m)\cap L^2(\mathbb R^m)$
for almost every $x$. Applying \eqref{e:atomic-area} on these slices and commuting
the Riesz transforms, we obtain
$$
\sum_{i=0}^d\|\mathsf T_2R_i^{(1)}g\|_{L^1(\mathcal H_2)}
\leq C\sum_{i=0}^d\sum_{j=0}^m\|R_j^{(2)}R_i^{(1)}g\|_1
\leq C\|g\|_{H_F^1}.
$$

Let $(e_l)$ be an orthonormal basis of $\mathcal H_2$ and
$h_l=\langle\mathsf T_2g,e_l\rangle_{\mathcal H_2}$. Since the convolutions commute
on $L^2(\mathbb R^d)$,
$$
R_i^{(1)}h_l=\langle\mathsf T_2R_i^{(1)}g,e_l\rangle_{\mathcal H_2}
\in L^1(\mathbb R^d)\cap L^2(\mathbb R^d),\qquad 0\leq i\leq d.
$$
Thus every finite sum $\sum_{l=1}^N\epsilon_lh_l$ belongs to $H^1(\mathbb R^d)\cap L^2(\mathbb R^d)$.
For independent Rademacher signs $\epsilon_l$, Khintchine's inequality and
\eqref{e:atomic-area} imply
\begin{align*}
\int_{\mathbb R^d}\bigg(\sum_{l=1}^N\|\mathsf T_1h_l(z)\|_{\mathcal H_1}^2\bigg)^{1/2}dz
&\leq C\mathbb E_\epsilon\bigg\|\mathsf T_1
       \Big(\sum_{l=1}^N\epsilon_lh_l\Big)\bigg\|_{L^1(\mathcal H_1)}
\leq C\sum_{i=0}^d\mathbb E_\epsilon
       \bigg\|\sum_{l=1}^N\epsilon_lR_i^{(1)}h_l\bigg\|_1\\
&\leq C\sum_{i=0}^d\int_{\mathbb R^d}
       \bigg(\sum_{l=1}^N
       |\langle\mathsf T_2R_i^{(1)}g(z),e_l\rangle_{\mathcal H_2}|^2\bigg)^{1/2}dz\\
&\leq C\sum_{i=0}^d\|\mathsf T_2R_i^{(1)}g\|_{L^1(\mathcal H_2)}
\leq C\|g\|_{H_F^1}.
\end{align*}
The $e_l$-coordinate of $\mathsf T_1\mathsf T_2g$ is $\mathsf T_1h_l$.
By Step~1, Parseval's identity and monotone convergence,
\begin{align*}
\|\mathcal A_{\lambda,F}(G_g)\|_1
&\leq CL^{-m/2}\|\mathsf T_1\mathsf T_2g\|_{L^1(\mathcal H_1\otimes\mathcal H_2)}\\
&=CL^{-m/2}\lim_{N\to\infty}
\int_{\mathbb R^d}\bigg(\sum_{l=1}^N\|\mathsf T_1h_l(z)\|_{\mathcal H_1}^2\bigg)^{1/2}dz
\leq C\|g\|_{H_F^1}.
\end{align*}
That is, we obtain \eqref{e:hardy-area}.
\end{proof}


\subsection{An open-set Carleson bound for the flag BMO norm}\label{s:deleted-bmo}

\begin{proposition}\label{p:carleson-bmo}
Fix $\kappa$ and a finite family $\mathcal E$, and write
$\lambda_{I,J}=\lambda_{\kappa,I,J}^{\mathrm{str}}(\mathcal E)$.
For every $\mathcal D\subseteq\{(I,J):\lambda_{I,J}>0\}$,
$$
\bigg\|\sum_{(I,J)\in\mathcal D}A_{\kappa,I,J}^{\mathrm{bi}}(\mathcal E)\bigg\|_{\mathrm{BMO}_F}
\leq C\bigg(\sup_{\substack{\Omega\ \mathrm{open}\\0<|\Omega|<\infty}}
\frac1{|\Omega|}\sum_{\substack{(I,J)\in\mathcal D\\I\times J\subseteq\Omega}}\lambda_{I,J}^2\bigg)^{1/2}.
$$
The constant is independent of $\kappa,\mathcal E,\mathcal D$ and the ratios $\ell(J)/\ell(I)$.
\end{proposition}

\begin{proof}
By homogeneity, assume that for every open $\Omega$ of finite measure,
$$
\sum_{\substack{(I,J)\in\mathcal D\\I\times J\subseteq\Omega}}\lambda_{I,J}^2\leq|\Omega|.
$$
Identify $(I,J)$ with $R=I\times J$, and write $\lambda_R=\lambda_{I,J}$ and
$\mathcal T_R=\mathcal T_{\kappa,I,J}^{\mathrm{str}}$.
Since $\lambda_R>0$ for $R\in\mathcal D$, we have $m_{\mathcal E}=1$ on each $\mathcal T_R$.
Let $H_{\mathcal D}=H_f^\sigma\sum_{R\in\mathcal D}\mathbf1_{\mathcal T_R}$.
Omitting $\kappa,\mathcal E$ from the summands, we have by \eqref{e:kernel-split},
$$
\mathcal R_\sigma H_{\mathcal D}
=\sum_{(I,J)\in\mathcal D}\big(A_{I,J}^{\mathrm{bi}}+A_{I,J}^{\mathrm{rem}}\big).
$$
We estimate this sum by duality and then subtract the remainders.

\smallskip
\noindent{\it Step 1.}
Fix $g\in H_F^1\cap L^2(\mathbb R^{n+m})$ and write
$G_g(u,v,t,s)=((K_{t,s}^\sigma)^* *g)(u,v)$.
By Plancherel, as in Proposition~\ref{p:reconstruction},
$$
\|G_g\|_{L^2_{\mathrm{par}}}
=\bigg(\int_0^\infty|a_\sigma(r)|^2\,\frac{dr}{r}\bigg)\|g\|_{L^2(\mathbb R^{n+m})}.
$$
For $R\in\mathcal D$, let $\gamma_R^2=\int_{\mathcal T_R}|G_g|^2\,d\nu$ and write
$$
S_\gamma(z)^2=\sum_{R\in\mathcal D}\frac{\gamma_R^2}{|R|}\mathbf1_R(z).
$$

We use the open-set argument of Chang and Fefferman~\cite{ChangFefferman}.
Let $M_{\mathrm s}$ denote the maximal operator over products of open axis-parallel cubes,
$$
M_{\mathrm s}h(x,y)=\sup_{\substack{I\ni x\\J\ni y}}\frac1{|I||J|}\int_{I\times J}|h|.
$$
For $\Omega_l=\{S_\gamma>2^l\}$, the set
$\widetilde\Omega_l=\{M_{\mathrm s}\mathbf1_{\Omega_l}>1/2\}$ is open, and
$$
|\widetilde\Omega_l|
\leq4\|M_{\mathrm s}\mathbf1_{\Omega_l}\|_2^2
\leq4\|M_xM_y\mathbf1_{\Omega_l}\|_2^2
\leq C|\Omega_l|.
$$
For each $R$ with $\gamma_R>0$, there is a unique $l$ such that
$$
|R\cap\Omega_l|>|R|/2\geq|R\cap\Omega_{l+1}|,
$$
since $\gamma_R/|R|^{1/2}\leq S_\gamma\leq\|S_\gamma\|_\infty$ on $R$.
Denote these rectangles by $\mathcal R_l$. The strict inequality persists on
a small open enlargement of $R$, so $\overline R\subset\widetilde\Omega_l$.
By Cauchy--Schwarz and the Carleson condition on $\widetilde\Omega_l$,
\begin{align*}
\sum_{R\in\mathcal R_l}\lambda_R\gamma_R
&\leq\bigg(\sum_{R\in\mathcal R_l}\lambda_R^2\bigg)^{1/2}
     \bigg(\sum_{R\in\mathcal R_l}\gamma_R^2\bigg)^{1/2}
\leq|\widetilde\Omega_l|^{1/2}
     \bigg(2\sum_{R\in\mathcal R_l}\frac{\gamma_R^2}{|R|}|R\setminus\Omega_{l+1}|\bigg)^{1/2}\\
&\leq|\widetilde\Omega_l|^{1/2}
     \bigg(2\int_{\widetilde\Omega_l\setminus\Omega_{l+1}}S_\gamma(z)^2\,dz\bigg)^{1/2}
\leq C2^l|\widetilde\Omega_l|.
\end{align*}
Summing in $l$, we obtain
\begin{equation}\label{e:coefficient-pairing}
\sum_{R\in\mathcal D}\lambda_R\gamma_R
\leq C\sum_l2^l|\widetilde\Omega_l|
\leq C\int_{\mathbb R^{n+m}}\sum_{2^l<S_\gamma(z)}2^l\,dz
\leq C\|S_\gamma\|_1.
\end{equation}

\smallskip
\noindent{\it Step 2.}
Applying the first inequality in \eqref{e:coefficient-control} to $G_g$ on the
regions $\mathcal T_R$, $R\in\mathcal D$, we have
$$
S_\gamma(z)\leq C\mathcal A_{\lambda_0,F}(G_g)(z).
$$
By the adjoint formula for $\mathcal R_\sigma$, Cauchy--Schwarz on each $\mathcal T_R$,
\eqref{e:coefficient-pairing} and \eqref{e:hardy-area},
\begin{align*}
|\langle g,\mathcal R_\sigma H_{\mathcal D}\rangle|
&=|\langle G_g,H_{\mathcal D}\rangle_{L^2_{\mathrm{par}}}|
=\bigg|\sum_{R\in\mathcal D}\int_{\mathcal T_R}G_g\,\overline{H_f^\sigma}\,d\nu\bigg|\\
&\leq\sum_{R\in\mathcal D}\lambda_R\gamma_R
\leq C\|S_\gamma\|_1
\leq C\|\mathcal A_{\lambda_0,F}(G_g)\|_1
\leq C\|g\|_{H_F^1}.
\end{align*}
Taking the supremum over $\|g\|_{H_F^1}\leq1$, we obtain
\begin{equation}\label{e:deleted-bmo}
\bigg\|\sum_{(I,J)\in\mathcal D}\big(A_{I,J}^{\mathrm{bi}}+A_{I,J}^{\mathrm{rem}}\big)\bigg\|_{\mathrm{BMO}_F}
=\|\mathcal R_\sigma H_{\mathcal D}\|_{\mathrm{BMO}_F}\leq C.
\end{equation}

\smallskip
\noindent{\it Step 3.}
Fix $k\geq0$ and consider the rectangles with $\ell(J)=2^k\ell(I)$.
For $T_k(x,y)=(x,2^ky)$, the cubes $Q=T_k^{-1}(I\times J)=I\times 2^{-k}J$
belong to a single dyadic system on $\mathbb R^{n+m}$. Write
$$
u_Q=A_{I,J}^{\mathrm{rem}}\circ T_k,\qquad
d_Q=2^{-2k}\frac{\lambda_{I,J}}{(|I||J|)^{1/2}}.
$$
Since $\ell(Q)=\ell(I)=2^{-k}\ell(J)$,
$$
\ell(Q)\nabla u_Q
=\big(\ell(I)\nabla_xA_{I,J}^{\mathrm{rem}},\ell(J)\nabla_yA_{I,J}^{\mathrm{rem}}\big)\circ T_k.
$$
Thus \eqref{e:rem-deriv} implies
\begin{equation}\label{e:rem-bounds}
|u_Q|+\ell(Q)|\nabla u_Q|\leq Cd_Q\mathbf1_Q.
\end{equation}
The support lies at distance at least $\ell(Q)/24$ from $\partial Q$, and
$\int_Q u_Q=2^{-km}\int_{I\times J}A_{I,J}^{\mathrm{rem}}=0$.

For a cube $U$ in this dyadic system, apply the Carleson condition to open
rectangles decreasing to $\overline{T_kU}$. Since $|T_kU|=2^{km}|U|$,
\begin{equation}\label{e:rem-carleson}
\frac1{|U|}\sum_{Q\subseteq U}d_Q^2|Q|
=\frac{2^{-4k}}{|T_kU|}
\sum_{\substack{(I,J)\in\mathcal D,\ \ell(J)=2^k\ell(I)\\I\times J\subseteq T_kU}}\lambda_{I,J}^2
\leq2^{-4k}.
\end{equation}
In particular, $d_Q\leq2^{-2k}$. By \eqref{e:rem-bounds}, cancellation and
\eqref{e:row-orthogonality}, for every subfamily $\mathcal F$ we have
\begin{equation}\label{e:rem-l2}
\bigg\|\sum_{Q\in\mathcal F}u_Q\bigg\|_{L^2(\mathbb R^{n+m})}^2
\leq C\sum_{Q\in\mathcal F}d_Q^2|Q|.
\end{equation}

Let $h_k=\sum_Q u_Q$. For a test cube $P$ of side $r$, let $U$ range over the dyadic
cubes meeting $P$ with $r\leq\ell(U)<2r$. As in Proposition~\ref{p:vertical-selection},
use \eqref{e:rem-l2} for $\ell(Q)\leq r$ and the gradient bound \eqref{e:rem-bounds}
for $\ell(Q)>r$. By \eqref{e:rem-carleson},
\begin{align*}
\bigg(\frac1{|P|}\int_P|h_k-(h_k)_P|^2\bigg)^{1/2}
&\leq C\bigg(\frac1{|P|}\sum_U\sum_{Q\subseteq U}d_Q^2|Q|\bigg)^{1/2}
+Cr\sum_{\substack{\ell(Q)>r\\Q\cap P\ne\emptyset}}\frac{d_Q}{\ell(Q)}\\
&\leq C2^{-2k}\bigg[\bigg(\frac{\sum_U|U|}{|P|}\bigg)^{1/2}
+\sum_{j\geq0}2^{-j}\bigg]\leq C2^{-2k},
\end{align*}
where at most $C_{n,m}$ cubes meet $P$ at each dyadic side length at least $r$.

Applying \eqref{e:dilated-bmo} with $\rho=2^k$, we obtain
$$
\bigg\|\sum_{\substack{(I,J)\in\mathcal D\\\ell(J)=2^k\ell(I)}}A_{I,J}^{\mathrm{rem}}\bigg\|_{\mathrm{BMO}_F}
\leq C\|h_k\|_{\mathrm{BMO}(\mathbb R^{n+m})}\leq C2^{-2k}.
$$
Subtracting the remainders from \eqref{e:deleted-bmo}, we conclude that
$$
\bigg\|\sum_{(I,J)\in\mathcal D}A_{I,J}^{\mathrm{bi}}\bigg\|_{\mathrm{BMO}_F}
\leq\|\mathcal R_\sigma H_{\mathcal D}\|_{\mathrm{BMO}_F}
+\sum_{k\geq0}\bigg\|\sum_{\substack{(I,J)\in\mathcal D\\\ell(J)=2^k\ell(I)}}
A_{I,J}^{\mathrm{rem}}\bigg\|_{\mathrm{BMO}_F}
\leq C+C\sum_{k\geq0}2^{-2k}\leq C.
$$
By rescaling, we obtain  Proposition~\ref{p:carleson-bmo}.
\end{proof}


\subsection{Construction of the bounded approximant}\label{s:approximation-proof}

\begin{proof}[Proof of Theorem~\ref{t:approximation}]
If $\|S_F(f)\|_\infty=0$, then $f=0$ by \eqref{e:flag-plancherel}.
Suppose $0<\varepsilon\leq\|S_F(f)\|_\infty$. Fix a finite family $\mathcal E$
and a pair of shifts $\kappa$.

For $0<\vartheta_x\leq1/2$, let $B_J=\sum_{I\in\mathcal K_J}A_{I,J}$ be the rows
selected in Proposition~\ref{p:horizontal-selection}, and write $G=\sum_JB_J$.
By Proposition~\ref{p:carleson-bmo}, \eqref{e:deleted-carleson} and \eqref{e:row-coefficients},
\begin{align*}
\bigg\|\sum_{I,J}A_{I,J}-G\bigg\|_{\mathrm{BMO}_F}
&=\bigg\|\sum_J\sum_{I\in\mathcal D_J}A_{I,J}\bigg\|_{\mathrm{BMO}_F}\\
&\leq C\bigg(\sup_{\substack{\Omega\ \mathrm{open}\\0<|\Omega|<\infty}}
\frac1{|\Omega|}\sum_J\sum_{\substack{I\in\mathcal D_J\\I\times J\subseteq\Omega}}
c_{I,J}^2|I||J|\bigg)^{1/2}
\leq C\vartheta_x\|S_F(f)\|_\infty.
\end{align*}

Each $B_J(x,\cdot)$ is supported in $J$, has integral zero, and by \eqref{e:scaled-sobolev},
$$
|B_J(x,y)|+\ell(J)|\nabla_yB_J(x,y)|
\leq C\|B_J(x,\cdot)\|_J\mathbf1_J(y).
$$
For the selection in $y$, \eqref{e:selected-row} and \eqref{e:row-coefficients} imply
$$
\sum_J\|B_J(x,\cdot)\|_J^2\mathbf1_J(y)
\leq C\frac{\log(e/\vartheta_x)}{\vartheta_x^2}
\sum_{I,J}c_{I,J}^2\mathbf1_{I\times J}(x,y)
\leq C\frac{\log(e/\vartheta_x)}{\vartheta_x^2}S_F(f)(x,y)^2.
$$
With these row coefficients, the proof of Proposition~\ref{p:vertical-selection}
applies at each $x$. For $0<\vartheta_y\leq1/2$, the selected sum
$b^{\mathrm{bi}}(x,y)=\sum_J\eta_J(x)B_J(x,y)$ is measurable and satisfies
$$
\|b^{\mathrm{bi}}\|_\infty
\leq C\frac{\sqrt{\log(e/\vartheta_x)\log(e/\vartheta_y)}}{\vartheta_x\vartheta_y}
\|S_F(f)\|_\infty.
$$
Since $G-b^{\mathrm{bi}}$ is bounded and supported on a finite union of rectangles,
we may apply \eqref{e:vertical-bmo} to the discarded rows. Thus
\begin{align*}
\bigg\|\sum_{I,J}A_{I,J}-b^{\mathrm{bi}}\bigg\|_{\mathrm{BMO}_F}
&\leq \bigg\|\sum_{I,J}A_{I,J}-G\bigg\|_{\mathrm{BMO}_F}
+C\operatorname*{ess\,sup}_x
\left \|G(x,\cdot)-b^{\mathrm{bi}}(x,\cdot) \right \|_{\mathrm{BMO}(\mathbb R^m)}\\
&\leq C\bigg(\vartheta_x+
\frac{\vartheta_y}{\vartheta_x}\sqrt{\log(e/\vartheta_x)}\bigg)  \,  \|S_F(f)\|_\infty.
\end{align*}

For the remainder and cubic rows, \eqref{e:coefficient-control} implies
$$
\sum_J\big[(c_J^{\mathrm{rem}}(x))^2+(c_J^{\mathrm{cub}}(x))^2\big]\mathbf1_J(y)
\leq C S_F(f)(x,y)^2.
$$
Apply Proposition~\ref{p:vertical-selection} to each family with $\vartheta=\vartheta_x$.
The selected sums satisfy
$$
\|b^{\mathrm{rem}}\|_\infty+\|b^{\mathrm{cub}}\|_\infty
\leq C\frac{\sqrt{\log(e/\vartheta_x)}}{\vartheta_x}\|S_F(f)\|_\infty,
$$
and by \eqref{e:vertical-bmo},
$$
\|G^{\mathrm{rem}}-b^{\mathrm{rem}}\|_{\mathrm{BMO}_F}
+\|G^{\mathrm{cub}}-b^{\mathrm{cub}}\|_{\mathrm{BMO}_F}
\leq C\vartheta_x\|S_F(f)\|_\infty.
$$

Let $b_{\varepsilon,\mathcal E}$ be the sum of
$b^{\mathrm{bi}}+b^{\mathrm{rem}}+b^{\mathrm{cub}}$ over the $3^{n+m}$ pairs of shifts.
For some $C_1\geq1$ depending only on $n,m$,
$$
\|f_{\mathcal E}-b_{\varepsilon,\mathcal E}\|_{\mathrm{BMO}_F}
\leq C_1\bigg(\vartheta_x+
\frac{\vartheta_y}{\vartheta_x}\sqrt{\log(e/\vartheta_x)}\bigg)\|S_F(f)\|_\infty.
$$
Choose the parameters so that the two terms on the right are equal:
$$
\vartheta_x=\frac{\varepsilon}{2C_1\|S_F(f)\|_\infty},\qquad
\vartheta_y=\frac{\vartheta_x^2}{\sqrt{\log(e/\vartheta_x)}}.
$$
Then $0<\vartheta_y\leq\vartheta_x\leq1/2$ and
$$
\|f_{\mathcal E}-b_{\varepsilon,\mathcal E}\|_{\mathrm{BMO}_F}
\leq2C_1\vartheta_x\|S_F(f)\|_\infty=\varepsilon.
$$
Since $\log(e/\vartheta_y)\leq\tfrac52\log(e/\vartheta_x)$, we also have
\begin{align*}
\|b_{\varepsilon,\mathcal E}\|_\infty
&\leq C\bigg(
\frac{\sqrt{\log(e/\vartheta_x)\log(e/\vartheta_y)}}{\vartheta_x\vartheta_y}
+\frac{\sqrt{\log(e/\vartheta_x)}}{\vartheta_x}\bigg)  \,  \|S_F(f)\|_\infty\\
&\leq C\frac{[\log(e/\vartheta_x)]^{3/2}}{\vartheta_x^3}\|S_F(f)\|_\infty
\leq C\frac{\|S_F(f)\|_\infty^4}{\varepsilon^3}
\bigg[\log\bigg(e+\frac{\|S_F(f)\|_\infty}{\varepsilon}\bigg)\bigg]^{3/2}.
\end{align*}

Choose finite families $\mathcal E_N\uparrow\mathcal I_{\rm par}$.
By \eqref{e:finite-recon-limit}, $f_{\mathcal E_N}\to f$ in $L^2(\mathbb R^{n+m})$.
By weak-* compactness and the separability of $L^1(\mathbb R^{n+m})$,
a subsequence of $b_{\varepsilon,\mathcal E_N}$ converges weak-* to $b_\varepsilon$
with the same $L^\infty$ bound. For $g\in H_F^1\cap L^2(\mathbb R^{n+m})$,
we have $g\in L^1\cap L^2$, and hence
$$
|\langle g,f-b_\varepsilon\rangle|
=\lim_{j\to\infty}|\langle g,f_{\mathcal E_{N_j}}-b_{\varepsilon,\mathcal E_{N_j}}\rangle|
\leq\varepsilon\|g\|_{H_F^1}.
$$
Taking the supremum over $\|g\|_{H_F^1}\leq1$, we obtain
$\|f-b_\varepsilon\|_{\mathrm{BMO}_F}\leq\varepsilon$, as desired.
\end{proof}


\subsection{Rectangular moments and the size of bounded approximants}\label{s:rect-approx}

\begin{proposition}\label{p:rect-moments}
Let $f\in L^2(\mathbb R^{n+m})+L^\infty$ represent a bounded functional on $H_F^1$.
For $2\leq p<\infty$,
\begin{equation}\label{e:rect-bmo-moments}
\sup_{R=I\times J}\bigg(\frac1{|R|}\int_R|Q_Rf|^p\bigg)^{1/p}
\leq Cp^2\|f\|_{\mathrm{BMO}_F},
\end{equation}
where $I,J$ range over axis-parallel cubes in $\mathbb R^n,\mathbb R^m$, respectively.
If $f\in\mathcal X_F\cap L^2(\mathbb R^{n+m})$, the left-hand side is $o(p^2)$ as $p\to\infty$, and
\begin{equation}\label{e:closure-half}
\sup_R\frac1{|R|}\int_R\exp\big(\lambda|Q_Rf|^{1/2}\big)\,dz<\infty,\qquad \lambda>0.
\end{equation}
\end{proposition}

\begin{proof}
Fix $R=I\times J$, let $p'=p/(p-1)$, and extend $h=Q_Ru$ by zero off $R$ for
$u\in L^{p'}(R)\cap L^2(R)$. Then $h$ has integral zero in each variable and
$\|h\|_{p'}\leq4\|u\|_{L^{p'}(R)}$. Write
$$
g_I(y)=\int_I|h(v,y)|\,dv,\qquad g_J(x)=\int_J|h(x,w)|\,dw.
$$

Let $K_t^{(q)}$ be the heat kernel on $\mathbb R^q$. For a cube $D$ with centre $z_D$ and $z\notin3D$,
$$
\sup_{\substack{t>0\\v\in D}}|K_t^{(q)}(z-v)-K_t^{(q)}(z-z_D)|
\leq C_q\ell(D)\sup_{\substack{t>0\\w\in D}}|\nabla K_t^{(q)}(z-w)|
\leq C_q\frac{\ell(D)}{|z-z_D|^{q+1}}=:k_D(z),
$$
where $\int_{(3D)^c}k_D\leq C_q$. Subtracting the horizontal kernel at the centre of $I$
when $x\notin3I$, and the vertical kernel at the centre of $J$ when $y\notin3J$, we obtain
\begin{align*}
\mathcal M_Fh(x,y)
&=\sup_{t,s>0}|e^{t^2\Delta_x}e^{(t^2+s^2)\Delta_y}h(x,y)|
\leq\sup_{a,b>0}|e^{a\Delta_x}e^{b\Delta_y}h(x,y)|\\
&\leq C
\begin{cases}
M_xM_y|h|(x,y),       &(x,y)\in3I\times3J,\\[2pt]
k_I(x)M_yg_I(y),      &(x,y)\in(3I)^c\times3J,\\[2pt]
k_J(y)M_xg_J(x),      &(x,y)\in3I\times(3J)^c,\\[2pt]
k_I(x)k_J(y)\|h\|_1,  &(x,y)\in(3I)^c\times(3J)^c.
\end{cases}
\end{align*}
By H\"older's inequality and the $L^{p'}$ maximal bound,
\begin{align*}
\|\mathcal M_Fh\|_1
&\leq C\Big(|R|^{1/p}\|M_xM_y|h|\|_{p'}
  +|J|^{1/p}\|M_yg_I\|_{p'}+|I|^{1/p}\|M_xg_J\|_{p'}+\|h\|_1\Big)\\
&\leq C\Big(p^2|R|^{1/p}\|h\|_{p'}
  +p|J|^{1/p}\|g_I\|_{p'}+p|I|^{1/p}\|g_J\|_{p'}+\|h\|_1\Big)\\
&\leq C(p^2+2p+1)|R|^{1/p}\|h\|_{p'}
\leq Cp^2|R|^{1/p}\|u\|_{L^{p'}(R)},
\end{align*}
where $\|g_I\|_{p'}\leq|I|^{1/p}\|h\|_{p'}$ and $\|g_J\|_{p'}\leq|J|^{1/p}\|h\|_{p'}$.
Thus $h\in H_F^1\cap L^2$, and self-adjointness of $Q_R$ on $L^2(R)$ implies
$$
\bigg|\int_R u\,\overline{Q_Rf}\bigg|
=|\langle h,f\rangle|
\leq\|\mathcal M_Fh\|_1\|f\|_{\mathrm{BMO}_F}
\leq Cp^2|R|^{1/p}\|u\|_{L^{p'}(R)}\|f\|_{\mathrm{BMO}_F}.
$$
By density and duality, this proves \eqref{e:rect-bmo-moments}.

Write $m_p(f)$ for the left-hand side of \eqref{e:rect-bmo-moments}.
For $f\in\mathcal X_F\cap L^2$ and $b\in L^\infty$,
$$
m_p(f)\leq m_p(b)+m_p(f-b)\leq4\|b\|_\infty+Cp^2\|f-b\|_{\mathrm{BMO}_F}.
$$
Hence
$$
0\leq\limsup_{p\to\infty}\frac{m_p(f)}{p^2}
\leq\inf_{b\in L^\infty}\limsup_{p\to\infty}
\bigg(\frac{4\|b\|_\infty}{p^2}+C\|f-b\|_{\mathrm{BMO}_F}\bigg)
=C\inf_{b\in L^\infty}\|f-b\|_{\mathrm{BMO}_F}=0.
$$
For every $\lambda>0$ and integer $j\geq4$,
$$
\sup_R\bigg(\frac{\lambda^j}{j!|R|}\int_R|Q_Rf|^{j/2}\bigg)^{1/j}
=\frac{\lambda\,m_{j/2}(f)^{1/2}}{(j!)^{1/j}}
\leq\frac{e\lambda}{2}\bigg(\frac{m_{j/2}(f)}{(j/2)^2}\bigg)^{1/2}\to0.
$$
The $p=2$ estimate bounds the first four terms uniformly in $R$.
Summing the exponential series, we obtain \eqref{e:closure-half}.
\end{proof}

\begin{remark}\label{r:approx-lower}
For every $t>0$ and $0<\varepsilon\leq c_0t$, there is a real Schwartz function $F$
with $\|S_F(F)\|_\infty=t$ such that
$$
\inf_{\substack{b\in L^\infty\\ \|F-b\|_{\mathrm{BMO}_F}\leq\varepsilon}}\|b\|_\infty
\geq c\,t^2/\varepsilon.
$$
Indeed, let $f_L,R_L$ be as in the proof of Proposition~\ref{p:two-orders} and write
$F_L=tf_L/\|S_F(f_L)\|_\infty$. Since $Q_{R_L}F_L=F_L$ on $R_L$,
\eqref{e:trace-free-sup} and \eqref{e:zero-trace-values} imply
$$
\bigg(\frac1{|R_L|}\int_{R_L}|Q_{R_L}F_L|^L\bigg)^{1/L}
\geq\frac{tL^2}{4\|S_F(f_L)\|_\infty}\big(c2^{-(n+m)L}\big)^{1/L}
\geq c_1tL.
$$
If $\|F_L-b\|_{\mathrm{BMO}_F}\leq\varepsilon$, then \eqref{e:rect-bmo-moments} gives
$$
c_1tL\leq m_L(F_L-b)+4\|b\|_\infty\leq C\varepsilon L^2+4\|b\|_\infty.
$$
Choose $L=\lfloor c_1t/(2C\varepsilon)\rfloor$. For $c_0$ sufficiently small,
$\|b\|_\infty\geq c_1tL/8\geq ct^2/\varepsilon$.
\end{remark}

\begin{remark}\label{r:scalar-approx}
For the one-parameter Poisson area,
$$
\sup_{\substack{f\in L^2(\mathbb R^d)\ \mathrm{real}\\ \|S_{\mathrm{cl}}^{(d)}(f)\|_\infty=1}}
\inf_{\substack{b\in L^\infty\\ \|f-b\|_{\mathrm{BMO}(\mathbb R^d)}\leq\varepsilon}}
\|b\|_\infty\simeq_d\varepsilon^{-1},\qquad \varepsilon\to0+.
$$
Fix such an $f$ and let $w=e^{\lambda f}$, $\lambda\geq1$. By the exponential-square estimate
of Chang, Wilson and Wolff~\cite[Theorem~3.2]{CWW},
$$
[w]_{A_2}
=\sup_Q\langle e^{\lambda(f-f_Q)}\rangle_Q\langle e^{-\lambda(f-f_Q)}\rangle_Q
\leq e^{\lambda^2/(2c_d)}\sup_Q\langle e^{c_d|f-f_Q|^2}\rangle_Q^2
\leq C_de^{C_d\lambda^2},
$$
where $\lambda|f-f_Q|\leq c_d|f-f_Q|^2+\lambda^2/(4c_d)$.

Let $M$ be the Hardy--Littlewood maximal operator. Following Coifman, Jones and Rubio de Francia~\cite{CJR}, write
$$
Tg=w^{-1/2}M(w^{1/2}g)+w^{1/2}M(w^{-1/2}g).
$$
By Buckley's estimate~\cite{Buckley},
$$
A:=2\|T\|_{2\to2}
\leq2\big(\|M\|_{L^2(w^{-1})\to L^2(w^{-1})}+\|M\|_{L^2(w)\to L^2(w)}\big)
\leq C_d[w]_{A_2}\leq C_de^{C_d\lambda^2}.
$$
Choose $h\in L^2(\mathbb R^d)$ with $h>0$ almost everywhere and let
$$
H=\sum_{j=0}^\infty A^{-j}T^jh,\qquad u=w^{1/2}H,\qquad v=w^{-1/2}H.
$$
Then $\|H\|_2\leq2\|h\|_2$, $w=u/v$, and
$$
\frac{Mu}{u}+\frac{Mv}{v}
=\frac{TH}{H}
\leq\frac1H\sum_{j=0}^\infty A^{-j}T^{j+1}h
=A\frac{H-h}{H}\leq A.
$$
Since $Mu/u,Mv/v\geq1$ almost everywhere, the function
$$
b=\frac1\lambda\log\frac{u\,Mv}{v\,Mu}
=f-\frac1\lambda\big(\log Mu-\log Mv\big)
$$
satisfies
$$
\|b\|_\infty
=\frac1\lambda\bigg\|\log\frac{Mv}{v}-\log\frac{Mu}{u}\bigg\|_\infty
\leq\frac{\log A}{\lambda}\leq C_d\lambda.
$$
By the Coifman--Rochberg estimate~\cite{CoifmanRochberg},
$$
\|f-b\|_{\mathrm{BMO}}
=\lambda^{-1}\|\log Mu-\log Mv\|_{\mathrm{BMO}}
\leq\lambda^{-1}\big(\|\log Mu\|_{\mathrm{BMO}}+\|\log Mv\|_{\mathrm{BMO}}\big)
\leq C_d/\lambda.
$$
Taking $\lambda=C_d/\varepsilon$ with $C_d$ sufficiently large proves the upper bound for $0<\varepsilon\leq1$.

For the lower bound, take $k=1$ and $M=\lceil p\rceil$ in the construction of
Proposition~\ref{p:flag-reverse-sharp}, with $p$ sufficiently large.
Write $g_M$ for the resulting real Schwartz function, $Q$ for its cube and
$f_p=g_M/\|S_{\mathrm{cl}}^{(d)}(g_M)\|_\infty$. Since $(f_p)_Q=0$,
\eqref{e:flag-sharp-sup}, \eqref{e:flag-large-values} and the John--Nirenberg inequality imply
\begin{align*}
c_d\sqrt p
&\leq\frac{M(c2^{-dM})^{1/p}}{2\|S_{\mathrm{cl}}^{(d)}(g_M)\|_\infty}
\leq\langle|f_p|^p\rangle_Q^{1/p}\\
&\leq\langle|(f_p-b)-(f_p-b)_Q|^p\rangle_Q^{1/p}
     +\|b-b_Q\|_{L^\infty(Q)}
\leq C_dp\|f_p-b\|_{\mathrm{BMO}}+2\|b\|_\infty.
\end{align*}
For $p=(c_d/(2C_d\varepsilon))^2$ and $\varepsilon$ sufficiently small, every $b$ with
$\|f_p-b\|_{\mathrm{BMO}}\leq\varepsilon$ therefore satisfies
$$
2\|b\|_\infty\geq c_d\sqrt p-C_dp\varepsilon
=\frac{c_d^2}{4C_d\varepsilon}.
$$

For the two-step flag, the functions in \eqref{e:two-weight-growth} satisfy
$\|S_F(f_M)\|_\infty=1$ and, on one fixed rectangle $R_0$,
$$
\langle e^{\lambda f_M}\rangle_{R_0}\langle e^{-\lambda f_M}\rangle_{R_0}
\geq1+c\lambda^2M\to\infty,\qquad \lambda\ne0.
$$
\end{remark}


\subsection{Bounded approximation and exponential integrability on rectangles}\label{s:strict-classes}

The logarithmic example below belongs to $\mathscr E_F$ and has positive $\mathrm{BMO}_F$ distance from $L^\infty$.
In the other direction, Proposition~\ref{p:exp-divergence} uses two Riesz transforms along the flag to construct
$F\in\mathcal X_F\cap L^2$ for which, on one flag rectangle $R_0$,
$$
\int_{R_0}\exp\big(\lambda|Q_{R_0}F|^\gamma\big)\,dz=\infty,\qquad \lambda>0,\quad \gamma>1/2.
$$

\begin{proposition}\label{p:log-separation}
For every $n,m\geq1$, there is a real $F\in\mathscr E_F$ such that
$$
\inf_{b\in L^\infty}\|F-b\|_{\mathrm{BMO}_F}>0.
$$
Moreover, for some $c>0$,
$$
\sup_{\ell(I)\leq\ell(J)}\frac1{|I||J|}\int_{I\times J}e^{c|F-F_{I\times J}|}\,dz<\infty.
$$
\end{proposition}

\begin{proof}
Choose a nonzero real $a_0\in C_c^\infty(\mathbb R^n)$ with $\int a_0=0$, and let
$$
F(x,y)=a_0(x)L_0(y),\qquad L_0(y)=\log_+(1/|y|),\qquad \log_+t=\max\{0,\log t\}.
$$
For a cube $J\subset\mathbb R^m$, write $s=\ell(J)$, $D_J=\sup_{y\in J}|y|$ and $a_J=\log_+(1/D_J)$.
Since $s/2\leq D_J\leq |y|+\sqrt m\,s$ on $J$, for $\gamma=m/2$ we have
\begin{align*}
\frac1{|J|}\int_J e^{\gamma(L_0-a_J)}\,dy
&\leq\frac1{|J|}\int_J\bigg(\frac{D_J}{|y|}\bigg)^\gamma\,dy
\leq\frac1{|J|}\int_J\bigg(1+\frac{\sqrt m\,s}{|y|}\bigg)^\gamma\,dy\\
&\leq C_m\bigg(1+s^{\gamma-m}\int_{|y|<s}|y|^{-\gamma}\,dy\bigg)\leq C_m.
\end{align*}
By Jensen's inequality,
$$
(L_0)_J\leq a_J+\gamma^{-1}\log\langle e^{\gamma(L_0-a_J)}\rangle_J
\leq C_m+\log_+(1/s).
$$
Since $L_0\geq a_J$ on $J$, another application of Jensen's inequality gives
\begin{equation}\label{e:log-oscillation}
\frac1{|J|}\int_J e^{\gamma|L_0-(L_0)_J|}\,dy
\leq e^{\gamma((L_0)_J-a_J)}\langle e^{\gamma(L_0-a_J)}\rangle_J
\leq\langle e^{\gamma(L_0-a_J)}\rangle_J^2\leq C_m.
\end{equation}
Thus $L_0\in\mathrm{BMO}(\mathbb R^m)\cap L^2(\mathbb R^m)$, since
$\int_0^1r^{m-1}\log^2(1/r)\,dr<\infty$. Hence $F\in L^2$ and, by \eqref{e:vertical-bmo},
$$
\|F\|_{\mathrm{BMO}_F}\leq C\|a_0\|_\infty\|L_0\|_{\mathrm{BMO}(\mathbb R^m)}<\infty.
$$

On every $R=I\times J$, we have $Q_RF=(a_0-(a_0)_I)(L_0-(L_0)_J)$.
Choose $c>0$ with $2c\|a_0\|_\infty\leq\gamma$. By \eqref{e:log-oscillation},
$$
\sup_R\frac1{|R|}\int_R e^{c|Q_RF|}\,dz
\leq\sup_J\frac1{|J|}\int_J e^{2c\|a_0\|_\infty|L_0-(L_0)_J|}\,dy
\leq C_m,
$$
so $F\in\mathscr E_F$. If $r=\ell(I)\leq s=\ell(J)$, then
$$
\|a_0-(a_0)_I\|_{L^\infty(I)}(L_0)_J
\leq C\min\{r,1\}\big(1+\log_+(1/s)\big)
\leq C\min\{s,1\}\big(1+\log_+(1/s)\big)\leq C.
$$
Using $F-F_R=a_0(L_0-(L_0)_J)+(a_0-(a_0)_I)(L_0)_J$, we obtain
$$
\frac1{|R|}\int_R e^{c|F-F_R|}\,dz
\leq e^{c\|a_0-(a_0)_I\|_{L^\infty(I)}(L_0)_J}
\frac1{|J|}\int_J e^{c\|a_0\|_\infty|L_0-(L_0)_J|}\,dy\leq C.
$$

Choose $0\leq\rho\in C_c^\infty(B(0,1))$ on $\mathbb R^m$ with $\int\rho=1$.
For $\rho_r(y)=r^{-m}\rho(y/r)$ and $h_j=\rho_{2^{-j}}-\rho_{2^{-j+1}}$, let
$$
g_N=\frac1N a_0\otimes(\rho_{2^{-N}}-\rho)
=\frac1N\sum_{j=1}^N a_0\otimes h_j,\qquad N\geq1.
$$
The $h_j$ are dilates of the fixed mean-zero function $\rho_{1/2}-\rho$.
By the classical heat maximal estimate for atoms~\cite{FS,Latter},
\begin{align*}
\|g_N\|_{H_F^1}
&\leq\frac1N\sum_{j=1}^N
\bigg\|\sup_{t,s>0}\big|(e^{t^2\Delta_x}a_0)(x)(e^{(t^2+s^2)\Delta_y}h_j)(y)\big|\bigg\|_{L^1(dx\,dy)}\\
&\leq\frac C N\|a_0\|_{H^1(\mathbb R^n)}\sum_{j=1}^N\|h_j\|_{H^1(\mathbb R^m)}
=C\|a_0\|_{H^1(\mathbb R^n)}\|\rho_{1/2}-\rho\|_{H^1(\mathbb R^m)}\leq C.
\end{align*}
Since $L_0(2^{-N}y)-L_0(y)=N\log2$ for $0<|y|<1$, for every $b\in L^\infty$ we have
\begin{align*}
\|a_0\|_2^2\log2
&=\frac{\|a_0\|_2^2}{N}\int_{\mathbb R^m}\rho(y)\big(L_0(2^{-N}y)-L_0(y)\big)\,dy
=\langle g_N,F\rangle\\
&\leq|\langle g_N,F-b\rangle|+|\langle g_N,b\rangle|
\leq\|g_N\|_{H_F^1}\|F-b\|_{\mathrm{BMO}_F}+\|g_N\|_1\|b\|_\infty\\
&\leq C\|F-b\|_{\mathrm{BMO}_F}+\frac{2\|a_0\|_1}{N}\|b\|_\infty.
\end{align*}
Letting $N\to\infty$ and taking the infimum over $b\in L^\infty$, we obtain the claimed positive distance.
\end{proof}



For the function $F$ in Proposition~\ref{p:log-separation}, choose $c>0$ as in the definition of $\mathscr E_F$.
Since $\lambda\sqrt t\leq ct+\lambda^2/(4c)$ for $t\geq0$,
$$
\sup_R\frac1{|R|}\int_R e^{\lambda|Q_RF|^{1/2}}\,dz
\leq e^{\lambda^2/(4c)}\sup_R\frac1{|R|}\int_R e^{c|Q_RF|}\,dz<\infty,\qquad \lambda>0.
$$
Thus \eqref{e:closure-half} holds for a function at positive $\mathrm{BMO}_F$ distance from $L^\infty$.

\begin{proposition}\label{p:exp-divergence}
For every $n,m\geq1$, there are a real $F\in\mathcal X_F\cap L^2(\mathbb R^{n+m})$ and a fixed flag rectangle $R_0$ such that
$$
\int_{R_0}\exp\big(\lambda|Q_{R_0}F|^\gamma\big)\,dz=\infty,\qquad \gamma>1/2,\quad \lambda>0.
$$
\end{proposition}

\begin{proof}
Write $x=(x_1,x')$, $y=(y_1,y')$ and let $\mathcal T=R_1^{(1)}R_1^{(2)}$, with the full-space transform in $x_1$
and the vertical transform in $y_1$. Since $\mathcal T^*=\mathcal T$, \eqref{e:hardy-riesz} implies
\begin{equation}\label{e:riesz-bmo}
\|\mathcal T B\|_{\mathrm{BMO}_F}
=\sup_{\substack{g\in H_F^1\cap L^2\\ \|g\|_{H_F^1}\leq1}}|\langle\mathcal T g,B\rangle|
\leq C\|B\|_\infty,\qquad B\in L^\infty\cap L^2.
\end{equation}

Choose $0\leq\chi\leq1$ in $C_c^\infty((-1,1))$, equal to one on $[0,1/2]$, and let
$$
v(u)=
\begin{cases}
\chi(u)/\log(e+\log(1/u)),&0<u<1,\\[2pt]
0,&\text{otherwise}.
\end{cases}
$$
Choose smooth compactly supported $0\leq\phi,\psi\leq1$ in $x',y'$, equal to one on
$[-2,2]^{n-1}$ and $[-2,2]^{m-1}$, respectively. Let
$$
B(x,y)=v(x_1)\phi(x')v(y_1)\psi(y'),\qquad F=\mathcal T B,
$$
omitting the transverse factor in dimension one. Then $F$ is real and $\|F\|_2\leq\|B\|_2$.
For $B_j\in C_c^\infty$ converging uniformly to $B\in C_c$, Fourier inversion gives
$\|\mathcal T B_j\|_\infty\leq C\|\widehat B_j\|_1$. Hence
$$
\inf_{b\in L^\infty}\|F-b\|_{\mathrm{BMO}_F}
\leq\|F-\mathcal T B_j\|_{\mathrm{BMO}_F}
=\|\mathcal T(B-B_j)\|_{\mathrm{BMO}_F}
\leq C\|B-B_j\|_\infty\to0,
$$
so $F\in\mathcal X_F\cap L^2$.

Let $K_t(u)=(4\pi t)^{-1/2}e^{-u^2/(4t)}$. For $z,t>0$,
\begin{align*}
0\leq V_t(z)&:=(K_t*v)'(-z)
=\int_0^1\frac{z+u}{2t}K_t(z+u)v(u)\,du\\
&\leq\min\{K_t(z),\|K_t'\|_\infty\|v\|_1\}
\leq C\min\{t^{-1/2}e^{-z^2/(4t)},t^{-1}\}.
\end{align*}
For $x_1=-r_x<0$ and $y_1=-r_y<0$, the heat representation and
$e^{t\Delta_{x,y}}e^{s\Delta_y}=e^{t\Delta_x}e^{(t+s)\Delta_y}$ imply
\begin{equation}\label{e:riesz-integral}
F(x,y)=\frac1\pi\int_0^\infty\int_0^\infty
V_t(r_x)V_{t+s}(r_y)(e^{t\Delta_{x'}}\phi)(x')(e^{(t+s)\Delta_{y'}}\psi)(y')
\,\frac{ds\,dt}{\sqrt{st}}.
\end{equation}
Indeed, truncation to $[L^{-1},L]$ multiplies each Riesz multiplier by
$$
\frac{|\zeta|}{\sqrt\pi}\int_{L^{-1}}^L e^{-t|\zeta|^2}\frac{dt}{\sqrt t}
=\frac1{\sqrt\pi}\int_{|\zeta|^2/L}^{L|\zeta|^2}e^{-u}\frac{du}{\sqrt u}\to1,\qquad \zeta\ne0.
$$
Plancherel and monotone convergence of the nonnegative integrand prove \eqref{e:riesz-integral} almost everywhere.

Fix $0<r<1/16$, let $R_0=I\times J=(-r,0)^n\times(-r,0)^m$, and write $L(z)=1+\log(1/z)$ for $0<z<r$.
Splitting at $s=a$ and using $s\leq a+s\leq2s$ for $s\geq a$, we obtain
$$
\sup_{a\geq0}\int_0^\infty V_{a+s}(z)\frac{ds}{\sqrt s}
\leq C+C\int_0^1e^{-z^2/(8s)}\frac{ds}{s}+C\int_1^\infty\frac{ds}{s^{3/2}}
\leq C L(z).
$$
Taking $a=t$ and then $a=0$ in \eqref{e:riesz-integral}, we have on $R_0$
$$
0\leq F(x,y)
\leq C L(r_y)\int_0^\infty V_t(r_x)\frac{dt}{\sqrt t}
\leq C L(r_x)L(r_y).
$$
Averaging $|F(x,y)|\leq C L(r_x)L(r_y)$ in $y$ and $x$, respectively, and using
$r^{-1}\int_0^r L(u)\,du=2+\log(1/r)$, we obtain
\begin{equation}\label{e:separation-j}
|F_J(x)|\leq\frac1{|J|}\int_J|F(x,y)|\,dy
\leq\frac{C L(r_x)}r\int_0^r L(u)\,du\leq C_r L(r_x)
\end{equation}
and
\begin{equation}\label{e:separation-i}
|F_I(y)|\leq\frac1{|I|}\int_I|F(x,y)|\,dx
\leq\frac{C L(r_y)}r\int_0^r L(u)\,du\leq C_r L(r_y).
\end{equation}
Since $F_{R_0}=(F_J)_I$, a further average of \eqref{e:separation-j} over $I$ gives
$$
|F_{R_0}|\leq\frac1{|I|}\int_I|F_J(x)|\,dx
\leq\frac{C_r}{r}\int_0^r L(u)\,du\leq C_r.
$$

For all sufficiently large integers $N$, let
$$
D_N=\big\{(x,y)\in R_0:e^{-4N}<r_x<2e^{-4N},\ e^{-N}<r_y<2e^{-N}\big\},
\qquad |D_N|=r^{n+m-2}e^{-5N}.
$$
Restrict \eqref{e:riesz-integral} to $e^{-6N}<t<e^{-4N}$ and $e^{-N}<s<1/128$. On $D_N$,
$$
r_x<2e^{-4N}<e^{-3N}<\sqrt t,\qquad
r_y<2e^{-N}<\sqrt s\leq\sqrt{t+s}\leq\sqrt{2s}<1/8.
$$
For $\tau\in\{t,t+s\}$ and $0<z\leq\sqrt\tau$, the interval $[\sqrt\tau,2\sqrt\tau]$ lies where $\chi=1$. Hence
$$
V_\tau(z)\geq\frac c\tau\int_{\sqrt\tau}^{2\sqrt\tau}v(u)\,du
\geq\frac{c}{\sqrt\tau\,\log(e+\log(1/\sqrt\tau))}
\geq\frac{c'}{\sqrt\tau\,\log(e+N)}.
$$
At these scales the transverse heat averages are bounded below on $R_0$ by integration over
the cubes of radii $\sqrt t$ and $\sqrt{t+s}$ where $\phi$ and $\psi$ equal one. Thus, almost everywhere on $D_N$,
\begin{equation}\label{e:separation-growth}
F(x,y)\geq\frac{c}{\log^2(e+N)}
\int_{e^{-6N}}^{e^{-4N}}\frac{dt}{t}\int_{e^{-N}}^{1/128}\frac{ds}{s}
=\frac{2cN(N-\log128)}{\log^2(e+N)}
\geq\frac{c'N^2}{\log^2(e+N)}.
\end{equation}
By \eqref{e:separation-j} and \eqref{e:separation-i}, on the same set,
\begin{align}
\min\big\{F(x,y),(Q_{R_0}F)(x,y)\big\}
&\geq F(x,y)-|F_J(x)|-|F_I(y)|-|F_{R_0}|\notag\\
&\geq\frac{cN^2}{\log^2(e+N)}-C_r(1+N)
\geq\frac{c'N^2}{\log^2(e+N)}.\label{e:separation-full}
\end{align}
For every $\lambda>0$ and $\gamma>1/2$, it follows that
$$
\int_{R_0}e^{\lambda|Q_{R_0}F|^\gamma}\,dz
\geq |D_N|\exp\bigg[\lambda\bigg(\frac{c'N^2}{\log^2(e+N)}\bigg)^\gamma\bigg]
=r^{n+m-2}\exp\bigg(\frac{\lambda(c')^\gamma N^{2\gamma}}{\log^{2\gamma}(e+N)}-5N\bigg)\to\infty,
$$
since $2\gamma>1$.
\end{proof}

\begin{remark}\label{r:bmo-moment-sharp}
A single function on $R_0$ attains the order $p^2$ in \eqref{e:rect-bmo-moments}.
Replace $v$ by $v_0=\chi\mathbf1_{(0,1)}$ and let
$$
F_0=\mathcal T\big[v_0(x_1)\phi(x')v_0(y_1)\psi(y')\big].
$$
By \eqref{e:riesz-bmo} and the $L^2$ bound for $\mathcal T$, we have $F_0\in\mathrm{BMO}_F\cap L^2$.
The bounds for the coordinate averages used only $0\leq v\leq\mathbf1_{(0,1)}$, so their contribution
on $D_N$ is still at most $C_r(1+N)$. At the scales in \eqref{e:separation-growth},
$v_0=1$ on $[\sqrt\tau,2\sqrt\tau]$, whence
$$
(K_\tau*v_0)'(-z)\geq\frac c\tau\int_{\sqrt\tau}^{2\sqrt\tau}du
=\frac c{\sqrt\tau},\qquad 0<z\leq\sqrt\tau.
$$
Thus for almost every $(x,y)\in D_N$ and all sufficiently large $N$,
$$
(Q_{R_0}F_0)(x,y)\geq c\int_{e^{-6N}}^{e^{-4N}}\frac{dt}{t}
\int_{e^{-N}}^{1/128}\frac{ds}{s}-C_r(1+N)
=2cN(N-\log128)-C_r(1+N)\geq c_rN^2.
$$
Taking $N=\lceil p\rceil$ and using \eqref{e:rect-bmo-moments}, we obtain for all sufficiently large $p$,
$$
c_r'p^2\leq c_rN^2\big(r^{-2}e^{-5N}\big)^{1/p}
\leq\bigg(\frac1{|R_0|}\int_{R_0}|Q_{R_0}F_0|^p\,dz\bigg)^{1/p}
\leq C_{n,m}p^2\|F_0\|_{\mathrm{BMO}_F}.
$$
The $o(p^2)$ bound in Proposition~\ref{p:rect-moments} therefore implies $F_0\notin\mathcal X_F$.
\end{remark}

\begin{proof}[Proof of Theorem~\ref{t:classes}]
\noindent\textup{(i)} Proposition~\ref{p:log-separation} and Proposition~\ref{p:exp-divergence} with $\gamma=1$ prove, respectively,
$$
\mathscr E_F\setminus(\mathcal X_F\cap L^2)\ne\emptyset,\qquad
(\mathcal X_F\cap L^2)\setminus\mathscr E_F\ne\emptyset.
$$
Since $\mathscr E_F\subseteq\mathrm{BMO}_F\cap L^2$, the first example also proves
$\mathcal X_F\cap L^2\subsetneq\mathrm{BMO}_F\cap L^2$.

By Theorem~\ref{t:approximation} and \eqref{e:flag-exponential} with $k=2$,
$\mathscr S_F\subseteq(\mathcal X_F\cap L^2)\cap\mathscr E_F$.
The function $f=\mathbf1_{[0,1]^{n+m}}$ belongs to both classes, since $f$ is bounded and $|Q_Rf|\leq4$ on every $R$.
By Proposition~\ref{p:flag-area-sharp}, $S_F(f)\notin L^\infty$, so both inclusions are strict.

\smallskip
\noindent\textup{(ii)} Proposition~\ref{p:rect-moments} proves \eqref{e:closure-half} for every $\lambda>0$.
Proposition~\ref{p:exp-divergence} proves divergence on one fixed flag rectangle for every $\gamma>1/2$ and $\lambda>0$.
\end{proof}


\section{Weighted estimates, Riesz commutators and Poisson maximal functions}\label{s:consequences}

We first estimate the weighted moments of $f-f_R$ on a product of cubes $R$.
For two-step flags, bounds on the coordinate areas also control exponential weights and repeated Riesz commutators.
For finite flags, we compare the radial and non-tangential Poisson maximal functions with the flag area function.

\subsection{Weighted moments of the local oscillation}\label{s:change-measure}

\begin{corollary}\label{c:weighted-local}
Let $R=\prod_{j=1}^k E_j$, where $k\geq2$ and $E_j\subset X_j$ are axis-parallel cubes.
Let $w\geq0$ be locally integrable, with $w(\Omega)=\int_\Omega w\,dx$ and $0<w(R)<\infty$.
Suppose that for some $K\geq1$ and $0<\delta\leq1$,
$$
\frac{w(\Omega)}{w(R)}\leq K\bigg(\frac{|\Omega|}{|R|}\bigg)^\delta
\qquad\text{for every measurable }\Omega\subseteq R.
$$
If $f\in L^u(\mathbb R^N)$ for some $1\leq u<\infty$ and $\mathcal E_{\ell,R}(f)\in L^\infty(R)$ for $1\leq\ell\leq k$, then
$$
\bigg(\frac1{w(R)}\int_R|f-f_R|^p w\,dx\bigg)^{1/p}
\leq C\sum_{\ell=1}^k\bigg(\frac{p+\log K}{\delta}\bigg)^{\ell/2}
\|\mathcal E_{\ell,R}(f)\|_{L^\infty(R)},\qquad 2\leq p<\infty.
$$
The constant depends only on $k,d_1,\ldots,d_k$. For fixed $p$ and $\delta=1$, the power $k/2$ of $\log K$
in the term $\ell=k$ is optimal as $K\to\infty$.
\end{corollary}

\begin{proof}
Write $d\mu=w\,dx/w(R)$ on $R$ and $a_\ell=\|\mathcal E_{\ell,R}(f)\|_{L^\infty(R)}$.
For $a_\ell>0$, let $Y=c\delta(|F_R^{(\ell)}|/a_\ell)^{2/\ell}$, with $c$ as in \eqref{e:coord-order-exp}.
By the weight hypothesis and \eqref{e:coord-order-exp}, for $s\geq0$,
$$
\mu\{Y>\log(2K)+s\}
\leq K\bigg(e^{-(\log(2K)+s)/\delta}\frac1{|R|}\int_R e^{Y/\delta}\,dx\bigg)^\delta
\leq 2^{\delta-1}e^{-s}\leq e^{-s}.
$$
For $r=\ell p/2\geq1$, it follows that
\begin{align*}
c\delta\bigg(\frac{\|F_R^{(\ell)}\|_{L^p(\mu)}}{a_\ell}\bigg)^{2/\ell}
&=\|Y\|_{L^r(\mu)}
\leq\log(2K)+\|(Y-\log(2K))_+\|_{L^r(\mu)}\\
&=\log(2K)+\bigg(r\int_0^\infty s^{r-1}\mu\{Y>\log(2K)+s\}\,ds\bigg)^{1/r}\\
&\leq\log(2K)+\Gamma(1+r)^{1/r}\leq C_k(p+\log K).
\end{align*}
If $a_\ell=0$, the second moment estimate for $F_R^{(\ell)}$ implies $F_R^{(\ell)}=0$ almost everywhere on $R$.
Thus, by \eqref{e:coord-expansion},
$$
\|f-f_R\|_{L^p(\mu)}
\leq\sum_{\ell=1}^k\|F_R^{(\ell)}\|_{L^p(\mu)}
\leq C\sum_{\ell=1}^k\bigg(\frac{p+\log K}{\delta}\bigg)^{\ell/2}a_\ell.
$$

For sharpness, take $f_M,R_M$ from Proposition~\ref{p:flag-reverse-sharp} and let
$w_M=\mathbf1_{\Omega_M}$, where $\Omega_M$ is the product of balls in \eqref{e:flag-large-values}. Then
$$
\frac{w_M(\Omega)}{w_M(R_M)}
=\frac{|\Omega\cap\Omega_M|}{|\Omega_M|}
\leq K_M\frac{|\Omega|}{|R_M|},\qquad
K_M=\frac{|R_M|}{|\Omega_M|}\simeq2^{NM},
$$
with equality for $\Omega=\Omega_M$. All coordinate averages of $f_M$ vanish, so $(f_M)_{R_M}=0$ and
$\mathcal E_{\ell,R_M}(f_M)=0$ for $\ell<k$. By \eqref{e:flag-large-values} and \eqref{e:flag-sharp-sup},
$$
\frac{\big(w_M(R_M)^{-1}\int_{R_M}|f_M|^p w_M\,dx\big)^{1/p}}
     {\|S_F^{(k)}(f_M)\|_{L^\infty(R_M)}}
\geq\frac{\big(|\Omega_M|^{-1}\int_{\Omega_M}|f_M|^p\,dx\big)^{1/p}}
          {\|S_F^{(k)}(f_M)\|_\infty}
\geq\frac{(M/2)^k}{CM^{k/2}}
\geq c(\log K_M)^{k/2}.
$$
This proves sharpness for $\delta=1$ and each fixed $p\geq2$.
\end{proof}


\subsection{Exponential weights from the coordinate areas}\label{s:coordinate-weights}

Return to the two-step flag, with $d=n+m$. For a product $R=I\times J$ of axis-parallel cubes, write
$\langle g\rangle_R=|R|^{-1}\int_R g$. The flag rectangles and their weight characteristic are
$$
\mathcal R_F=\{I\times J:\ell(I)\leq\ell(J)\},\qquad
[w]_{A_{F,2}}=\sup_{R\in\mathcal R_F}\langle w\rangle_R\langle w^{-1}\rangle_R,
$$
where $w>0$ almost everywhere and $w,w^{-1}\in L^1_{\mathrm{loc}}$. We write $w\in A_{F,2}$ if $[w]_{A_{F,2}}<\infty$.

For real $f$ with $\sup_{R\in\mathcal R_F}\langle|f-f_R|\rangle_R<\infty$,
Duong, Li, Ou, Pipher and Wick~\cite[Proposition~4.2]{DLOPW} proved that $e^{\lambda f}\in A_{F,2}$ for small $|\lambda|$.
The function in Proposition~\ref{p:log-separation} satisfies this little flag BMO condition and has positive
$\mathrm{BMO}_F$ distance from $L^\infty$. A bound on the coordinate areas permits every real $\lambda$.
For flag weights in Hardy space estimates, see Han, Lin and X.~Wu~\cite{HLW}.

For $f\in L^2(\mathbb R^d)$, define the slice areas almost everywhere by
$$
\mathcal S_{P,x}f(x,y)=\mathcal S_P^{(n)}(f(\cdot,y))(x),\qquad
\mathcal S_{P,y}f(x,y)=\mathcal S_P^{(m)}(f(x,\cdot))(y).
$$
On $R=I\times J$ the coordinate area satisfies
$$
\mathcal E_{P,R}(f)(x,y)^2=|\mathcal S_P^{(n)}f_J(x)|^2+|\mathcal S_P^{(m)}f_I(y)|^2.
$$
Let $\mathcal D$ be the standard dyadic system in $\mathbb R^d$, and write
$A(f)=\sup_{Q\in\mathcal D}\|\mathcal E_{P,Q}(f)\|_{L^\infty(Q)}$.
Here and in Subsection~\ref{s:commutators}, the constants depend only on $n,m$.

\begin{proposition}\label{p:coordinate-weights}
For real $f\in L^2(\mathbb R^d)$,
$$
\max\{\|\mathcal S_{P,x}f\|_\infty,\|\mathcal S_{P,y}f\|_\infty\}
\leq A(f)\leq\big(\|\mathcal S_{P,x}f\|_\infty^2+\|\mathcal S_{P,y}f\|_\infty^2\big)^{1/2}.
$$
If $0<A(f)<\infty$, then
$$
\sup_R\big\langle\exp\big(c|f-f_R|^2/A(f)^2\big)\big\rangle_R\leq2.
$$
For $A(f)<\infty$ and every $\lambda\in\mathbb R$,
\begin{equation}\label{e:exp-weights}
[e^{\lambda f}]_{A_{F,2}}
\leq\sup_R\langle e^{\lambda f}\rangle_R\langle e^{-\lambda f}\rangle_R
\leq C\exp\big(C\lambda^2 A(f)^2\big).
\end{equation}
Both suprema over $R$ range over all products of axis-parallel cubes. On one fixed flag rectangle $R_*$,
for each sufficiently large $|\lambda|$ there is a real Schwartz function $f$ with $A(f)\leq1$ and
$$
\langle e^{\lambda f}\rangle_{R_*}\langle e^{-\lambda f}\rangle_{R_*}\geq e^{c\lambda^2}.
$$
\end{proposition}

\begin{proof}
By Minkowski's inequality, on every $I\times J$,
\begin{align*}
\mathcal E_{P,I\times J}(f)(x,y)^2
&\leq\bigg(\frac1{|J|}\int_J\mathcal S_{P,x}f(x,v)\,dv\bigg)^2
    +\bigg(\frac1{|I|}\int_I\mathcal S_{P,y}f(u,y)\,du\bigg)^2\\
&\leq\|\mathcal S_{P,x}f\|_\infty^2+\|\mathcal S_{P,y}f\|_\infty^2.
\end{align*}
For the reverse bound assume $A(f)<\infty$. For each dyadic $J$, the cubes $I$ of side $\ell(J)$
tile $\mathbb R^n$ and $I\times J\in\mathcal D$. Using these cubes and then interchanging $I,J$, we obtain
\begin{equation}\label{e:coordinate-average-areas}
\max\{\|\mathcal S_P^{(n)}f_J\|_{L^\infty(\mathbb R^n)},\|\mathcal S_P^{(m)}f_I\|_{L^\infty(\mathbb R^m)}\}\leq A(f)
\qquad(I,J\text{ dyadic}).
\end{equation}
Let $(\mathbb E_j^y f)(x,y)=\sum_{\ell(J)=2^{-j}}\mathbf1_J(y) f_J(x)$ be the dyadic average in $y$.
By \eqref{e:coordinate-average-areas} and Plancherel,
$$
\|(\mathcal S_{P,x}f-A(f))_+\|_2^2
\leq\|\mathcal S_{P,x}(f-\mathbb E_j^y f)\|_2^2
=\frac{3|B_{\mathbb R^n}(0,1)|}{8}\|f-\mathbb E_j^y f\|_2^2\to0.
$$
Dyadic averaging in $x$ proves $\|\mathcal S_{P,y}f\|_\infty\leq A(f)$ as well.

On every $R=I\times J$, apply \eqref{e:one-local} to $f(x,\cdot)$ on $J$ and $f_J$ on $I$. For $2\leq p<\infty$,
\begin{align*}
\langle|f-f_R|^p\rangle_R^{1/p}
&\leq\langle|f-f_J|^p\rangle_R^{1/p}+\langle|f_J-f_R|^p\rangle_I^{1/p}\\
&\leq C\sqrt p\big(\langle|\mathcal S_{P,y}f|^p\rangle_R^{1/p}
      +\langle|\mathcal S_P^{(n)}f_J|^p\rangle_I^{1/p}\big)\\
&\leq C\sqrt p\big(\|\mathcal S_{P,y}f\|_\infty+\|\mathcal S_{P,x}f\|_\infty\big)
\leq C\sqrt p\,A(f).
\end{align*}
If $A(f)=0$, then $f=f_R$ on every $R$, hence $f=0$ almost everywhere. Otherwise, summing the moments as in
\eqref{e:coord-exponential}, we obtain $\langle\exp(c|f-f_R|^2/A(f)^2)\rangle_R\leq2$. Consequently,
\begin{align*}
\langle e^{\lambda f}\rangle_R\langle e^{-\lambda f}\rangle_R
&=\langle e^{\lambda(f-f_R)}\rangle_R\langle e^{-\lambda(f-f_R)}\rangle_R\\
&\leq e^{\lambda^2 A(f)^2/(2c)}
\big\langle e^{c|f-f_R|^2/A(f)^2}\big\rangle_R^2
\leq4e^{\lambda^2 A(f)^2/(2c)},
\end{align*}
which proves \eqref{e:exp-weights}.

Take $R_*=E\times J_0$, $h_L$ and $H_L$ from Proposition~\ref{p:two-orders} and its proof.
Their definitions give $\|h_L\|_\infty+\|H_L\|_\infty\leq CL$. By Lemma~\ref{l:one-lacunary},
$\|\mathcal S_P^{(n)}h_L\|_\infty\leq C\sqrt L$, and the same calculation at scales $\rho_j^{-1}=\rho^{-1}2^{1-j}$ gives
$$
\|\mathcal S_P^{(m)}H_L\|_\infty^2
\leq C\sum_{i,j=1}^L\frac{\rho_i^{-2}\rho_j^{-2}}{(\rho_i^{-1}+\rho_j^{-1})^4}
\leq C\sum_{i,j=1}^L2^{-2|i-j|}\leq CL.
$$
For $g_L(x,y)=h_L(x)H_L(y)/(C_0L^{3/2})$ and $C_0$ sufficiently large,
$$
A(g_L)\leq\frac{\|H_L\|_\infty\|\mathcal S_P^{(n)}h_L\|_\infty
+\|h_L\|_\infty\|\mathcal S_P^{(m)}H_L\|_\infty}{C_0L^{3/2}}\leq1.
$$
Since $(h_L)_E=0$, we have $(g_L)_{R_*}=0$. By Lemma~\ref{l:one-lacunary} and \eqref{e:coord-flat},
$g_L\geq c_1\sqrt L$ on a subset of $R_*$ of relative measure at least $c2^{-nL}$. Hence by Jensen's inequality,
$$
\langle e^{\lambda g_L}\rangle_{R_*}\langle e^{-\lambda g_L}\rangle_{R_*}
\geq\langle e^{|\lambda|g_L}\rangle_{R_*}e^{-|\lambda|(g_L)_{R_*}}
\geq c\exp\big(c_1|\lambda|\sqrt L-nL\log2\big).
$$
Take $L=\lfloor a\lambda^2\rfloor$ with $0<a<(c_1/(n\log2))^2$. For sufficiently large $|\lambda|$, the product is at least $e^{c\lambda^2}$.
\end{proof}

For $F,R_0,D_N$ in Proposition~\ref{p:exp-divergence}, \eqref{e:separation-full} implies, for every $\lambda>0$,
$$
\int_{R_0}e^{\lambda F}\,dz\geq\int_{D_N}e^{\lambda F}\,dz
\geq r^{n+m-2}\exp\bigg(\frac{c\lambda N^2}{\log^2(e+N)}-5N\bigg)\to\infty.
$$
Applying this bound to the reciprocal weight when $\lambda<0$, we obtain $e^{\lambda F}\notin A_{F,2}$ for every real $\lambda\ne0$.

\subsection{Iterated commutators of flag Riesz transforms}\label{s:commutators}

Fix $T=R_i^{(1)}R_j^{(2)}$, $1\leq i\leq d$, $1\leq j\leq m$, as in Subsection~\ref{s:hardy-estimates}.
For bounded $f$, write
$$
[f,S]g=fSg-S(fg),\qquad \operatorname{ad}_f^0T=T,\qquad
\operatorname{ad}_f^{r+1}T=[f,\operatorname{ad}_f^rT],\quad r\geq0.
$$
Duong, Li, Ou, Pipher and Wick~\cite[Theorem~4.3 and Corollary~4.6]{DLOPW} used analytic conjugation for little flag BMO symbols.
Here \eqref{e:exp-weights} controls the conjugation on every circle in the complex plane.

\begin{corollary}\label{c:commutators}
Let $f\in L^2(\mathbb R^d)$ be real with $A(f)<\infty$. Then
\begin{equation}\label{e:commutators}
\|\operatorname{ad}_f^r(R_i^{(1)}R_j^{(2)})\|_{L^2(\mathbb R^d)\to L^2(\mathbb R^d)}
\leq(C\sqrt r\,A(f))^r,\qquad r=1,2,\ldots.
\end{equation}
For unbounded $f$, this is the unique bounded operator with
$$
\langle\operatorname{ad}_f^rTg,h\rangle
=\sum_{\ell=0}^r(-1)^\ell\binom r\ell\langle T(f^\ell g),f^{r-\ell}h\rangle,
\qquad g,h\in C_c^\infty(\mathbb R^d).
$$
\end{corollary}

\begin{proof}
For $w\in A_{F,2}$, every cube in $\mathbb R^d$ is a flag rectangle, so $[w]_{A_2(\mathbb R^d)}\leq[w]_{A_{F,2}}$.
For a cube $J\subset\mathbb R^m$ and the cube $I_\varepsilon(x)$ of side $\varepsilon$ centred at $x$, differentiation in $x$ gives
$$
\langle w(x,\cdot)\rangle_J\langle w(x,\cdot)^{-1}\rangle_J
=\lim_{\varepsilon\to0+}\langle w\rangle_{I_\varepsilon(x)\times J}\langle w^{-1}\rangle_{I_\varepsilon(x)\times J}
\leq[w]_{A_{F,2}}.
$$
Taking first cubes $J$ with rational centres and side lengths, then approximating by such cubes, we obtain
$[w(x,\cdot)]_{A_2(\mathbb R^m)}\leq[w]_{A_{F,2}}$ for almost every $x$.
By Petermichl's estimate~\cite[Theorem~2.1]{PetRiesz} on $\mathbb R^d$ and on the vertical slices,
\begin{equation}\label{e:weighted-riesz}
\|R_i^{(1)}R_j^{(2)}\|_{L^2(w)\to L^2(w)}
\leq\|R_i^{(1)}\|_{L^2(w)\to L^2(w)}
\operatorname*{ess\,sup}_{x}\|R_j^{(2)}\|_{L^2(w(x,\cdot))\to L^2(w(x,\cdot))}
\leq C[w]_{A_{F,2}}^2.
\end{equation}
Density of $C_c^\infty$ in $L^2((1+w)\,dx\,dy)$ identifies each weighted transform with its unweighted extension on $L^2\cap L^2(w)$.

Put $A=A(f)$ and $f_L=\max\{-L,\min\{f,L\}\}$. Since truncation is $1$-Lipschitz, \eqref{e:exp-weights} implies
\begin{align*}
[e^{\lambda f_L}]_{A_{F,2}}
&=\sup_{R\in\mathcal R_F}\frac1{|R|^2}\int_R\int_R\cosh\big(\lambda(f_L(u)-f_L(v))\big)\,du\,dv\\
&\leq\sup_{R\in\mathcal R_F}\frac1{|R|^2}\int_R\int_R\cosh\big(\lambda(f(u)-f(v))\big)\,du\,dv
=[e^{\lambda f}]_{A_{F,2}}\leq C e^{C\lambda^2 A^2}.
\end{align*}
For $F_L(z)=e^{zf_L}Te^{-zf_L}$ and $w=e^{2\operatorname{Re}(z)f_L}$, \eqref{e:weighted-riesz} implies
$$
\|F_L(z)g\|_2=\|T(e^{-zf_L}g)\|_{L^2(w)}
\leq C[w]_{A_{F,2}}^2\|e^{-zf_L}g\|_{L^2(w)}
\leq C e^{C|z|^2 A^2}\|g\|_2.
$$
The family $F_L$ is entire. For $A>0$, Cauchy's formula on $|z|=\sqrt r/A$ gives
$$
\|\operatorname{ad}_{f_L}^rT\|_{2\to2}
=\|F_L^{(r)}(0)\|_{2\to2}
\leq r!\bigg(\frac{A}{\sqrt r}\bigg)^r\sup_{|z|=\sqrt r/A}\|F_L(z)\|_{2\to2}
\leq C r!\bigg(\frac{A}{\sqrt r}\bigg)^r e^{Cr}\leq(C\sqrt r\,A)^r.
$$
For $A=0$ the same formula on circles of arbitrarily large radius gives $F_L^{(r)}(0)=0$.

By \eqref{e:exp-weights}, $e^{|f|}\leq e^f+e^{-f}\in L^1_{\mathrm{loc}}$.
Thus $f_L^\ell g\to f^\ell g$ in $L^2$ for $g\in C_c^\infty$ and $0\leq\ell\leq r$, by dominated convergence.
For $g,h\in C_c^\infty$,
\begin{align*}
\bigg|\sum_{\ell=0}^r(-1)^\ell\binom r\ell\langle T(f^\ell g),f^{r-\ell}h\rangle\bigg|
&=\lim_{L\to\infty}\bigg|\sum_{\ell=0}^r(-1)^\ell\binom r\ell\langle T(f_L^\ell g),f_L^{r-\ell}h\rangle\bigg|\\
&=\lim_{L\to\infty}|\langle\operatorname{ad}_{f_L}^rTg,h\rangle|
\leq(C\sqrt r\,A)^r\|g\|_2\|h\|_2.
\end{align*}
The form extends uniquely by density, proving \eqref{e:commutators}.
\end{proof}

\subsection{Poisson maximal functions and the flag area function}\label{s:poisson-maximal}

For the Poisson extension $U_f$ of a finite flag, the radial and non-tangential maximal functions are
$$
U_f^+(x)=\sup_{\mathbf t>0}|U_f(x,\mathbf t)|,\qquad
U_f^*(x)=\sup_{\substack{y\in\mathbb R^N,\ \mathbf t>0\\ |x_i-y_i|\leq\sum_{j=1}^it_j,\ 1\leq i\leq k}}|U_f(y,\mathbf t)|.
$$

\begin{corollary}\label{c:poisson-maximal}
Let $f\in L^u(\mathbb R^N)$ for some $1\leq u<\infty$, and let $2\leq p<\infty$.
If either $f$ or $S_F^{(k)}(f)$ belongs to $L^p(\mathbb R^N)$, then
\begin{equation}\label{e:poisson-maximal}
\|f\|_p\leq\|U_f^+\|_p\leq\|U_f^*\|_p
\leq C\|f\|_p\leq Cp^{k/2}\|S_F^{(k)}(f)\|_p.
\end{equation}
Moreover,
$$
cp^{k/2}
\leq\sup_{0\ne f\in\mathcal S(\mathbb R^N)}\frac{\|U_f^+\|_p}{\|S_F^{(k)}(f)\|_p}
\leq\sup_{0\ne f\in\mathcal S(\mathbb R^N)}\frac{\|U_f^*\|_p}{\|S_F^{(k)}(f)\|_p}
\leq Cp^{k/2}.
$$
The constants depend only on $k,d_1,\ldots,d_k$.
\end{corollary}

\begin{proof}
By Theorem~\ref{t:concentration}\textup{(ii)}, $f\in L^p$.
At positive scales the translated Poisson kernels depend continuously on the parameters in $L^{p'}$,
so both suprema may be taken over countable dense sets.
For the normalized cube densities $q_{r,D}$, summation over dyadic annuli gives
$$
P_t^{(D)}\leq C_D\sum_{\nu\geq0}2^{-\nu}q_{2^\nu t,D}.
$$
For $r_j=2^{\nu_j}t_j$ and $r_i^*=\max_{j\leq i}r_j$, the convolution of the cube densities factors as
$$
q_{r_1,D_1}*_{V_2}q_{r_2,D_2}*\cdots*_{V_k}q_{r_k,D_k}(u)
=\prod_{i=1}^k(q_{r_1,d_i}*\cdots*q_{r_i,d_i})(u_i).
$$
The $i$th factor has supremum at most $(2r_i^*)^{-d_i}$ and support in $|u_i|_\infty\leq ir_i^*$.
If $|y_i-x_i|\leq\sum_{j\leq i}t_j\leq ir_i^*$, it follows that
$$
(q_{r_1,d_i}*\cdots*q_{r_i,d_i})(y_i-u_i)\leq (2i)^{d_i}q_{2ir_i^*,d_i}(x_i-u_i)
\quad\text{for almost every }u_i.
$$
Writing $M_{X_i}$ for the Hardy--Littlewood maximal operator on $X_i$, we obtain
\begin{align*}
|U_f(y,\mathbf t)|
&\leq C\sum_{\boldsymbol\nu\in\mathbb N_0^k}2^{-\sum_j\nu_j}
\int_{\mathbb R^N}|f(u)|\prod_{i=1}^k q_{2ir_i^*,d_i}(x_i-u_i)\,du\\
&\leq C\sum_{\boldsymbol\nu\in\mathbb N_0^k}2^{-\sum_j\nu_j}
(M_{X_1}\cdots M_{X_k}f)(x)\leq C(M_{X_1}\cdots M_{X_k}f)(x).
\end{align*}
Since $U_f(\cdot,\varepsilon,\ldots,\varepsilon)\to f$ in $L^p$ as $\varepsilon\to0+$,
we have $|f|\leq U_f^+$ almost everywhere. The maximal bounds on the $X_i$ are uniform for $p\geq2$, so
$$
\|f\|_p\leq\|U_f^+\|_p\leq\|U_f^*\|_p\leq C\|M_{X_1}\cdots M_{X_k}f\|_p
\leq C\|f\|_p\leq Cp^{k/2}\|S_F^{(k)}(f)\|_p
$$
by Theorem~\ref{t:flag-reverse}. Testing $|f|\leq U_f^+$ on the functions of
Proposition~\ref{p:flag-reverse-sharp} proves sharpness.
\end{proof}


\section{Other spectral functions and the heat area}\label{s:multipliers-heat}


For the classical one-parameter Poisson and heat multipliers $re^{-r}$ and $r^2e^{-r^2}$,
the sharp reverse factor is $\sqrt p$ at every fixed positive aperture. For other smooth radial
mean-zero kernels the factor can be $p$. Even for one fixed kernel the order can change from
$p$ at small apertures to $\sqrt p$ at large ones, as Li, Pipher and Wu show
in~\cite[Proposition~3.7 and Theorems~B and~C]{Zygmund}.

We keep the reconstruction multiplier $a_\sigma$ normalized in \eqref{e:dual-normalization} and
replace the Poisson multipliers in the coefficients by functions in $L^2((0,\infty),dr/r)$
with nonzero pairing with $a_\sigma$. We estimate their coefficient area at the aperture
$\lambda_0$ of \eqref{e:rectangle-cone}. For the heat multipliers, a change of scale
bounds this area by the flag heat area at aperture one.

\subsection{Changing the spectral functions in the reconstruction}\label{s:spectral-reconstruction}

For $b\in L^2((0,\infty),dr/r)$, write
$$
m_\sigma(b)=\int_0^\infty a_\sigma(r)b(r)\,\frac{dr}{r}.
$$
Choose Borel functions $b_1,b_2$ in this space with $m_\sigma(b_1)m_\sigma(b_2)\ne0$.
For $\mathbf b=(b_1,b_2)$ and $f\in L^2(\mathbb R^{n+m})$, define $H_f^{\mathbf b}$ by its Fourier transform in $(u,v)$,
$$
\widehat H_f^{\mathbf b}(\xi,\eta,r_1,r_2)
=\frac{b_1(r_1|(\xi,\eta)|)b_2(r_2|\eta|)}{m_\sigma(b_1)m_\sigma(b_2)}\widehat f(\xi,\eta).
$$
The right-hand side is jointly square integrable for the scale measure $dr_1\,dr_2/(r_1r_2)$ and Lebesgue measure in $(\xi,\eta)$,
as the norm calculation below shows. Partial Fourier inversion therefore defines $H_f^{\mathbf b}$ in $L^2_{\mathrm{par}}$.
We measure these coefficients by the area $\mathcal A_{\lambda_0,F}$ in \eqref{e:aux-area}.

\begin{proposition}\label{p:spectral-multipliers}
Let $f\in L^2(\mathbb R^{n+m})$. Then $H_f^{\mathbf b}\in L^2_{\mathrm{par}}$ and $\mathcal R_\sigma H_f^{\mathbf b}=f$ in $L^2(\mathbb R^{n+m})$.
For $2\leq p<\infty$,
$$
\|f\|_p\leq Cp\|\mathcal A_{\lambda_0,F}(H_f^{\mathbf b})\|_p,
$$
and for every product $R=I\times J$ of axis-parallel cubes,
\begin{equation}\label{e:spectral-local}
\left \|Q_Rf  \right \|_{L^p(R)}\leq Cp  \left  \|\mathcal A_{\lambda_0,F}(H_f^{\mathbf b})\right \|_{L^p(R)}.
\end{equation}
For the deviation from the mean,
$$
\|f-f_R\|_{L^p(R)}
\leq Cp\|\mathcal A_{\lambda_0,F}(H_f^{\mathbf b})\|_{L^p(R)}
+C\sqrt p\,\|\mathcal E_{P,R}(f)\|_{L^p(R)}.
$$
The exponential estimate \eqref{e:two-step-exp} and the tail bound \eqref{e:two-tails} hold with
$\mathcal A_{\lambda_0,F}(H_f^{\mathbf b})$ in place of $S_F(f)$. The coordinate term in \eqref{e:two-tails}
uses the Poisson areas of $f_J$ and $f_I$. With the normalization in $H_f^{\mathbf b}$,
the constants are independent of $\mathbf b,f,p,R$. For the unnormalized coefficients,
the factor multiplying their area is $1/|m_\sigma(b_1)m_\sigma(b_2)|$ times the stated constant.
\end{proposition}

\begin{proof}
For $\eta\ne0$, rescale $r_1$ by $|(\xi,\eta)|$ and $r_2$ by $|\eta|$ in the two scale integrals.
By Plancherel's identity and Tonelli's theorem,
$$
\|H_f^{\mathbf b}\|_{L^2_{\mathrm{par}}}^2
=\prod_{i=1}^2\frac{\|b_i\|_{L^2(dr/r)}^2}{|m_\sigma(b_i)|^2}\,\|f\|_{L^2(\mathbb R^{n+m})}^2<\infty.
$$
Integrating first in the cone centres, we also have
$$
\left \|\mathcal A_{\lambda_0,F}(H_f^{\mathbf b})\right \|_{L^2(\mathbb R^{n+m})}^2
=|B_{\mathbb R^n}(0,\lambda_0)|\,|B_{\mathbb R^m}(0,\lambda_0)|
\|H_f^{\mathbf b}\|_{L^2_{\mathrm{par}}}^2<\infty.
$$
For each $\rho>0$ and $i=1,2$, we use $v=r\rho$ in the truncated scale integral,
$$
\int_{1/L}^L\frac{a_\sigma(r\rho)b_i(r\rho)}{m_\sigma(b_i)}\,\frac{dr}{r}
=\frac1{m_\sigma(b_i)}\int_{\rho/L}^{L\rho}a_\sigma(v)b_i(v)\,\frac{dv}{v}\to  1.
$$
By Cauchy--Schwarz the absolute integrals are bounded by
$\|a_\sigma\|_{L^2(dr/r)}\|b_i\|_{L^2(dr/r)}/|m_\sigma(b_i)|$.
With $\rho=|(\xi,\eta)|$ and $\rho=|\eta|$ the product tends to one almost everywhere.
By Plancherel's identity and dominated convergence, as in Proposition~\ref{p:reconstruction},
we obtain $\mathcal R_\sigma H_f^{\mathbf b}=f$ in $L^2(\mathbb R^{n+m})$.

For any finite selection $\mathcal E$, the atom bounds in Proposition~\ref{p:two-atoms} and
\eqref{e:cube-sum} use only the $L^2(d\nu)$ norms of the coefficients on the parameter regions.
Proposition~\ref{p:coeff-control} bounds their square sum by
$C\mathcal A_{\lambda_0,F}(m_{\mathcal E}H_f^{\mathbf b})^2$.
Applying Propositions~\ref{p:local-product} and~\ref{p:local-vertical} to the terms with cancellation in
both variables and the two sums with vertical cancellation, respectively, we obtain
$$
\|Q_R\mathcal R_\sigma(m_{\mathcal E}H_f^{\mathbf b})\|_{L^p(R)}
\leq Cp\|\mathcal A_{\lambda_0,F}(H_f^{\mathbf b})\|_{L^p(R)}.
$$
Taking $H=H_f^{\mathbf b}$ in \eqref{e:finite-recon-limit}, we obtain
$\mathcal R_\sigma(m_{\mathcal E_N}H_f^{\mathbf b})\to f$ in $L^2(\mathbb R^{n+m})$.
After centring on $R$, a subsequence converges almost everywhere on $R$ to $Q_Rf$, as in Theorem~\ref{t:two-local}.
By Fatou's lemma, we obtain \eqref{e:spectral-local}. The global bound follows by exhausting $\mathbb R^{n+m}$ with rectangles as in Theorem~\ref{t:flag-reverse}.
Inserting \eqref{e:spectral-local} and \eqref{e:coord-mean-lp} into \eqref{e:coord-splitting},
we obtain the bound for $f-f_R$. The exponential estimates follow from the same moment
calculation as \eqref{e:exponential-series} and the level-set decomposition in \eqref{e:two-tails}.
\end{proof}

The pairing condition in Proposition~\ref{p:spectral-multipliers} depends on the support of the dual kernel.
To see this, choose a nonzero real even $h\in C_c^\infty(\{x\in\mathbb R:2<|x|<3\})$
and let $b(r)=r\widehat h(r)$ for $r>0$. Then $0\ne b\in L^2(dr/r)$. For every
$\theta_0$ as in Subsection~\ref{s:compact-kernels}, we have by Plancherel's identity
$$
2c_\sigma m_\sigma(b)
=\int_{\mathbb R}r^2\widehat\theta_0(r)\widehat h(r)\,dr
=2\pi\int_{\mathbb R}\theta_0(x)(-h''(x))\,dx=0,
$$
since $\operatorname{supp}\theta_0\subset(-1,1)$.

For $h_\Lambda(r)=r^\alpha e^{-(\Lambda r)^\beta}$, $\alpha,\beta,\Lambda>0$, the scale integral is explicit,
$$
\|h_\Lambda\|_{L^2(dr/r)}^2
=\frac{\Gamma(2\alpha/\beta)}{\beta\,2^{2\alpha/\beta}\Lambda^{2\alpha}}.
$$
Substituting $u=\Lambda r$ in $m_\sigma(h_\Lambda)$, we obtain
\begin{equation}\label{e:spectral-normalization}
c_\sigma\Lambda^{\alpha+2}m_\sigma(h_\Lambda)
=\int_0^\infty u^{\alpha+1}\widehat\theta_0(u/\Lambda)e^{-u^\beta}\,du 
\to\frac1\beta\Gamma\Big(\frac{\alpha+2}{\beta}\Big)>0,\qquad (\Lambda\to\infty).
\end{equation}
The limit follows by dominated convergence, with majorant
$\|\widehat\theta_0\|_\infty u^{\alpha+1}e^{-u^\beta}$.
It follows that $m_\sigma(h_\Lambda)>0$ for all sufficiently large $\Lambda$. For the original Poisson
multipliers $b_i(r)=re^{-\sigma r}$, $m_\sigma(b_i)=1$ by construction. Then $H_f^{\mathbf b}=H_f^\sigma$
and \eqref{e:coeff-area} bounds $\mathcal A_{\lambda_0,F}(H_f^{\mathbf b})$ by $CS_F(f)$.

\subsection{Heat areas and coordinate averages}\label{s:heat-area}

For $1\leq u<\infty$ and $f\in L^u(\mathbb R^{n+m})$, consider the heat extension
$$
W_f(\cdot,t,s)=e^{t^2\Delta_{x,y}}e^{s^2\Delta_y}f .
$$
The first semigroup acts on $\mathbb R^{n+m}$ and the second in the $y$ variable.
With the cone weight $\chi_{t,s}$ from \eqref{e:poisson-area-weight}, define the flag heat area by
$$
S_{F,\mathrm h}(f)(x,y)^2
=\int_0^\infty\int_0^\infty\int_{\mathbb R^{n+m}}
\chi_{t,s}(x-u,y-v)|t^2(-\Delta_{x,y})s^2(-\Delta_y)W_f(u,v,t,s)|^2
\frac{du\,dv}{t^{n+m}s^m}\frac{dt}{t}\frac{ds}{s}.
$$
For $g\in L^u(\mathbb R^q)$ in the same range of $u$, define the one-parameter heat area by
$$
\mathcal S_{\mathrm h}^{(q)}g(x)^2
=\int_0^\infty\int_{|x-y|<t}
|t^2(-\Delta)e^{t^2\Delta}g(y)|^2\frac{dy}{t^q}\frac{dt}{t}.
$$
For $R=I\times J$, define the area of its coordinate averages by
$$
\mathcal E_{\mathrm h,R}(f)(x,y)^2
=\mathbf1_J(y)|\mathcal S_{\mathrm h}^{(n)}f_J(x)|^2
+\mathbf1_I(x)|\mathcal S_{\mathrm h}^{(m)}f_I(y)|^2.
$$

\begin{theorem}\label{t:heat-area}
Let $f\in L^u(\mathbb R^{n+m})$ for some $1\leq u<\infty$, and let $2\leq p<\infty$.
For every product $R=I\times J$ of axis-parallel cubes,
\begin{equation}\label{e:heat-centred}
\left \|Q_Rf  \right \|_{L^p(R)}\leq Cp  \left  \|S_{F,\mathrm h}(f)  \right \|_{L^p(R)},
\end{equation}
and
\begin{equation}\label{e:heat-means}
 \left  \|f_J+f_I-2f_R \right \|_{L^p(R)}\leq C\sqrt p\,  \left  \|\mathcal E_{\mathrm h,R}(f)\right \|_{L^p(R)}.
\end{equation}
For the full oscillation, we therefore have
\begin{align}
\|f-f_R\|_{L^p(R)}
\leq &Cp\|S_{F,\mathrm h}(f)\|_{L^p(R)}\notag\\
&\quad+C\sqrt p\,\bigg\|\Big(|\mathcal S_{\mathrm h}^{(n)}f_J(x)|^2
+|\mathcal S_{\mathrm h}^{(m)}f_I(y)|^2\Big)^{1/2}\bigg\|_{L^p(R)}.\label{e:heat-scalar}
\end{align}
If $0<\|S_{F,\mathrm h}(f)\|_{L^\infty(R)}<\infty$, then
$$
\frac1{|R|}\int_R\exp\bigg(c\frac{|Q_Rf|}{\|S_{F,\mathrm h}(f)\|_{L^\infty(R)}}\bigg)\,dz\leq2.
$$
For $0<\|\mathcal E_{\mathrm h,R}(f)\|_{L^\infty(R)}<\infty$, the coordinate terms satisfy
$$
\frac1{|R|}\int_R\exp\bigg(c\frac{|f_J+f_I-2f_R|^2}{\|\mathcal E_{\mathrm h,R}(f)\|_{L^\infty(R)}^2}\bigg)\,dz\leq2.
$$
When both areas are bounded on $R$, for $\tau>0$ we have
$$
\left |  \left  \{z\in R:|f(z)-f_R|>\tau  \right\}  \right |
\leq 2|R|\exp\bigg(-\frac{c\tau}{\|S_{F,\mathrm h}(f)\|_{L^\infty(R)}}\bigg) 
 +2|R|\exp\bigg(-\frac{c\tau^2}{\|\mathcal E_{\mathrm h,R}(f)\|_{L^\infty(R)}^2}\bigg),
$$
where terms with zero denominator are omitted.
On $\mathbb R^{n+m}$, $f\in L^p$ if and only if $S_{F,\mathrm h}(f)\in L^p$. In this case,
$$
(Cp)^{-1}\|f\|_p\leq\|S_{F,\mathrm h}(f)\|_p\leq Cp\|f\|_p .
$$
The order $p$ is sharp in the global reverse estimate and in \eqref{e:heat-centred}.
More precisely, if $a_p,b_p$ are coefficients for the flag and coordinate areas in \eqref{e:heat-scalar},
respectively, then $a_p\geq cp$ and $b_p\geq c\sqrt p$.
The exponential powers $1$ and $2$ are also optimal.
The constants depend only on $n,m$ and the fixed reconstruction functions.
\end{theorem}

\begin{proof}
\noindent\emph{Step 1.} First suppose $f\in L^2(\mathbb R^{n+m})$. We compare the coefficient area
with the heat area.
Let $h_\Lambda(r)=r^2e^{-(\Lambda r)^2}$. By \eqref{e:spectral-normalization} with $\alpha=\beta=2$,
$\Lambda^4m_\sigma(h_\Lambda)\to(2c_\sigma)^{-1}>0$.
Fix $\Lambda$ large enough that $m_\sigma(h_\Lambda)>0$ and $\lambda_0/\Lambda\leq1/8$.
For $\mathbf b=(h_\Lambda,h_\Lambda)$, $t=\Lambda r_1$ and $s=\Lambda r_2$,
$$
h_\Lambda(r_1|(\xi,\eta)|)h_\Lambda(r_2|\eta|)
=\Lambda^{-4}t^2|(\xi,\eta)|^2s^2|\eta|^2
 e^{-t^2|(\xi,\eta)|^2-s^2|\eta|^2}.
$$
Consequently,
$$
H_f^{\mathbf b}(u,v,r_1,r_2)
=m_\sigma(h_\Lambda)^{-2}\Lambda^{-4}t^2(-\Delta_{x,y})s^2(-\Delta_y)W_f(u,v,t,s).
$$
The scale measures are invariant under this substitution and
$r_1^n\max\{r_1,r_2\}^m=\Lambda^{-d}t^n\max\{t,s\}^m$. Hence
\begin{align}
&\mathcal A_{\lambda_0,F}(H_f^{\mathbf b})(x,y)^2\notag\\
 = &\frac{\Lambda^{d-8}}{m_\sigma(h_\Lambda)^4}\int_0^\infty\int_0^\infty
\int_{\substack{|x-u|<(\lambda_0/\Lambda)t\\|y-v|<(\lambda_0/\Lambda)\max\{t,s\}}}
\left  |t^2(-\Delta_{x,y})s^2(-\Delta_y)W_f(u,v,t,s)  \right |^2
\frac{du\,dv}{t^n\max\{t,s\}^m}\frac{dt}{t}\frac{ds}{s}\notag\\
 \leq &C\frac{\Lambda^{d-8}}{m_\sigma(h_\Lambda)^4}S_{F,\mathrm h}(f)(x,y)^2,\label{e:heat-comparison}
\end{align}
where the last inequality follows from $\lambda_0/\Lambda\leq1/8$ and \eqref{e:cone-lower}.
Combining \eqref{e:heat-comparison} with the local estimate \eqref{e:spectral-local}, we obtain \eqref{e:heat-centred}
for $f\in L^2(\mathbb R^{n+m})$.

\smallskip
\noindent\emph{Step 2.} For the coordinate averages we need the local heat estimate on a single cube.
Let $q\in\{n,m\}$, $g\in L^2(\mathbb R^q)$ and $E\subset\mathbb R^q$ be an axis-parallel cube.
Write
$$
H_t=m_\sigma(h_\Lambda)^{-1}t^2(-\Delta)e^{(\Lambda t)^2\Delta}g.
$$
For this $H_t$, Plancherel's identity implies
\begin{equation}\label{e:heat-energy}
\int_0^\infty\|H_t\|_{L^2(\mathbb R^q)}^2\frac{dt}{t}=\frac{\|g\|_{L^2(\mathbb R^q)}^2}{8\Lambda^4m_\sigma(h_\Lambda)^2},
\end{equation}
and the normalization of $h_\Lambda$ implies the reproducing formula
$$
g=\int_0^\infty\varphi_{q,t}^\sigma*H_t\frac{dt}{t}\quad\text{in }L^2(\mathbb R^q),
$$
where $\varphi_{q,t}^\sigma(z)=t^{-q}\varphi_q^\sigma(z/t)$.
For $\omega\in\{0,1/3,2/3\}^q$ and $Q\in\mathcal D_{\mathbb R^q}^\omega$, use \eqref{e:dyadic-choice} to define
$$
\mathcal T_{\omega,Q}=\{(y,t)\in\mathbb R^q\times(0,\infty):\kappa_q(y,t)=\omega,\ Q_q(y,t)=Q\}.
$$
The nonempty sets partition $\mathbb R^q\times(0,\infty)$. The contribution from one set is
$$
A_{\omega,Q}(z)=\int_{\mathcal T_{\omega,Q}}\varphi_{q,t}^\sigma(z-y)H_t(y)\,dy\,\frac{dt}{t}.
$$
Write $c_{\omega,Q}=|Q|^{-1/2}\|\mathbf1_{\mathcal T_{\omega,Q}}H_t(y)\|_{L^2(dy\,dt/t)}$.
For $(y,t)\in\mathcal T_{\omega,Q}$ we have $\ell(Q)/48\leq t<\ell(Q)/24$ and
$\overline{B_\infty(y,4t)}\subset Q^\circ$. The kernel is supported in $\overline{B(y,t)}$ and has integral
zero. Thus $A_{\omega,Q}$ has integral zero and its support lies at distance at least $\ell(Q)/16$ from
$\partial Q$. By Cauchy--Schwarz, for $|\alpha|\leq1$,
$$
\|\partial^\alpha A_{\omega,Q}\|_\infty
\leq C c_{\omega,Q}|Q|^{1/2}
\bigg(\int_{\ell(Q)/48}^{\ell(Q)/24}t^{-q-2|\alpha|}\,\frac{dt}{t}\bigg)^{1/2}
\leq C_q c_{\omega,Q}\ell(Q)^{-|\alpha|}.
$$
If $x\in Q$, then $|x-y|\leq48\sqrt q\,t<\Lambda t$ on $\mathcal T_{\omega,Q}$.
By disjointness of the parameter sets and the substitution $r=\Lambda t$, we have
\begin{equation}\label{e:heat-coeff}
\sum_{\omega,Q}c_{\omega,Q}^2\mathbf1_Q(x)
\leq C_q\int_0^\infty\int_{|x-y|<\Lambda t}|H_t(y)|^2\,\frac{dy}{t^q}\frac{dt}{t} 
 = C_q\frac{\Lambda^{q-4}}{m_\sigma(h_\Lambda)^2}\mathcal S_{\mathrm h}^{(q)}g(x)^2.
\end{equation}
By Plancherel's identity and Cauchy--Schwarz, for every finite selection $\mathcal E$ we have
$$
\bigg\|g-\sum_{(\omega,Q)\in\mathcal E}A_{\omega,Q}\bigg\|_{L^2(\mathbb R^q)}
\leq\|a_\sigma\|_{L^2(dr/r)}
\bigg(\sum_{(\omega,Q)\notin\mathcal E}\int_{\mathcal T_{\omega,Q}}|H_t(y)|^2\,dy\,\frac{dt}{t}\bigg)^{1/2}.
$$
By \eqref{e:heat-energy}, the right-hand side tends to zero as the finite selections exhaust the partition.
Since $Q_E$ is a contraction on $L^2(E)$, a subsequence of the centred sums converges almost everywhere on $E$ to $Q_Eg$.
Applying \eqref{e:one-coeff-local} in each of the $3^q$ adjacent dyadic systems and then using \eqref{e:heat-coeff}
and Fatou's lemma, we obtain
\begin{equation}\label{e:heat-coordinate}
\|Q_Eg\|_{L^p(E)}\leq C_q\sqrt p\,\|\mathcal S_{\mathrm h}^{(q)}g\|_{L^p(E)}.
\end{equation}

\smallskip
\noindent\emph{Step 3.} We extend both local estimates to $L^u$.
For the upper estimate on $\mathbb R^q$ we use the mean-zero Schwartz kernel
$(-\Delta)[(4\pi)^{-q/2}e^{-|\cdot|^2/4}]$. By Chanillo and Wheeden~\cite[Theorem~1(a) and the paragraph
following it, p.~278]{ChanilloWheeden}, for Schwartz $g$ and bounded $w\geq0$,
$$
\int_{\mathbb R^q}|\mathcal S_{\mathrm h}^{(q)}g|^2 w
\leq C_q\int_{\mathbb R^q}|g|^2M_q w.
$$
Let $\mathcal S_{\mathrm h,y}f(x,y)=\mathcal S_{\mathrm h}^{(m)}(f(x,\cdot))(y)$.
Here $M_d$ acts on $\mathbb R^{n+m}$ and $M_y$ on each vertical slice.
Applying the same estimate first on $\mathbb R^{n+m}$ and then on each $y$-slice,
as in Proposition~\ref{p:two-step-upper}, we have for Schwartz $f$
$$
\int_{\mathbb R^{n+m}} S_{F,\mathrm h}(f)^2w
\leq C\int_{\mathbb R^{n+m}}|\mathcal S_{\mathrm h,y}f|^2 M_dw
\leq C\int_{\mathbb R^{n+m}}|f|^2M_yM_dw .
$$
For $p>2$, each maximal operator has norm at most $Cp$ on $L^{(p/2)'}$.
By duality this yields $\|S_{F,\mathrm h}(f)\|_p\leq Cp\|f\|_p$
and, in one parameter, $\|\mathcal S_{\mathrm h}^{(q)}g\|_p\leq C_q\sqrt p\,\|g\|_p$.
For $p=2$ take $w=1$. The spatial kernels at fixed positive scales are Schwartz, so H\"older's inequality in the convolution and Fatou's lemma extend these bounds to $L^p$.

When $2\leq p\leq u$, approximate $f$ in $L^u$ by Schwartz functions.
The upper bound at exponent $u$ ensures convergence of the areas in $L^u$ and hence in $L^p(R)$.
Centring on $R$ is bounded on $L^p(R)$, so \eqref{e:heat-centred} follows.
For $1\leq u<p$, choose inner rectangles $R_i\uparrow R^\circ$ and mollifiers $\rho_i$ as in
Corollary~\ref{c:two-extension}. Then $\rho_i*f\in L^p$, and $R_i-\operatorname{supp}\rho_i\subset R$.
By Minkowski's inequality in the heat area integral,
\begin{align*}
\|Q_{R_i}(\rho_i*f)\|_{L^p(R_i)}
\leq &Cp\|\rho_i*(\mathbf1_R S_{F,\mathrm h}(f))\|_{L^p(R_i)}
\leq Cp\|S_{F,\mathrm h}(f)\|_{L^p(R)} .
\end{align*}
Extend $Q_{R_i}(\rho_i*f)$ by zero to $R$. As in Corollary~\ref{c:two-extension}, a subsequence
converges almost everywhere on $R$ to $Q_Rf$. Fatou's lemma proves \eqref{e:heat-centred}.
The same two approximations on a cube $E$ prove \eqref{e:heat-coordinate} for $g\in L^u$.
On expanding rectangles the coordinate averages tend to zero almost everywhere, as in Theorem~\ref{t:flag-reverse}.
Applying Fatou's lemma in \eqref{e:heat-centred}, we obtain $\|f\|_p\leq Cp\|S_{F,\mathrm h}(f)\|_p$ whenever
$S_{F,\mathrm h}(f)\in L^p$. The global upper bound proves the converse.

By Jensen's inequality, $f_J\in L^u(\mathbb R^n)$ and $f_I\in L^u(\mathbb R^m)$.
Applying \eqref{e:heat-coordinate} to $f_J$ on $I$ and $f_I$ on $J$, we have
\begin{align*}
\|f_J+f_I-2f_R\|_{L^p(R)}
\leq &|J|^{1/p}\|f_J-f_R\|_{L^p(I)}+|I|^{1/p}\|f_I-f_R\|_{L^p(J)}\\
\leq &C\sqrt p\,\|\mathcal E_{\mathrm h,R}(f)\|_{L^p(R)}.
\end{align*}
Combining \eqref{e:heat-means} and \eqref{e:heat-centred} in \eqref{e:coord-splitting}, we obtain \eqref{e:heat-scalar}.
Summing the moments as in \eqref{e:flag-exp-series}, we obtain the exponential estimates for the doubly centred term
and the coordinate averages. The tail bound follows by splitting the level set as in \eqref{e:two-tails}.
A term with zero area vanishes by the case $p=2$.

\smallskip
\noindent\emph{Step 4.} We test sharpness with the functions used for the Poisson area.
Write $\beta_{\mathrm h}(u)=u^2e^{-cu^2}$ with $c>0$ sufficiently small.
On a frequency region with $|(\xi,\eta)|\simeq\lambda$ and $|\eta|\simeq\mu$,
$$
(t|(\xi,\eta)|)^2(s|\eta|)^2e^{-t^2|(\xi,\eta)|^2-s^2|\eta|^2}
\leq C\beta_{\mathrm h}(t\lambda)\beta_{\mathrm h}(s\mu).
$$
For $\lambda,\mu>0$,
$$
\int_0^\infty\beta_{\mathrm h}(t\lambda)\beta_{\mathrm h}(t\mu)\frac{dt}{t}
=\frac{\lambda^2\mu^2}{2c^2(\lambda^2+\mu^2)^2}
\leq C\min\{\lambda/\mu,\mu/\lambda\}^2.
$$
In particular, for $\lambda>0$ and $L\in\mathbb N_0$,
\begin{equation}\label{e:heat-scale-sum}
\int_0^\infty\bigg(\sum_{j=0}^L\beta_{\mathrm h}(2^j\lambda t)\bigg)^2\frac{dt}{t}
\leq C\sum_{j,r=0}^L2^{-2|j-r|}\leq C(L+1).
\end{equation}
Use the functions $h_L$, $f_L$ and $F_{L,M}$ from Lemma~\ref{l:one-lacunary} and
Proposition~\ref{p:two-orders}. The frequencies of $h_L$ have sizes $2^j$, $0\leq j\leq L$,
and each piece has uniformly bounded Fourier $L^1$ norm. By Fourier inversion and \eqref{e:heat-scale-sum},
\begin{equation}\label{e:heat-one-sup}
\|\mathcal S_{\mathrm h}^{(q)}h_L\|_\infty^2
\leq C\int_0^\infty\bigg(\sum_{j=0}^L\beta_{\mathrm h}(2^jt)\bigg)^2\frac{dt}{t}\leq CL.
\end{equation}
For $f_L$ the two frequency sizes are comparable to $2^{L+2+j}$ and $2^r$, $0\leq j,r\leq L$.
Since $\int_{\mathbb R^{n+m}}\chi_{t,s}=|B_{\mathbb R^{n+m}}(0,1)|\,|B_{\mathbb R^m}(0,1)|\,t^{n+m}s^m$, we obtain
\begin{equation}\label{e:heat-flag-sup}
\|S_{F,\mathrm h}(f_L)\|_\infty^2
\leq C\int_0^\infty\bigg(\sum_{j=0}^L\beta_{\mathrm h}(2^{L+2+j}t)\bigg)^2\frac{dt}{t}
\int_0^\infty\bigg(\sum_{r=0}^L\beta_{\mathrm h}(2^rs)\bigg)^2\frac{ds}{s}
\leq CL^2 .
\end{equation}
For $F_{L,M}=h_L\otimes H_M/M$, the two sizes are comparable to $2^j$ and
$\rho_r^{-1}$, $1\leq r\leq M$, where $\rho_r=\rho 2^{r-1}$. Summing in the reverse order
in the second integral, we have
\begin{equation}\label{e:heat-mixed-sup}
\|S_{F,\mathrm h}(F_{L,M})\|_\infty^2
\leq \frac C{M^2}\int_0^\infty\bigg(\sum_{j=0}^L\beta_{\mathrm h}(2^jt)\bigg)^2\frac{dt}{t}
\int_0^\infty\bigg(\sum_{r=1}^M\beta_{\mathrm h}(s/\rho_r)\bigg)^2\frac{ds}{s}
\leq C L/M .
\end{equation}

Since $\int_0^\infty r^4e^{-2r^2}\,dr/r=1/8$, the one-parameter factor in the squared area norm is $\omega_q/8$.
For the flag area, integrating in both scales yields
$$
\|S_{F,\mathrm h}(f)\|_2^2=\frac{\omega_d\omega_m}{64}\|f\|_2^2.
$$
By Lemma~\ref{l:one-lacunary} and the horizontal dilation in $f_L$,
$\|h_L\|_2\leq C$ and $\|f_L\|_2\leq C2^{-(L+2)n/2}$. Interpolating with \eqref{e:heat-flag-sup},
we obtain
\begin{equation}\label{e:heat-example-lp}
\|S_{F,\mathrm h}(f_L)\|_p\leq C L^{1-2/p}2^{-(L+2)n/p}.
\end{equation}
Take $L=\lceil Kp\rceil$ as in Proposition~\ref{p:two-orders}. Both averages vanish on $R_L$.
Using $|R_L|=2^{-(L+2)n}|E_1||E_2|$ in \eqref{e:zero-trace-values} and applying \eqref{e:heat-example-lp}, we obtain
$$
\|Q_{R_L}f_L\|_{L^p(R_L)}
=\|f_L\|_{L^p(R_L)}
\geq (L^2/4)(c|R_L|2^{-dL})^{1/p}
\geq c p\|S_{F,\mathrm h}(f_L)\|_p .
$$
This proves sharpness in the global and local reverse estimates and gives $a_p\geq cp$.
For $F_{L,M}$ the only nonzero coordinate average on $R_*=E\times J_0$ is
$(H_M)_{J_0}h_L/M$, and $(H_M)_{J_0}/M\to1$. For fixed $p,L$, the flag term tends to zero
by \eqref{e:heat-mixed-sup}. Letting $M\to\infty$ in \eqref{e:heat-scalar} with coefficients $a_p,b_p$,
we obtain, by \eqref{e:heat-one-sup} and Lemma~\ref{l:one-lacunary},
$$
cp\leq\|h_L\|_{L^p(E)}
\leq b_p\|\mathcal S_{\mathrm h}^{(n)}h_L\|_{L^p(E)}
\leq C b_p\sqrt p .
$$

For the exponential powers the same two level sets as in Proposition~\ref{p:two-orders} suffice.
On the set in \eqref{e:zero-trace-values}, the flag oscillation is at least $L^2/4$ and its area
is bounded by $CL$. Since the relative measure is at least $c2^{-dL}$, the flag power is at most $1$
in both the exponential and the tail estimates. For $F_{L,M}$, the factor
$(H_M)_{J_0}/M$ cancels in the coordinate exponential estimate. On the ball in Lemma~\ref{l:one-lacunary},
$h_L\geq L/2$ and $\mathcal S_{\mathrm h}^{(n)}h_L\leq C\sqrt L$, so the coordinate power is at most $2$.
For the two-term tail, fix $\alpha>0$ and apply \eqref{e:trace-values} with
$M=\lceil L^{2/\alpha}\rceil$ and $\tau=L/4$. A coordinate power $\beta>2$
would imply, by \eqref{e:heat-mixed-sup} and \eqref{e:heat-one-sup},
$$
c2^{-nL}\leq C e^{-c_1L^{1+\alpha/2}}+C e^{-c_2L^{\beta/2}}=o(2^{-nL}).
$$
Thus the coordinate power in the two-term tail is at most $2$, completing the proof.
\end{proof}

The upper bounds for exponential weights and commutators in Proposition~\ref{p:coordinate-weights}
and Corollary~\ref{c:commutators} also hold for the coordinate heat areas. Indeed, the
heat operators commute with averaging in the other coordinate. The dyadic averages in the proof of
Proposition~\ref{p:coordinate-weights} converge to $f$ in $L^2(\mathbb R^{n+m})$. By the identity
$\|\mathcal S_{\mathrm h}^{(q)}g\|_2^2=(\omega_q/8)\|g\|_2^2$, the dyadic supremum bounds both slice areas. Applying \eqref{e:heat-coordinate}
on each block of any $I\times J$, we obtain for real $f\in L^2(\mathbb R^{n+m})$,
$$
\langle|f-f_{I\times J}|^p\rangle_{I\times J}^{1/p}
\leq C\sqrt p\sup_{Q\in\mathcal D}\|\mathcal E_{\mathrm h,Q}(f)\|_{L^\infty(Q)},\qquad 2\leq p<\infty .
$$
Thus \eqref{e:exp-weights} and \eqref{e:commutators} hold with
$A(f)=\sup_{Q\in\mathcal D}\|\mathcal E_{\mathrm h,Q}(f)\|_{L^\infty(Q)}$
by the same moment and conjugation estimates.

\medskip
\medskip
\medskip
\medskip
\noindent\textbf{Acknowledgements:}\par\nobreak
\noindent Li is supported by the Australian Research Council (DP220100285 and DP260100485).
Wu is supported by the National Natural Science Foundation of China (No.~12201002).

This work began in 2023 when Wu was a postdoctoral fellow at Macquarie University.
It continued during Wick's visits to Macquarie University in 2023 and 2024,
Li's visits to Auburn University in 2023 and 2024, and Li's visit to Anhui University in 2025. It was finalized in Feb 2026 when Wu visited Macquarie University.


\end{document}